\documentclass[a4paper,11pt,parskip=half]{scrartcl}

\usepackage[T1]{fontenc}
\usepackage{amsthm,amsmath,amsfonts,amssymb,mathtools,amscd}
\usepackage{aliascnt}
\usepackage{hyperref}
\hypersetup{
  colorlinks=true,
  linkcolor=black,
  citecolor=blue,
  urlcolor=blue,
  breaklinks,
  bookmarksopen
}
\usepackage{xurl}
\usepackage{graphicx,xcolor,color}

\usepackage{enumitem}
\usepackage{float,caption,subfig,grffile,wrapfig}
\usepackage{etoolbox}
\usepackage{array,multirow,tabularx,tabulary,booktabs,diagbox,makecell}
\usepackage{tikzsymbols}
\usepackage{tikz}
\usetikzlibrary{arrows,petri,topaths}
\usepackage{tkz-berge,tkz-graph}
\usepackage{adjustbox}
\allowdisplaybreaks[3]
\usepackage[nameinlink,capitalize]{cleveref}
\numberwithin{equation}{section}

\newtheoremstyle{thmstyleone}%
  {6pt}{6pt}{\itshape}{}{}{.}{0.5em}{}
\newtheoremstyle{thmstyletwo}%
  {6pt}{6pt}{\normalfont}{}{}{.}{0.5em}{}

\theoremstyle{thmstyleone}
\newtheorem{theorem}{Theorem}[section]

\newaliascnt{assumption}{theorem}

\aliascntresetthe{assumption}
\newaliascnt{proposition}{theorem}
\newtheorem{proposition}[proposition]{Proposition}
\aliascntresetthe{proposition}
\newaliascnt{lemma}{theorem}
\newtheorem{lemma}[lemma]{Lemma}
\aliascntresetthe{lemma}
\newaliascnt{definition}{theorem}
\newtheorem{definition}[definition]{Definition}
\aliascntresetthe{definition}
\newaliascnt{corollary}{theorem}
\newtheorem{corollary}[corollary]{Corollary}
\aliascntresetthe{corollary}
\newaliascnt{conjecture}{theorem}

\aliascntresetthe{conjecture}

\theoremstyle{thmstyletwo}
\newaliascnt{example}{theorem}

\aliascntresetthe{example}

\newtheorem*{claim}{Claim}

\crefname{theorem}{Theorem}{Theorems}
\crefname{section}{Section}{Sections}
\crefname{figure}{Figure}{Figures}
\crefname{subfigure}{Figure}{Figures}
\crefname{table}{Table}{Tables}
\crefname{lemma}{Lemma}{Lemmas}
\crefname{assumption}{Assumption}{Assumptions}
\crefname{corollary}{Corollary}{Corollaries}
\crefname{definition}{Definition}{Definitions}
\crefname{example}{Example}{Examples}
\crefname{remark}{Remark}{Remarks}
\crefname{proposition}{Proposition}{Propositions}
\crefname{conjecture}{Conjecture}{Conjectures}
\crefname{appendix}{Appendix}{Appendices}
\crefname{subappendix}{Appendix}{Appendices}

\providecommand{\ManuscriptAppendixHeadingSetup}{}
\newcommand{\ManuscriptAppendixReferenceSetup}{%
  \renewcommand{\thesection}{\Alph{section}}%
  \renewcommand{\thesubsection}{\thesection.\arabic{subsection}}%
  \renewcommand{\thesubsubsection}{\thesubsection.\arabic{subsubsection}}%
  \crefalias{section}{appendix}%
  \crefalias{subsection}{appendix}%
  \crefalias{subsubsection}{appendix}%
  \ManuscriptAppendixHeadingSetup
}

\newcommand{\ManuscriptAppendixTheoremNumbering}{%
  \renewcommand{\thetheorem}{\Alph{section}.\arabic{theorem}}%
  \aliascntresetthe{assumption}%
  \aliascntresetthe{proposition}%
  \aliascntresetthe{lemma}%
  \aliascntresetthe{definition}%
  \aliascntresetthe{corollary}%
  \aliascntresetthe{conjecture}%
  \aliascntresetthe{example}%
}

\newcommand*{\dif}{\mathop{}\!\mathrm{d}}
\newcommand*{\comp}{\mathbin{\vcenter{\hbox{\scalebox{.6}{$\circ$}}}}}
\newcommand*{\bigcomp}{\operatornamewithlimits{\kern .4em\prod\nolimits^{\circ}}}
\newcommand*{\iid}{\ifmmode\mathrel{\overset{\scalebox{.6}{ i.i.d.}}{\scalebox{1.2}[1]{$\sim$}}}\else i.i.d. \fi}

\newcommand*{\Parallel}{\mathrel{\rlap{$/$}{\kern2pt/}}}
\providecommand*{\Perp}{\mathrel{\rlap{$\perp$}{\kern2pt\perp}}}

\newcommand\extralabel[4][0mm]{\node[label={[label distance=#1]#2:#3}] at (#4){};}

\definecolor{RED}{RGB}{255,0,0}

\setkomafont{paragraph}{\normalfont\bfseries}
\makeatletter
\renewcommand{\title}[1]{\gdef\@title{\Large #1}}
\makeatother
\setkomafont{title}{\fontfamily{ptm}\selectfont\bfseries}
\setkomafont{author}{\large}

\newcommand{\ArxivAuthor}[3]{%
  {\bfseries #1}\thanks{#2.%
    \ifdefempty{#3}{}{\ Email: \href{mailto:#3}{\nolinkurl{#3}}.}}%
}

\usepackage{placeins}
\usepackage{setspace}
\usepackage[margin=0.85in,headsep=10pt]{geometry}
\usepackage{times}
\usepackage{newpxmath}
\usepackage[authoryear]{natbib}
\setcitestyle{square,citesep={;},aysep={},yysep={,}}

\renewenvironment{abstract}
  {\vspace{10pt}\begin{quote}\small\textbf{Abstract.\quad}}
  {\end{quote}}
\newenvironment{keywords}{\par\textbf{Keywords: }}{}

\usepackage{microtype}

\makeatletter
\renewcommand\@seccntformat[1]{\csname the#1\endcsname.\quad}
\renewcommand\thesubsection{\thesection.\arabic{subsection}}
\renewcommand\thesubsubsection{\thesubsection.\arabic{subsubsection}}
\renewcommand\section{\@startsection{section}{1}{0pt}%
  {-2.5ex plus -1ex minus -.2ex}{1.3ex plus .2ex}{\normalsize\bfseries}}
\renewcommand\subsection{\@startsection{subsection}{2}{0pt}%
  {-2.25ex plus -1ex minus -.2ex}{1ex plus .2ex}{\normalsize\itshape}}
\renewcommand\subsubsection{\@startsection{subsubsection}{3}{0pt}%
  {-2.25ex plus -1ex minus -.2ex}{1ex plus .2ex}{\small\upshape}}
\makeatother

\newcommand{\PaperTitleFirst}{Wasserstein-$p$ Bounds via Cumulant-Based Edgeworth Expansion}
\newcommand{\PaperTitleSecond}{for $\alpha$-Mixing Random Fields}

\newcommand{\PaperPDFTitle}{Wasserstein-p Bounds via Cumulant-Based Edgeworth Expansion for alpha-Mixing Random Fields}

\newcommand{\PaperAuthorOne}{Tianle Liu}
\newcommand{\PaperAuthorTwo}{Morgane Austern}
\newcommand{\PaperAuthorOneEmail}{tianletl@unc.edu}
\newcommand{\PaperAuthorTwoEmail}{maustern@fas.harvard.edu}

\newcommand{\PaperResearchAffiliationFull}{Department of Statistics, Harvard University, 1 Oxford Street, Cambridge, MA 02138, USA}

\newcommand{\PaperAuthorOneCurrentAffiliationFull}{Department of Statistics and Operations Research, University of North Carolina at Chapel Hill, Chapel Hill, North Carolina 27599, USA}

\newcommand{\PaperKeywordOne}{$\alpha$-mixing random field}
\newcommand{\PaperKeywordTwo}{Wasserstein-$p$ distance}
\newcommand{\PaperKeywordThree}{Stein's method}
\newcommand{\PaperKeywordFour}{Edgeworth expansion}
\newcommand{\PaperPDFKeywords}{alpha-mixing random field, Wasserstein-p distance, Stein's method, Edgeworth expansion}
 
\hypersetup{
  pdftitle={\PaperPDFTitle},
  pdfauthor={\PaperAuthorOne; \PaperAuthorTwo},
  pdfkeywords={\PaperPDFKeywords}
}

\title{\PaperTitleFirst\\\PaperTitleSecond}
\date{}

\author{%
  \ArxivAuthor{\PaperAuthorOne}{\PaperAuthorOneCurrentAffiliationFull}{\PaperAuthorOneEmail}
  \and
  \ArxivAuthor{\PaperAuthorTwo}{\PaperResearchAffiliationFull}{\PaperAuthorTwoEmail}%
}

\newcommand{\ManuscriptBeforeAppendices}{%
  \newpage
\section*{Declaration of generative AI and AI-assisted technologies in the manuscript preparation process}
During the preparation of this work, the authors used OpenAI tools, including
GPT-5.6 Sol, for proof auditing and Lean-code assistance. The authors reviewed
and verified the outputs, edited them as needed, and take full responsibility
for the content of the published article.
  \bibliographystyle{chicago}%
  \bibliography{reference}%
  \newpage
}
\newcommand{\ManuscriptAfterAppendices}{}

\begin{document}

\maketitle
\vspace*{-55pt}

\begin{abstract}
We establish normal-approximation bounds in the Wasserstein-$p$ distance for $\alpha$-mixing random fields, extending results previously available mainly for independent or locally dependent variables. For standardized sums indexed by $T\subset\mathbb Z^d$, we obtain rates $\mathcal O(|T|^{-\beta})$, where $\beta\in(0,1/2]$ depends on $p$, the moment order, the dimension, and the mixing decay; the optimal square-root rate holds under sufficiently fast polynomial mixing. For integer $p\geq2$, a refined bound retains the stronger mixing power produced by the expansion, covers the endpoint of finite $(p+2)$-th moments, and identifies exact logarithmic corrections at critical decay. Under finite-range dependence, the endpoint moment condition gives the optimal rate for every real $p\geq1$. The proof combines a cumulant-based Edgeworth expansion with Stein's method. Its main technical tool is a constructive graph calculus that organizes high-order dependence and preserves cancellations before mixing bounds are applied. The resulting Wasserstein estimates also imply two-regime non-uniform normal-approximation bounds with sharper sample-size decay in the far tail. A Lean~4/Mathlib companion machine-checks the principal theorem interfaces and the complete paper-internal genogram argument for the expansions, relative to seven explicitly identified literature and standard inputs.
\end{abstract}

\setlength{\parskip}{6pt}

\section{Introduction}

Over the last decade, considerable efforts have been made to control the normal approximation error in the central limit theorem (CLT) using transport metrics \citep{rio2009upper,ledoux2015stein,bobkov2018berry,bonis2020stein}.
In particular, the Wasserstein-$p$ distance \citep{villani2009optimal} between two probability measures $\nu$ and $\mu$ of random variables is defined as 
\begin{equation}
\mathcal{W}_p(\nu,\mu):=\inf_{\gamma\in \Gamma(\nu,\mu)}\Bigl(\mathbb{E}_{(X,Y)\sim \gamma}[|X-Y|^p]\Bigr)^{1/p},
\end{equation}
where $\Gamma(\nu,\mu)$ denotes the set of all couplings of $\nu$ and $\mu$. 
If we set $W_n:=\frac{1}{\sqrt{n}}\sum_{i=1}^nX_i$ where $(X_{i})$ are \emph{i.i.d.} random variables with mean zero, variance one, and finite third moment, it is well known that $W_n$ converges in the Wasserstein-$1$ distance to $N(0,1)$ at the optimal rate $\mathcal O(1/\sqrt{n})$ \citep{petrov1975sums}.
More recently, \cite{ledoux2015stein,bonis2020stein,bobkov2018berry} extended this result to the Wasserstein-$p$ distances for $p\ge 1$ and demonstrated that under moment conditions we have
\begin{equation}\label{eq:wrateiid}
\mathcal{W}_p(\mathcal{L}(W_n),N(0,1))=\mathcal O(1/\sqrt{n}).
\end{equation} 
Their results have in turn been used to establish moderate deviation bounds \citep{fang2023p}, efficient concentration inequalities \citep{austern2022efficient}, and to formalize neural network approximations by Gaussian processes \citep{eldan2021non}. 

However, it is important to note that these results rely heavily on the assumption of independence whereas in many applications the data of interest is often temporally or spatially dependent. A common scenario is that such dependence gradually decays over time and space. For example, in quantitative finance interest rates are usually heavily dependent over short horizons but the dependence diminishes over longer time horizons \citep{nze2004weak,luo2020nonparametric}. In climate science or environmental modeling, the dependence between measurements taken for precipitation, air pollution, or soil properties at distant locations weakens with increasing spatial separation \citep{yan2021conditional,franzke2022stochastic,thakur2023spatial}. 

The dependence structure in such cases is often modeled using $\alpha$-mixing random sequences or, more generally, $\alpha$-mixing random fields. For a given random field, a sequence of $\alpha$-mixing coefficients $(\alpha_{\ell})_{\ell \geq 1}$ ($\alpha_{\ell} \geq 0$) is associated to quantify the dependence. These coefficients provide a standard measure of the maximum dependence between two groups of variables separated by a distance of $\ell$ \citep{kolmogorov1960strong,doukhan1994mixing,bradley2005basic}. The decay of $\alpha_{\ell}$ as $\ell$ increases reflects how dependence diminishes over time or spatial distance.

For $p = 1$, extensions of \eqref{eq:wrateiid} to locally dependent or $\alpha$-mixing variables have been established using Stein's method \citep{baldi1989normal,barbour1989central,sunklodas2007normal,bentkus2007normal,ross2011fundamentals}. For $p > 1$, however, obtaining Wasserstein-$p$ bounds as in \eqref{eq:wrateiid} under dependence proves to be significantly more challenging. While recent works have addressed this for locally dependent variables and U-statistics \citep{fang2019wasserstein,fang2023p,liu2023wassersteinp}, no previous results extend to general $\alpha$-mixing random fields.

To bridge this gap, we establish Wasserstein-$p$ bounds in the CLT for general $p \geq 1$ for $\alpha$-mixing random fields with $\alpha$-mixing coefficients $\alpha_{\ell}$ that decay polynomially in $\ell$. Consider, for example, a $d$-dimensional random field $(X_{i})_{i \in T_{n} \subset \mathbb{Z}^{d}}$ with $\alpha$-mixing coefficients $\alpha_{\ell}$, and define $W_n := \sigma_{n}^{-1} \sum_{i \in T_n} X_{i}$, where $\sigma^2_n := \operatorname{Var}\bigl(\sum_{i \in T_n} X_{i}\bigr)$. The CLT for $\alpha$-mixing random fields \citep{bolthausen1982central} states that $W_n$ is asymptotically $N(0,1)$ if $\sigma^{2}_{n}$ is non-degenerate (see \cref{sec:mainressetup}) and there exists $\delta > 0$ such that 
\begin{equation}\label{eq:alphamixclt}
\sum_{\ell=1}^{\infty} \ell^{d-1} \alpha_{\ell}^{\delta/(2+\delta)} < \infty, \quad \sup_{i, n} \|X_i\|_{2+\delta} < \infty,
\end{equation}
where $\lVert X \rVert_{p} := \bigl(\mathbb{E} [\lvert X \rvert^{p}]\bigr)^{1/p}$ denotes the $p$-norm of $X$.

In this paper, we significantly strengthen this result by providing Wasserstein-$p$ bounds for any $p \geq 1$. Specifically, we show in \cref{thm:mixingconditions} that under stronger conditions of a form similar to \eqref{eq:alphamixclt}, we have 
\begin{equation}\label{eq:ratefirsttimebeta}
\mathcal{W}_{p}(\mathcal{L}(W_n), N(0,1)) = \mathcal{O}(\lvert T_{n} \rvert^{-\beta}),
\end{equation}
where $0 < \beta \leq 1/2$ is a constant that depends on the decay rate of $\alpha_{\ell}$. In particular, for $p \in \mathbb{N}_{+}$, if the random field has moment order $r\in(p+2,\infty]$ and $\alpha_{\ell} = \mathcal{O}(\ell^{-v})$ with $v > \frac{d(p+1)r}{r-p-2}$, we achieve $\beta = 1/2$. This convergence rate is sharp, matching the rate for \emph{i.i.d.} variables.

For integer $p\geq2$, our refined expansion retains two different powers of
the mixing coefficient.  Exploiting this distinction gives convergence under
the weaker sufficient condition
\[
  v>\frac{d(p+1)r}{2(r-p-1)},
\]
including the endpoint $r=p+2$, and produces sharper polynomial rates than
those obtained by treating the two mixing tails alike.  At the critical decay
rates we retain the exact logarithmic factors instead of absorbing them into
an arbitrary loss in the power of $|T_n|$.  For noninteger $p$, the endpoint
$r=p+2$ is also available for fields with a fixed finite dependence range,
and again gives the sharp $\mathcal O(|T_n|^{-1/2})$ rate.

The contributions of this paper are threefold.
\begin{enumerate}[label=(\roman*),leftmargin=*,itemsep=1pt,topsep=3pt]
  \item We obtain quantitative Wasserstein-$p$ bounds for every $p\geq 1$
  for normalized sums of $\alpha$-mixing random fields, with explicit
  dependence on the polynomial mixing rate. Under sufficiently fast mixing
  and the stated moment conditions, the rate is
  $\mathcal{O}(|T_n|^{-1/2})$, matching the optimal independent rate.  For
  integer $p\geq2$, the refined branch reaches the endpoint of finite
  $(p+2)$-th moments and improves both the sufficient mixing threshold and the
  rates in the slower-mixing regimes.  For noninteger $p$, the same endpoint
  and sharp rate hold under fixed finite-range dependence.
  \item We establish cumulant-based Edgeworth expansions for H\"older-smooth
  test functions at any prescribed finite order, subject to the corresponding
  moment and mixing assumptions. The principal terms involve only genuine
  cumulants of $W_n$; the geometry and decay of the field enter through the
  remainder estimates.  The resulting cumulant inequalities are quantitative
  finite-sample bounds and may also be used independently of the normal
  approximation conclusion.
  \item We develop a constructive graph calculus for producing these
  expansions. Its finite-sum identities are independent of the particular
  dependence assumption: rooted trees encode factorization, integer
  identifiers encode nested spatial constraints, and two graph operations
  represent Taylor expansion and telescoping. This separates the exact
  expansion from the probabilistic estimates used to control its remainder.
\end{enumerate}
The Wasserstein bounds also yield the concentration and non-uniform
Berry--Esseen inequalities presented in \cref{sec:applications}.  In the far
tail, combining the transport and Gaussian estimates gives a faster
sample-size exponent than the uniform smoothing bound near the center.

\subsection{Machine-checked companion}\label{sec:machinecheckedcompanion}

A Lean~4/Mathlib companion checks the principal theorem interfaces and the
complete paper-internal genogram proof chain behind
\cref{thm:wfwgraph3,thm:wfwgraph4}.  In particular, the finite-range endpoint
\cref{eq:wfwgraph3endpoint} has no admitted mathematical dependency, while the
strict-moment and refined graph bounds use only the Doukhan covariance
inequality \citep{doukhan1994mixing}.

The formal dependency boundary is explicit.  The main Wasserstein theorems use
seven admitted literature or standard inputs: (i) Rio's Wasserstein--Zolotarev
comparison, with Kantorovich--Rubinstein duality at $p=1$
\citep{rio2009upper}; (ii) Barbour's independent-sum expansion
\citep{barbour1986asymptotic}; (iii) Bobkov's independent-sum Wasserstein bound
\citep{bobkov2018berry}; (iv) Liu--Austern cumulant matching
\citep{liu2023wassersteinp}; (v) Stein-solution regularity
\citep{barbour1986asymptotic}; (vi) the Doukhan covariance inequality
\citep{doukhan1994mixing}; and (vii) the Wasserstein gluing/triangle inequality
\citep{villani2009optimal}.  The software archive
  \nolinkurl{wasserstein-p-clt-support.7z} contains the source, pinned toolchain,
theorem map, dependency audit, and reproduction instructions.  OpenAI tools
assisted with the Lean code; the authors reviewed the files and checked them
with the pinned toolchain.

\subsection{Cumulant-based Edgeworth expansion}\label{sec:cumulantbasedee}

To derive the desired Wasserstein-$p$ bounds for general $p \geq 1$, \cite{liu2023wassersteinp} introduced a generalization of Stein's method by connecting the Wasserstein-$p$ distance to cumulant-based Edgeworth expansions \eqref{eq:smoothedg}. In this paper, we demonstrate that their approach is not only effective for locally dependent variables but also applicable to $\alpha$-mixing random fields. The primary challenge lies in establishing a high-order expansion for $\alpha$-mixing random fields that satisfies all the technical requirements of the method, which is considerably more difficult compared to locally dependent settings.

Specifically, given an integer $k \geq 2$ and $\omega \in (0,1]$, for any Hölder smooth function $h \in \mathcal{C}^{k-1,\omega}(\mathbb{R})$ (see \cref{thm:defholder}), we establish that under certain conditions on $(X_{i})$, there exist coefficients $h_{r,s_{1:r}}$, depending only on $h$, such that  
\begin{equation}\label{eq:smoothedg}
\mathbb{E}[h(W_n)] - \mathcal{N}h = \sum_{(r,s_{1:r}) \in \Gamma(k-1)} (-1)^{r} \prod_{j=1}^{r} \frac{\kappa_{s_{j}+2}(W_n) h_{r,s_{1:r}}}{(s_{j}+1)!} + \mathcal{O}(\lvert T_n \rvert^{-(k-1+\gamma)/2}),
\end{equation}
where $\mathcal{N}h := \mathbb{E}[h(Z)]$ with $Z \sim N(0,1)$, and $\Gamma(k-1) := \bigl\{(r,s_{1},\ldots,s_{r}) \in \mathbb{N}_{+} : \sum_{j=1}^{r}s_{j} \leq k-1\bigr\}$. For simplicity, $s_{1:r}$ is shorthand for $s_{1},\ldots,s_{r}$, and $\kappa_{\ell}(\cdot)$ denotes the $\ell$-th cumulant of $W_n$. Importantly, $\gamma$ satisfies $0 < \gamma \leq \omega$ and depends on the decay rate of the $\alpha$-mixing coefficients. Notably, $\gamma = \omega$ can be achieved if $\alpha_{\ell}$ decreases at a sufficiently large polynomial rate. As shown in \cite{barbour1986asymptotic,liu2023wassersteinp}, $\gamma = \omega$ yields the optimal rate for \emph{i.i.d.} or locally dependent variables. Moreover, \cite{rinott2003edgeworth,rotar2008edgeworth} demonstrated that \eqref{eq:smoothedg} immediately implies an expansion in the form of \eqref{eq:asympexpan}, commonly referred to as the Edgeworth expansion in the literature. Thus, we term \eqref{eq:smoothedg} the \emph{cumulant-based Edgeworth expansion}, as the expansion terms involve only the cumulants of $W_n$ and the function $h$.

As a historical note, \eqref{eq:smoothedg} with $\gamma = \omega$ was first derived for \emph{i.i.d.} random variables in \cite{barbour1986asymptotic} and later extended to locally dependent settings by \cite{rinott2003edgeworth,liu2023wassersteinp}. Under local dependence, \eqref{eq:smoothedg} was achieved through the analysis of chains of dependency neighborhoods, as demonstrated in \cite{rotar2008edgeworth,fang2019wasserstein}. However, for $\alpha$-mixing random fields, this method is not applicable. In fact, \cite{rotar2008edgeworth} showed that similar techniques could be used to derive Edgeworth expansions but not a cumulant-based expansion in the form of \eqref{eq:smoothedg}. Specifically, their approach required replacing cumulants with quantities defined using local neighborhoods. 

We, however, demonstrate that a cumulant-based expansion \eqref{eq:smoothedg} can also be obtained for $\alpha$-mixing random fields. This result, as we will show, is a crucial step in deriving Wasserstein-$p$ bounds. Additionally, we argue that the cumulant-based expansion offers deeper insights. Since a distribution on $\mathbb{R}$ is uniquely determined by its cumulants under mild conditions \citep{reed1975ii}, one might intuitively expect that $\mathbb{E}[h(W_{n})] - \mathcal{N}h$ would be small as long as $\kappa_{\ell}(W_{n})$ are sufficiently close to the cumulants of $N(0,1)$, regardless of the dependence structure of $(X_{i})$. Ideally, the dependence structure should primarily influence the control of the remainder terms, and quantities relying on local neighborhoods should not appear in the main expansion.

A natural approach to analyzing $\mathbb{E}[h(W_{n})] - \mathcal{N}h$ is to consider the so-called Stein's equation:  
\begin{equation}\label{eq:stein}
    xf(x) - f'(x) = \mathcal{N}h - h(x), \quad \forall x \in \mathbb{R},
\end{equation}  
where $h \in \mathcal{C}^{k-1,\omega}(\mathbb{R})$, owing to its resemblance in form. Specifically, if $f = \Theta h$ is chosen as the solution of \eqref{eq:stein}, it suffices to derive an expansion for $\mathbb{E}[W_{n}f(W_{n})] - \mathbb{E}[f'(W_{n})]$. Here, $\Theta$ denotes the operator mapping the function $h$ to the solution of Stein's equation associated with $h$.

Indeed, to derive \eqref{eq:smoothedg}, where the expansion depends on $W_{n}$ only through its cumulants, we require the following type of result: For $f \in \mathcal{C}^{k,\omega}(\mathbb{R})$ with $k \in \mathbb{N}_{+}$, we establish the expansion
\begin{equation}\label{eq:mysterious}
\mathbb{E}[W_n f(W_n)]=\sum_{j=1}^{k}\frac{\kappa_{j+1}(W_n)}{j!}\mathbb{E}[ f^{(j)}(W_n)]+\text{remainder terms}.
\end{equation}

Since $f$ is the solution of \eqref{eq:stein}, we have
\begin{align}
    \mathbb{E} [h(W_{n})]-\mathcal{N}h
    &=-\sum_{j=2}^{k}\frac{\kappa_{j+1}(W_{n})}{j!}\mathbb{E} [f^{(j)}(W_{n})]+\text{remainder terms}\nonumber\\
    &=-\sum_{j=2}^{k}\frac{\kappa_{j+1}(W_{n})}{j!}\mathbb{E} [\partial^{j}\Theta h(W_{n})]+\text{remainder terms},\label{eq:parttoapply}
\end{align}
where $\partial$ is the differential operator. Notice that when $k\geq 2$, this expansion depends not only on the cumulants but also on $\mathbb{E} [\partial^{j}\Theta h(W_{n})]$, which is not desirable for our purpose. To address this issue, we re-express each term $\mathbb{E} [\partial^{j}\Theta h(W_{n})]$ by exploiting \eqref{eq:parttoapply} again. This time we replace $h$ with $\partial^{j}\Theta h$ and replace $k$ with $k-j+1$ in \eqref{eq:parttoapply} for $j=2,\cdots,k$, resulting in
\begin{equation}\label{eq:subsubstitute}
    \mathbb{E} [\partial^{j}\Theta h(W_{n})]-\mathcal{N}\partial^{j}\Theta h=-\sum_{\ell=2}^{k-j+1}\frac{\kappa_{\ell+1}(W_{n})}{\ell!}\mathbb{E} [\partial^{\ell}\Theta\partial^{j}\Theta h(W_{n})]+\text{remainder terms}.
\end{equation}
Here the sum in \eqref{eq:subsubstitute} is over $\ell$ ranging from $2$ to $k-j+1$. Importantly, the upper limit of $\ell$ is set to ensure that all remainder terms are of the same order. Specifically, we associate an integer $s$ to any composition of $\partial$ and $\Theta$, defined as $2$ plus a sum of $\pm 1$, where each $\partial$ contributes $1$ and each $\Theta$ contributes $-1$. Here the constant $2$ is chosen to align with the definition of ``order'' in the constructive graph approach introduced in \cref{sec:illconst}. For example,
\begin{equation}\label{eq:increasebyone}
3\leq s(\partial^{j}\Theta)=2+j-1\leq k+1,\quad 4\leq s(\partial^{\ell}\Theta\partial^{j}\Theta)=2+\ell-1+j-1=\ell+j\leq k+1.
\end{equation}
Intuitively, $s$ represents the order of each non-cumulant factor related to $W_{n}$ in the expansion. Furthermore, we observe that for $j=k$, i.e., $s(\partial^{j}\Theta)=k+1$, the summation in \eqref{eq:subsubstitute} vanishes. In other words, if a non-cumulant factor is of order $k+1$, it is approximated by some quantity that depends on $h$ but does not depend on $W_{n}$ with an error of order $k+2$.

By substituting the equations in \eqref{eq:subsubstitute} into \eqref{eq:parttoapply}, the smallest order of non-cumulant factors in the expansion of $\mathbb{E}[h(W_{n})] - \mathcal{N}h$ increases by $1$, as demonstrated in \eqref{eq:increasebyone}. This incremental increase in order allows us to establish a proof by induction, as the substitution process can be repeated iteratively, raising the order of non-cumulant factors step by step, until all such factors reach order $k+1$. At this stage, they are approximated by quantities determined solely by $h$, thereby achieving an expansion in the form of \eqref{eq:smoothedg}.

\subsection{Constructive graph approach: role and scope}\label{sec:illconst}

The main methodological problem is to establish \eqref{eq:mysterious} at any
prescribed finite order. Although its principal part resembles the
moment--cumulant identity below, a direct expansion under long-range
dependence produces three coupled difficulties. The summation ranges are
nested and depend on previously chosen indices; Taylor and telescoping
operations generate products of generalized covariance terms; and these
terms must be recombined exactly to recover genuine cumulants. At the same
time, near and separated spatial configurations require different estimates.
In particular, taking absolute values term by term would discard
cancellations needed for the stated remainder bounds.

We call the resulting integer-decorated rooted trees \emph{spatial-expansion
genograms}, or simply \emph{genograms}. A vertex represents a summation index, the
rooted-tree structure records traversal and factorization, and the integer
identifiers record the associated spatial constraints. Adding a path
represents a Taylor expansion, whereas attaching a child with varying
identifiers represents a telescoping operation. The genealogy, identifiers,
and signed coefficients are therefore handled jointly: related genograms are
first grouped according to the exact algebraic identities, and only then are
moment, neighborhood-volume, and mixing estimates applied. This order of
operations makes it possible to use the same dependence-free expansion at
each prescribed order while confining the probabilistic input to the
remainder estimates.

The polynomial case illustrates why the desired principal terms should be
cumulants.
\begin{lemma}\label{thm:identitycumulant}
  Let $X$ be a mean-zero random variable with finite moments
  $\mu_{1},\ldots,\mu_{k+1}$ and cumulants
  $\kappa_{1},\ldots,\kappa_{k+1}$. Then
  \begin{equation}
    \mu_{k+1}=\sum_{j=1}^{k-1}\frac{k!}{j!(k-j)!}\kappa_{j+1}\mu_{k-j}+\kappa_{k+1}.
  \end{equation}
\end{lemma}

Thus, \eqref{eq:mysterious} is exact with vanishing remainder when $f$ is a
polynomial.  For a general H\"older-smooth $f$, the expansion must instead be
constructed from the normalized field sum and its remainder controlled by the
$\alpha$-mixing coefficients.

The first-order case already contains the two operations used at every order.
For a stationary mean-zero sequence, set
\begin{equation}\label{eq:definewim}
  W_{i,m}:=\frac{1}{\sigma_n}\sum_{\ell:\lvert\ell-i\rvert\geq m+1}X_\ell,
  \qquad
  W_{i,m}^{*}:=W_{i,m}+\frac{X_{i+m}}{\sigma_n}.
\end{equation}
The near part $W_n-W_{i,m}$ is treated by Taylor expansion.  For the separated
part, the exact finite telescope gives
\begin{align}
  E_2
  &:=\frac{1}{\sigma_n}\sum_{i=1}^n\mathbb E[X_i f(W_{i,m})]
   \notag\\
  &\quad-\frac{1}{\sigma_n^2}\sum_{i=1}^n
     \sum_{j:\lvert j-i\rvert\geq m+1}
       \mathbb E[X_iX_j]\,\mathbb E[f'(W_n)] \notag\\
  &=\frac{1}{\sigma_n}\sum_{i=1}^n\sum_{j\geq m+1}
       \mathbb E\!\left[X_i\{f(W_{i,j-1})-f(W_{i,j})\}\right]
   \notag\\
  &\quad-\frac{1}{\sigma_n^2}\sum_{i=1}^n
     \sum_{j:\lvert j-i\rvert\geq m+1}
       \mathbb E[X_iX_j]\,\mathbb E[f'(W_n)] \notag\\
  &=E_4+E_5+E_6.\label{eq:split2}
\end{align}
Here $E_4$ is the Taylor remainder, $E_5$ is a separated generalized-covariance
term, and $E_6$ records the centering correction.  Their full formulas, together
with the analogous near terms, are retained in the Technical Companion.

A genogram turns this calculation into a finite construction at arbitrary
order.  Vertices record summation indices, traversal and integer identifiers
record their nested spatial domains, and branches record generalized-covariance
factors.  Path insertion implements Taylor expansion, while child extensions
with varying identifiers implement telescoping.  Related genograms are kept in
their complete signed blocks until the exact cumulant and telescoping identities
have been formed; only then are moment, neighborhood-volume, and mixing bounds
applied.  The resulting finite-sum identities do not use mixing, so the geometry
and decay of the random field enter only through the remainder estimates.

Conceptually, the calculus follows four steps: encode each nested sum as an
integer-decorated rooted tree; expand it by the two graph operations; regroup
the resulting signed families exactly into cumulants and telescoping blocks;
and only then estimate the remaining local and separated configurations.
Table~\ref{tab:genogram-concept} gives a practical dictionary for this process.

\begin{table}[t]
  \centering
  \caption{Conceptual roles of the genogram construction.}
  \label{tab:genogram-concept}
  \resizebox{\linewidth}{!}{%
    \begin{minipage}{1.20\linewidth}
      \small
      \setlength{\tabcolsep}{4pt}
      \renewcommand{\arraystretch}{1.08}
      \begin{tabularx}{\linewidth}{
        @{}
        >{\raggedright\arraybackslash}p{0.20\linewidth}
        >{\raggedright\arraybackslash}X
        >{\raggedright\arraybackslash}p{0.26\linewidth}
        @{}}
        \toprule
        Difficulty & Genogram mechanism & Consequence \\
        \midrule
        Nested relative sums
          & Traversal and ancestor-based identifiers encode each successive
            index domain.
          & A finite systematic representation replaces term-by-term
            enumeration. \\
        Generalized-covariance products
          & Rooted-tree branches and identifier signs encode successive factors.
          & Signed factorizations remain visible and regroup into genuine
            cumulants. \\
        Taylor expansion and telescoping
          & Path insertion represents Taylor expansion; child extensions with
            varying identifiers represent telescoping.
          & Exact recursive identities reach any prescribed finite order. \\
        Near and separated configurations
          & Identifier magnitudes and neighborhood shells locate the first block
            separated from its ancestors.
          & Local terms use volume counts; separated terms use mixing estimates. \\
        Cancellation among related terms
          & Coefficients and signed identifier blocks remain intact before absolute
            values are taken.
          & Telescoping cancellation yields the powers in the remainder bounds. \\
        \bottomrule
      \end{tabularx}
    \end{minipage}%
  }
\end{table}

The proof architecture is summarized in \cref{sec:genoroadmap}.  The complete
definitions, bijections, low-order formulas, and remainder bookkeeping are given
in the Technical Companion, beginning with \cref{sec:mixingmainpart}.  The full
extended version of the present overview is preserved there as
\cref{sec:technical-overview}.

\subsection{Related literature}\label{sec:relatedlit}

For \emph{i.i.d.} random variables \cite{bartfai1970entfernung} first obtained the following result: For any $p\geq 1$, assuming finite exponential moments, it holds that
$$\mathcal{W}_p(\mathcal{L}(W_n),N(0,1))=\mathcal O\bigl(n^{-\frac{1}{2}+\frac{1}{p}}\bigr).$$
\cite{sakhanenko1985estimates} weakened the condition to finite $p$-th moments. However, in both works the obtained rate is sub-optimal. \cite{rio1998distances,rio2009upper} showed that in order to get $\mathcal{O}(n^{-1/2})$, it is necessary for variables to have finite ($p$+$2$)-th moments. They proved such optimal results for $p\le 2$ and conjectured that it should also be valid for any arbitrary $p> 2$, which was demonstrated to be true by \cite{bobkov2018berry,bonis2020stein}. \cite{bobkov2018berry} proved by exploiting a version of Edgeworth expansion together with Rosenthal inequalities. And \cite{bonis2020stein} improved a technique called the Ornstein--Uhlenbeck interpolation introduced by \cite{ledoux2015stein}, and used it along with Stein exchangeable pair arguments. Notably \cite{bonis2020stein} was able to show their results for \emph{i.i.d.} random vectors. Moreover, we mention that for the special case $p=2$ \eqref{eq:wrateiid} can also be shown via its connection with the Kullback-Leibler divergence implied by the HWI inequality \citep{otto2000generalization} and Talagrand quadratic transport inequality \citep{talagrand1996transportation}.

Unlike the independent case, much less is known for Wasserstein-$p$ CLT of dependent variables. \cite{baldi1989normal,barbour1989central} obtained Wasserstein-$1$ bounds under local dependence, and \cite{sunklodas2007normal,bentkus2007normal} derived first order cumulant-based Edgeworth expansions that further lead to Wasserstein-$1$ bounds for $\alpha$-mixing random fields. \cite{fang2019wasserstein} showed sharp Wasserstein-$2$ bounds for locally dependent variables by relating it to the second order cumulant-based Edgeworth expansions using the main result of \cite{rio2009upper}. \cite{liu2023wassersteinp} extended this technique by relating general Wasserstein-$p$ bounds to cumulant-based Edgeworth expansions of ($\lceil p\rceil\!\!-\!\!1$)-th order. As a result, they were able to achieve sharp Wasserstein-$p$ bounds for locally dependent variables. 
From another line of research \cite{fang2023p} modified the approach of \cite{bonis2020stein} and obtained the rate $\mathcal{O}\bigl(\lvert T_n \rvert^{-1/2}\log\, \lvert T_n \rvert\bigr)$ in Wasserstein-$p$ CLT for locally dependent random vectors.

The term Edgeworth expansion originally refers to an asymptotic result that expands $\mathbb{P}(W_n\geq t)-\Phi^{c}(t)$ in terms of $n$ and $t$ as derived in \cite{hsu1945approximate}, and was generalized or refined by \cite{bhattacharya1972recent,gotze1978asymptotic,bhattacharya1986normal,zhilova2022new}. In particular, \cite{gotze1978asymptotic} showed that for independent random vectors with mean zero and covariance matrix $\mathbf{I}$, under proper moment conditions, there exists a signed measure $\Psi_{n,k}$ such that for any general function $h$ we have
\begin{equation}\label{eq:asympexpan}
    \mathbb{E} [h(W_{n})]-\mathcal{N}h=\int h\dif \Psi_{n,k}+o(\lvert T_n\rvert^{-k/2}),
\end{equation}
where $\Psi_{n,k}$ only depends on $k,n$ and the joint distribution of $(X_{i})$, and can be implicitly defined using Edgeworth polynomials and inverse Fourier transform. In particular, if we take $h$ to be indicator functions, \eqref{eq:asympexpan} reduces to the original result in \cite{hsu1945approximate}. 

Equation \eqref{eq:asympexpan} is traditionally obtained by carefully examining the characteristic functions. It has been extended using similar approaches to various dependent settings including Markov chains \citep{statulevivcius1969limit,hipp1985asymptotic,jensen1989asymptotic,malinovskii1987limit,fernando2021edgeworth,dolgopyat2023berry,gouezel2009local}, martingales \citep{mykland1993asymptotic}, U-statistics \citep{bickel1986edgeworth,loh1996edgeworth,bentkus1997edgeworth}, $m$-dependent random fields \citep{rhee1985edgeworth,heinbich1987asymptotic,heinrich1990asymptotic}, and mixing random fields \citep{gotze1983asymptotic,lahiri1993refinements,lahiri1996asymptotic,jensen1993note,lahiri2010edgeworth}. We point out that most works from this line of research adopted Cramér's condition for the characteristic function, which is satisfied if the variables have non-degenerate absolutely continuous components.

Interestingly Cramér's condition can be removed if we require the Hölder smoothness of $h$, the studies of which \citep{jirak2021sharp,jirak2023berry} provide sharp connections between Berry--Esseen bounds and second-order Edgeworth expansions for dependent variables. For \emph{i.i.d.} variables and second-order Edgeworth expansions, \cite{fang2022edgeworth} showed that the smooth class $\mathcal{C}^{1,1}(\mathbb{R})$ could further be relaxed via Stein's method so that it includes indicator functions, and hence the classical Edgeworth result \citep{hsu1945approximate} can be recovered. For a high-order smooth function $h$, high-order expansions in the form of \eqref{eq:asympexpan} could also be derived from cumulant-based expansions \eqref{eq:smoothedg} without Cramér's condition as demonstrated in \cite{rinott2003edgeworth}, where the authors also extended the \emph{i.i.d.} result of \cite{barbour1986asymptotic} to local dependence. 

\cite{rotar2008edgeworth} followed the approach of \cite{rinott2003edgeworth} and considered $\alpha$-mixing on graphs but were unable to obtain results in exactly the culumant-based form of \eqref{eq:smoothedg} due to the limitation of their dependency neighborhood technique. Instead, they obtained a similar expansion by replacing the cumulants $\kappa_{\ell}(W_{n})$ with other quantities defined using local neighborhoods, which is sufficient for obtaining \eqref{eq:asympexpan}. However, with this modification the technique of \cite{liu2023wassersteinp} for obtaining Wasserstein-$p$ bounds does not apply any more. Moreover, such modification potentially leads to larger error rates. For $\omega=1$ and random fields with $\alpha$-mixing coefficients $\alpha_{\ell}$ decaying exponentially, the error term in \cite{rotar2008edgeworth} is at a sub-optimal rate $\mathcal{O}(\lvert T_n\rvert^{-(k+1)/2} \log^s \lvert T_n\rvert)$ where $s$ depends on $k$ and the dimension $d$. In contrast, our error term is $\mathcal{O}(\lvert T_{n} \rvert^{-(k+1)/2})$ with the same setup, and we even allow $\alpha_{\ell}$ to decay polynomially.

\subsection{Paper outline}

The separately submitted \emph{Technical Companion} is an integral reference
for this article; cross-references to its results are linked directly in the
electronic files.  In \cref{sec:mainressetup} we introduce the random-field and
mixing framework.
Sections~\ref{sec:edgeworthexpansion}--\ref{sec:applications} state the
cumulant-based expansions, Wasserstein bounds, cumulant estimates, and
concentration and Berry--Esseen consequences.  Section~\ref{sec:proofoutline}
gives the proof architecture: \cref{sec:pfoutlinewp} explains the
cumulant-matching reduction and \cref{sec:genoroadmap} the exact genogram
expansion and global remainder control.  Complete proofs and auxiliary
calculations appear in the Technical Companion, Appendices~A--I; Appendix~I
also preserves the complete version of the constructive overview condensed in
\cref{sec:illconst}.
\section{Main results}
\subsection{Setup}\label{sec:mainressetup}
A \emph{\(d\)-dimensional (discrete) random field} \((X_{i})_{i \in T}\) is a collection of random variables indexed by a subset of the lattice \(T \subset \mathbb{Z}^{d}\). The random field considered in our main results is not assumed to be stationary or independent, which means the central limit theorem does not hold without additional restrictions on the dependence structure. To address this, we use \(\alpha\)-mixing coefficients to quantify the dependence between groups of variables and to characterize how quickly the dependence decays in the random field.

\begin{definition}\label{thm:amix}
Let \((\Omega, \mathcal{F}, \mathbb{P})\) be a probability space. For two sub-\(\sigma\)-algebras \(\mathcal{A}, \mathcal{B} \subseteq \mathcal{F}\), the \emph{\(\alpha\)-mixing coefficient} between \(\mathcal{A}\) and \(\mathcal{B}\) is defined as
\begin{equation}\label{thm:deffieldalpha}
\alpha (\mathcal{A}, \mathcal{B}) := \sup_{A \in \mathcal{A}, B \in \mathcal{B}} \bigl| \mathbb{P}(A \cap B) - \mathbb{P}(A)\mathbb{P}(B) \bigr|.
\end{equation}

Given an index set \(T \subseteq \mathbb{Z}^d\), not necessarily finite, let \((X_i)_{i \in T}\) be a random field indexed by \(T\). For positive integers \(\ell\) and \(k\), the \emph{\(\alpha\)-mixing coefficient of} \((X_i)_{i \in T}\) is defined as
\begin{equation}\label{eq:deffieldalpha}
\alpha_{\ell}^{[k]} := 0 \vee \sup \bigl\{ \alpha (\mathcal{F}_{U_1}, \mathcal{F}_{U_2}) : U_1, U_2 \subseteq T, \lvert U_1 \rvert \leq k, d(U_1, U_2) \geq \ell \bigr\},
\end{equation}
where:
\begin{itemize}
    \item \(\mathcal{F}_U := \sigma(X_i : i \in U)\) for any subset \(U \subseteq T\),
    \item \(\lvert U_1 \rvert\) denotes the cardinality of \(U_1\),
    \item \(d(U_1, U_2) := \min \bigl\{ \lVert i_1 - i_2 \rVert : i_1 \in U_1, i_2 \in U_2 \bigr\}\), and
    \item \(\lVert \cdot \rVert\) denotes the maximum norm.
\end{itemize}

An infinite random field indexed by \(\mathbb{Z}^d\) is said to be \emph{\(\alpha\)-mixing of order \(k\)} if its coefficients \(\alpha_\ell^{[k]}\) decay to zero as \(\ell \to \infty\).
\end{definition}

For a sequence of finite random fields \(\bigl(X_i^{\scalebox{0.6}{$(n)$}}\bigr)_{i\in T_n}\), let \(\alpha_{\ell,n}^{[k]}\) denote the coefficient in \eqref{eq:deffieldalpha} for the \(n\)-th field. We call the array \emph{uniformly \(\alpha\)-mixing of order \(k\)} if
\[
  \sup_n \alpha_{\ell,n}^{[k]}\longrightarrow 0\quad\text{as }\ell\to\infty.
\]
This uniform condition, rather than the automatic eventual vanishing of the coefficients for each fixed finite \(T_n\), is relevant to asymptotic statements. Throughout this paper, we fix \(k = \lceil p \rceil + 1\), where \(p>0\) is the smoothness order in the Edgeworth expansion and, when \(p\geq1\), the order of the Wasserstein distance considered below. For simplicity, we will drop the superscript \([k]\) from the notation.

The coefficients in \eqref{eq:deffieldalpha} are standard for random fields \citep{bolthausen1982central,doukhan1994mixing}. For a stationary field indexed by \(\mathbb{Z}^d\), decay of the usual unrestricted strong-mixing coefficients implies ergodicity of the lattice-shift action \citep{bradley2005basic}; no ergodicity assertion is made for a finite or nonstationary field. Under appropriate moment and mixing conditions, laws of large numbers and central limit theorems are available \citep{rosenblatt1972central,bolthausen1982central}. Alternative notions of mixing coefficients have also been proposed, including \(\beta\)-mixing \citep{volkonskii1959some}, \(\rho\)-mixing \citep{kolmogorov1960strong}, and \(\phi\)-mixing \citep{ibragimov1959some,cogburn1960asymptotic}, which are typically more restrictive than \(\alpha\)-mixing, as illustrated in \cite{doukhan1994mixing,bradley2005basic}. 

Notably, the notion of $\alpha$-mixing is closely tied to other key concepts used to capture dependence in statistical mechanics, such as the Dobrushin conditions \citep{dobruschin1968description,doukhan1994mixing} and spatial mixing \citep{dyer2004mixing}. Under appropriate regularity constraints, \(\alpha\)-mixing conditions have been shown to hold for various classes of random fields, including Gaussian random fields \citep{rosenblatt1972central}, \(m\)-Markovian random fields \citep{dobruschin1968description}, and Gibbs random fields \citep{georgii2011gibbs,chen2023strong}. A specific example is \(m\)-dependent random fields, where the dependence is limited to a fixed range such that \(\alpha_{\ell} = 0\) for all \(\ell \geq m+1\). Another example is Markov random fields, whose \(\alpha\)-mixing coefficients often exhibit exponential decay under general conditions \citep{doukhan1994mixing,bradley2005basic}.

For \(d = 1\), our definition of \(\alpha\)-mixing coefficients for random sequences \((X_{i})_{i=1}^{n}\) differs slightly from the more common definition:  
\begin{equation}
\widetilde{\alpha}_{\ell} := \sup_{t \geq 0} \alpha\bigl(\sigma\left(X_i: i \leq t\right), \sigma\left(X_i: i \geq t+\ell\right)\bigr),
\end{equation}
where the supremum is taken over two slices of the sequence rather than all constrained pairs of subsets of $T$. For random sequences, the coefficients defined in \eqref{eq:deffieldalpha} are sometimes referred to as the \emph{interlaced \(\alpha\)-mixing coefficients}, as discussed in \cite{bradley1993equivalent,bradley2005basic}.

For \(d \geq 2\), the finite cardinality requirement in \eqref{eq:deffieldalpha} plays a crucial role. Without this constraint, an alternative definition
\[
\overline{\alpha}_{\ell} := \sup \bigl\{ \alpha (\mathcal{F}_{U_1}, \mathcal{F}_{U_2}) : U_1, U_2 \subseteq T, d(U_1, U_2) \geq \ell \bigr\}
\]  
may fail to decay to zero, even for certain well-known two-dimensional Markov random fields, as shown in \cite{dobruschin1968description}.

In this paper, we consider a sequence of random fields \(\bigl(X^{\scalebox{0.6}{$(n)$}}_{i}\bigr)_{i \in T_{n}}\) with increasing index sets \(T_1 \subseteq T_2 \subseteq \cdots \subset \mathbb{Z}^d\). We assume throughout that \(\lvert T_{n} \rvert \to \infty\) as \(n \to \infty\). While \(T_n\) is often chosen as \(\{1, \ldots, n\}^{d}\) for simplicity in previous works, our analysis is not limited to this specific choice. 

We assume that \(X^{\scalebox{0.6}{$(n)$}}_{i}\) are mean-zero random variables with finite second moments, and define the rescaled average \(W_n\) as:
\begin{equation}\label{eq:firstdefinewn}
W_n := \frac{1}{\sigma_n} \sum_{\ell \in T_n} X^{\scalebox{0.6}{$(n)$}}_{\ell}, \quad \text{where } \sigma_n^2 := \operatorname{Var}\Bigl(\sum_{\ell \in T_n} X^{\scalebox{0.6}{$(n)$}}_{\ell}\Bigr).
\end{equation}

Under general moment and mixing conditions (e.g., \eqref{moment} and \eqref{mixing}, \cite{doukhan1994mixing}), the variance \(\sigma_n^{2}\) grows at the rate \(\sigma_n^{2} = \mathcal{O}(|T_n|)\). For instance, if \(T_n = \{1, \ldots, n\}^{d}\), then \(\sigma_n^{2} = \mathcal{O}(n^{d})\). However, if the variance becomes too small, the central limit theorem may not hold, even for independent variables. To avoid this issue, we follow \cite{fang2019wasserstein} and impose a non-degeneracy condition on \(\sigma_n^{2}\), requiring:
\begin{equation}\label{rr}
    \liminf_{n \to \infty} \frac{\sigma_{n}^{2}}{\lvert T_{n} \rvert} > 0.
\end{equation}

\subsection[Cumulant-based Edgeworth expansion and Wasserstein-p bounds]{Cumulant-based Edgeworth expansion and Wasserstein-$p$ bounds}\label{sec:edgeworthexpansion}

For any real number \(p>0\), we establish an Edgeworth expansion of order \((\lceil p \rceil - 1)\) for \(\mathbb{E}[h(W_n)]\), where \(h\) belongs to the Hölder smoothness class \(\Lambda_p\). When \(p\geq1\), this expansion is then used to derive upper bounds on \(\mathcal{W}_p(\mathcal{L}(W_n), N(0,1))\), the Wasserstein-\(p\) distance between the law of \(W_n\) and the standard normal distribution. To proceed, we begin by providing a precise definition of the function class \(\Lambda_p\).

\begin{definition}\label{thm:defholder}
Let \(k \in \mathbb{N}\), \(f \in \mathcal{C}^{k}(\mathbb{R})\), and \(\omega \in (0,1]\). We say that \(f\) belongs to the \emph{H\"older space} \(\mathcal{C}^{k,\omega}(\mathbb{R})\) with \emph{exponent} \(\omega\) if and only if 
\begin{equation}
    \lvert f \rvert_{k, \omega} := \sup_{x \neq y \in \mathbb{R}} \frac{\lvert \partial^{k} f(x) - \partial^{k} f(y) \rvert}{\lvert x - y \rvert^{\omega}} < \infty,
\end{equation}
where \(\lvert f \rvert_{k, \omega}\) is referred to as the \emph{H\"older coefficient} of \(f\). 

For \(p>0\), we further define \(\Lambda_{p} := \{ f \in \mathcal{C}^{\lceil p \rceil - 1, \omega}(\mathbb{R}) : \lvert f \rvert_{\lceil p \rceil - 1, \omega} \leq 1 \}\), where \(\omega = p + 1 - \lceil p \rceil\). In particular, \(\Lambda_p=\{f\in\mathcal C^{0,p}(\mathbb R):\lvert f\rvert_{0,p}\leq1\}\) when \(0<p\leq1\).
\end{definition}

To ensure the convergence \(W_{n} \to N(0,1)\) in the Wasserstein metrics, it is well established that certain moments of \(\bigl(X^{\scalebox{0.6}{$(n)$}}_{i}\bigr)_{i \in T_n}\) must exist, and the mixing coefficients \((\alpha_{\ell,n})\) must be carefully controlled \citep{rio2013inequalities}. For our purposes, these conditions must depend on \(p\), reflecting the order of the Edgeworth expansions or the Wasserstein metrics. Specifically, for the cumulant-based Edgeworth expansion below, we assume the following for some moment order \(r\in(p+2,\infty]\):

1. The random variables satisfy the moment condition:
\begin{equation}\label{moment}
  \sup_{n,\, i \in T_{n}} \bigl\lVert X^{\scalebox{0.6}{$(n)$}}_{i}\bigr\rVert_r < \infty,
\end{equation}
where \(\lVert\cdot\rVert_\infty\) is the essential-supremum norm. When \(r=\infty\), ratios of the form \((r-a)/r\) and \(r/(r-a)\), for finite \(a\), are interpreted by their limits and equal \(1\).

2. The $\alpha$-mixing coefficients satisfy the summability condition up to the largest scale used by the expansion:
\begin{equation}\label{mixing}
  \sup_{n} \sum_{\ell=1}^{\lfloor\lvert T_n\rvert^{1/d}\rfloor} \ell^{d-1} \alpha_{\ell,n}^{(r-p-2)/r} < \infty.
\end{equation}

Notably, the mixing condition \eqref{mixing} becomes weaker as the number of moments \(X^{\scalebox{0.6}{$(n)$}}_{i}\) admits increases. For example, when the random variables are uniformly bounded, i.e., \(\sup_{n,\, i \in T_{n}} \bigl\| X^{\scalebox{0.6}{$(n)$}}_{i} \bigr\|_{\infty} < \infty\), condition \eqref{mixing} simplifies to:
\[
\sup_{n} \sum_{\ell=1}^{\lfloor\lvert T_n\rvert^{1/d}\rfloor} \ell^{d-1} \alpha_{\ell,n} < \infty,
\]
which aligns closely with conditions commonly imposed in the literature \citep{rio2013inequalities,austern2022limit}.

Moreover, the required rate of mixing depends on the dimension \(d\): higher dimensions demand faster mixing. For instance, a uniform polynomial bound of the form \(\sup_n\alpha_{\ell,n}\leq C\ell^{-v}\) implies \eqref{mixing} whenever \(v>dr/(r-p-2)\).

\begin{theorem}\label{thm:hugehugethm}
Let $\bigl(X^{\scalebox{0.6}{$(n)$}}_{i}\bigr)_{i\in T_n}$ be a sequence of mean-zero random fields indexed by $T_{n}\subset \mathbb{Z}^{d}$ with $\alpha$-mixing coefficients $(\alpha_{\ell,n} )_{\ell\geq 1}$. Given $p>0$ with $\omega:=p+1-\lceil p\rceil\in (0,1]$, suppose that there exists $r\in(p+2,\infty]$ such that conditions \eqref{rr}, \eqref{moment} and \eqref{mixing} hold.
   Define
  \begin{equation}\label{eq:defofm1}
    M_{n}:= \lvert T_n \rvert^{-p/2}\sum_{\ell=1}^{\lfloor\lvert T_n\rvert^{1/d}\rfloor}\ell^{d(p+1)-\omega}\alpha_{\ell,n}^{(r-p-2)/r}.
  \end{equation}
  \begin{enumerate}[label=(\alph*)]
  \item\label{THM:BARBOURGRAPH} The following Edgeworth expansion holds uniformly for $h \in \Lambda_{p}$:
  \begin{equation}\label{eq:hugeexpan3}
  \begin{aligned}
    \mathbb{E} [h(W_n)]-\mathcal{N}h
    &=\sum_{(r,s_{1:r})\in\Gamma( \lceil p\rceil -1)}(-1)^{r} \\
    &\quad{}\times\prod_{j=1}^{r}\frac{\kappa _{s _{j}+2}(W_n)}{(s _{j}+1)!}
    \mathcal{N}\ \Bigl[\operatornamewithlimits{\kern .4em\prod\nolimits^{\circ}}_{j=1}^{r}
    (\partial ^{s _{j}+1}\Theta)\ h\Bigr] \\
    &\quad{}+\mathcal{O}\bigl(\lvert T_n \rvert^{-p/2}\bigr)+\mathcal{O}(M_{n}).
  \end{aligned}
  \end{equation}
  where $\Gamma(\lceil p\rceil-1):=\bigl\{ r,s_{1:r}\in\mathbb{N}_{+}:\sum_{j=1}^{r}s _{j}\leq \lceil p\rceil-1\bigr\}$, with $\Gamma(0)=\varnothing$, $\Theta$ is the operator such that $\Theta h$ is the solution of the Stein's equation \eqref{eq:stein} for $h$, $\partial$ is the differential operator, and $\bigcomp$ indicates the composition of operators. Moreover, for every $1\leq j\leq\lceil p\rceil-1$,
  \begin{equation}\label{eq:lowercumulantsmain}
    \bigl\lvert\kappa_{j+2}(W_n)\bigr\rvert
    =\mathcal O\bigl(\lvert T_n\rvert^{-j/2}\bigr)
     +\mathcal O\bigl(M_n^{j/p}\bigr).
  \end{equation}
  This assertion is vacuous when $0<p\leq1$.
  \item\label{thm:amixingmain} If $p\geq1$ and $M_n\to0$ as $n\to\infty$, then $W_{n}$ converges to $N(0,1)$ in the Wasserstein-$p$ distance, and
  \begin{equation}\label{eq:amixingmain1}
    \mathcal{W}_{p}(\mathcal{L}(W_{n}),N(0,1))= \mathcal{O}(\lvert T_n \rvert^{-1/2})+
    \mathcal{O}(M_{n}^{1/p}).
  \end{equation}
  \end{enumerate}
\end{theorem}

\begin{corollary}[Finite-range endpoint]\label{thm:finite-range-noninteger}
Let $p\geq1$ be noninteger and set $k=\lceil p\rceil$.  Let
$\bigl(X_i^{\scalebox{0.6}{$(n)$}}\bigr)_{i\in T_n}$ be a sequence of
mean-zero random fields, with $W_n$ and $\sigma_n$ normalized as in
\eqref{eq:firstdefinewn}.  Suppose that \eqref{rr} holds, that
\begin{equation*}
  \sup_{n,\,i\in T_n}
  \bigl\lVert X_i^{\scalebox{0.6}{$(n)$}}\bigr\rVert_{p+2}<\infty,
\end{equation*}
and that there is a fixed $R\in\mathbb N_+$ such that
$\alpha_{R+1,n}^{[k+1]}=0$ for every $n$.  Then \eqref{eq:hugeexpan3}
holds with its two displayed error terms replaced by the single remainder
$\mathcal O(\lvert T_n\rvert^{-p/2})$, and, for
$1\leq j\leq k-1$,
\begin{equation*}
  \bigl\lvert\kappa_{j+2}(W_n)\bigr\rvert
  =\mathcal O(\lvert T_n\rvert^{-j/2}).
\end{equation*}
Moreover,
\begin{equation*}
  \mathcal W_p\bigl(\mathcal L(W_n),N(0,1)\bigr)
  =\mathcal O(\lvert T_n\rvert^{-1/2}).
\end{equation*}
\end{corollary}

This endpoint is not obtained by assigning a value to
$\alpha_{\ell,n}^{0}$.  Exact independence instead makes each complete
generalized-covariance block associated with a separated genogram vanish;
only the local terms at radius $R$ remain.

We clarify that the mixing condition \eqref{mixing} with \(p = 0\) or \(p = 1\) is sufficient to establish a central limit theorem or Wasserstein-\(1\) bounds \citep{austern2022limit, rio2013inequalities}, but it is not sufficient for Wasserstein-\(p\) convergence or cumulant-based Edgeworth expansions up to order \(p\). For our purposes, faster mixing rates and stricter moment conditions are required. 

As shown by \cite{rio2009upper}, achieving Wasserstein-\(p\) bounds at the rate \(\mathcal{O}(\lvert T_n \rvert^{-1/2})\) for \emph{i.i.d.} random variables requires the existence of \((p + 2)\)-th moments. The general mixing condition in \cref{thm:hugehugethm} imposes a slightly stricter moment requirement; \cref{thm:finite-range-noninteger} shows that this excess is unnecessary under exact finite-range dependence. The Edgeworth expansion in \cref{THM:BARBOURGRAPH} is quantitative without a convergence assumption on \(M_n\); the additional condition \(M_n\to0\) is imposed only for the Wasserstein conclusion in \cref{thm:amixingmain}. Both \eqref{mixing} and \(M_{n}\) depend on \(p\), and the error terms in \eqref{eq:hugeexpan3} and \eqref{eq:amixingmain1} also vary with \(p\) through \(M_{n}\). Their rates are made explicit in the following proposition when \(\alpha_{\ell,n}\) decreases polynomially with \(\ell\).

\begin{proposition}\label{thm:mixingconditions}
Let $p, \omega, T_{n}, \alpha_{\ell,n}, M_{n}$ be as specified in \cref{thm:hugehugethm}.
  Suppose that there are constants $C>0$ and
  \begin{equation*}
    v>\frac{(d(p+2)+2(1-\omega))r}{2(r-p-2)},
  \end{equation*}
  independent of $n$, such that $\alpha_{\ell,n}\leq C\ell^{-v}$ whenever
  $1\leq\ell\leq\lfloor\lvert T_n\rvert^{1/d}\rfloor$.
Then
  \begin{equation}\label{eq:exactMnpolyrate}
    M_{n}=
    \begin{cases}
      \mathcal O(\lvert T_n\rvert^{-p/2})
        & u>d(p+1),\\
      \mathcal O(\lvert T_n\rvert^{-p/2}\log\lvert T_n\rvert)
        & u=d(p+1),\\
      \mathcal O(\lvert T_n\rvert^{-p/2+p+1-u/d})
        & \frac{d(p+2)}2<u<d(p+1),
    \end{cases},
  \end{equation}
  where $u=(r-p-2)v/r-1+\omega$.
  In particular, if $p\geq1$, then
  \begin{equation}\label{eq:exactWppolyrate}
    \mathcal W_p\bigl(\mathcal L(W_n),N(0,1)\bigr)=
    \begin{cases}
      \mathcal O(\lvert T_n\rvert^{-1/2})
        & u>d(p+1),\\
      \mathcal O\bigl(\lvert T_n\rvert^{-1/2}
        (\log\lvert T_n\rvert)^{1/p}\bigr)
        & u=d(p+1),\\
      \mathcal O(\lvert T_n\rvert^{-1/2+(p+1)/p-u/(dp)})
        & \frac{d(p+2)}2<u<d(p+1).
    \end{cases}
  \end{equation}
  At the critical value $u=d(p+1)$, the remainder in
  \eqref{eq:hugeexpan3} is
  $\mathcal O(\lvert T_n\rvert^{-p/2}\log\lvert T_n\rvert)$.
\end{proposition} 

For sufficiently large \(u\), \cref{thm:hugehugethm} provides a sharp convergence rate, matching the rate for \emph{i.i.d.} random variables. However, the result is not optimal for \(u \leq d(p + 1)\). For instance, as shown in \cite{austern2022limit}, convergence in the Wasserstein-\(1\) metric is guaranteed for \(v > \frac{rd}{r-2}\), whereas \cref{thm:mixingconditions} only ensures this convergence for \(v > \frac{3rd}{2(r-3)}\). When \(p \in \mathbb{N}\) (i.e., \(p\) is an integer), we perform a more refined analysis and obtain a sharper rate, recovering the known convergence result for \(p = 1\).

\begin{theorem}\label{thm:amixingmain2}
Let $\bigl(X^{\scalebox{0.6}{$(n)$}}_{i}\bigr)_{i\in T_n}$ be a sequence of mean-zero random fields with $\alpha$-mixing coefficients $(\alpha_{\ell,n} )_{\ell\geq 1}$. For $p\in\mathbb{N}_{+}$, suppose that there exists $r\in[p+2,\infty]$ such that \eqref{rr} and \eqref{moment} hold, and that 
 \begin{equation}\label{eq:mixinginthm3.5}
  \sup_{n}\sum_{\ell=1}^{\infty}\ell^{d-1}\alpha_{\ell,n}^{(r-p-1)/r}<\infty, \quad M_{n,m,\delta}^{(i)}\to 0,\quad  M_{n,m}^{(ii)}\to 0\quad\textup{ as }n\to\infty
 \end{equation}
 for some positive integer sequence $m=o(\lvert T_{n} \rvert^{\frac{1}{4d}})$ and fixed $\delta\in [0,1]\cap[0,r-p-1)$,
 where
 \begin{align}
    M_{n,m,\delta}^{(i)}:=&\textstyle\lvert T_n \rvert^{-(p-1+\delta)/2}m^{dp}\sum_{\ell=m+1}^{m+1+\lfloor\frac{\lvert T_n\rvert^{1/d}}{2}\rfloor}\ell^{d\delta-\delta}{\alpha_{\ell,n}}^{(r-p-1-\delta)/r}, \label{eq:mnmidef}\\
    M_{n,m}^{(ii)}:=&\textstyle\lvert T_n \rvert^{-(p-1)/2}\sum_{\ell=m+1}^{m+1+\lfloor\frac{\lvert T_n\rvert^{1/d}}{2}\rfloor}\ell^{dp-1}{\alpha_{\ell,n}}^{(r-p-1)/r}.\label{eq:mnmiidef}
 \end{align}
 For $1\leq j\leq p-1$, set
 \begin{equation}\label{eq:lowerorderradius}
   \mathcal R_{j,n}:=\min_{1\leq b\leq m}\lvert T_n\rvert^{-j/2}
   \left(b^{d(j+1)}+
   \sum_{\ell=b+1}^{b+1+\lfloor\lvert T_n\rvert^{1/d}/2\rfloor}
   \ell^{d(j+1)-1}\alpha_{\ell,n}^{(r-j-2)/r}\right).
 \end{equation}
 Then $W_{n}$ converges to $N(0,1)$ in the Wasserstein-$p$ distance, and
 \begin{equation}\label{eq:amixingmain2exact}
   \begin{aligned}
   \mathcal{W}_{p}(\mathcal{L}(W_n),N(0,1))
   =\mathcal O\Biggl(\Biggl\{
   &\sum_{j=1}^{p-1}\mathcal R_{j,n}^{p/j}
   +\lvert T_n\rvert^{-p/2}m^{d(p+1)}\\
   &{}+M_{n,m,\delta}^{(i)}+M_{n,m}^{(ii)}\Biggr\}^{1/p}\Biggr).
   \end{aligned}
 \end{equation}
 where the sum is empty when $p=1$. In particular,
 \begin{equation}\label{eq:amixingmain2}
   \mathcal{W}_{p}(\mathcal{L}(W_n),N(0,1))
   = \mathcal{O}\bigl(\lvert T_n\rvert^{-1/2}m^{2d}\bigr)
   +\mathcal{O}\bigl({M_{n,m,\delta}^{(i)}}^{1/p}\bigr)
   +\mathcal{O}({M_{n,m}^{(ii)}}^{1/p}).
 \end{equation}
\end{theorem}

\cref{thm:amixingmain2} implies that convergence in the Wasserstein-\(1\) distance holds whenever \(v > \frac{rd}{r-2}\), as detailed in \cref{thm:w1result}. By combining \cref{thm:hugehugethm,thm:amixingmain2,thm:mixingconditions}, we obtain the following rate bounds when \(\alpha_{\ell,n}\) decays polynomially. This result will be used for applications in \cref{sec:applications}.

\begin{corollary}\label{thm:alphapoly}
Let $p, T_{n}, \alpha_{\ell,n}, W_{n}$ be as specified in \cref{thm:amixingmain2}.
Suppose that there exists $r\in[p+2,\infty]$ such that \eqref{rr} and
\eqref{moment} hold, and that
\begin{equation*}
  \alpha_{\ell,n}\leq C\ell^{-v}
\end{equation*}
for constants $C,v>0$ not depending on $n$. Put
\begin{equation*}
  A:=\frac{r-p-1}{r}v,
\end{equation*}
with $A=v$ when $r=\infty$, and put
$u:=(r-p-2)v/r$ (with $u=v$ when $r=\infty$). Whenever
$A>d(p+1)/2$ and $N\geq2$, define
\begin{align}
 \mathfrak R_{p,A}(N)&:=
 \begin{cases}
  N^{-\left\{A/(dp)-(p+1)/(2p)\right\}},
    &\dfrac{d(p+1)}2<A<dp,\\
  N^{-(p-1)/(2p)}(\log N)^{1/p},&A=dp,\\
  N^{-\left\{1/2-d/(A+dp)\right\}},&A>dp,
 \end{cases}\label{eq:refined-rate-envelope}\\
 \mathfrak S_{p,A}(N)&:=
 \begin{cases}
  \mathfrak R_{p,A}(N),&A\leq dp,\\
  N^{-\left\{1/2-d(p+1)/(2p(A+d))\right\}},&A>dp,
 \end{cases}\label{eq:separate-radius-envelope}
\end{align}
and, for $u>d(p+2)/2$, define
\begin{equation}\label{eq:general-rate-envelope}
 \mathfrak G_{p,u}(N):=
 \begin{cases}
  N^{-\left\{u/(dp)-(p+2)/(2p)\right\}},
    &\dfrac{d(p+2)}2<u<d(p+1),\\
  N^{-1/2}(\log N)^{1/p},&u=d(p+1),\\
  N^{-1/2},&u>d(p+1).
 \end{cases}
\end{equation}
Then the following assertions hold.
\begin{enumerate}[label=(\alph*)]
 \item If $p\geq2$ and $A>d(p+1)/2$, then
 \begin{equation}\label{eq:refined-poly-rate}
  \mathcal W_p\bigl(\mathcal L(W_n),N(0,1)\bigr)
  =\mathcal O\bigl(\mathfrak R_{p,A}(\lvert T_n\rvert)\bigr).
 \end{equation}
 This assertion includes the endpoint $r=p+2$.
 \item Suppose that $r>p+2$ and $u>d(p+1)/2$. Then
 \begin{equation}\label{eq:separate-radius-poly-rate}
  \mathcal W_p\bigl(\mathcal L(W_n),N(0,1)\bigr)
  =\mathcal O\bigl(\mathfrak S_{p,A}(\lvert T_n\rvert)\bigr).
 \end{equation}
 If, in addition, $u>d(p+2)/2$, then
 \begin{equation}\label{eq:combined-poly-rate}
  \mathcal W_p\bigl(\mathcal L(W_n),N(0,1)\bigr)
  =\mathcal O\left(\min\left\{
    \mathfrak S_{p,A}(\lvert T_n\rvert),
    \mathfrak G_{p,u}(\lvert T_n\rvert)\right\}\right).
 \end{equation}
\end{enumerate}
\end{corollary}

\begin{corollary}\label{thm:finite-range-integer-rate}
In the setting and notation of \cref{thm:amixingmain2}, suppose that
\eqref{rr} and \eqref{moment} hold with $r=p+2$, and that there is a fixed
integer $R\geq1$ such that
$\alpha_{\ell,n}=0$ for every $\ell>R$ and every $n$. Then
\begin{equation*}
 \mathcal W_p\bigl(\mathcal L(W_n),N(0,1)\bigr)
 =\mathcal O(\lvert T_n\rvert^{-1/2}).
\end{equation*}
\end{corollary}

\subsection[Bounds on the cumulants of Wn]{Bounds on the cumulants of $W_{n}$}\label{sec:cumulantbounds}

As a by-product of our analysis, we derive a bound on the cumulants of \(W_n\), which is of independent interest. For \emph{i.i.d.} random variables, it is well known that for integer \(p \geq 1\), \(\kappa_{p+2}(W_n) = \mathcal{O}(|T_n|^{-p/2})\). In this proposition, we extend this result to mixing random fields and observe that if the random fields exhibit sufficiently fast mixing, the obtained rates match those for \emph{i.i.d.} variables:

\begin{proposition}\label{thm:cumuctrl1}
Let $\bigl(X^{\scalebox{0.6}{$(n)$}}_{i}\bigr)$ be a sequence of
mean-zero random fields with $\alpha$-mixing coefficients
$(\alpha_{\ell,n})_{\ell\geq1}$. For $p\in\mathbb N_+$, suppose that
there exists $r\in(p+2,\infty]$ such that \eqref{rr} and
\eqref{moment} hold. Then
\begin{equation}\label{eq:printed-cumulant-bound}
  \bigl\lvert\kappa_{p+2}(W_n)\bigr\rvert
  \lesssim \lvert T_n\rvert^{-p/2}
  \left(
    1+\sum_{\ell=1}^{\lfloor\lvert T_n\rvert^{1/d}\rfloor}
    \ell^{d(p+1)-1}
    \alpha_{\ell,n}^{(r-p-2)/r}
  \right).
\end{equation}
The implicit constant is independent of $n$. In particular, if, uniformly
in $n$ and $1\leq\ell\leq\lfloor\lvert T_n\rvert^{1/d}\rfloor$,
\begin{equation*}
  \alpha_{\ell,n}^{(r-p-2)/r}\lesssim\ell^{-q}
\end{equation*}
for some $q>0$, then
\begin{equation*}
  \bigl\lvert\kappa_{p+2}(W_n)\bigr\rvert=
  \begin{cases}
    \mathcal O(\lvert T_n\rvert^{-p/2}),
      &q>d(p+1),\\
    \mathcal O(\lvert T_n\rvert^{-p/2}\log\lvert T_n\rvert),
      &q=d(p+1),\\
    \mathcal O(\lvert T_n\rvert^{p/2+1-q/d}),
      &q<d(p+1).
  \end{cases}
\end{equation*}
\end{proposition}

Similar bounds have been established in \cite{janson1988normal,heinrich1990some,gotze1995m,liu2023wassersteinp} under local dependence, and in \cite{gotze1983asymptotic,lahiri1993refinements,lahiri1996asymptotic} for \(\alpha\)-mixing random sequences. In particular, \cite{lahiri1996asymptotic} demonstrated that the rate \(\mathcal{O}(n^{-p/2})\) can be achieved for random sequences with polynomially decaying \(\alpha_{\ell,n}\), which aligns with our \cref{thm:cumuctrl1}. \cite{lahiri1996asymptotic} further extended these results to random vectors. 

While our focus is on real-valued random variables, our results generalize to random fields indexed by \(T_n \subset \mathbb{Z}^d\) for any \(d \geq 1\), unlike the sequence-based setting in \cite{lahiri1996asymptotic}. Moreover, \cite{doring2013moderate,doring2022method} highlighted that bounds on cumulants of this type are instrumental in deriving moderate deviation bounds and enhancing the understanding of normal approximation.
\subsection{Concentration and Berry--Esseen inequalities}\label{sec:applications}

Let \(T_n\), \(W_n\), and \(\alpha_{\ell,n}\) be as defined in the previous section. We apply our results to derive new concentration inequalities and non-uniform Berry--Esseen bounds for mixing random fields with polynomially decaying \(\alpha\)-mixing coefficients.

We first focus on concentration results by providing an upper bound for the tail probability \(\mathbb{P}(W_n \ge t)\) for \(t \geq 0\). It is well known that sub-Gaussian concentration inequalities can hold not only under independence but also for locally dependent random variables \citep{janson2004large,zhang2019mcdiarmid}, \(\phi\)-mixing sequences \citep{yu1994rates,samson2000concentration}, or random variables constrained by the Dobrushin interdependence coefficients \citep{chatterjee2005concentration,stroock1992logarithmic,kontorovich2008concentration}. 

For \(\alpha\)-mixing sequences, no comparable result exists in the literature. One of the most significant prior works is by \cite{merlevede2009bernstein}, where the authors demonstrated almost sub-Gaussian concentration under the assumption that the \(\alpha\)-mixing coefficients decay exponentially. More concretely, let \((X_i)_{i=1}^{n}\) be a stationary sequence of bounded, mean-zero random variables with \(\alpha\)-mixing coefficients defined as  
\begin{equation*}
  \alpha_{\ell} := 0 \vee \sup_{t \geq 0} \alpha\bigl(\sigma(X_i: i \leq t), \sigma(X_i: i \geq t + \ell)\bigr).
\end{equation*}
It was shown that if \(\alpha_{\ell} \leq c_1 e^{-c_2 \ell}\) for some \(c_1, c_2 > 0\), then there exist constants \(K_1, K_2 > 0\) such that
\begin{equation}\label{eq:oldconcentrationalpha}
\mathbb{P}\Bigl(\frac{1}{\sqrt{n}} \sum_{i=1}^n X_i \geq t\Bigr) \leq e^{-K_1 t^2} + e^{-\frac{K_2 t \sqrt{n}}{\log n \log \log n}}, \quad \forall t > 0.
\end{equation}
This inequality is \emph{almost sub-Gaussian} in \(t\), with an additional non-sub-Gaussian term that diminishes as \(n\) increases. 

However, this result requires exponentially fast decay of \(\alpha_{\ell}\) and is limited to random sequences. In contrast, we establish similar concentration inequalities for \emph{\(d\)-dimensional random fields} with \emph{polynomially} decaying \(\alpha\)-mixing coefficients.

Another well-known result related to the quantity \(\mathbb{P}(W_{n} \geq t)\) is the Berry--Esseen theorem. It guarantees that for \emph{i.i.d.} mean-zero, variance-one variables, there exists a constant \(C > 0\) such that, as long as \(\|X_i\|_3 < \infty\), we have:
\begin{equation}\label{eq:berryesseenoriginal}
\sup_{t \in \mathbb{R}} \bigl|\mathbb{P}(W_n \geq t) - \Phi^{c}(t)\bigr| \leq \frac{C \lVert X_{i} \rVert _{3}^{3}}{\sqrt{n}},
\end{equation}
where \(\Phi^{c}(t) := \mathbb{P}(Z \geq t)\) with \(Z \sim N(0,1)\) is the tail probability of the standard normal. This result has been extended to various dependent structures, including locally dependent random variables \citep{chen2004normal}, random sequences satisfying martingale approximation conditions \citep{jirak2016berry}, \(\phi\)-mixing sequences \citep{rio1996berry}, Markov chains \citep{bolthausen1982berry}, and \(\alpha\)-mixing sequences with exponentially decaying coefficients \citep{tik}.

One limitation of \eqref{eq:berryesseenoriginal}, however, is that the bound does not depend on \(t\). For large \(t\), it is natural to expect tighter bounds for \(\bigl|\mathbb{P}(W_n \geq t) - \Phi^{c}(t)\bigr|\), a refinement known as non-uniform Berry--Esseen bounds. For instance, Theorem 2.5 of \cite{chen2004normal} demonstrated that for mean-zero, variance-one locally dependent variables, under general conditions, there exists a constant \(C'\) such that:
\begin{equation}\label{eq:nuberryesseen}
\bigl|\mathbb{P}(W_n \geq t) - \Phi^{c}(t)\bigr| \leq \frac{C' \lVert X_{i} \rVert _{3}^{3}}{(1 + |t|^3) \sqrt{n}}, \quad \forall t \in \mathbb{R}.
\end{equation}
This bound decreases as \(|t|\) increases and has been extended to Bernoulli shifts \citep{jirak2015non,jirak2016berry,jirak2023berry} and associated sequences \citep{dewan2005non}. 

The next theorem isolates the general mechanism behind our applications. It
requires only a Wasserstein bound and therefore applies, in particular, to
the standardized sums of the nonstationary triangular arrays considered in
this paper.

\begin{theorem}\label{concentration}
Let $p\in[1,\infty)$, let $(W_n)$ be real-valued random variables, and let
$(a_n)$ be a deterministic sequence of nonnegative numbers such that
$a_n\to0$ and
\begin{equation}\label{eq:applicationwpinput}
  \mathcal W_p\bigl(\mathcal L(W_n),N(0,1)\bigr)\leq a_n.
\end{equation}
Then the following assertions hold.
\begin{enumerate}[label=(\alph*)]
  \item\label{itm:concent} For every $n$ and $t>0$,
  \begin{align}
    \mathbb P(W_n\geq t)
    &\leq \Phi^c(t/2)+\frac{2^p a_n^p}{t^p},
      \label{eq:concentrationgood}\\
    \mathbb P(\lvert W_n\rvert\geq t)
    &\leq 2\Phi^c(t/2)+\frac{2^{p+1}a_n^p}{t^p}.
      \label{eq:twosidedconcentration}
  \end{align}
  \item\label{itm:berry1} There are constants $C_p,c_p>0$, depending only
  on $p$, and an integer $n_0$, depending on $(a_n)$, such that, for every
  $n\geq n_0$ and $t\in\mathbb R$,
  \begin{equation}\label{eq:newberryesseengood10}
    \bigl\lvert\mathbb P(W_n\geq t)-\Phi^c(t)\bigr\rvert
    \leq \frac{C_p a_n^{p/(p+1)}}{1+\lvert t\rvert^p}.
  \end{equation}
  Moreover, for $\lvert t\rvert\geq1$ this can be combined with a sharper
  far-tail estimate to give
  \begin{equation}\label{eq:newberryesseentworegime}
    \bigl\lvert\mathbb P(W_n\geq t)-\Phi^c(t)\bigr\rvert
    \leq C_p\min\left\{
      \frac{a_n^{p/(p+1)}}{1+\lvert t\rvert^p},
      e^{-c_p t^2}+\frac{a_n^p}{\lvert t\rvert^p}
    \right\}.
  \end{equation}
\end{enumerate}
\end{theorem}

No independence, stationarity, or random-field structure is used in
\cref{concentration}; all such information enters only through the
Wasserstein bound in \cref{eq:applicationwpinput}. For an integer $p$ in the
corresponding regimes of \cref{thm:alphapoly}, write
$N=\lvert T_n\rvert$ and let $b_n$ be a fixed constant multiple of any one of
\begin{equation*}
  \mathfrak R_{p,A}(N),\qquad
  \mathfrak S_{p,A}(N),\qquad
  \min\{\mathfrak S_{p,A}(N),\mathfrak G_{p,u}(N)\},
\end{equation*}
according to the applicable conclusion there. Increasing the constant and,
if necessary, the finitely many initial values makes
\eqref{eq:applicationwpinput} hold with $a_n=b_n$. If, for example,
$b_n$ is of pure-power order $N^{-\gamma}$, then
\begin{align}
  \mathbb P(W_n\geq t)
  &\leq \Phi^c(t/2)+\frac{C}{t^pN^{p\gamma}},
    \label{eq:concentrationrateenvelope}\\
  \mathbb P(\lvert W_n\rvert\geq t)
  &\leq 2\Phi^c(t/2)+\frac{C}{t^pN^{p\gamma}},
    \label{eq:twosidedconcentrationrateenvelope}\\
  \bigl\lvert\mathbb P(W_n\geq t)-\Phi^c(t)\bigr\rvert
  &\leq C\min\left\{
    \frac{N^{-p\gamma/(p+1)}}{1+\lvert t\rvert^p},
    e^{-c t^2}+\frac{N^{-p\gamma}}{\lvert t\rvert^p}
  \right\},\qquad \lvert t\rvert\geq1.
  \label{eq:berryrateenvelope}
\end{align}
Thus the sample-size exponent in the far tail is $p\gamma$, rather than
$p\gamma/(p+1)$, once $\lvert t\rvert$ is a sufficiently large constant
multiple of $\sqrt{\log N}$.

At the critical value $u=d(p+1)$ in \cref{thm:alphapoly}, the harmonic-sum
calculation gives the exact rate
\begin{equation*}
  b_n=N^{-1/2}(\log N)^{1/p}.
\end{equation*}
Consequently, without replacing the logarithm by an arbitrary power loss,
\begin{align}
  \mathbb P(W_n\geq t)
  &\leq \Phi^c(t/2)+\frac{C\log N}{t^pN^{p/2}},
    \label{eq:criticalconcentration}\\
  \bigl\lvert\mathbb P(W_n\geq t)-\Phi^c(t)\bigr\rvert
  &\leq C\min\Biggl\{
    \frac{(\log N)^{1/(p+1)}}
      {(1+\lvert t\rvert^p)N^{p/[2(p+1)]}},\notag\\[-1mm]
  &\hspace{42mm}
    e^{-c t^2}+\frac{\log N}{\lvert t\rvert^pN^{p/2}}
  \Biggr\},\qquad \lvert t\rvert\geq1.
  \label{eq:criticalberry}
\end{align}

The concentration estimates remain almost sub-Gaussian: their only
non-Gaussian term vanishes as $a_n\to0$. When the bound in
\cref{eq:applicationwpinput} is obtained from \cref{thm:alphapoly}, this yields
concentration and non-uniform Berry--Esseen bounds for polynomially mixing
random fields, without a stationarity assumption. The power of
$\lvert t\rvert$ can be increased by using a higher Wasserstein order whenever
the corresponding moment and mixing conditions hold.
\section{Proof architecture}\label{sec:proofoutline}

This section records the logical route from the expansion machinery to the
main results. The separately submitted Technical Companion is an integral
reference for this article: it contains every proof and auxiliary calculation,
and its Appendix~A gives a result-by-result location guide.

\subsection[From cumulant expansions to Wasserstein-p bounds]
{From cumulant expansions to Wasserstein-$p$ bounds}\label{sec:pfoutlinewp}

The comparison step uses the Zolotarev ideal metric. For probability measures
$\mu$ and $\nu$ on $\mathbb R$, let
\begin{equation*}
  \mathcal Z_p(\mu,\nu)
  :=\sup_{h\in\Lambda_p}
  \left(\int h\,\dif\mu-\int h\,\dif\nu\right),
\end{equation*}
where $\Lambda_p$ is the class in the following definition.
\begin{definition}\label{thm:defzolo}
  For $p>0$, put $\omega=p+1-\lceil p\rceil\in(0,1]$ and
  \[
    \Lambda_p:=\left\{h\in\mathcal C^{\lceil p\rceil-1,\omega}(\mathbb R):
      |h|_{\lceil p\rceil-1,\omega}\leq1\right\}.
  \]
\end{definition}

\begin{lemma}[Rio's comparison theorem \citep{rio2009upper}]\label{thm:lemzolo}
  For every $p\geq1$ there is a constant $C_p>0$ such that, whenever the
  displayed quantities are finite,
  \[
    \mathcal W_p(\mu,\nu)
    \leq C_p\,\mathcal Z_p(\mu,\nu)^{1/p}.
  \]
  For $p=1$, this is the Kantorovich--Rubinstein duality.
\end{lemma}

Fix $p\geq1$ and write $k=\lceil p\rceil$. Lemma~B.3 of the Technical
Companion constructs, from the lower cumulants of $W_n$, integers $q_n\to
\infty$ and identically distributed independent variables
$\xi_1^{(n)},\ldots,\xi_{q_n}^{(n)}$ such that
\begin{equation*}
  V_n:=q_n^{-1/2}\sum_{i=1}^{q_n}\xi_i^{(n)},
  \qquad
  \kappa_j(V_n)=\kappa_j(W_n),\quad 1\leq j\leq k+1.
\end{equation*}
The required decay of these cumulants is supplied by
\eqref{eq:lowercumulantsmain}. The $p=1$ case is the corresponding base case,
where there are no higher cumulants to match.

The cumulant-based expansion in \cref{thm:hugehugethm} and Barbour's
independent-sum expansion have identical principal terms after this matching.
Their difference therefore consists only of their two controlled remainders.
Taking the supremum over $h\in\Lambda_p$ and applying
\cref{thm:lemzolo} bounds
$\mathcal W_p(\mathcal L(W_n),\mathcal L(V_n))$. Bobkov's independent-sum
bound controls $\mathcal W_p(\mathcal L(V_n),N(0,1))$, and the Wasserstein
triangle inequality gives
\begin{equation*}
  \mathcal W_p(\mathcal L(W_n),N(0,1))
  \leq \mathcal W_p(\mathcal L(W_n),\mathcal L(V_n))
     +\mathcal W_p(\mathcal L(V_n),N(0,1)).
\end{equation*}
This proves the general mechanism behind Theorem~2.3(b). The same comparison,
using the refined terminal expansion in Theorem~B.4 of the Technical
Companion, proves Theorem~2.6. The common-radius and separate-radius choices
then give the rate envelopes in Propositions~2.5 and Corollary~2.7. Under
fixed finite-range dependence, the complete separated block vanishes before
absolute values are taken, which yields Corollaries~2.4 and~2.8 at the moment
endpoint.

\subsection{Genogram expansion and remainder control}\label{sec:genoroadmap}

The engine for the preceding comparison is the expansion of
$\mathbb E[W_nf(W_n)]$. Its only external probabilistic estimate is the
following separated-block covariance inequality.

\begin{lemma}[Doukhan's covariance inequality \citep{doukhan1994mixing}]
  \label{thm:covineq}
  Suppose that $X$ and $Y$ are measurable with respect to sigma-algebras
  $\mathcal A$ and $\mathcal B$, respectively. If $a,b,c\geq1$ and
  $a^{-1}+b^{-1}+c^{-1}=1$, then
  \begin{equation}\label{eq:covineqalphap}
    |\operatorname{Cov}(X,Y)|
    \leq 8\,\alpha(\mathcal A,\mathcal B)^{1/c}
      \|X\|_a\|Y\|_b.
  \end{equation}
\end{lemma}

The algebraic part is independent of mixing. The Technical Companion defines
three evaluations of a genogram $G$: $\mathcal T_f(G)$ is a term still to be
expanded, $\mathcal S(G)$ is a coefficient independent of $f$, and
$\mathcal U_f(G)$ is an integral-form remainder. Its exact expansion theorem
has the following schematic form; all index classes and coefficients are
defined in Appendix~E of the Technical Companion.

\begin{theorem}[Genogram expansion]\label{thm:expansion1graphrestate}
  Let $G$ be a genogram and let $k\geq|G|$. For
  $f\in\mathcal C^{k-1}(\mathbb R)$,
  \begin{equation}\label{eq:expansion1graphrestate}
    \mathcal T_f(G)
    =\sum_{\substack{H\in\mathcal A_G\\ |H|\leq k}}
       a_{H,G}\,\mathcal S(H)\,
       \mathbb E[\partial^{|H|-1}f(W_n)]
     +\sum_{\substack{H\in\mathcal B_G\\ |H|=k+1}}
       b_{H,G}\,\mathcal U_f(H).
  \end{equation}
\end{theorem}

The proof is an induction using two exact operations. A varying-identifier
child extension forms a complete finite telescope, and path insertion applies
Taylor's formula with its integral remainder. The signed terms generated by
these operations are assembled before any estimate is made. Taking $G$ to be
the one-vertex genogram and identifying the resulting $\mathcal S$-sum by the
moment--cumulant identity gives the following consequence.

\begin{corollary}[Cumulant form]\label{thm:wfwgraph12restate}
  For $k\geq1$ and $f\in\mathcal C^k(\mathbb R)$,
  \begin{equation}\label{eq:wfwgraph1restate}
    \mathbb E[W_nf(W_n)]
    =\sum_{j=1}^k\frac{\kappa_{j+1}(W_n)}{j!}
       \mathbb E[\partial^jf(W_n)]
     +\sum_{|H|=k+2}b_H\mathcal U_f(H).
  \end{equation}
\end{corollary}

The remaining task is to estimate the complete signed remainder in
\eqref{eq:wfwgraph1restate}. Genograms are first regrouped with their
genealogy, orientations, identifiers, and coefficients intact. Positive
terminal-identifier ranges are telescoped across both maximum sectors, with
ties assigned to the terminal range. Only after that exact regrouping are
H\"older bounds, moment estimates, lattice-shell counts, and
\cref{thm:covineq} applied. This produces the four global estimates in
Lemma~F.5 of the Technical Companion. Combining those estimates with the
exact expansion yields the general expansion in Proposition~E.5 and the
refined integer-order expansion in Proposition~E.6.

The deductions of the main results are then as follows.
\begin{enumerate}[label=(\roman*),leftmargin=*,itemsep=2pt,topsep=3pt]
  \item Proposition~E.5, Stein-solution regularity, and induction over the
  Stein operator give Theorem~2.3(a); its coefficient estimates also give
  Proposition~2.9.
  \item Cumulant matching and the comparison argument of
  \cref{sec:pfoutlinewp} give Theorem~2.3(b), Proposition~2.5, and the
  noninteger finite-range Corollary~2.4.
  \item Proposition~E.6 gives the refined expansion in Theorem~B.4; the same
  comparison gives Theorem~2.6 and Corollaries~2.7--2.8.
  \item A Wasserstein coupling, Markov's inequality, Gaussian tail bounds,
  smoothing, and reflection give the distribution-free concentration and
  nonuniform Berry--Esseen statements in Theorem~2.10.
\end{enumerate}

The Technical Companion retains the complete proof chain: Appendices~B--D
prove the analytic deductions and applications, Appendices~E--H contain the
genogram construction and all remainder estimates, and Appendix~I records the
low-order notation concordance together with the complete overview
corresponding to \cref{sec:illconst}.
  
\ManuscriptBeforeAppendices
\appendix
\ManuscriptAppendixReferenceSetup
\ManuscriptAppendixTheoremNumbering

\section{Proof location guide}

In \cref{SEC:LOCALthmpf} we prove all the results related to the Wasserstein-$p$ bounds, i.e., \cref{thm:amixingmain} of \cref{thm:hugehugethm}, \cref{thm:mixingconditions}, \cref{thm:amixingmain2} and \cref{thm:alphapoly} using the cumulant-based Edgeworth expansions as stated in \cref{THM:BARBOURGRAPH} of \cref{thm:hugehugethm} and \cref{THM:BARBOURGRAPH2}. 
In \cref{sec:finalpflemma} we establish the bounds on cumulants, i.e., \cref{thm:cumuctrl1} and the cumulant-based Edgeworth expansions, i.e., \cref{THM:BARBOURGRAPH} of \cref{thm:hugehugethm} and \cref{THM:BARBOURGRAPH2} using \cref{thm:wfwgraph3,thm:wfwgraph4}. Then we prove \cref{concentration} in \cref{sec:concentration}.

In \cref{sec:mixingmainpart} we present the details on genograms as well as specific results in performing the constructive graph approach. In \cref{sec:tech1} we prove all the results from \cref{sec:mixingmainpart} along with formal versions of \cref{thm:expansion1graphrestate,thm:wfwgraph12restate} after introducing three key technical lemmas, i.e., \cref{THM:STEP1GRAPH,THM:STEP2GRAPH,THM:REMAINDERCTRL1234}. Finally we provide the proofs for \cref{THM:STEP1GRAPH,THM:STEP2GRAPH} in \cref{sec:lemma4} and the proof for \cref{THM:REMAINDERCTRL1234} in \cref{sec:lemma5}.

\section[Proofs of Wasserstein-p bounds]{Proofs of Wasserstein-$p$ bounds}\label{SEC:LOCALthmpf}

First, we present two lemmas on the normal approximation for independent random variables. Let $W_n$ and $\sigma_{n}$ be as defined in \eqref{eq:firstdefinewn}. \cref{thm:barbour} provides the cumulant-based Edgeworth expansion for $\mathbb{E}[h(W_n)]-\mathcal{N}h$ while \cref{thm:lemiidwp} gives an upper bound on the Wasserstein-$p$ distance between the distribution of $W_{n}$ and the standard normal.

\begin{lemma}[Theorem 1 of \cite{barbour1986asymptotic}]\label{thm:barbour}
  For any $p>0$, let $(X_{i})$ be independent variables with mean zero and finite ($p+2$)-th moment. Suppose that
  \begin{equation}\label{eq:lyapunovcond}
    \frac{1}{\sigma_{n}^{p+2}}\sum_{i=1}^{n}\mathbb{E} [\lvert X_{i} \rvert ^{p+2}]\to 0\quad \textup{as }n\to\infty.
  \end{equation}
  Then, uniformly over $h \in \Lambda_{p}$,
  \begin{align}
      &\mathbb{E}[h(W_{n})]-\mathcal{N}h\notag\\
      &={}\!\!\sum_{\substack{(r,s_{1:r})\in\\ \Gamma (\lceil p\rceil-1)}}\!\!(-1)^{r}
      \prod_{j=1}^{r}\frac{\kappa _{s _{j}+2}(W_{n})}{(s _{j}+1)!}
      \mathcal{N}\!\Bigl[\bigcomp_{j=1}^{r}(\partial ^{s _{j}+1}\Theta)\ h\Bigr]\notag\\
      &\quad{}+\mathcal{O}\biggl(\frac{1}{\sigma_{n}^{p+2}}
        \sum_{i=1}^{n}\mathbb{E} [\lvert X_{i} \rvert ^{p+2}]\biggr).\label{eq:barbour}
  \end{align}
 where $\Gamma (\lceil p\rceil -1):=\bigl\{ r, s _{1:r}\in \mathbb{N}_{+}:\sum_{j=1}^{r}s _{j}\leq \lceil p\rceil-1\bigr\}$.
\end{lemma}
Here \eqref{eq:lyapunovcond} is also known as the Lyapunov condition. If we require uniformly bounded $(p+2)$-th moments and the non-degeneracy of $\sigma_{n}$ as stated in \eqref{rr}, then the remainder term in \eqref{eq:barbour} can be rewritten as $\mathcal{O}(n^{-p/2})$.
We can see that \cref{thm:barbour} looks quite similar to \cref{THM:BARBOURGRAPH} of \cref{thm:hugehugethm}. The differences arise from the dependence structures of $\bigl(X^{\scalebox{0.6}{$(n)$}}_{i}\bigr)$ and the remainder terms in the expansions. This similarity inspires the proof of \cref{thm:amixingmain} of \cref{thm:hugehugethm}, as illustrated in \cref{sec:pfoutlinewp}.

\begin{lemma}[Theorem 1.1 of \cite{bobkov2018berry}]\label{thm:lemiidwp}
  For any $p\geq 1$, let $(X_{i})$ be independent variables with mean zero and finite ($p+2$)-th moment. If \eqref{eq:lyapunovcond} holds, then we have
  \begin{equation}\label{eq:iidwp}
    \mathcal{W}_{p}(\mathcal{L}(W_{n}), N(0,1)) \leq C_{p}\biggl(\frac{1}{\sigma_{n}^{p+2}}\sum_{i=1}^{n} \mathbb{E}[\lvert X_{i}\rvert ^{p+2}]\biggr)^{1/p},
  \end{equation}
  where $C_{p}$ depends continuously on $p$.
\end{lemma}

We shall use the following cumulant-matching lemma. For $p>1$, this is
\cite[Lemma A.8]{liu2023wassersteinp}.
\begin{lemma}\label{THM:EXISTENCEXI}
  Let $p\geq 1$ and denote $k=\lceil p\rceil$. If $p>1$, let
  $\bigl(u_{j}^{\scalebox{0.6}{$(n)$}}\bigr)_{j=1}^{k-1}$ be a sequence
  of real numbers such that $u_{j}^{\scalebox{0.6}{$(n)$}}\to 0$ as
  $n\to \infty$ for every $1\leq j\leq k-1$. No sequence $(u_j^{(n)})$
  is specified when $p=1$.
  Then there exist constants $C_{p}, C_{p}'$, depending only on $p$, and a threshold $N>0$ such that for any $n >N$, there exists $q_{n}\in \mathbb{N}_{+}$ and a random variable $\xi^{\scalebox{0.6}{$(n)$}}$ such that
  \begin{enumerate}[label=(\alph*)]
    \item \label{itm:match12} $\mathbb{E} [\xi^{\scalebox{0.6}{$(n)$}}]=0$,\quad $\mathbb{E} [(\xi^{\scalebox{0.6}{$(n)$}})^{2}]=1$.
    \item \label{itm:match3more} If $p>1$, then
    $\kappa_{j+2}(\xi^{\scalebox{0.6}{$(n)$}})
    =q_{n}^{j/2}u_{j}^{\scalebox{0.6}{$(n)$}}$ for every
    $1\leq j\leq k-1$.
    \item \label{itm:boundedaway} If $p>1$, then either all the indexed
    cumulants in part~\ref{itm:match3more} vanish, or
    $\max_{1\leq j\leq k-1}
    \bigl\lvert \kappa_{j+2}(\xi^{\scalebox{0.6}{$(n)$}}) \bigr\rvert
    \geq C_{p}>0$.
    \item \label{itm:momentbound} $\mathbb{E} [\lvert \xi^{\scalebox{0.6}{$(n)$}}\rvert^{p+2}]\leq C_{p}'$.
  \end{enumerate}
  Furthermore, $q_{n}$ can be chosen to be such that $q_{n}\to\infty$ as $n\to\infty$.
\end{lemma}
For $p>1$, the proof is given in Appendix C of the cited paper. The threshold
$N$ may depend on the sequences $\bigl(u_j^{\scalebox{0.6}{$(n)$}}\bigr)$.
For $p=1$, take $q_n=n$ and let $\xi^{\scalebox{0.6}{$(n)$}}$ be standard
normal; Conditions \ref{itm:match3more} and \ref{itm:boundedaway} then make
no assertion. We point out that
\cref{THM:EXISTENCEXI} is an asymptotic statement in the sense that, for a
given $n\leq N$, $q_n$ and $\xi^{\scalebox{0.6}{$(n)$}}$ need not be supplied.
When $p>1$, Condition \ref{itm:match3more} determines the indexed cumulants
of $\xi^{\scalebox{0.6}{$(n)$}}$, while Condition
\ref{itm:boundedaway} says that
$\max_{1\leq j\leq k-1}\lvert\kappa_{j+2}(\xi^{\scalebox{0.6}{$(n)$}})\rvert$
is either zero or bounded away from zero, while Condition
\ref{itm:momentbound} gives the required uniform absolute-moment bound.

The proofs of our main results work in three stages:
\begin{enumerate}[label=\arabic*.]
  \item By \cref{THM:EXISTENCEXI} we find a sequence of \emph{i.i.d.} random variables $\bigl(\xi^{\scalebox{0.6}{$(n)$}}_{\ell}\bigr)_{\ell}$ and a sample size $q_{n}$ such that the first $k$+$1$ cumulants of $W_n$ match the first $k$+$1$ cumulants of $V_n=q_{n}^{-1/2}\sum_{i=1}^{q_{n}}\xi^{\scalebox{0.6}{$(n)$}}_i$.
  \item \cref{thm:lemzolo} helps bound the Wasserstein-$p$ distance between the distributions of $W_n$ and $V_{n}$ in terms of $\bigl\lvert \mathbb{E} [h(W_n)]-\mathbb{E} [h(V_n)] \bigr\rvert$ for Hölder smooth functions $h$. We do so by exploiting \cref{thm:barbour} and \cref{THM:BARBOURGRAPH} of \cref{thm:hugehugethm}.
  \item \cref{thm:lemiidwp} provides us with the bound on the Wasserstein-$p$ distance between the distribution of $V_{n}$ and $N(0,1)$.
\end{enumerate} Then \cref{thm:amixingmain} of \cref{thm:hugehugethm} follows from the triangle inequality of the Wasserstein metric: $$\mathcal{W}_p(W_n,N(0,1))\le \mathcal{W}_p(\mathcal{L}(W_n),\mathcal{L}(V_n))+\mathcal{W}_p(\mathcal{L}(V_n),N(0,1)).$$

Notably in Step 2 we need lemmas on the high-order cumulant-based Edgeworth expansions for the $\alpha$-mixing fields, for which purpose we utilize \cref{THM:BARBOURGRAPH} of \cref{thm:hugehugethm} for the proof of \cref{thm:amixingmain} of \cref{thm:hugehugethm} and introduce \cref{THM:BARBOURGRAPH2} for the proof of \cref{thm:amixingmain2}. Furthermore, \cref{thm:mixingconditions,thm:alphapoly} are directly applications of \cref{thm:hugehugethm,thm:amixingmain2} for random fields with $\alpha$-mixing coefficients converging at a polynomial rate.

\begin{proof}[Proof of \cref{thm:amixingmain} of \cref{thm:hugehugethm}.]
  Let $k:=\lceil p\rceil$.
  If $p>1$, \eqref{eq:lowercumulantsmain} and the assumption $M_n\to0$
  imply that $\kappa_{j+2}(W_n)\to0$ for every $1\leq j\leq k-1$.

  If $p>1$, apply \cref{THM:EXISTENCEXI} with
  $u_{j}^{\scalebox{0.6}{$(n)$}}=\kappa_{j+2}(W_n)$ for
  $1\leq j\leq k-1$; for $p=1$, use its Gaussian construction.
  For $n$ large enough, there exist constants $C_{p}$ and $C_{p}'$, positive integer sequence $(q_{n})$ such that $q_{n}\to\infty$ as $n\to\infty$, and random variables $(\xi^{\scalebox{0.6}{$(n)$}})$ such that
  \begin{enumerate}[label=(\alph*)]
    \item $\mathbb{E} [\xi^{\scalebox{0.6}{$(n)$}}]=0$,\quad $\mathbb{E} [(\xi^{\scalebox{0.6}{$(n)$}})^{2}]=1$.
    \item \label{itm:bla2} If $p>1$, then
    $\kappa_{j+2}(\xi^{\scalebox{0.6}{$(n)$}})
    =q_{n}^{j/2}\kappa_{j+2} (W_{n})$ for every $1\leq j\leq k-1$.
    \item \label{itm:lowerbound2} If $p>1$, then either all the indexed
    cumulants in part~\ref{itm:bla2} vanish, or
    $\max_{1\leq j\leq k-1}
    \bigl\lvert \kappa_{j+2}(\xi^{\scalebox{0.6}{$(n)$}}) \bigr\rvert
    \geq C_{p}>0$.
    \item \label{itm:uppermombound2} $\mathbb{E} [ \lvert \xi^{\scalebox{0.6}{$(n)$}}\rvert^{p+2}]\leq C_{p}'$.
  \end{enumerate}

We use this construction to compare $W_n$ with a normalized sum of at least
$q_n$ \emph{i.i.d.} variables. Introduce an integer sequence
$(\widetilde q_n)$ by setting
$\widetilde q_n:=\lceil\lvert T_n\rvert^{2(p+1)/p}\rceil\vee q_n$ if
$p=1$, or if $p>1$ and
$\kappa_3(W_n)=\cdots=\kappa_{k+1}(W_n)=0$, and
$\widetilde q_n:=q_n$ otherwise. Then $\widetilde q_n\to\infty$ as
$n\to\infty$.

  Let $\xi_{1}^{\scalebox{0.6}{$(n)$}},\cdots,\xi_{\widetilde{q}_{n}}^{\scalebox{0.6}{$(n)$}}$ be \emph{i.i.d.} copies of $\xi^{\scalebox{0.6}{$(n)$}}$. Define $V_{n}:=\widetilde{q_{n}}^{-1/2}\sum_{i=1}^{\widetilde{q}_{n}}\xi_{i}^{\scalebox{0.6}{$(n)$}}$. When $p>1$, by construction, for every $1\leq j\leq k-1$ we have
  \begin{equation*}
    \kappa_{j+2}(V_{n})=\widetilde{q_{n}}^{-\frac{j+2}{2}} \sum_{i=1}^{\widetilde{q}_{n}}\kappa_{j+2}(\xi_{i}^{\scalebox{0.6}{$(n)$}})=\widetilde{q_{n}}^{-j/2}\kappa_{j+2}(\xi^{\scalebox{0.6}{$(n)$}})=\kappa_{j+2}(W_{n}).
  \end{equation*}
  For $p=1$, $V_n$ is standard normal, and the two Edgeworth expansions
  below contain no cumulant-correction term.
  Before taking absolute values, the complete signed sums indexed by
  $\Gamma(\lceil p\rceil-1)$ cancel term by term because the matched
  cumulants of $W_n$ and $V_n$ agree.
  Thus, by \cref{thm:barbour} and \cref{THM:BARBOURGRAPH} of \cref{thm:hugehugethm}, uniformly over $h\in\Lambda_{p}$ we have
  \begin{equation}\label{eq:random56}
    \bigl\lvert \mathbb{E} [h(W_{n})]-\mathbb{E} [h(V_{n})] \bigr\rvert \lesssim \lvert T_{n}\rvert^{-p/2} +M_{n}+\widetilde{q_{n}}^{-\frac{p+2}{2}}\sum_{i=1}^{\widetilde{q}_{n}}\mathbb{E} \bigl[\bigl\lvert \xi_{i}^{\scalebox{0.6}{$(n)$}} \bigr\rvert^{p+2}\bigr].
  \end{equation}

  To remove the dependency on $\xi^{\scalebox{0.6}{$(n)$}}$ in the inequality \eqref{eq:random56} we need to upper bound the following quantity in terms of $M_{n}$ as defined in \cref{thm:hugehugethm}:
  $$\widetilde{q_{n}}^{-\frac{p+2}{2}}\sum_{i=1}^{\widetilde{q}_{n}}\mathbb{E} \bigl[\bigl\lvert \xi_{i}^{\scalebox{0.6}{$(n)$}} \bigr\rvert^{p+2}\bigr].$$ 
  To do so we use the lower bounds on $(\widetilde{q}_{n})$ implied by their choices. If $p>1$ and
  $\max_{1\leq j\leq k-1}\bigl\lvert \kappa_{j+2}(W_{n})\bigr\rvert>0$,
  \cref{itm:lowerbound2} implies
  \begin{equation*}
    \begin{aligned}
      C_{p}
      &\leq \max_{1\leq j\leq k-1}
        \bigl\lvert \kappa_{j+2}(\xi_{1}^{\scalebox{0.6}{$(n)$}}) \bigr\rvert\\
      &\overset{(*)}{=}\max_{1\leq j\leq k-1}
        \bigl\{\widetilde{q_{n}}^{j/2}\bigl\lvert \kappa_{j+2}(W_{n})\bigr\rvert\bigr\}\\
      &\overset{(**)}{\lesssim}\max_{1\leq j\leq k-1}
        \bigl\{\widetilde{q_{n}}^{j/2}
        \bigl(\lvert T_{n}\rvert^{-j/2}+M_{n}^{j/p}\bigr)\bigr\}.
    \end{aligned}
  \end{equation*}
  where $(*)$ follows from \cref{itm:bla2}, while $(**)$ follows from
  \eqref{eq:lowercumulantsmain}.
  Thus, the following holds
  \begin{equation*}
    \widetilde{q_{n}}^{-p/2}=(\widetilde{q_{n}}^{-j_{0}/2})^{p/j_{0}}\lesssim \lvert T_{n}\rvert^{-p/2}+M_{n},
  \end{equation*}
  where $j_{0}$ is the integer satisfying that $\bigl\lvert \kappa_{j_{0}+2}(\xi^{\scalebox{0.6}{$(n)$}}) \bigr\rvert=\max\limits_{1\leq j\leq k-1}\bigl\lvert \kappa_{j+2}(\xi^{\scalebox{0.6}{$(n)$}}) \bigr\rvert$. $M_{n}$ does not depend on the value of $j_{0}$ anymore.

  In the complementary case---namely, if $p=1$, or if $p>1$ and
  $\kappa_3(W_n)=\cdots=\kappa_{k+1}(W_n)=0$---the definition gives
  $\widetilde{q}_{n}\geq \lvert T_{n}\rvert^{2(p+1)/p}$.
  Moreover, by H\"older's inequality we obtain that
  \begin{equation}\label{eq:random57}
    \sum_{i\in T_{n}}\mathbb{E} \bigl[\bigl\lvert X^{\scalebox{0.6}{$(n)$}}_{i} \bigr\rvert^{2}\bigr]\leq \lvert T_{n}\rvert^{p/(p+2)}\Bigl(\sum_{i\in T_{n}}\mathbb{E} \bigl[\bigl\lvert X^{\scalebox{0.6}{$(n)$}}_{i} \bigr\rvert^{p+2}\bigr]\Bigr)^{2/(p+2)},
  \end{equation}
  and that
  \begin{equation}\label{eq:random58}
    \Bigl(\sum_{i\in T_{n}}X^{\scalebox{0.6}{$(n)$}}_{i}\Bigr)^{2}\leq \lvert T_{n}\rvert \sum_{i\in T_{n}}\bigl\lvert X^{\scalebox{0.6}{$(n)$}}_{i} \bigr\rvert^{2}.    
  \end{equation}
  Since $\sigma_{n}^{2}=\mathbb{E} \Bigl[\Bigl(\sum_{i\in T_{n}}X^{\scalebox{0.6}{$(n)$}}_{i}\Bigr)^{2}\Bigr]$, we have
  \begin{align*}
    \widetilde{q_{n}}^{-p/2}
    \leq    & \lvert T_{n}\rvert^{-(p+1)}\sigma_{n}^{-(p+2)}\Bigl(\mathbb{E} \Bigl[\Bigl(\sum_{i\in T_{n}}X^{\scalebox{0.6}{$(n)$}}_{i}\Bigr)^{2}\Bigr]\Bigr)^{(p+2)/2} \\
    \overset{(*)}{\leq} & \sigma_{n}^{-(p+2)}\lvert T_{n}\rvert^{-p/2}\Bigl(\sum_{i\in T_{n}}\mathbb{E} \bigl[\bigl\lvert X^{\scalebox{0.6}{$(n)$}}_{i} \bigr\rvert^{2}\bigr]\Bigr)^{(p+2)/2}           \\
    \overset{(**)}{\leq} & \sigma_{n}^{-(p+2)}\sum_{i\in T_{n}}\mathbb{E} \bigl[\bigl\lvert X^{\scalebox{0.6}{$(n)$}}_{i} \bigr\rvert^{p+2}\bigr]\lesssim \lvert T_{n}\rvert^{-p/2},
  \end{align*}
  where to obtain $(*)$ we use \eqref{eq:random58} and to obtain $(**)$ we use \eqref{eq:random57}.

  Thus, using \cref{itm:uppermombound2} and the fact that $\xi_{1}^{\scalebox{0.6}{$(n)$}},\cdots,\xi_{\widetilde{q}_{n}}^{\scalebox{0.6}{$(n)$}}$ are \emph{i.i.d.}, we obtain
  \begin{equation}\label{eq:compareconnect21}
    \widetilde{q_{n}}^{-\frac{p+2}{2}}\sum_{i=1}^{\widetilde{q}_{n}}\mathbb{E} \bigl[\bigl\lvert \xi_{i}^{\scalebox{0.6}{$(n)$}} \bigr\rvert^{p+2}\bigr]\leq C_{p}'\widetilde{q_{n}}^{-p/2}\lesssim \lvert T_{n}\rvert^{-p/2}+M_{n}.
  \end{equation}
  We combine this with \eqref{eq:random56} to conclude that there is a constant $K>0$ such that it does not depend on $h$ and we have
  \begin{equation*}
    \bigl\lvert \mathbb{E} [h(W_{n})] -\mathbb{E} [h(V_{n})]\bigr\rvert\leq K(\lvert T_{n}\rvert^{-p/2}+ M_{n}).
  \end{equation*}
  Taking supremum over $h\in \Lambda_{p}$, by \cref{thm:lemzolo} we obtain that
  \begin{align*}
    &\mathcal{W}_{p}(\mathcal{L}(W_{n}),\mathcal{L}(V_{n}))\lesssim \sup_{h\in\Lambda_{p}}\bigl\lvert \mathbb{E}[h(W_{n})]-\mathbb{E} [h(V_{n})] \bigr\rvert^{1/p}\\
    \lesssim &\Bigl(\lvert T_{n}\rvert^{-p/2}+M_{n}+\widetilde{q_{n}}^{-\frac{p+2}{2}}\sum_{i=1}^{\widetilde{q}_{n}}\mathbb{E} \bigl[\bigl\lvert \xi_{i}^{\scalebox{0.6}{$(n)$}} \bigr\rvert^{p+2}\bigr]\Bigr)^{1/p}\lesssim \lvert T_{n}\rvert^{-1/2}+M_{n}^{1/p}.
  \end{align*}
  The condition $M_n\to0$ was used above to obtain the vanishing cumulants
  required by \cref{THM:EXISTENCEXI}; it also makes the final bound tend to
  zero. Combining \cref{thm:lemiidwp} and \eqref{eq:compareconnect21} we have
  \begin{equation*}
    \mathcal{W}_{p}(\mathcal{L}(V_{n}),N(0,1))\lesssim \Bigl(\widetilde{q_{n}}^{-\frac{p+2}{2}}\sum_{i=1}^{\widetilde{q}_{n}}\mathbb{E} \bigl[\bigl\lvert \xi_{i}^{\scalebox{0.6}{$(n)$}} \bigr\rvert^{p+2}\bigr]\Bigr)^{1/p}\lesssim \lvert T_{n}\rvert^{-1/2}+M_n^{1/p}.
  \end{equation*}
  By triangle inequality we conclude that
  \begin{align*}
    \mathcal{W}_{p}(\mathcal{L}(W_{n}),N(0,1))\leq & \mathcal{W}_{p}(\mathcal{L}(W_{n}),\mathcal{L}(V_{n}))+\mathcal{W}_{p}(\mathcal{L}(V_{n}),N(0,1))                                                                                      \\
    \lesssim                                         \lvert T_{n}\rvert^{-1/2}+    M_{n}^{1/p}\lesssim &\lvert T_{n}\rvert^{-1/2}+\lvert T_{n}\rvert^{-1/2}\Bigl(\sum_{\ell=1}^{\lfloor\lvert T_{n}\rvert^{1/d}\rfloor}\ell^{d( p+1)-\omega}\alpha_{\ell,n}^{(r-p-2)/r}\Bigr)^{1/p}.
  \end{align*}
  \hfill
\end{proof}

\begin{proof}[Proof of \cref{thm:finite-range-noninteger}.]
  Write $k=\lceil p\rceil$ and $\omega=p+1-k\in(0,1)$, and take $m=R$.
  At smoothness $p$, the endpoint version \eqref{eq:wfwgraph3endpoint}
  gives a one-step remainder of order
  $\lvert T_n\rvert^{-p/2}R^{d(p+1)}=
  \mathcal O(\lvert T_n\rvert^{-p/2})$.  For every
  $1\leq j\leq k-1$, the recursive smoothness order is $p-j$, and the
  available moment order satisfies
  \[
    p+2>(p-j)+2.
  \]
  Hence the ordinary, strict-moment version of \cref{thm:wfwgraph3}, again
  with $m=R$, gives a remainder
  $\mathcal O(\lvert T_n\rvert^{-(p-j)/2})$, because every mixing term in
  its sum has index at least $R+1$ and positive exponent.  Similarly,
  \cref{thm:cumuctrl}, applied at the integer order $j$ with $m=R$, gives
  \begin{equation}\label{eq:finite-range-lower-cumulants}
    \bigl\lvert\kappa_{j+2}(W_n)\bigr\rvert
    =\mathcal O(\lvert T_n\rvert^{-j/2}),
    \qquad 1\leq j\leq k-1.
  \end{equation}

  We may now repeat the exact-composition induction in the proof of
  \cref{THM:BARBOURGRAPH}.  Its terminal step uses the endpoint estimate
  above, while all lower-smoothness steps use the ordinary strict-moment
  estimate.  Each error indexed by an exact composition
  $s_1+\cdots+s_a=k$ is of order
  \begin{equation*}
    \lvert T_n\rvert^{-s_1/2}\cdots
    \lvert T_n\rvert^{-s_{a-1}/2}
    \lvert T_n\rvert^{-(s_a+\omega-1)/2}
    =\lvert T_n\rvert^{-p/2}.
  \end{equation*}
  There are only finitely many such compositions, with their number
  depending on $p$ alone.  This proves the asserted Edgeworth remainder;
  \eqref{eq:finite-range-lower-cumulants} proves the cumulant assertion.

  Finally use the cumulant-matching construction from the proof of
  \cref{thm:amixingmain}.  The cumulants in
  \eqref{eq:finite-range-lower-cumulants} tend to zero, so that construction
  supplies the comparison sum $V_n$.  If at least one matched cumulant is
  nonzero, its lower-bound property and
  \eqref{eq:finite-range-lower-cumulants} imply
  $\widetilde q_n^{-p/2}=\mathcal O(\lvert T_n\rvert^{-p/2})$; if all of
  them vanish, this follows directly from the defining lower bound on
  $\widetilde q_n$.  By cumulant matching, the complete signed
  $\Gamma(k-1)$ correction sums in the two expansions agree term by term and
  cancel before either remainder is bounded.  The endpoint expansion just
  proved and the i.i.d. expansion therefore yield
  \begin{equation*}
    \sup_{h\in\Lambda_p}
    \bigl\lvert\mathbb Eh(W_n)-\mathbb Eh(V_n)\bigr\rvert
    =\mathcal O(\lvert T_n\rvert^{-p/2}).
  \end{equation*}
  Applying \cref{thm:lemzolo} and \cref{thm:lemiidwp}, and then the
  triangle inequality, gives
  $\mathcal W_p(\mathcal L(W_n),N(0,1))
  =\mathcal O(\lvert T_n\rvert^{-1/2})$.
\end{proof}

Next we show \cref{thm:mixingconditions} through direct calculation.
\begin{proof}[Proof of \cref{thm:mixingconditions}.]
  By assumption we know $\alpha_{\ell,n}^{(r-p-2)/r}\leq C\ell^{-(u-\omega+1)}$. Thus, we have
  \begin{equation*}
    \sum_{\ell=1}^{\lfloor\lvert T_{n}\rvert^{1/d}\rfloor}\ell^{d(p+1)-\omega}\alpha_{\ell,n}^{(r-p-2)/r}\lesssim \sum_{\ell=1}^{\lfloor\lvert T_{n}\rvert^{1/d}\rfloor}\ell^{d(p+1)-\omega}\ell^{-(u-\omega+1)}=\sum_{\ell=1}^{\lfloor\lvert T_{n}\rvert^{1/d}\rfloor}\ell^{d(p+1)-u-1}.
  \end{equation*}
  If $u>d(p+1)$, then we have $d(p+1)-u-1< -1$. Thus, the sum is finite and does not depend on $n$, which implies that $M_{n}=\mathcal{O}(\lvert T_{n}\rvert^{-p/2})$. 
  If $u=d(p+1)$, then similarly we have
  \[
    M_{n}=\mathcal{O}\biggl(\lvert T_{n}\rvert^{-p/2}\sum_{\ell=1}^{\lfloor\lvert T_{n}\rvert^{1/d}\rfloor}\ell^{-1}\biggr)=\mathcal{O}\bigl(\lvert T_{n}\rvert^{-p/2}\log \lvert T_{n}\rvert\bigr).
  \]
  Lastly if $d(p/2+1)<u<d(p+1)$, we derive that
  \begin{align*}
    M_{n}
    = & \mathcal{O}\biggl(\lvert T_{n}\rvert^{-p/2}\sum_{\ell=1}^{\lfloor\lvert T_{n}\rvert^{1/d}\rfloor}\ell^{d(p+1)-u-1}\biggr) \\
    = & \mathcal{O}\bigl(\lvert T_{n}\rvert^{-p/2}\lfloor\lvert T_{n}\rvert^{1/d}\rfloor^{d(p+1)-u}\bigr)=\mathcal{O}\bigl(\lvert T_{n}\rvert^{-p/2-u/d+(p+1)}\bigr),
  \end{align*}
  The three Wasserstein rates in \eqref{eq:exactWppolyrate} now follow by
  substituting these estimates into \eqref{eq:amixingmain1}. At the critical
  value, substituting the middle estimate into \eqref{eq:hugeexpan3} also
  gives the stated Edgeworth remainder. This concludes the proof.\hfill
\end{proof}

In order to prove \cref{thm:amixingmain2}, we need the following result, which is a refined version of \cref{THM:BARBOURGRAPH} of \cref{thm:hugehugethm} for the case where $p$ is an integer.
\begin{theorem}\label{THM:BARBOURGRAPH2}
Let $\bigl(X^{\scalebox{0.6}{$(n)$}}_{i}\bigr)_{i\in T_n}$ be a sequence of mean-zero random fields with $\alpha$-mixing coefficients $(\alpha_{\ell,n} )_{\ell\geq 1}$. For $p\in\mathbb{N}_{+}$, suppose that there exists $r\in[p+2,\infty]$ such that \eqref{rr} and \eqref{moment} hold, and that 
 \begin{equation}\label{eq:mixinginthmB}
  \sup_{n}\sum_{\ell=1}^{\infty}\ell^{d-1}\alpha_{\ell,n}^{(r-p-1)/r}<\infty, \quad M_{n,m,\delta}^{(i)}\to 0,\quad  M_{n,m}^{(ii)}\to 0\quad\textup{ as }n\to\infty
 \end{equation}
 for some positive integer sequence $m=o(\lvert T_{n} \rvert^{\frac{1}{4d}})$ and fixed $\delta\in [0,1]\cap[0,r-p-1)$,
 where $M_{n,m,\delta}^{(i)}$ and $M_{n,m}^{(ii)}$ are as defined in \eqref{eq:mnmidef} and \eqref{eq:mnmiidef}. For $1\leq j\leq p-1$, define
 \begin{equation*}
   \mathcal R_{j,n}:=\min_{1\leq b\leq m}\lvert T_n\rvert^{-j/2}
   \left(b^{d(j+1)}+
   \sum_{\ell=b+1}^{b+1+\lfloor\lvert T_n\rvert^{1/d}/2\rfloor}
   \ell^{d(j+1)-1}\alpha_{\ell,n}^{(r-j-2)/r}\right).
 \end{equation*}
 Then we have
  \begin{enumerate}[label=(\alph*)]
    \item For any $j=1,\cdots,p-1$, we have 
    \begin{equation}
    \bigl\lvert\kappa_{j+2}(W_{n})\bigr\rvert=\mathcal O(\mathcal R_{j,n}).
    \end{equation}
    \item If $p\geq 2$, there exists a quantity
    $\widetilde{\kappa}_{p+1,n}$, depending on $p$ and the joint distribution
    of $(X_i)_{i\in T_n}$ and on the terminal radius $m$, such that
    \begin{equation}
    \bigl\lvert\widetilde{\kappa}_{p+1,n}\bigr\rvert^{p/(p-1)}
    =\mathcal O\bigl(\mathcal R_{p-1,n}^{p/(p-1)}\bigr)
      +\mathcal O\bigl((M_{n,m}^{(ii)})^{p/(p-1)}\bigr),
    \end{equation}
    and, uniformly over $h \in \Lambda_{p}$, the following expansion holds:
  \begin{equation}\label{eq:hugeexpan4}
    \begin{aligned}
      \mathbb{E} [h(W_{n})]-\mathcal{N}h
      ={}&\sum_{(r,s_{1:r})\in \Gamma(p -1)}(-1)^{r}\\
      &\quad{}\times\prod_{j=1}^{r}
        \frac{\kappa _{s _{j}+2}(W_{n})}{(s _{j}+1)!}\\
      &\quad{}\times\mathcal{N}\ \Bigl[\bigcomp_{j=1}^{r}
        (\partial ^{s _{j}+1}\Theta)\ h\Bigr]\\
      &\quad{}+\frac{\kappa_{p +1}(W_{n})-\widetilde{\kappa}_{p+1,n}}{p !}
        \mathcal{N}\ [\partial ^{p}\Theta\ h]\\
      &\quad{}+\mathcal O\Biggl(
        \sum_{j=1}^{p-1}\mathcal R_{j,n}^{p/j}
        +\lvert T_n\rvert^{-p/2}m^{d(p+1)}\\
      &\hspace{37mm}{}+M_{n,m,\delta}^{(i)}+M_{n,m}^{(ii)}\Biggr).
    \end{aligned}
  \end{equation}
  \end{enumerate}
\end{theorem}

  \cref{THM:BARBOURGRAPH2} is different from \cref{thm:hugehugethm} in the following ways:
  \begin{enumerate}[label=(\alph*)]
    \item $p$ is required to be an integer due to the proof technique we use.
    \item The remainder is controlled using $M_{n,m,\delta}^{(i)}$ and $M_{n,m}^{(ii)}$ instead of $M_{n}$, which will lead to different convergence rates in \cref{thm:hugehugethm,thm:amixingmain2}.
    \item $\kappa_{p+1}(W_{n})$ is replaced by $\widetilde{\kappa}_{p+1,n}$.
  \end{enumerate}
  In general $M_{n}$ does not dominate $M_{n,m,\delta}^{(i)}$ or $M_{n,m}^{(ii)}$, and vice versa. They lead to different conditions and convergence rates for the Wasserstein-$p$ distance in \cref{thm:hugehugethm,thm:amixingmain2}.  

\begin{proof}[Proof of \cref{thm:amixingmain2}.]
  The case $p=1$ follows directly from \eqref{eq:barbourgraph2base} and
  $\mathcal{W}_{1}=\mathcal{Z}_{1}$, with the sum in
  \eqref{eq:amixingmain2exact} empty. Hence below we assume $p\geq2$ and put
  \begin{equation*}
    E_n:=\sum_{j=1}^{p-1}\mathcal R_{j,n}^{p/j}
      +\lvert T_n\rvert^{-p/2}m^{d(p+1)}
      +M_{n,m,\delta}^{(i)}+M_{n,m}^{(ii)}.
  \end{equation*}
  Taking $b=m$ in \eqref{eq:lowerorderradius} and applying H\"older's
  inequality between the terminal mixing tail and the summable
  $\ell^{d-1}$ tail gives
  \begin{equation}\label{eq:lower-radius-coarsening}
    \mathcal R_{j,n}^{p/j}\lesssim
      \lvert T_n\rvert^{-p/2}m^{2dp}+M_{n,m}^{(ii)},
      \qquad 1\leq j\leq p-1.
  \end{equation}
  Thus $E_n\to0$ under \eqref{eq:mixinginthm3.5}.

  Define the target cumulants
  \begin{equation*}
    a_{j,n}:=
    \begin{cases}
      \kappa_{j+2}(W_n),&j<p-1,\\
      \widetilde\kappa_{p+1,n},&j=p-1,
    \end{cases}
    \qquad
    A_n:=\max_{1\leq j\leq p-1}\lvert a_{j,n}\rvert.
  \end{equation*}
  By \cref{THM:BARBOURGRAPH2},
  \begin{equation}\label{eq:target-cumulants-exact}
    \lvert a_{j,n}\rvert^{p/j}\lesssim E_n,
    \qquad 1\leq j\leq p-1,
  \end{equation}
  where we also used $(M_{n,m}^{(ii)})^{p/(p-1)}\lesssim
  M_{n,m}^{(ii)}$ for large $n$. In particular all the target cumulants
  converge to zero.

  Apply \cref{THM:EXISTENCEXI} to $(a_{j,n})_{j=1}^{p-1}$. Let $q_n$ and
  $\xi^{\scalebox{0.6}{$(n)$}}$ be the resulting integer and random variable,
  and let
  \begin{equation*}
    \widetilde q_n:=
    \begin{cases}
      \lceil\lvert T_n\rvert^{2(p+1)/p}\rceil\vee q_n,&A_n=0,\\
      q_n,&A_n>0.
    \end{cases}
  \end{equation*}
  For \emph{i.i.d.} copies $(\xi_i^{\scalebox{0.6}{$(n)$}})$, set
  $V_n:=\widetilde q_n^{-1/2}\sum_{i=1}^{\widetilde q_n}
  \xi_i^{\scalebox{0.6}{$(n)$}}$. Then
  $\kappa_{j+2}(V_n)=a_{j,n}$ for $1\leq j\leq p-1$.

  We next bound the \emph{i.i.d.} remainder without losing the refined
  lower-order scales. If $A_n>0$, the lower bound in
  \cref{THM:EXISTENCEXI} and \eqref{eq:target-cumulants-exact} give, for
  some $1\leq j_0\leq p-1$,
  \begin{equation*}
    1\lesssim q_n^{j_0/2}\lvert a_{j_0,n}\rvert,
    \qquad\text{and hence}\qquad
    \widetilde q_n^{-p/2}\lesssim E_n.
  \end{equation*}
  If $A_n=0$, the definition of $\widetilde q_n$ gives instead
  \begin{equation*}
    \widetilde q_n^{-p/2}\leq\lvert T_n\rvert^{-(p+1)}
      \leq\lvert T_n\rvert^{-p/2}m^{d(p+1)}\leq E_n.
  \end{equation*}
  In \eqref{eq:hugeexpan4}, the correction term changes the terminal
  singleton from $-\kappa_{p+1}(W_n)/p!$ to
  $-\widetilde\kappa_{p+1,n}/p!$. Thus the complete signed cumulant sums in
  the two Edgeworth expansions agree exactly.
  Therefore, \cref{thm:barbour,THM:BARBOURGRAPH2} and the uniform moment
  bound in \cref{THM:EXISTENCEXI} yield, uniformly over $h\in\Lambda_p$,
  \begin{equation*}
    \bigl\lvert\mathbb E[h(W_n)]-\mathbb E[h(V_n)]\bigr\rvert
    \lesssim E_n+\widetilde q_n^{-p/2}\lesssim E_n.
  \end{equation*}
  It follows from \cref{thm:lemzolo} that
  $\mathcal W_p(\mathcal L(W_n),\mathcal L(V_n))\lesssim E_n^{1/p}$.
  Moreover, \cref{thm:lemiidwp} gives
  \begin{equation*}
    \mathcal W_p(\mathcal L(V_n),N(0,1))
    \lesssim\widetilde q_n^{-1/2}\lesssim E_n^{1/p}.
  \end{equation*}
  The triangle inequality proves \eqref{eq:amixingmain2exact}.

  Finally, \eqref{eq:lower-radius-coarsening} and
  $m^{d(p+1)}\leq m^{2dp}$ reduce the exact bound to
  \eqref{eq:amixingmain2}. This proves the theorem.\hfill
\end{proof}

  As an application of \cref{thm:amixingmain2}, we show \cref{thm:w1result}.

  \begin{corollary}\label{thm:w1result}
    For $p=1$, the Wasserstein distance
    \[
      \mathcal{W}_{1}(\mathcal{L}(W_n),N(0,1))
    \]
    converges to $0$ as long as the moment and non-degeneracy conditions in \cref{thm:amixingmain2} hold with $r\geq 3$ and $\alpha_{\ell,n}$ satisfies that
  \begin{align*}
    &\sup_n\sum_{\ell=1}^{\infty}\ell^{d-1}\alpha_{\ell,n}^{(r-2)/r}<\infty;\\
    &\sup_n\sum_{\ell=m}^{\infty}\ell^{d-1}\alpha_{\ell,n}^{(r-2-\epsilon)/r}\to 0\quad\textup{as }m\to \infty \quad\textup{for some }\epsilon>0.
  \end{align*}
  \end{corollary}

\begin{proof}[Proof of \cref{thm:w1result}.]
  Set $\epsilon_0:=\frac12\min\{\epsilon,1,r-2\}$ and
  $m:=\lceil\lvert T_n\rvert^{(\epsilon_0\wedge(1/3))/(2d)}\rceil$.
  We apply \cref{thm:amixingmain2} with $p=1$ and
  $\delta=\epsilon_0$. Then $m\to\infty$ as $n\to\infty$. Since
  $$\sum_{\ell=1}^{\infty}\ell^{d-1}\alpha_{\ell,n}^{(r-2)/r}\leq \sum_{\ell=1}^{\infty}\ell^{d-1}\alpha_{\ell,n}^{(r-2-\epsilon)/r},$$
  the first part of the  condition \eqref{eq:mixinginthm3.5} is satisfied. Now we check that
  \begin{align*}
    & \lvert T_{n}\rvert^{-1/2}m^{2d}\lesssim \lvert T_{n}\rvert^{-1/2+(\epsilon_0\wedge 1/3)}\xrightarrow{n\to\infty} 0,\\
    &\lvert T_{n}\rvert^{-\epsilon_0/2}m^{d}\sum_{\ell=m+1}^{m+1+\lfloor\frac{\lvert T_{n}\rvert^{1/d}}{2}\rfloor}\ell^{d\epsilon_0-\epsilon_0}\alpha_{\ell,n}^{(r-2-\epsilon_0)/r}\lesssim \sum_{\ell=m+1}^{\infty}\ell^{d-1}\alpha_{\ell,n}^{(r-2-\epsilon)/r},\\
    &\sum_{\ell=m+1}^{m+1+\lfloor\frac{\lvert T_{n}\rvert^{1/d}}{2}\rfloor}\ell^{d-1}\alpha_{\ell,n}^{(r-2)/r}\leq\sum_{\ell=m+1}^{\infty}\ell^{d-1}\alpha_{\ell,n}^{(r-2-\epsilon)/r}.
  \end{align*}
Since $m\to\infty$ as $n\to\infty$, we have that by assumption $\sum_{\ell=m+1}^{\infty}\ell^{d-1}\alpha_{\ell,n}^{(r-2-\epsilon)/r}$ converges to zero. Thus, \eqref{eq:mixinginthm3.5} of \cref{thm:amixingmain2} is fully satisfied and the result follows. \hfill
\end{proof}

Lastly, we prove \cref{thm:alphapoly} by applying \cref{thm:amixingmain2} to random fields with $\alpha$-mixing coefficients that converge at a polynomial rate, and combining the results of \cref{thm:mixingconditions}.
\begin{proof}[Proof of \cref{thm:alphapoly}.]
  Write $N:=\lvert T_n\rvert$. We repeatedly use the elementary one-sided
  estimates
  \begin{equation}\label{eq:power-tail-estimates}
    \sum_{\ell=m+1}^{m+1+\lfloor N^{1/d}/2\rfloor}\ell^{a-1}
    \lesssim
    \begin{cases}
      N^{a/d},&a>0,\\
      \log N,&a=0,\\
      m^a,&a<0.
    \end{cases}
  \end{equation}
  The first mixing sum in \eqref{eq:mixinginthm3.5} is uniformly bounded
  whenever $A>d$, since it is bounded by a constant multiple of
  $\sum_{\ell\geq1}\ell^{d-1-A}$.

  We first prove part~(a). Thus $p\geq2$ and $A>d(p+1)/2>d$.
  If $d(p+1)/2<A<dp$, take $m=1$ and $\delta=0$. Then
  \begin{equation*}
    M_{n,1,0}^{(i)}\lesssim N^{-(p-1)/2},\qquad
    M_{n,1}^{(ii)}\lesssim
      N^{-(p-1)/2+(dp-A)/d}.
  \end{equation*}
  The second term is the slowest term in \eqref{eq:amixingmain2}, and hence
  \begin{equation*}
    \mathcal W_p(\mathcal L(W_n),N(0,1))
    \lesssim N^{-\{A/(dp)-(p+1)/(2p)\}}.
  \end{equation*}
  At $A=dp$, the same choices and the critical case of
  \eqref{eq:power-tail-estimates} give instead
  \begin{equation*}
    \mathcal W_p(\mathcal L(W_n),N(0,1))
    \lesssim N^{-(p-1)/(2p)}(\log N)^{1/p}.
  \end{equation*}

  Now suppose $A>dp$, set $x:=1/(2(A+dp))$, and take
  $m=\lceil N^x\rceil$. Then
  \begin{equation}\label{eq:A-terminal-balance}
    N^{-1/2}m^{2d}\lesssim N^{-1/2+d/(A+dp)},\qquad
    (M_{n,m}^{(ii)})^{1/p}\lesssim N^{-1/2+d/(A+dp)}.
  \end{equation}
  It remains to control $M_{n,m,\delta}^{(i)}$. If $u>d$, take $\delta=1$;
  then
  \begin{equation*}
    (M_{n,m,1}^{(i)})^{1/p}\lesssim
      N^{-1/2}m^{d+(d-u)/p}\leq N^{-1/2}m^{2d}.
  \end{equation*}
  If $u=d$, the same choice gives
  $(M_{n,m,1}^{(i)})^{1/p}\lesssim
  N^{-1/2}m^d(\log N)^{1/p}=o(N^{-1/2}m^{2d})$.
  Finally, if $0\leq u<d$, take
  \begin{equation*}
    \delta:=\frac{A-1}{A-u+d-1}\in[0,1).
  \end{equation*}
  The power of $\ell$ in $M_{n,m,\delta}^{(i)}$ is then $-1$, so
  \begin{equation*}
    (M_{n,m,\delta}^{(i)})^{1/p}\lesssim
      N^{-1/2+(1-\delta)/(2p)}m^d(\log N)^{1/p}.
  \end{equation*}
  This is $o(N^{-1/2}m^{2d})$: indeed,
  \begin{equation*}
    \frac{(d-u)(A+dp)}{A-u+d-1}<dp,
  \end{equation*}
  because the difference after clearing the positive denominator is
  $A\{d(p-1)+u\}-dp>0$. Combining this with
  \eqref{eq:A-terminal-balance} proves the third case in
  \eqref{eq:refined-rate-envelope}. These choices satisfy the convergence
  conditions of \cref{thm:amixingmain2}; in particular
  $m=o(N^{1/(4d)})$. This proves part~(a), including $r=p+2$, for which
  $u=0$ and the last choice of $\delta$ applies.

  We next prove part~(b). If $p\geq2$ and $A\leq dp$, the assertion for
  $\mathfrak S_{p,A}$ is already part~(a). Suppose, therefore, that $A>dp$;
  this is automatic when $p=1$, since then $A>u>d$. For
  $1\leq j\leq p-1$, put $A_j:=(r-j-2)v/r$, with $A_j=v$ when
  $r=\infty$. Since
  \begin{equation*}
    A_j\geq A>dp\geq d(j+1),
  \end{equation*}
  choosing $b=1$ in \eqref{eq:lowerorderradius} gives
  \begin{equation}\label{eq:lower-radii-iid-order}
    \mathcal R_{j,n}\lesssim N^{-j/2},
    \qquad 1\leq j\leq p-1.
  \end{equation}
  (There is no such term when $p=1$.) Take $\delta=1$ and
  $m=\lceil N^{1/(2(A+d))}\rceil$. The assumption
  $r>p+2$ makes $\delta=1$ admissible, and $A>d$ gives
  $m=o(N^{1/(4d)})$. The assumption $u>d(p+1)/2$ implies $u>d$, and
  \eqref{eq:power-tail-estimates} yields
  \begin{align*}
    \left(N^{-p/2}m^{d(p+1)}\right)^{1/p}
      &\lesssim N^{-1/2+d(p+1)/(2p(A+d))},\\
    (M_{n,m}^{(ii)})^{1/p}
      &\lesssim N^{-1/2+d(p+1)/(2p(A+d))},\\
    (M_{n,m,1}^{(i)})^{1/p}
      &\lesssim N^{-1/2}m^{d+(d-u)/p}
       \leq N^{-1/2}m^{d(p+1)/p}.
  \end{align*}
  Together with \eqref{eq:lower-radii-iid-order}, the exact bound
  \eqref{eq:amixingmain2exact} proves
  \eqref{eq:separate-radius-poly-rate}. Notice that the lower cumulants use
  radius $b=1$, while only the terminal remainder uses the growing radius
  $m$; no graph estimate is altered.

  Finally, if $u>d(p+2)/2$, apply \cref{thm:mixingconditions} with
  $\omega=1$. Its three exact alternatives in \eqref{eq:exactWppolyrate}
  give the bound $\mathfrak G_{p,u}$. Taking the better of that bound and
  \eqref{eq:separate-radius-poly-rate} proves
  \eqref{eq:combined-poly-rate}, and completes the proof.\hfill
\end{proof}

\begin{proof}[Proof of \cref{thm:finite-range-integer-rate}.]
  Take $m=R$ and $\delta=0$ in \cref{thm:amixingmain2}. Both tail terms
  $M_{n,R,0}^{(i)}$ and $M_{n,R}^{(ii)}$ vanish, while the first mixing sum
  in \eqref{eq:mixinginthm3.5} has only finitely many nonzero terms. Hence
  \eqref{eq:amixingmain2} gives
  \begin{equation*}
    \mathcal W_p(\mathcal L(W_n),N(0,1))
    \lesssim N^{-1/2}R^{2d}=\mathcal O(N^{-1/2}),
  \end{equation*}
  since $R$ is fixed.\hfill
\end{proof}
\section{Proofs of cumulant-based Edgeworth expansions}\label{sec:finalpflemma}

We use $W$, $\sigma$, $T$, $\alpha_{\ell}$ and $M$ instead of $W_{n}$, $\sigma_{n}$, $T_{n}$, $\alpha_{\ell,n}$ and $M_{n}$ to simplify notation. First we prove a more general version of the cumulant bound \cref{thm:cumuctrl1}.
\begin{proposition}\label{thm:cumuctrl}
Let $\bigl(X^{\scalebox{0.6}{$(n)$}}_{i}\bigr)$ be a sequence of mean-zero random fields with $\alpha$-mixing coefficients $(\alpha_{\ell} )_{\ell\geq 1}$. For $p\in\mathbb{N}_{+}$, suppose that there exists $r\in(p+2,\infty]$ such that \eqref{rr} and \eqref{moment} hold. Then, for every choice of positive integers $m=m_n$ and all sufficiently large $n$,
  \begin{equation*}
    \bigl\lvert\kappa_{p+2}(W)\bigr\rvert
    \lesssim \lvert T\rvert^{-p/2}\biggl(m^{d(p+1)}
    +\sum_{\ell=m+1}^{m+1+\lfloor\frac{\lvert T\rvert^{1/d}}{2}\rfloor}
      \ell^{d(p+1)-1}\alpha_{\ell}^{(r-p-2)/r}\biggr).
  \end{equation*}
The implicit constant is uniform in $n$ and $m$. In particular, if, uniformly
in $n$, $\alpha_{\ell}^{(r-p-2)/r}\lesssim\ell^{-q}$ for some $q>0$, then
taking $m=1$ gives
\begin{equation*}
  \bigl\lvert\kappa_{p+2}(W)\bigr\rvert=
  \begin{cases}
    \mathcal O(\lvert T\rvert^{-p/2}),&q>d(p+1),\\
    \mathcal O(\lvert T\rvert^{-p/2}\log\lvert T\rvert),&q=d(p+1),\\
    \mathcal O(\lvert T\rvert^{p/2+1-q/d}),&q<d(p+1).
  \end{cases}
\end{equation*}
\end{proposition}

\begin{proof}[Proof of \cref{thm:cumuctrl}.]
  Apply \cref{thm:wfwgraph3} with its parameters $k=p$, $\omega=1$, and
  $f(x)=x^{p+1}/(p+1)!$. This function belongs to
  $\mathcal C^{p,1}(\mathbb R)$ and has H\"older coefficient one, so
  \begin{align*}
    \mathbb{E} [Wf(W)]
    ={}& \sum_{j=1}^{p}\frac{\kappa_{j+1}(W)}{j !}
      \mathbb{E}[\partial^{j}f(W)]\\
    &\quad{}+\mathcal{O}\biggl(\lvert T \rvert^{-p/2}
      \Bigl(m^{d(p+1)}+\!\!\sum_{\ell=m+1}^{m+1+\lfloor\frac{\lvert T\rvert^{1/d}}{2}\rfloor}
      \ell^{d(p+1) -1}\alpha_{\ell}^{(r-p-2)/r}\Bigr)\biggr).
  \end{align*}
  On the other hand, by \cref{thm:identitycumulant} we have
  \begin{align*}
    \mathbb{E} [Wf(W)]
    ={}&\frac{1}{(p+1)!}\mu_{p+2}(W)
    =\frac{1}{(p+1)!}\sum_{j=1}^{p+1}\binom{p+1}{j}\kappa_{j+1}(W)\mu_{p+1-j}(W)\\
    = &\sum_{j=1}^{p}\frac{\kappa_{j+1}(W)}{j !}\mathbb{E} [\partial^{j} f(W)]+\frac{\kappa_{p+2}(W)}{(p+1)!}.
  \end{align*}
  Thus, we conclude that
  \begin{equation*}
    \bigl\lvert\kappa_{p+2}(W)\bigr\rvert \lesssim \lvert T \rvert^{-p/2}\Bigl(m^{d(p+1)}+\sum_{\ell=m+1}^{m+1+\lfloor\frac{\lvert T\rvert^{1/d}}{2}\rfloor}\ell^{d(p+1)-1 }\alpha_{\ell}^{(r-p-2)/r}\Bigr).
  \end{equation*}\hfill
\end{proof}

Before showing \cref{THM:BARBOURGRAPH} of \cref{thm:hugehugethm}, we need the following lemma:
\begin{lemma}[Consequence of Lemmas 5 and 6 of \cite{barbour1986asymptotic}]\label{thm:lemsteinsol}
  For any nonnegative integer $k$, $0<\omega\leq 1$ and
  $h \in \mathcal{C}^{k,\omega}(\mathbb{R})$, the solution to
  \eqref{eq:stein}, i.e., $f=\Theta h$, satisfies
  $f\in \mathcal{C}^{k,\omega}(\mathbb{R})\cap
  \mathcal{C}^{k+1,\omega }(\mathbb{R})$ and the H\"older coefficients
  (see \cref{thm:defholder}) $\lvert f \rvert_{k,\omega}$ and
  $\lvert f \rvert_{k+1,\omega }$ are bounded by a constant that depends
  only on $k,\omega$ and $\lvert h\rvert_{k,\omega}$.
\end{lemma}

We shall also use the following fractional base estimate. It is stated for
$0<q\leq 1$ because such exponents arise after differentiating a test
function of noninteger smoothness.
\begin{lemma}\label{thm:fractionalsteinbase}
  Let $0<q\leq 1$, and suppose that the assumptions of
  \cref{thm:wfwgraph3} hold with $p=q$, so that $k=1$ and $\omega=q$.
  Then, uniformly over $h\in\mathcal C^{0,q}(\mathbb R)$ with
  $\lvert h\rvert_{0,q}\leq1$,
  \begin{equation*}
    \bigl\lvert\mathbb E[h(W)]-\mathcal Nh\bigr\rvert
    \lesssim
    \lvert T\rvert^{-q/2}
    \Bigl(1+\sum_{\ell=1}^{\lfloor\lvert T\rvert^{1/d}\rfloor}
      \ell^{d(1+q)-q}\alpha_{\ell}^{(r-q-2)/r}\Bigr).
  \end{equation*}
\end{lemma}

\begin{proof}
  Let $f=\Theta h$. By \cref{thm:lemsteinsol} with $k=0$, we have
  $f\in\mathcal C^{0,q}(\mathbb R)\cap\mathcal C^{1,q}(\mathbb R)$,
  with $\lvert f\rvert_{1,q}$ bounded uniformly over the stated class of
  test functions. Applying \cref{thm:wfwgraph3} with $k=1$, $\omega=q$
  and $m=1$, and using $\kappa_2(W)=1$, gives the claimed range of
  summation once $\lfloor\lvert T\rvert^{1/d}\rfloor\geq3$; indeed,
  $2+\lfloor\lvert T\rvert^{1/d}/2\rfloor
  \leq\lfloor\lvert T\rvert^{1/d}\rfloor$. In the remaining cases where
  $\lfloor\lvert T\rvert^{1/d}\rfloor<3$, the finite sums in
  \cref{thm:wfwgraph3} are uniformly bounded and are absorbed into the
  implicit constant. Therefore,
  \begin{align*}
    \mathbb E[h(W)]-\mathcal Nh
    &=\mathbb E[f'(W)]-\mathbb E[Wf(W)]\\
    &=\mathcal O\!\left(
      \lvert T\rvert^{-q/2}
      \Bigl(1+\sum_{\ell=1}^{\lfloor\lvert T\rvert^{1/d}\rfloor}
        \ell^{d(1+q)-q}\alpha_{\ell}^{(r-q-2)/r}\Bigr)\right),
  \end{align*}
  as claimed.
\end{proof}

\begin{proof}[Proof of \cref{THM:BARBOURGRAPH} of \cref{thm:hugehugethm}.]
  Let $k:=\lceil p\rceil$. For convenience, denote
  \begin{equation*}
    L:=\sup_{n}\sum_{\ell=1}^{\lfloor\lvert T_n\rvert^{1/d}\rfloor}
      \ell^{d-1}\alpha_{\ell,n}^{(r-p-2)/r},
  \end{equation*}
  and set
  \begin{equation*}
    \widehat{R}_{j,\eta}:=\lvert T \rvert^{-(j+\eta-1)/2}\Bigl(1+\!\!\sum_{\ell=1}^{\lfloor\lvert T\rvert^{1/d}\rfloor}\ell^{d(j+\eta)-\eta}\alpha_{\ell}^{(r-j-1-\eta)/r}\Bigr)
  \end{equation*}    
  only for the pairs $(j,\eta)$ such that either
  $1\leq j\leq k-1$ and $\eta=1$, or $1\leq j\leq k$ and
  $\eta=\omega$.
  Then we have that $\widehat{R}_{k,\omega}=\lvert T \rvert^{-p/2}+M$. 

  For a positive integer $a$, define the exact-composition set
  \begin{equation*}
    \Gamma_{=}(a):=\bigl\{(r,s_{1:r}):r,s_1,\ldots,s_r\in\mathbb N_+,
      \ s_1+\cdots+s_r=a\bigr\}.
  \end{equation*}
  We prove the following auxiliary assertion for every $p>0$, by induction
  on $k=\lceil p\rceil$:
  \begin{equation}
  \begin{aligned}
    \mathbb{E} [h(W)]-\mathcal{N}h
    =&\sum_{(r,s_{1:r})\in\Gamma(k-1)}(-1)^{r}\prod_{j=1}^{r}\frac{\kappa _{s _{j}+2}(W)}{(s _{j}+1)!}\mathcal{N}\ \Bigl[\bigcomp_{j=1}^{r}(\partial ^{s _{j}+1}\Theta)\ h\Bigr]\\*
    &\ +\mathcal{O}\biggl(\sum_{(r,s_{1:r})\in\Gamma_{=}(k)}
      \widehat{R}_{s_{1},1}\cdots
      \widehat{R}_{s_{r-1},1}\widehat{R}_{s_{r},\omega}\biggr),\label{eq:hugeexpan0}
  \end{aligned}
\end{equation}
  where $\Gamma(k-1)=\{ r,s_{1:r}\in\mathbb{N}_{+}:s_{1}+\cdots+s_{r}\leq k -1\}$,
  with the convention that $\Gamma(0)$ is empty.

  If $k=1$, then $0<p=\omega\leq1$, and \cref{thm:fractionalsteinbase}
  gives \eqref{eq:hugeexpan0}, since $\Gamma_{=}(1)$ consists only of
  $(1,1)$. Suppose now that the assertion holds for all smaller ceiling
  orders, and consider $k\geq2$ and $p=k+\omega-1$. By
  \cref{thm:lemsteinsol}, $f=\Theta h\in\mathcal C^{k-1,\omega}(\mathbb R)
  \cap\mathcal C^{k,\omega}(\mathbb R)$, and $\lvert f\rvert_{k,\omega}$
  is bounded by a constant depending only on $p$. Thus, by
  \cref{thm:wfwgraph3} with $m=1$, we have
  \begin{align*}
    &\mathbb{E} [h(W)]-\mathcal{N}h=  \mathbb{E} [ f'(W)]-\mathbb{E} [W f(W)]\\
    =& -\sum_{j=2}^{k}\frac{\kappa_{j+1}(W)}{j !}\mathbb{E} [\partial^{j} f(W)] + \mathcal{O}\biggl(\lvert T \rvert^{-p/2}\Bigl(1+\sum_{\ell=1}^{\lfloor\lvert T\rvert^{1/d}\rfloor}\ell^{d(p+1) -\omega}\alpha_{\ell}^{(r-p-2)/r}\Bigr)\biggr)\\
    =& -\sum_{j=1}^{k-1}\frac{\kappa_{j+2}(W)}{(j+1) !}\mathbb{E} [\partial^{j+1} \Theta h(W)] + \mathcal{O}(\widehat{R}_{k,\omega}).
  \end{align*}

  For $1\leq j\leq k-1$, the function
  $\partial^{j+1}\Theta h$ belongs to
  $\mathcal C^{k-j-1,\omega}(\mathbb R)$ and has a uniformly bounded
  H\"older coefficient. Its smoothness order is $p-j$, whose ceiling is
  $k-j$. The inductive hypothesis therefore gives
  \begin{align*}
    &\mathbb{E} [\partial^{j+1}\Theta h(W)]-\mathcal{N}[\partial^{j+1}\Theta h]\\
    =  & \sum_{(r,s_{1:r})\in\Gamma(k-j-1)}(-1)^{r}\prod_{\ell=1}^{r}\frac{\kappa _{s _{\ell}+2}(W)}{(s _{\ell}+1)!}\mathcal{N}\ \Bigl[\bigcomp_{\ell=1}^{r}(\partial ^{s _{\ell}+1}\Theta)\comp \partial^{j+1}\Theta \ h\Bigr] \\
     &\ +\mathcal{O}\biggl(\sum_{(r,s_{1:r})\in\Gamma_{=}(k-j)}
       \widehat{R}_{s_{1},1}\cdots
       \widehat{R}_{s_{r-1},1}\widehat{R}_{s_{r},\omega}\biggr).
  \end{align*}

For $1\leq j\leq k-1$, \cref{thm:cumuctrl}, applied at order $j$ with
$m=1$, gives
\begin{equation*}
  \bigl\lvert\kappa_{j+2}(W)\bigr\rvert
  \lesssim \widehat R_{j,1}.
\end{equation*}
Indeed, its terminal sum is contained in the sum defining
$\widehat R_{j,1}$ once $\lfloor\lvert T\rvert^{1/d}\rfloor\geq3$, and the
remaining finitely many cases are absorbed into the implicit constant.
H\"older interpolation between the weights
$d(p+1)-1$ and $d-1$, together with
$d(p+1)-1\leq d(p+1)-\omega$, gives
\begin{equation*}
  \widehat{R}_{j,1}
  \lesssim
  \bigl(\lvert T\rvert^{-j/2}+M^{j/p}\bigr)(1+L)^{1-j/p}.
\end{equation*}
Since $L<\infty$ by \eqref{mixing}, this proves
\eqref{eq:lowercumulantsmain}, without requiring $M=M_n$ to converge to zero.
As usual, a product over an empty index set equals one.
Thus, we have
  \begin{align*}
    \mathbb{E} [h(W)]-\mathcal{N}h
    ={}&-\sum_{j=1}^{k-1}\frac{\kappa_{j+2}(W)}{(j+1)!}
      \mathbb{E} [\partial^{j+1}\Theta h(W)]
      +\mathcal{O} (\widehat{R}_{k,\omega})\\
    ={}&-\sum_{j=1}^{k-1}\frac{\kappa_{j+2}(W)}{(j+1)!}
      \mathcal{N} [\partial^{j+1}\Theta h]\\
    &{}-\sum_{j=1}^{k-1}\frac{\kappa_{j+2}(W)}{(j+1)!}
      \sum_{(r,s_{1:r})\in\Gamma(k-j-1)}(-1)^{r}
      \prod_{\ell=1}^{r}\frac{\kappa _{s _{\ell}+2}(W)}{(s _{\ell}+1)!}\\
    &\qquad{}\times\mathcal{N}\ \Bigl[
      \bigcomp_{\ell=1}^{r}(\partial ^{s _{\ell}+1}\Theta)
      \comp \partial^{j+1}\Theta \ h\Bigr] \\
    &{}+\mathcal{O}\Biggl(
      \begin{aligned}
      \widehat{R}_{k,\omega}
      &+\sum_{j=1}^{k-1}\widehat{R}_{j,1}
       \sum_{(r,s_{1:r})\in\Gamma_{=}(k-j)}\\[-1mm]
      &\qquad{}\times
       \left(\prod_{\ell=1}^{r-1}\widehat{R}_{s_{\ell},1}\right)
       \widehat{R}_{s_{r},\omega}
      \end{aligned}
      \Biggr)\\
    ={}&\sum_{(r,s_{1:r})\in\Gamma(k-1)}(-1)^{r}
      \prod_{\ell=1}^{r}\frac{\kappa _{s _{\ell}+2}(W)}{(s _{\ell}+1)!}\\
    &\qquad{}\times\mathcal{N}\ \Bigl[
      \bigcomp_{\ell=1}^{r}(\partial ^{s _{\ell}+1}\Theta)\ h\Bigr] \\
    &{}+\mathcal{O}\Biggl(\sum_{(r,s_{1:r})\in\Gamma_{=}(k)}
      \left(\prod_{\ell=1}^{r-1}\widehat{R}_{s_{\ell},1}\right)
      \widehat{R}_{s_{r},\omega}\Biggr).
  \end{align*}
  This completes the induction.

Next we prove
\begin{equation}\label{eq:needslabel4}
  \begin{gathered}
  \widehat{R}_{s_{1},1}\cdots \widehat{R}_{s_{r-1},1}
  \widehat{R}_{s_{r},\omega}\leq \widehat{R}_{k,\omega}(1+ L)^{k},\\[-1mm]
  s_{1}+\cdots+s_{r}=k,\qquad s_{j}\geq 1\quad(1\leq j\leq r).
  \end{gathered}
\end{equation}
In fact, by H\"older's inequality, we get
\begin{align*}
  \widehat{R}_{j,1}\leq &\lvert T \rvert^{-j/2}\Bigl(1+\sum_{\ell=1}^{\lfloor\lvert T\rvert^{1/d}\rfloor}\ell^{d(j+1)-1}\alpha_{\ell}^{(r-p-2)/r}\Bigr)\\
  \leq &\lvert T \rvert^{-j/2}\Bigl(1+\sum_{\ell=1}^{\lfloor\lvert T\rvert^{1/d}\rfloor}\ell^{d(k+\omega)-\omega}\alpha_{\ell}^{(r-p-2)/r}\Bigr)^{\frac{jd}{kd-(d-1)(1-\omega)}}\cdot\\
  &\ \Bigl(1+\sum_{\ell=1}^{\lfloor\lvert T\rvert^{1/d}\rfloor}\ell^{d-1}\alpha_{\ell}^{(r-p-2)/r}\Bigr)^{\frac{(k-j)d-(d-1)(1-\omega)}{kd-(d-1)(1-\omega)}},\\
  \widehat{R}_{j,\omega}\leq &\lvert T \rvert^{-(j+\omega-1)/2}\Bigl(1+\sum_{\ell=1}^{\lfloor\lvert T\rvert^{1/d}\rfloor}\ell^{d(j+\omega)-\omega}\alpha_{\ell}^{(r-p-2)/r}\Bigr)\\
  \leq &\lvert T \rvert^{-(j+\omega-1)/2}\Bigl(1+\sum_{\ell=1}^{\lfloor\lvert T\rvert^{1/d}\rfloor}\ell^{d(k+\omega)-\omega}\alpha_{\ell}^{(r-p-2)/r}\Bigr)^{\frac{jd-(d-1)(1-\omega)}{kd-(d-1)(1-\omega)}}\cdot\\
  &\ \Bigl(1+\sum_{\ell=1}^{\lfloor\lvert T\rvert^{1/d}\rfloor}\ell^{d-1}\alpha_{\ell}^{(r-p-2)/r}\Bigr)^{\frac{(k-j)d}{kd-(d-1)(1-\omega)}}.
\end{align*}
Substituting them into \eqref{eq:needslabel4}, we get
\begin{align*}
  &\widehat{R}_{s_{1},1}\cdots \widehat{R}_{s_{r-1},1}\widehat{R}_{s_{r},\omega}\\
  \leq &\lvert T \rvert^{-(k+\omega-1)/2}\Bigl(1+\sum_{\ell=1}^{\lfloor\lvert T\rvert^{1/d}\rfloor}\ell^{d(k+\omega)-\omega}\alpha_{\ell}^{(r-p-2)/r}\Bigr)\Bigl(1+\sum_{\ell=1}^{\lfloor\lvert T\rvert^{1/d}\rfloor}\ell^{d -1}\alpha_{\ell}^{(r-p-2)/r}\Bigr)^{k}\\
  \leq &\widehat{R}_{k,\omega}
  \Bigl(1+\sum_{\ell=1}^{\lfloor\lvert T\rvert^{1/d}\rfloor}
    \ell^{d -1}\alpha_{\ell}^{(r-p-2)/r}\Bigr)^{k}
  \leq \widehat{R}_{k,\omega}(1+ L)^{k}.
\end{align*}
Note that $L<\infty$ by assumption and $\widehat{R}_{k,\omega}=\lvert T \rvert^{-p/2}+M$. \eqref{eq:hugeexpan3} is proven. \hfill
\end{proof}

We now prove \cref{THM:BARBOURGRAPH2}.

\begin{proof}[Proof of \cref{THM:BARBOURGRAPH2}.]
  For $p=1$, let $f=\Theta h$. By \cref{thm:lemsteinsol}, both
  $\lvert f\rvert_{0,1}$ and $\lvert f\rvert_{1,1}$ are bounded uniformly
  over $h\in\Lambda_1$. The refined expansion \cref{thm:wfwgraph4}, with
  $k=1$, and its bound on $\lvert\widetilde\kappa_2-\kappa_2(W)\rvert$
  give
  \begin{align}
    \sup_{h\in\Lambda_1}\bigl\lvert\mathbb E[h(W)]-\mathcal Nh\bigr\rvert
    \lesssim{}& \lvert T\rvert^{-1/2}m^{2d}
      +M_{n,m,\delta}^{(i)}+M_{n,m}^{(ii)}.
    \label{eq:barbourgraph2base}
  \end{align}
  This supplies the base estimate used for Wasserstein--$1$ below. The
  assertions in parts (a) and (b) of \cref{THM:BARBOURGRAPH2} are otherwise
  empty when $p=1$.

  Henceforth assume $p\geq2$. For $1\leq j\leq p-1$ and $1\leq b\leq m$,
  put
  \begin{align*}
    R_j(b)&:=\lvert T\rvert^{-j/2}
    \left(b^{d(j+1)}+
      \sum_{\ell=b+1}^{b+1+\lfloor\lvert T\rvert^{1/d}/2\rfloor}
      \ell^{d(j+1)-1}\alpha_\ell^{(r-j-2)/r}\right),\\
    \mathcal R_j&:=\min_{1\leq b\leq m}R_j(b).
  \end{align*}
  Thus $\mathcal R_j$ is $\mathcal R_{j,n}$ with the index $n$ suppressed.
  Choose a minimizer $b_j$. Notice that no quantity $\mathcal R_p$ is
  introduced. This is important
  when $r=p+2$: the full order-$p$ expansion requires a strict moment
  inequality, whereas every lower order $j\leq p-1$ satisfies $r>j+2$.

  We first compare the lower-order scales with the common-radius bound used
  previously. By
  H\"older's inequality and the assumed summability, for $1\leq j\leq p-1$,
  \begin{equation}
    R_j(m)^{1/j}\lesssim \lvert T\rvert^{-1/2}
    \Bigl(m^{2d(p-1)}+
      \sum_{\ell=m+1}^{m+1+\lfloor\lvert T\rvert^{1/d}/2\rfloor}
      \ell^{dp-1}\alpha_\ell^{(r-p-1)/r}\Bigr)^{1/(p-1)}.
    \label{eq:needslabel3}
  \end{equation}
  Indeed, after replacing $\alpha_\ell^{(r-j-2)/r}$ by the larger
  $\alpha_\ell^{(r-p-1)/r}$, the weight $d(j+1)-1$ is the interpolation,
  with parameter $j/(p-1)$, between $dp-1$ and $d-1$. Also
  $m^{d(j+1)}\leq m^{2dj}$. Consequently,
  \begin{equation}
    \mathcal R_j^{p/j}\leq R_j(m)^{p/j}\lesssim
      \lvert T\rvert^{-p/2}m^{2dp}+M_{n,m}^{(ii)},
      \qquad 1\leq j\leq p-1,
    \label{eq:lowerorderscale}
  \end{equation}
  where we used $M_{n,m}^{(ii)}\to0$. In particular
  $\mathcal R_j\to0$. Applying \cref{thm:cumuctrl} with order $j$ and
  radius $b_j$ gives
  \begin{equation}
    \bigl\lvert\kappa_{j+2}(W)\bigr\rvert\lesssim\mathcal R_j,
    \qquad 1\leq j\leq p-1.
    \label{eq:lowercumulants}
  \end{equation}

  We next record the ordinary expansion only at the lower orders needed in
  the terminal argument. Since $r>q+2$ for every $1\leq q\leq p-1$, an
  induction using \cref{thm:wfwgraph3} gives, uniformly over $g\in\Lambda_q$,
  \begin{align}
    \mathbb E[g(W)]-\mathcal Ng
    ={}&\sum_{(a,s_{1:a})\in\Gamma(q-1)}(-1)^a
       \prod_{b=1}^a\frac{\kappa_{s_b+2}(W)}{(s_b+1)!}\notag\\
       &\quad{}\times\mathcal N\!\left[\bigcomp_{b=1}^a
       (\partial^{s_b+1}\Theta)g\right]
       +\mathcal O\!\left(\sum_{j=1}^q\mathcal R_j^{q/j}\right).
    \label{eq:hugeexpan1}
  \end{align}
  For completeness, at order $q$ the Stein identity and
  \cref{thm:wfwgraph3}, used with radius $b_q$, give
  \begin{equation*}
    \mathbb E[g(W)]-\mathcal Ng
    =-\sum_{j=1}^{q-1}\frac{\kappa_{j+2}(W)}{(j+1)!}
      \mathbb E[(\partial^{j+1}\Theta)g(W)]+\mathcal O(\mathcal R_q).
  \end{equation*}
  Substitution of the lower-order identities produces the sign
  $(-1)^{a+1}$, as required after adjoining the outer index $j$. The error
  closes because \eqref{eq:lowercumulants} and Young's inequality imply
  \begin{equation*}
    \mathcal R_j\mathcal R_\ell^{(q-j)/\ell}
    \lesssim \mathcal R_j^{q/j}+\mathcal R_\ell^{q/\ell}.
  \end{equation*}

  It remains to use the refined expansion at order $p$. Set
  \begin{align*}
    \widetilde R_p:={}&\lvert T\rvert^{-p/2}m^{d(p+1)}
      +\lvert T\rvert^{-(p-1+\delta)/2}m^{dp}
       \sum_{\ell=m+1}^{m+1+\lfloor\lvert T\rvert^{1/d}/2\rfloor}
       \ell^{d\delta-\delta}\alpha_\ell^{(r-p-1-\delta)/r}\\
      &+\lvert T\rvert^{-(p-1)/2}
       \sum_{\ell=m+1}^{m+1+\lfloor\lvert T\rvert^{1/d}/2\rfloor}
       \ell^{dp-1}\alpha_\ell^{(r-p-1)/r}\\
    \leq{}&\lvert T\rvert^{-p/2}m^{d(p+1)}
       +M_{n,m,\delta}^{(i)}+M_{n,m}^{(ii)}.
  \end{align*}
  By \cref{thm:wfwgraph4}, for $f=\Theta h$ and $h\in\Lambda_p$,
  \begin{align*}
    \mathbb E[h(W)]-\mathcal Nh
    ={}&-\sum_{j=1}^{p-2}\frac{\kappa_{j+2}(W)}{(j+1)!}
       \mathbb E[(\partial^{j+1}\Theta)h(W)]\\
      &-\frac{\widetilde\kappa_{p+1}}{p!}
       \mathbb E[(\partial^p\Theta)h(W)]
       +\mathcal O(\widetilde R_p),
  \end{align*}
  where, writing $D_p:=\lvert\widetilde\kappa_{p+1}-
  \kappa_{p+1}(W)\rvert$, we have
  \begin{equation*}
    D_p\lesssim \lvert T\rvert^{-(p-1)/2}
      \sum_{\ell=m+1}^{m+1+\lfloor\lvert T\rvert^{1/d}/2\rfloor}
      \ell^{dp-1}\alpha_\ell^{(r-p-1)/r}
    =M_{n,m}^{(ii)}.
  \end{equation*}

  Apply \eqref{eq:hugeexpan1} with $q=p-j$ to the terms in the sum and
  with $q=1$ to $(\partial^p\Theta)h$. The outer minus sign changes an
  inner coefficient $(-1)^a$ to $(-1)^{a+1}$. Hence the resulting terms
  are exactly those indexed by $\Gamma(p-1)$, apart from the terminal
  singleton, and
  \begin{align}
    \mathbb E[h(W)]-\mathcal Nh
    ={}&\sum_{(a,s_{1:a})\in\Gamma(p-1)}(-1)^a
       \prod_{b=1}^a\frac{\kappa_{s_b+2}(W)}{(s_b+1)!}
       \mathcal N\!\left[\bigcomp_{b=1}^a
       (\partial^{s_b+1}\Theta)h\right]\\
      &+\frac{\kappa_{p+1}(W)-\widetilde\kappa_{p+1}}{p!}
       \mathcal N[(\partial^p\Theta)h]
       +\mathcal O\!\left(\sum_{j=1}^{p-1}\mathcal R_j^{p/j}
       +\widetilde R_p\right).
    \label{eq:hugeexpan2}
  \end{align}
  Here the error arising from the last expectation is bounded using
  \begin{equation*}
    (\mathcal R_{p-1}+D_p)\mathcal R_1
    \lesssim \mathcal R_{p-1}^{p/(p-1)}+D_p^{p/(p-1)}+\mathcal R_1^p,
  \end{equation*}
  and $D_p^{p/(p-1)}\lesssim D_p$ for all sufficiently large $n$.

  Combining \eqref{eq:hugeexpan2} with the definition of
  $\widetilde R_p$ proves \eqref{eq:hugeexpan4}, while
  \eqref{eq:lowercumulants} proves part (a). Finally,
  \begin{align*}
    \lvert\widetilde\kappa_{p+1}\rvert^{p/(p-1)}
    &\lesssim \lvert\kappa_{p+1}(W)\rvert^{p/(p-1)}
      +D_p^{p/(p-1)}\\
    &\lesssim \mathcal R_{p-1}^{p/(p-1)}
      +(M_{n,m}^{(ii)})^{p/(p-1)},
  \end{align*}
  which proves part (b). Thus, \cref{THM:BARBOURGRAPH2} follows.\hfill
\end{proof}

\section{Proofs of concentration and Berry--Esseen inequalities}\label{sec:concentration}
In this section we prove the general Wasserstein-to-tail transfer in
\cref{concentration}.
\begin{proof}[Proof of \cref{concentration}.]
Fix $n$ and write
$\omega_{p,n}:=\mathcal W_p(\mathcal L(W_n),N(0,1))$. For every
$\epsilon>0$ there is a coupling of $W_n$ and $Z\sim N(0,1)$ such that
$\lVert W_n-Z\rVert_p\leq\omega_{p,n}+\epsilon$. Hence, for $t>0$,
\begin{align*}
  \mathbb P(W_n\geq t)
  &\leq \mathbb P(Z\geq t/2)
    +\mathbb P(W_n-Z\geq t/2)\\
  &\leq \Phi^c(t/2)+\frac{2^p(\omega_{p,n}+\epsilon)^p}{t^p}.
\end{align*}
Letting $\epsilon\downarrow0$ and using $\omega_{p,n}\leq a_n$ proves
\eqref{eq:concentrationgood}. Reflection preserves the Wasserstein distance,
so
\begin{equation*}
  \mathcal W_p(\mathcal L(-W_n),N(0,1))=\omega_{p,n}.
\end{equation*}
Applying
\eqref{eq:concentrationgood} to both $W_n$ and $-W_n$, and adding the two
bounds, gives \eqref{eq:twosidedconcentration}.

We next establish a uniform smoothing estimate. For any $x\in\mathbb R$ and
$\eta>0$, the same coupling gives
\begin{align*}
  \mathbb P(W_n\geq x)
  &\leq \Phi^c(x-\eta)
    +\frac{(\omega_{p,n}+\epsilon)^p}{\eta^p},\\
  \mathbb P(W_n\geq x)
  &\geq \Phi^c(x+\eta)
    -\frac{(\omega_{p,n}+\epsilon)^p}{\eta^p}.
\end{align*}
Since the standard normal density is bounded by $(2\pi)^{-1/2}$, it follows
that
\begin{equation*}
  \sup_{x\in\mathbb R}
  \bigl\lvert\mathbb P(W_n\geq x)-\Phi^c(x)\bigr\rvert
  \leq \frac{\eta}{\sqrt{2\pi}}
    +\frac{(\omega_{p,n}+\epsilon)^p}{\eta^p}.
\end{equation*}
Letting $\epsilon\downarrow0$, using $\omega_{p,n}\leq a_n$, and choosing
$\eta=a_n^{p/(p+1)}$ when $a_n>0$, we obtain
\begin{equation}\label{eq:wpkolmogorov}
  \sup_{x\in\mathbb R}
  \bigl\lvert\mathbb P(W_n\geq x)-\Phi^c(x)\bigr\rvert
  \leq C_p a_n^{p/(p+1)}.
\end{equation}
If $a_n=0$, the same conclusion follows because the Wasserstein metric then
implies $\mathcal L(W_n)=N(0,1)$.

Put $q=p/(p+1)$. Since $a_n\to0$, choose $n_0$ so that
$a_n^q\leq1/2$ for $n\geq n_0$. Fix such an $n$ with $a_n>0$, let
$t\geq1$, and first take any $h\in(0,t)$. After letting
$\epsilon\downarrow0$ in the coupling above, we have
\begin{align}
  \mathbb P(W_n\geq t)-\Phi^c(t)
  &\leq \Phi^c(t-h)-\Phi^c(t)+\frac{a_n^p}{h^p},
    \label{eq:weighteduppercoupling}\\
  \Phi^c(t)-\mathbb P(W_n\geq t)
  &\leq \Phi^c(t)-\Phi^c(t+h)+\frac{a_n^p}{h^p}.
    \label{eq:weightedlowercoupling}
\end{align}
Indeed, the upper and lower comparisons follow, respectively, from
\begin{align*}
 \{W_n\geq t\}&\subseteq
   \{Z\geq t-h\}\cup\{\lvert W_n-Z\rvert>h\},\\
 \{Z\geq t+h\}&\subseteq
   \{W_n\geq t\}\cup\{\lvert W_n-Z\rvert>h\}.
\end{align*}
Both inclusions remain valid with all displayed tail events strict, and the
normal probabilities are unchanged because $Z$ has no atoms. Thus the same
two inequalities hold with $\{W_n>t\}$ throughout. Retaining both directions
is essential for the later reflection.

Now set $h=a_n^q t$. Since
$0<h\leq t/2$ and the normal density $\varphi$ decreases on $[0,\infty)$,
\begin{align*}
  \Phi^c(t-h)-\Phi^c(t)&\leq h\varphi(t-h)
    \leq a_n^q t\varphi(t/2),\\
  \Phi^c(t)-\Phi^c(t+h)&\leq h\varphi(t)
    \leq a_n^q t\varphi(t).
\end{align*}
Moreover, $a_n^p/h^p=a_n^q/t^p$. Because
$\sup_{t\geq1}t^{p+1}\varphi(t/2)<\infty$, \eqref{eq:weighteduppercoupling}
and \eqref{eq:weightedlowercoupling}, together with their strict-event
versions, imply
\begin{equation}\label{eq:weightedpositivetail}
  \max\left\{
    \bigl\lvert\mathbb P(W_n\geq t)-\Phi^c(t)\bigr\rvert,
    \bigl\lvert\mathbb P(W_n>t)-\Phi^c(t)\bigr\rvert
  \right\}
  \leq \frac{C_p a_n^q}{t^p},\qquad t\geq1.
\end{equation}
For negative arguments we retain the strict-event version. If $s\geq1$,
then
\begin{equation}\label{eq:negative-tail-reflection}
  \bigl\lvert\mathbb P(W_n\geq-s)-\Phi^c(-s)\bigr\rvert
  =\bigl\lvert\mathbb P(-W_n>s)-\Phi^c(s)\bigr\rvert.
\end{equation}
Reflection preserves the Wasserstein distance, so applying
\eqref{eq:weightedpositivetail} to $-W_n$ proves the same weighted estimate
for $t\leq-1$. Combining it with \eqref{eq:wpkolmogorov} on
$\lvert t\rvert\leq1$, and using
$1+\lvert t\rvert^p\asymp\lvert t\rvert^p$ when $\lvert t\rvert\geq1$,
proves \eqref{eq:newberryesseengood10}. Indices with $a_n=0$ are immediate.

It remains to record the stronger far-tail alternative. In the same two
coupling comparisons, choose $h=t/2$. For $t\geq1$, their non-strict and
strict versions give
\begin{equation}\label{eq:farpositivetail}
  \begin{aligned}
  &\max\left\{
    \bigl\lvert\mathbb P(W_n\geq t)-\Phi^c(t)\bigr\rvert,
    \bigl\lvert\mathbb P(W_n>t)-\Phi^c(t)\bigr\rvert
  \right\}\\
  &\qquad\leq \Phi^c(t/2)+\frac{2^pa_n^p}{t^p}
  \leq e^{-t^2/8}+\frac{2^pa_n^p}{t^p}.
  \end{aligned}
\end{equation}
Using the strict-event term in \eqref{eq:farpositivetail} after the
reflection \eqref{eq:negative-tail-reflection} proves the same estimate for
$t\leq-1$. Taking the better of this bound and
\eqref{eq:newberryesseengood10} proves
\eqref{eq:newberryesseentworegime}, with $c_p=1/8$ if desired. \hfill
\end{proof}
\section{Details on constructive graph approach}\label{sec:mixingmainpart}

For notational simplicity whenever it is not ambiguous we will write $\alpha_{\ell}$, $\sigma$, $W$, $X_i$ and $T$ for $\alpha_{\ell,n}$, $\sigma_n$, $W_n$, $X^{\scalebox{0.6}{$(n)$}}_{i}$ and $T_n$ respectively. The rest of the section is constructed as follows: In \cref{sec:genogram}, we formally define a genogram and related concepts. In \cref{sec:summationterms}, we define three types of sums corresponding to a genogram, which will be used later in the expansion. In \cref{sec:proofkeylemma2}, we show how to achieve the expansion by growing a class of genograms. Finally, in \cref{sec:remaindercontrol}, we control the remainders using the mixing coefficients.

\subsection{Graph theory background and definition of genograms}\label{sec:genogram}
A \emph{rooted tree} is a tree in which one vertex has been designated the root. In a rooted tree, the \emph{parent} of a vertex $v$ is the vertex connected to $v$ on the path from the root to $v$; every vertex has a unique parent except the root, which has no parent. A \emph{child} of a vertex $v$ is a vertex of which $v$ is the parent. An \emph{ancestor} of a vertex $v$ is any vertex other than $v$ which is on the path from the root to $v$. A \emph{sibling} to a vertex $v$ is any other vertex on the tree which has the same parent as $v$. A \emph{leaf} is a vertex with no children. See \cite{bender2010lists} for a detailed exposition.

An order $k$ \emph{genogram} is defined as the tuple $G:=\bigl(V,E, \{ s_{v} \}_{v\in V}\bigr)$, where $(V,E)$ is a rooted tree with a vertex set $V$ ($\lvert G \rvert:=\lvert V \rvert=k$) and an edge set $E$, and $s_{v}$ is an integer called the \emph{identifier} associated to each $v\in V$ that satisfies the requirements below:
\begin{enumerate}[label=(\alph*)]
  \item $s_{v}=0$ for the root $v$. $s_{v}\geq -1$ for any other vertex $v$;
  \item $s_{v}\geq 0$ if $v$ is a child of the root;
  \item If $v$ has more than one child, identifiers of $v$'s children must be non-negative and mutually different.
\end{enumerate}
Beware that the identifiers are part of the genograms by definition. We say that a vertex $v$ is \emph{negative} if and only if $s_{v}=-1$, \emph{nil} if and only if $s_{v}=0$, \emph{positive} if and only if $s_{v}\ge 1$. The requirements above implies that the identifier of each child of $v$ is different, and therefore, $v$'s children can be uniquely identified by their identifiers. The last requirement also suggests that if $v$ has more than one child, there is no negative and at most one nil among them. In other words, a negative vertex has no sibling, a nil vertex only has positive siblings, and any vertex must have an identifier different from all its siblings.

\begin{figure}
  \centering
  \subfloat[$G_{1}$]{%
    \centering
    \begin{tikzpicture}[scale=0.75,transform shape]
      \begin{scope}[VertexStyle/.append style = {fill=lightgray}]
        \Vertex[x=0,y=-1.2,L=$0$]{1}
      \end{scope}
      \Vertex[x=1.5,y=0,L=$2$]{2}
      \Vertex[x=2.5,y=-0.8,L=$1$]{3}
      \Vertex[x=2,y=-2,L=$0$]{4}
      \Vertex[x=4,y=-1.6,L=$-1$]{5}
      \Vertex[x=5.5,y=-0.8,L=$2$]{6}
      \Vertex[x=6,y=-2,L=$0$]{7}
      \tikzstyle{EdgeStyle}=[]
      \Edge(1)(2)
      \Edge(1)(3)
      \Edge(1)(4)
      \Edge(4)(5)
      \Edge(5)(6)
      \Edge(5)(7)
    \end{tikzpicture}
    \label{fig:ex1}
  }
  \hspace{15pt}
  \subfloat[$G_{2}$]{%
    \centering
    \begin{tikzpicture}[scale=0.75,transform shape]
      \begin{scope}[VertexStyle/.append style = {fill=lightgray}]
        \Vertex[x=0.5,y=-1.2,L=$0$]{1}
      \end{scope}
      \Vertex[x=2,y=-0.4,L=$0$]{2}
      \Vertex[x=4,y=-0.8,L=$3$]{3}
      \Vertex[x=5.5,y=0,L=$2$]{4}
      \Vertex[x=6,y=-1.6,L=$1$]{5}
      \Vertex[x=7.5,y=-0.8,L=$-1$]{6}
      \Vertex[x=3,y=-2,L=$5$]{7}
      \tikzstyle{EdgeStyle}=[]
      \Edge(1)(2)
      \Edge(2)(3)
      \Edge(3)(4)
      \Edge(3)(5)
      \Edge(5)(6)
      \Edge(2)(7)
    \end{tikzpicture}
    \label{fig:ex2}
  }
  \caption{Examples of order-$7$ genograms.}
  \label{fig:expofgeno}
\end{figure}
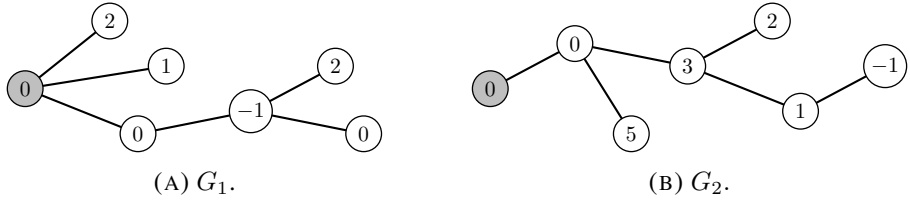

Furthermore, denote the set of all possible order-$k$ genograms by $\mathcal{G}(k)$. {\cref{fig:expofgeno} depicts two examples of genograms $G_{1},G_{2}\in\mathcal{G}(7)$, where each circle represents a vertex, the one representing the root is filled with gray, and the identifiers are marked inside the circles.}

We remark that the notion of genograms resembles the ordered trees in combinatorics. An \emph{ordered (rooted) tree} $(V,E,\prec)$ is a rooted tree $(V,E)$ where the children of every vertex are ordered (the order denoted by $\prec$) \citep{stanley2011enumerative}. Note that $\prec$ is a strict partial order on the vertex set $V$. By definition every genogram induces a unique ordered tree if we set $v_{1}\prec v_{2}\Leftrightarrow s_{v_{1}}>s_{v_{2}}$ whenever $v_{1},v_{2}\in V$ are siblings. However, an ordered tree corresponds to infinitely many genograms since the largest identifier is allowed to take any sufficiently large value in $\mathbb{N}\cup\{-1\}$.

\paragraph*{Compatible labeling, parent, progenitor, and ancestor}

Next we consider a labeling of the vertices of a genogram (or the induced ordered tree of a genogram). We say a labeling $V=\bigl\{ v[1],\cdots,v[k] \bigr\}$ is \emph{compatible} with $G$ (or $(V,E,\prec)$) if and only if
\begin{enumerate}[label=(\alph*)]
  \item \label{itm:labeling1} It follows from a depth-first traversal: $v[1]$ is the root, and for any $1\leq j\leq k-1$, the vertex $v[j+1]$ is chosen to be a child of the vertex with the largest label $\ell\le j$ that has an as-yet unlabeled child. In particular, $v[j+1]$ is $v[j]$'s child as long as $v[j]$ has an as-yet unlabeled child;
  \item \label{itm:labeling2} It respects the partial order $\prec$ induced by $G$: If $v[j]$ and $v[h]$ ($2\leq j, h\leq k$, $j\neq h$) are siblings, then we have $s_{j}>s_{h}\Leftrightarrow j< h$ (or equivalently, $v[j]\prec v[h]\Leftrightarrow j<h$). In other words, if a vertex has more than one child, a child with a larger identifier has a smaller label. In particular, if $s_{j}=0$, then $v[j]$ has the largest label $j$.
\end{enumerate}

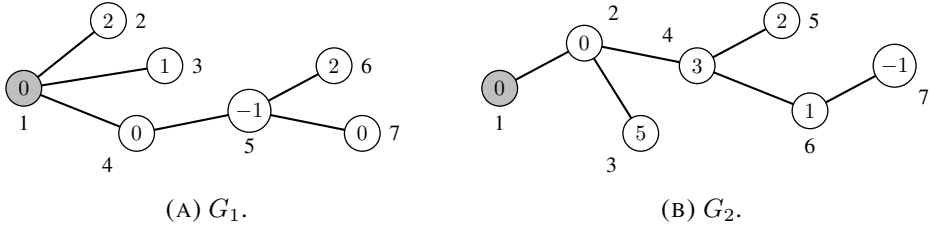
\begin{figure}
  \centering
  \subfloat[$G_{1}$]{%
    \centering
    \begin{tikzpicture}[scale=0.75,transform shape]
      \begin{scope}[VertexStyle/.append style = {fill=lightgray}]
        \Vertex[x=0,y=-1.2,L=$0$]{1}   \extralabel[6pt]{-90}{1}{1}
      \end{scope}
      \Vertex[x=1.5,y=0,L=$2$]{2}    \extralabel[6pt]{0}{2}{2}
      \Vertex[x=2.5,y=-0.8,L=$1$]{3}   \extralabel[6pt]{0}{3}{3}
      \Vertex[x=2,y=-2,L=$0$]{4}   \extralabel[6pt]{-135}{4}{4}
      \Vertex[x=4,y=-1.6,L=$-1$]{5}\extralabel[6pt]{-90}{5}{5}
      \Vertex[x=5.5,y=-0.8,L=$2$]{6}\extralabel[6pt]{0}{6}{6}
      \Vertex[x=6,y=-2,L=$0$]{7}\extralabel[6pt]{0}{7}{7}
      \tikzstyle{EdgeStyle}=[]
      \Edge(1)(2)
      \Edge(1)(3)
      \Edge(1)(4)
      \Edge(4)(5)
      \Edge(5)(6)
      \Edge(5)(7)
    \end{tikzpicture}
    \label{fig:ex1withlabel}
  }
  \hspace{15pt}
  \subfloat[$G_{2}$]{%
    \centering
    \begin{tikzpicture}[scale=0.75,transform shape]
      \begin{scope}[VertexStyle/.append style = {fill=lightgray}]
        \Vertex[x=0.5,y=-1.2,L=$0$]{1}   \extralabel[6pt]{-90}{1}{1}
      \end{scope}
      \Vertex[x=2,y=-0.4,L=$0$]{2}    \extralabel[6pt]{45}{2}{2}
      \Vertex[x=4,y=-0.8,L=$3$]{3}   \extralabel[6pt]{135}{4}{3}
      \Vertex[x=5.5,y=0,L=$2$]{4}   \extralabel[6pt]{0}{5}{4}
      \Vertex[x=6,y=-1.6,L=$1$]{5}\extralabel[6pt]{-90}{6}{5}
      \Vertex[x=7.5,y=-0.8,L=$-1$]{6}\extralabel[6pt]{-45}{7}{6}
      \Vertex[x=3,y=-2,L=$5$]{7}\extralabel[6pt]{-135}{3}{7}
      \tikzstyle{EdgeStyle}=[]
      \Edge(1)(2)
      \Edge(2)(3)
      \Edge(3)(4)
      \Edge(3)(5)
      \Edge(5)(6)
      \Edge(2)(7)
    \end{tikzpicture}
    \label{fig:ex2withlabel}
  }
  \caption{Examples of order-$7$ genograms with the compatible labeling.}
  \label{fig:expofgenowithlabel}
\end{figure}

\cref{fig:expofgenowithlabel} shows the compatible labelings of $G_{1}$ and $G_{2}$, where the labels are marked outside the circles that represent the vertices.

We remark that there is a unique compatible labeling given any genogram $G$ (or any ordered tree $(V,E,\prec)$). %

Now we introduce more notations in order to express the compatible requirements for identifiers and the labeling in a more concise manner. Let $G=(V,E,s_{1:k})$ be an order-$k$ genogram with vertices labelled as $V=\bigl\{ v[1],\cdots,v[k] \bigr\}$and  where $s_{j}$ is the identifier of $v[j]$ for $1\leq j\leq k$. We denote the label of $v[j]$'s \emph{parent} by $p(j,G)$ for $2\leq j\leq k$, and the label set of $v[j]$'s \emph{ancestors} by $A(j,G)$ (we set $A(1,G)=\emptyset$).
Moreover, we write
  \begin{equation}\label{eq:defgj}
    g(j,G):=\sup \{ \ell: \ell=1\text{ or }\ell\in A(j,G)~\&~s_{\ell}\geq 1 \},
  \end{equation}
and call $v[g(j,G)]$ the \emph{progenitor} of $v[j]$. In particular, we have that $g(1,G)= 1$. Intuitively, $v[g(j,G)]$ is the positive vertex closest to $v[j]$ in its ancestry if such vertex exists, in which case there is a path from $v[g(j,G)]$ to $v[p(j,G)]$, the parent of $v[j]$, such that $v[g(j,G)]$ is the only positive vertex along the path. Otherwise, $v[g(j,G)]$ is set to be the root. Note that $v[g(j,G)]\neq v[j]$ for $2\leq j\leq k$. Take the genograms $G_{1}$ and $G_{2}$ in \cref{fig:expofgenowithlabel} as examples, $g(j,G_{1})=1$ for all $1\leq j\leq 7$ while in $G_{2}$, $g(1,G_{2})=g(2,G_{2})=g(3,G_{2})=g(4,G_{2})=1$, $g(5, G_{2})=g(6,G_{2})=4$, and $g(7,G_{2})=6$.
We further denote
\begin{equation}\label{eq:defuj}
  u(j,G):=\sup\{ \ell\in \{ j \}\cup A(j,G): s_{\ell}\geq 0 \}.
\end{equation}
In other words, $v[u(j,G)]$ is the closest non-negative vertex in $v[j]$'s ancestry if $v[j]$ is negative, otherwise $u(j,G)=j$. In particular, $s_{j}=-1\Leftrightarrow u(j,G)<j$, $s_{j}\geq 0\Leftrightarrow u(j,G)=j$. For example, in the genogram $G_{1}$ shown in \cref{fig:ex1withlabel}, $u(5,G_{1})=4$ and $u(j,G_{1})=j$ for $j\neq 5$. {For ease of notation, when there is no ambiguity, we will abuse notations and write $p(j),A(j),g(j),u(j)$ to mean $p(j,G),A(j,G),g(j,G),u(j,G)$.}

We remark that the labeling has to respect the following properties:
\begin{proposition}\label{thm:rqmofgenogram}
  Let $k$ be a positive integer, $(V,E)$ be a rooted tree with the vertex set $V=\bigl\{ v[1],\cdots,v[k]\bigr\}$ and edge set $E$, and $s_{1},\cdots,s_{k}$ be $k$ integers. $\bigl(V,E,\{ s_{1:k} \}\bigr)$ is a genogram with the compatible labeling if and only if all the following statements are true:
  \begin{enumerate}[label=(\alph*)]
    \item \label{itm:dft} $p(j+1)=\max\{ p(\ell):\ell\geq j+1, p(\ell)\leq j \}$ for $1\leq j\leq k-1$;
    \item \label{itm:secondrq} $s_{1}=0$. $s_{j}\geq -1$ for $2\leq j\leq k$;
    \item \label{itm:thirdrq} If $s_{j}=-1$ ($2\leq j\leq k$), then (i) $p(j)\neq 1$, (ii) $p(j)=p(h)\Leftrightarrow j=h$ for $2\leq h\leq k$;
    \item \label{itm:fourthrq} If $p(j)=p(h)$ ($2\leq j,h\leq k$), then $s_{j}>s_{h}\Leftrightarrow j<h$, $s_{j}=s_{h}\Leftrightarrow j=h$.
  \end{enumerate}
\end{proposition}

\begin{proof}
  After the labels $1,\ldots,j$ have been assigned, the vertices that can be
  parents of an as-yet unlabeled vertex are precisely
  $\{v[p(\ell)]:\ell\geq j+1,\ p(\ell)\leq j\}$.  Thus the depth-first rule
  in \ref{itm:labeling1} is equivalent to \ref{itm:dft}; conversely,
  \ref{itm:dft} reproduces that rule at every step.  Condition
  \ref{itm:secondrq} gives the allowed identifier range and fixes the root
  identifier.  The two parts of \ref{itm:thirdrq} say respectively that a
  child of the root cannot be negative and that a negative vertex has no
  sibling.  Finally, \ref{itm:fourthrq} says that sibling identifiers are
  distinct and that their label order is the decreasing order of their
  identifiers.  These are exactly the remaining genogram requirements and
  the compatibility condition \ref{itm:labeling2}.
\end{proof}

\paragraph*{Induced sub-genograms}

Lastly, given $G=(V,E,s_{1:k})$ and $1\leq j\leq k$, we call an order-$j$ genogram $G[j]:=(V',E',s_{1:j})$ the induced \emph{sub-genogram} of $G$ the genograms by setting $V':=\bigl\{ v[1],\cdots,v[j] \bigr\}$ and $E'\subseteq E$ be the set of all edges between the vertices $V'$ in $G$. %
We further denote $H\subseteq G$ or $G\supseteq H$ if and only if a genogram $H$ is a sub-genogram of $G$. If $j<k$, we say $G[j]$ is a \emph{proper sub-genogram} of $G$ and write $G[j]\subsetneq G$ or $G\supsetneq G[j]$.

\subsection{Constructing sums from genograms}\label{sec:summationterms}
Consider a $d$-dimensional random field $(X_{i})_{i\in T}$ with the index set $T$ satisfying $T\subset \mathbb{Z}^{d}$ and $\lvert T \rvert<\infty$. We assume that $$\sigma^{2}:=\operatorname{Var} \bigl(\sum_{i\in T}X_{i}\bigr)>0\qquad\text{and write}\qquad W\!:=\sigma^{-1}\!\sum_{i\in T}X_{i}.$$ For any index subset $J\subseteq T$, we denote $$W(J):=\sigma^{-1}\sum_{i\in T\backslash  J}X_{i}.$$ All identities below are understood whenever the displayed expectations are finite. In this subsection, we build sums $\mathcal{S}(G),\mathcal{T}_{f}(G),\mathcal{U}_{f}(G)$ from a genogram $G$ and a given function $f\in \mathcal{C}^{k-1}(\mathbb{R})$ in four steps:
\begin{enumerate}[label=(\alph*)]
  \item Use the genogram $G$ to define the sets of values taken by running indices;
  \item Introduce the generalized covariance operator $\mathcal{D}^{*}$;
  \item Construct an operator $\mathcal{E}_{G}$ from $G$, which leads to the summand;
  \item Define $\mathcal{S}(G),\mathcal{T}_{f}(G),\mathcal{U}_{f}(G)$.
\end{enumerate}
Note that these sums will be used in the next subsection to track the expansion of the quantity $\mathbb{E} [Wf(W)]$.

Firstly, as we have pointed in the roadmap, we will construct from an order-$k$ genogram $G$ sums with $k$ running indices, where the $v[j]$ corresponds to the $j$-th running index, denoted by $i_{j}$. Since $i_{j}$ will appear in the subscript of $X_{i_{j}}$, the value of $i_{j}$ needs to be chosen from $T$. It is important to note that the vertices $v[1],\cdots,v[k]$ as well as the genogram do not represent specific values of $i_1,\cdots,i_k$. The genogram reflects the dependency structure of random variables appearing in the sum by encoding the distance structure between the running indices.

Setting $B_{1}:=T$ and $D_{1}:=\emptyset$, $i_{1}$ will be summed over $B_{1}\backslash D_{1}=T$. Next given the choice of the first $j-1$ running indices ($j\geq 2$), we aim to define two index sets, $B_{j}$ and $D_{j}$, using the chosen values $i_{1},\cdots,i_{j-1}$ and the order-$j$ sub-genogram $G[j]$. In the last step, we will take the sums over $i_{k}\in B_{k}\backslash D_{k}, T\backslash B_{k}\text{ or }T\backslash D_{k}$, and then $i_{k-1}\in B_{k-1}\backslash D_{k-1},\cdots,i_{1}\in B_{1}\backslash D_{1}$ in turn. We call $B_{j}$ the \emph{outer constraint (set)} of the running index $i_{j}$, and $D_{j}$, the \emph{inner constraint (set)} of $i_{j}$. For ease of notation, on most occasions we do not explicitly write out the dependencies on $G[j]$ and $i_{1},\cdots,i_{j-1}$ when referring to the constraint sets $B_{j}$ and $D_{j}$. However, if we are considering multiple genograms, we will use $B_{j}(G)$ and $D_{j}(G)$ to specify the constraint sets of $i_{j}$ with respect to $G$ to avoid ambiguity.

We will formally define $B_{j}$ and $D_{j}$ for $2\leq j\leq k$ later by induction. But first we consider the case $j=2$ to build intuition. When the first running index is set to be some specific element $i_{1}\in T$, we define $B_{2},D_{2}\subseteq T$ using $i_{1}$ and $s_{2}$, and $i_{2}$ will be summed over $B_{2}\backslash D_{2}$. If $s_{2}=0$, $i_{2}$ will be summed over all the indices of distance no greater than $m$ from $i_{1}$, in which case $B_{2}$ and $D_{2}$ are defined by
\begin{equation*}
  B_{2}:=\{ i\in T: \lVert i-i_{1} \rVert\leq m \},\quad D_{2}:=\emptyset,
\end{equation*}
where $\lVert \cdot \rVert$ is the maximum norm on $\mathbb{Z}^{d}$. Note that $m$ is a positive integer that we have fixed earlier.

Otherwise, $s_{2}\geq 1$, we let $i_{2}$ take the sum over a singleton $B_{2}\backslash D_{2}$. In other words, the second running index has only one possible choice in the summation. Different values of $1\leq s_{2}\leq \bigl\lvert \{ i\in T: \lVert i-i_{1} \rVert\geq m+1 \} \bigr\rvert=:s^{*}$ will correspond to different singletons of elements in $\{ i\in T: \lVert i-i_{1} \rVert\geq m+1 \}$. Let $\prec$ be a strict total order on $\mathbb{Z}^d$. With the first level comparing the value of $\lVert i-i_{1} \rVert$ and the second level using the strict order $\prec$, we perform a two-level sorting of all elements $i$ from $\{ i\in T: \lVert i-i_{1} \rVert\geq m+1\}$ and obtain an ascending sequence, $z_{1},\cdots,z_{s^{*}}$. Now $i_{2}$ is chosen to be $z_{s_{2}}$, and
\begin{align*}
  B_{2}:=&\{ i\in T: \lVert i-i_{1} \rVert\leq m \}\cup \{ z_{j}:1\leq j\leq s_{2} \},\\ 
  D_{2}:=&\{ i\in T: \lVert i-i_{1} \rVert\leq m \}\cup \{ z_{j}:1\leq j\leq s_{2}-1 \}.
\end{align*}

The motivation of using singletons arises from deriving \eqref{eq:split2}, where we have decomposed the quantity $\mathbb{E} [X_{i}f(W_{i,m})]$ into a telescoping sum:
\begin{equation*}
  \mathbb{E} [X_{i}f(W_{i,m})]=\sum_{j=m+1}^{n-1}\mathbb{E} \bigl[X_{i}\bigl((f(W_{i,j-1})-f(W_{i,j}^{*}))+(f(W_{i,j}^{*})-f(W_{i,j}))\bigr)\bigr].
\end{equation*}
In order to accurately approximate the differences $f(W_{i,j-1})-f(W_{i,j}^{*})$ and $f(W_{i,j}^{*})-f(W_{i,j})$ by the Taylor expansions, the differences between the inputs of $f$ need to be small enough. The best we can do is to put exactly one random variable in each of such differences.

In general, we introduce some new notations in order to define $B_{j}$ and $D_{j}$ ($2\leq j\leq k$) more conveniently. Denote the $m$-neighborhood of $J\subseteq T$ as
\begin{equation}\label{eq:defnj}
  N(J):=\{ i\in T: d(i,J)\leq m\},\text{ where }d(i,J):=\min_{j\in J} \lVert i-j\rVert.
\end{equation}
  We will treat each element of $T\backslash N(J)$ sequentially starting from the closest elements to $J$. To make this precise, for any positive integer $j$, we write $A^{(j)}(J):=\{i\in T: d(i,J)=j\}$ to be the set of indices in $T$ that are at distance of $j$ from $J$. We have $T\backslash N(J)=\bigcup_{j=m+1}^{\infty}A^{(j)}(J)$, with only finitely many nonempty sets in this union. We write $\operatorname{rk}(i,J)= \sum_{\ell=m+1}^{j} \lvert A^{(\ell)}(J)\rvert+r$ if $d(i,J)= j+1$ ($j\geq m$) and $i$ is the $r$-th smallest element of $A^{(j+1)}(J)$ with respect to the order $\prec$. Therefore, $\operatorname{rk}(\cdot,J)$ is a bijection between $T\backslash N(J)$ and $\{\ell\in \mathbb{Z}:1\leq \ell\leq \lvert T\backslash N(J)\rvert\}$. The smaller $\operatorname{rk}(i,J)$ is, the closer $i$ is from $J$. The value of $\operatorname{rk}(\cdot,J)$ does not have an intrinsic mathematical significance but will allow us to determine the order in which we will treat the indices in $T\backslash N(J)$. Now denote
\begin{equation}\label{eq:defnsj}
  N^{(s)}(J):=\{ i\in T: i\in N(J) \text{ or } \operatorname{rk}(i,J)\leq s \}.
\end{equation}
  We note in particular that $N^{(0)}(J)=N(J)$, and $N^{(s_L)}(J)=\{i\in T:d(i,J)\le L\}$ for $L\geq m+1$ and $s_L:=\sum_{\ell=m+1}^L\lvert A^{(\ell)}(J)\rvert$.

For any $2\leq j\leq k$, fix the genogram $G$ and indices
$i_{1}\in B_{1}\backslash D_{1},\ldots,$
$i_{j-1}\in B_{j-1}\backslash D_{j-1}$. We define
\begin{align}
  B_{j}:= &
  \begin{cases}
    N^{(s_{j})}\bigl( i_{\ell}: \ell\in A(j)\bigr)\cup D_{g(j)} & \quad\text{ if } s_{j}\geq 0 \\
    B_{u(j)}                                                    & \quad \text{ if } s_{j}=-1
  \end{cases}., \label{eq:defbj} \\
  D_{j}:= &
  \begin{cases}
    N^{(s_{j}-1)}\bigl( i_{\ell}: \ell\in A(j)\bigr)\cup D_{g(j)} & \ \text{ if } s_{j}\geq 1 \\
    D_{g(j)}                                                      & \ \text{ if } s_{j}\leq 0
  \end{cases}.\label{eq:defdj}
\end{align}
Here (by abuse of notation) $A(j)$ is the label set of $v[j]$'s ancestors, and $u(j)$ and $g(j)$ are defined in \eqref{eq:defuj} and \eqref{eq:defgj}. {Note that $B_{j}$ and $D_{j}$ depend on $G[j]$ through $s_{j},A(j),u(j),$ and $g(j)$.} Moreover, by definition, we have that $ D_{j}\subseteq B_{j}$ for any $1\leq j\leq k$. We remark that when the identifier $s_j=0$ then $B_j\backslash D_j\subseteq N(i_{\ell}: \ell\in A(j))$, and when $s_j\geq 1$ then $B_j\backslash D_j$ is either empty or a singleton with element the unique $i$ such that $\operatorname{rk}\bigl(i,\{i_{\ell}\in T: \ell\in A(j)\}\bigr)=s_j$. Finally, if $s_j=-1$, then $B_j\backslash D_j=B_{p(j)}\backslash D_{p(j)}$.

For instance, we consider the genograms $G_{1}$ and $G_{2}$ shown in \cref{fig:expofgenowithlabel}. The constraint sets of the running indices are presented in \cref{tab:sets}.

\begin{table}
  \centering
  \caption{The constraint sets with respect to $G_{1}$ and $G_{2}$.}
  \label{tab:sets}
  \small
  \begin{tabular}{@{}rrcccc@{}}
    \toprule
                             & $j$            & $1$                                                                       & $2$                                                   & $3$                                                   & $4$                    \\
    \midrule
    \multirow{2}{*}{$G_{1}$} & $B_{j}(G_{1})$ & $T$                                                                       & $N^{(2)}(i_{1})$                                      & $N^{(1)}(i_{1})$                                      & $N(i_{1})$             \\
                             & $D_{j}(G_{1})$ & $\emptyset$                                                               & $N^{(1)}(i_{1})$                                      & $N(i_{1})$                                            & $\emptyset$            \\
    \midrule
    \multirow{2}{*}{$G_{2}$} & $B_{j}(G_{2})$ & $T$                                                                       & $N(i_{1})$                                            & $N^{(5)}(i_{1},i_{2})$                                & $N^{(3)}(i_{1},i_{2})$ \\
                             & $D_{j}(G_{2})$ & $\emptyset$                                                               & $\emptyset$                                           & $N^{(4)}(i_{1},i_{2})$                                & $N^{(2)}(i_{1},i_{2})$ \\
    \bottomrule
  \end{tabular}

  \medskip
  \begin{tabular}{@{}rrccc@{}}
    \toprule
                             & $j$            & $5$ & $6$ & $7$ \\
    \midrule
    \multirow{2}{*}{$G_{1}$} & $B_{j}(G_{1})$ & $N(i_{1})$ & $N^{(2)}(i_{1},i_{4},i_{5})$ & $N(i_{1},i_{4},i_{5})$ \\
                             & $D_{j}(G_{1})$ & $\emptyset$ & $N^{(1)}(i_{1},i_{4},i_{5})$ & $\emptyset$ \\
    \midrule
    \multirow{2}{*}{$G_{2}$} & $B_{j}(G_{2})$
      & \makecell[c]{$N^{(2)}(i_{1},i_{2},i_{4})$\\$\cup N^{(2)}(i_{1},i_{2})$}
      & \makecell[c]{$N^{(1)}(i_{1},i_{2},i_{4})$\\$\cup N^{(2)}(i_{1},i_{2})$}
      & \makecell[c]{$N^{(1)}(i_{1},i_{2},i_{4})$\\$\cup N^{(2)}(i_{1},i_{2})$} \\
                             & $D_{j}(G_{2})$
      & \makecell[c]{$N^{(1)}(i_{1},i_{2},i_{4})$\\$\cup N^{(2)}(i_{1},i_{2})$}
      & \makecell[c]{$N(i_{1},i_{2},i_{4})$\\$\cup N^{(2)}(i_{1},i_{2})$}
      & \makecell[c]{$N(i_{1},i_{2},i_{4})$\\$\cup N^{(2)}(i_{1},i_{2})$} \\
    \bottomrule
  \end{tabular}
\end{table}

Secondly, as described in the roadmap, we define the \emph{generalized covariance operator} $\mathcal{D}^{*}$ on a finite sequence of random variables $(Y_{i})_{i\geq  1}$. To do so, we also need to inductively define another operator $\mathcal{D}$ that takes in a finite sequence of random variables and outputs a new random variable. For any random variable $Y$, define
\begin{align*}
  \mathcal{D}^{*}(Y):=\mathbb{E} [Y],   \quad \mathcal{D}(Y):=Y-\mathcal{D}^{*}(Y)=Y-\mathbb{E} [Y].
\end{align*}
Suppose $\mathcal{D}$ is already defined for a random sequence of length $t-1$. Then for any random variables $Y_{1},\cdots,Y_{t}$, let
\begin{align*}
  \mathcal{D}(Y_{1},\ldots,Y_{t})
    &:=\mathcal{D}\bigl(Y_{1}\mathcal{D}(Y_{2},\ldots,Y_{t})\bigr) \\
    &=Y_{1}\mathcal{D}(Y_{2},\ldots,Y_{t})
      -\mathbb{E}\bigl[Y_{1}\mathcal{D}(Y_{2},\ldots,Y_{t})\bigr], \\
  \mathcal{D}^{*}(Y_{1},\ldots,Y_{t})
    &:=\mathcal{D}^{*}\bigl(Y_{1}\mathcal{D}(Y_{2},\ldots,Y_{t})\bigr) \\
    &=\mathbb{E}\bigl[Y_{1}\mathcal{D}(Y_{2},\ldots,Y_{t})\bigr].
\end{align*}

In particular, for any two random variables $Y_{1}$ and $Y_{2}$, $\mathcal{D}^{*}(Y_{1},Y_{2})=\operatorname{Cov} (Y_{1},Y_{2})$ gives the covariance between $Y_{1}$ and $Y_{2}$. Here we remark that $\mathcal{D}(Y_{1},\cdots,Y_{t})$ and $\mathcal{D}^{*}(Y_{1},\cdots,Y_{t})$ are well-defined for a tuple of $t$ random variables $(Y_{i})_{i=1}^{t}$ supposing that for any $i,j\in\mathbb{N}_{+}$ such that $i\leq j\leq t$ we have $\mathbb{E} \bigl[\lvert Y_{i}Y_{i+1}\cdots Y_{j} \rvert\bigr]<\infty$.
It is straightforward to see from the definition that both operators are multilinear. We will show more properties of them in \cref{sec:lemma5}.

Thirdly, we construct from any genogram a new operator $\mathcal{E}_{G}$ that maps from $\lvert G \rvert$ random variables to a real number. Note that this $\mathcal{E}_{G}$ operator will provide us with the summands in $\mathcal{S}(G)$, $\mathcal{T}_{f}(G)$, and $\mathcal{U}_{f}(G)$. If $\lvert G \rvert=1$, for any random variable $Y$, define $\mathcal{E}_{G}(Y):=\mathcal{D}^{*}(Y)=\mathbb{E} [Y]$. Suppose $\mathcal{E}_{G}$ is already defined for $\lvert G \rvert\leq k-1$. Consider the case where $\lvert G \rvert=k$. Let
\begin{equation}\label{eq:defvq0}
  q_{0}:=\sup \{ j: j=1\text{ or }p(j)\neq j-1\text{ for }2\leq j\leq k \},
\end{equation}

  In other words, either $v[q_{0}-1]$ is the leaf with the largest label smaller than $k$, or alternatively $v[k]$ is the only leaf and $q_{0}=1$. Intuitively, $v[q_{0}]$ is the starting vertex of the last branch of $G$. Next we set $w:=\bigl\lvert \{t:q_{0}+1\leq t\leq k\ \&\ s_t\ge 0\}\bigr\rvert$ to be the number of all indices $q_{0}+1\leq t\leq k$ such that the identifier $s_{t}\geq 0$. If $w=0$, define
  \begin{equation}
    \mathcal{E}_{G}(Y_{1},\cdots,Y_{k}):=
    \begin{cases}
      \mathcal{D}^{*}\bigl(Y_{1}Y_{2}\cdots Y_{k}\bigr)                                                                             & \text{ if }q_{0}=1     \\
      \mathcal{E}_{G[q_{0}-1]}\bigl(Y_{1},\cdots,Y_{q_{0}-1}\bigr)\cdot \mathcal{D}^{*}\bigl(Y_{q_{0}}Y_{q_{0}+1}\cdots Y_{k}\bigr) & \text{ if }q_{0}\geq 2
    \end{cases},
  \end{equation}
  where $G[q_{0}-1]\subseteq G$ is the unique order-($q_{0}-1$) sub-genogram of $G$ as defined in \cref{sec:genogram}.

  Otherwise, we write $\{t:q_{0}+1\leq t\leq k\ \&\ s_t\ge 0\}=\{ q_{1},\cdots,q_{w} \}$. Without loss of generality, we suppose that  $q_{0}+1\leq q_{1}<\cdots<q_{w}\leq k$ is increasing. We define
  \begin{equation}\label{eq:defepsilong1}
    \mathcal{E}_{G}(Y_{1},\cdots,Y_{k}):=
    \begin{cases}
      \begin{aligned}
        \mathcal{D}^{*}\bigl(&Y_{1}\cdots Y_{q_{1}-1},
          Y_{q_{1}}\cdots Y_{q_{2}-1},\ldots,\\
          &Y_{q_{w}}\cdots Y_{k}\bigr)
      \end{aligned} & \text{ if }q_{0}=1 \\
      \begin{aligned}
      &\mathcal{E}_{G[q_{0}-1]}\bigl(Y_{1},\cdots,Y_{q_{0}-1}\bigr)\cdot \\
      &\quad \mathcal{D}^{*}\bigl(Y_{q_{0}}\cdots Y_{q_{1}-1},
        Y_{q_{1}}\cdots Y_{q_{2}-1},\ldots,\\
      &\qquad Y_{q_{w}}\cdots Y_{k}\bigr)
      \end{aligned}& \text{ if }q_{0}\geq 2
    \end{cases}.
  \end{equation}

By definition, we can see that $\mathcal{E}_{G}(Y_{1},\cdots,Y_{k})$ is either a $\mathcal{D}^{*}$ term or the product of multiple $\mathcal{D}^{*}$ terms, each of which corresponds to a branch of the rooted tree $(V,E)$.

Taking the genograms $G_{1}$ and $G_{2}$ shown in \cref{fig:expofgenowithlabel} as examples, $\mathcal{E}_{G_{1}}$ and $\mathcal{E}_{G_{2}}$ are provided by
\begin{align*}
  \mathcal{E}_{G_{1}}(Y_{1},\cdots,Y_{7}) & =\mathcal{D}^{*}(Y_{1},Y_{2})\ \mathcal{D}^{*}(Y_{3})\ \mathcal{D}^{*}(Y_{4}Y_{5}\,,\,Y_{6})\ \mathcal{D}^{*}(Y_{7}), \\
  \mathcal{E}_{G_{2}}(Y_{1},\cdots,Y_{7}) & =\mathcal{D}^{*}(Y_{1},Y_{2},Y_{3})\ \mathcal{D}^{*}(Y_{4},Y_{5})\ \mathcal{D}^{*}(Y_{6}Y_{7}).
\end{align*}

Finally, we define $\mathcal{S}(G)$ and $\mathcal{T}_{f}(G)$ with respect to any genogram $G=(V,E,s_{1:k})\in \mathcal{G}(k)$ and function $f\in \mathcal{C}^{k-1}(\mathbb{R})$.  When $k\geq 2$, we also define $\mathcal{U}_{f}(G)$ and $\Delta_f(G)$ by the third display and the following formula.
\begin{align}
  \mathcal{S}(G)
    &:=\sigma^{-k}\sum_{i_{1}\in B_{1}\backslash D_{1}}\cdots
       \sum_{i_{k}\in B_{k}\backslash D_{k}} \notag\\
    &\quad\times \mathcal{E}_{G}\bigl(X_{i_{1}},\ldots,X_{i_{k}}\bigr),
       \label{eq:graphsum1} \\
  \mathcal{T}_{f}(G)
    &:=\sigma^{-k}\sum_{i_{1}\in B_{1}\backslash D_{1}}\cdots
       \sum_{i_{k}\in B_{k}\backslash D_{k}} \notag\\
    &\quad\times \mathcal{E}_{G}\bigl(
       X_{i_{1}},\ldots,X_{i_{k-1}},
       X_{i_{k}}\partial^{k-1}f\bigl(W(D_{k})\bigr)\bigr),
       \label{eq:graphsum2} \\
  \mathcal{U}_{f}(G)
    &:=\sigma^{-(k-1)}\sum_{i_{1}\in B_{1}\backslash D_{1}}\cdots
       \sum_{i_{k-1}\in B_{k-1}\backslash D_{k-1}} \notag\\
    &\quad\times \mathcal{E}_{G}\bigl(
       X_{i_{1}},\ldots,X_{i_{k-1}},\Delta_{f}(G)\bigr).
       \label{eq:graphsum3}
\end{align}
where
\begin{align*}
   & \Delta_{f}(G) :=
  \begin{cases}
    \begin{aligned}
      &\partial^{k-2}f\bigl(W(B_{k})\bigr)\\[-2pt]
      &\quad-\partial^{k-2}f\bigl(W(D_{k})\bigr)
    \end{aligned} & \text{if }u(k)=k \\
    \begin{aligned}
    &\int_{0}^{1}(k-u(k))v^{k-1-u(k)}\Bigl(\\[-2pt]
    &\quad \partial^{k-2}f\bigl(vW(D_{k})+(1-v)W(B_{k})\bigr)\\[-2pt]
    &\qquad-\partial^{k-2}f\bigl(W(D_{k})\bigr)\Bigr)\dif v
    \end{aligned} & \text{if } u(k)\leq k-1
  \end{cases}\!\!
\end{align*}
is called the \emph{adjusted $f$-difference}. Note that $\Delta_{f}(G)$ depends on $G$ and $i_{1},\cdots,i_{k-1}$ through $B_{k},D_{k}$ and $u(k)$, where $k=\lvert G \rvert$. For ease of notation, we do not write out the other dependencies $i_{1},\cdots,i_{k-1}$.

Intuitively, $\mathcal{S}(G)$ is analogous to the $\mathcal{S}$-sums defined for the local dependence case while $\mathcal{T}_{f}(G)$ is analogous to the $\mathcal{T}$-sums. And $\mathcal{U}_{f}(G)$ is defined in a similar spirit to the $\mathcal{R}$-sums as they are both used to handle the remainders. $\Delta_{f}$ is obtained from the integral-form remainders of the Taylor expansions with integral forms of remainders (see \cref{THM:STEP2GRAPH}). Eventually, we would like to expand $\mathcal{T}_{f}(G)$ using $\mathcal{S}(H)$ and $\mathcal{U}_{f}(H)$ for some $H\supseteq G$ as shown in \cref{thm:expansion1graph}.

\subsection{Relating Edgeworth expansion to graphs}\label{sec:proofkeylemma2}
Firstly, we consider how to grow a genogram by following the compatible labeling order of the vertices. Initially, there is only the root and $\lvert G \rvert=1$. We would like $\lvert G \rvert$ to increase from $1$ to $k$ after repeatedly choosing a ``growing'' vertex and adding a child to that vertex. In order to obtain a genogram with the compatible labeling at each step, we observe that the growing vertex needs to be $v[\lvert G \rvert]$ or an ancestor of $v[\lvert G \rvert]$, which is formally stated in \cref{thm:growingvertex}. Moreover, a negative vertex can only be added to $v[\lvert G \rvert]$ because negative vertices do not have siblings as required in the definition of genograms (see \ref{itm:thirdrq} of \cref{thm:rqmofgenogram}).
\begin{lemma}\label{thm:growingvertex}
  Let $(V,E)$ be a rooted tree with $V:=\{ v[1],\cdots, v[k] \}$, whose vertex labels satisfy \ref{itm:dft} of \cref{thm:rqmofgenogram}. Then for any $1\leq j\leq k-1$ either $p(j+1)=j$ or $p(j+1)\in A(j)$ and $j$ is a leaf.
\end{lemma}

\subsection*{Graph operations}
In general, a genogram can be constructed by repeating the following two operations (not necessarily consecutively):
In Operation 1, the minimum over the children of a leaf is understood to be $+\infty$; in particular, any non-negative identifier is admissible when the growing vertex is $v[\lvert G\rvert]$.
We use the convention $\Lambda[0](G)=G$.
\begin{enumerate}
  \item[Operation 1.] \label{itm:grow} $G\rightsquigarrow \Omega [j,s_{\lvert G \rvert+1}](G)$: Fix the growing vertex $v[j]$ to be $v[\lvert G \rvert]$ or an ancestor of $v[\lvert G \rvert]$ that satisfies $\min_{2\leq t\leq \lvert G \rvert:p(t)=j}\{s_{t}\}\geq 1$. Add a non-negative child $v[\lvert G \rvert+1]$ to $v[j]$ and choose $s_{\lvert G \rvert+1}$ to satisfy $0\leq s_{\lvert G \rvert+1}< \min_{2\leq t\leq \lvert G \rvert:p(t)=j}\{s_{t}\}$;
  \item[Operation 2.] \label{itm:glue} $G\rightsquigarrow \Lambda[h](G)$: When $\lvert G\rvert\geq 2$ and $h\geq 1$, add a path of $h$ negative vertices to $v[\lvert G \rvert]$.  In other words, $v[\lvert G \rvert+t]$ is added as the single child of $v[\lvert G \rvert+t-1]$ for $t=1,\cdots,h$, and we set $s_{\lvert G \rvert+1}=\cdots =s_{\lvert G \rvert+h}=-1$.
\end{enumerate}

Take $G_{2}$ shown in \cref{fig:ex2withlabel} as an example. A negative vertex can only be added as a child of $v[7]$ while a non-negative vertex can be added as a child of $v[4]$ with the identifier $0,1$ or $2$, or as a child of $v[7]$ with any non-negative identifier.

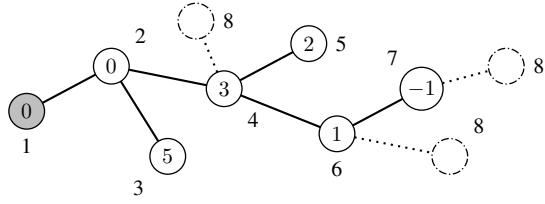
\begin{figure}
  \centering
  \begin{tikzpicture}[scale=0.75,transform shape]
    \begin{scope}[VertexStyle/.append style = {fill=lightgray}]
      \Vertex[x=0.5,y=-1.2,L=$0$]{1}   \extralabel[6pt]{-90}{1}{1}
    \end{scope}
    \Vertex[x=2,y=-0.4,L=$0$]{2}    \extralabel[6pt]{45}{2}{2}
    \Vertex[x=4,y=-0.8,L=$3$]{3}   \extralabel[6pt]{-45}{4}{3}
    \Vertex[x=5.5,y=0,L=$2$]{4}   \extralabel[6pt]{0}{5}{4}
    \Vertex[x=6,y=-1.6,L=$1$]{5}\extralabel[6pt]{-90}{6}{5}
    \Vertex[x=7.5,y=-0.8,L=$-1$]{6}\extralabel[6pt]{135}{7}{6}
    \Vertex[x=3,y=-2,L=$5$]{7}\extralabel[6pt]{-135}{3}{7}
    \begin{scope}[VertexStyle/.append style = {densely dash dot}]
      \Vertex[x=3.5,y=0.4,L=$ $]{8}\extralabel[6pt]{0}{8}{8}
      \Vertex[x=9,y=-0.4,L=$ $]{10}\extralabel[6pt]{0}{8}{10}
    \end{scope}
    \tikzstyle{EdgeStyle}=[]
    \Edge(1)(2)
    \Edge(2)(3)
    \Edge(3)(4)
    \Edge(3)(5)
    \Edge(5)(6)
    \Edge(2)(7)
    \begin{scope}[EdgeStyle/.append style = {dotted}]
      \Edge(3)(8)
      \Edge(6)(10)
    \end{scope}
  \end{tikzpicture}
  \caption{Adding a new vertex to $G_{2}$.}
  \label{fig:adding}
\end{figure}

For clarity we restate \cref{thm:expansion1graphrestate} below.
\begin{theorem}\label{thm:expansion1graph}
  Given a genogram $G$ and an integer $k\geq \lvert G \rvert$, the equation below holds for any $f\in \mathcal{C}^{k-1}(\mathbb{R})$ for which the displayed expectations are finite:
  \begin{equation}\label{eq:expansion1graph}
    \mathcal{T}_{f}(G)
    =\sum_{\substack{
        H\supseteq G:\\
        \lvert H \rvert\leq k,\\
        s_{\lvert G \rvert+1}\geq 0
      }}
    a_{H,G}\,\mathcal{S}(H)\ \mathbb{E} \bigl[\partial^{\lvert H \rvert-1}f(W)\bigr]
    +\sum_{\substack{
        H\in \mathcal{G}(k+1):\\
        H\supseteq G,\\
        s_{\lvert G \rvert+1}\geq 0
      }}
    b_{H,G}\,\mathcal{U}_{f}(H),
  \end{equation}
  where the coefficients $a_{H,G}$ and $b_{H,G}$ are provided by
  \begin{align}
    a_{H,G} & :=\begin{cases}1                                                                                   & \text{if } \lvert H\rvert =\lvert G\rvert      \\
             (-1)^{\gamma_{H}-\gamma_{G}+\tau_{H}-\tau_{G}}\prod_{j=\lvert G \rvert+1}^{\lvert H \rvert}\frac{1}{j+1-u(j,H)} & \text{if } \lvert H\rvert\geq \lvert G\rvert+1\end{cases}, \label{eq:ahg}      \\
    b_{H,G} & :=\begin{cases}(-1)^{\gamma_{H}-\gamma_{G}+\tau_{H}-\tau_{G}+1}                                                                    & \text{if } \lvert H\rvert =\lvert G\rvert +1     \\
             (-1)^{\gamma_{H}-\gamma_{G}+\tau_{H}-\tau_{G}+1}\prod_{j=\lvert G \rvert+1}^{\lvert H \rvert-1}\frac{1}{j+1-u(j,H)} & \text{if } \lvert H\rvert \geq \lvert G\rvert +2\end{cases}.\label{eq:bhg}
  \end{align}
  Here $\gamma_{G}$ denotes the number of leaves on $G$ and $\tau_{G}$ is the number of negative vertices on $G$. It will be used later that for all genograms $H\supseteq G$ with $|H|\ge |G|+1$ we have $a_{H,G}=-\frac{b_{H,G}}{|H|+1-u(|H|,H)}$.
  In the first sum in \eqref{eq:expansion1graph}, the condition $s_{\lvert G\rvert+1}\geq0$ is understood to be vacuous when $H=G$.
\end{theorem}

Let $\mathcal{P}_{0}(k):=\{G\in \mathcal{G}(k): s_{j}\leq 0,\text{ for any }1\leq j\leq k \}$ denote the set of order-$k$ genograms with no positive vertex. Let $\mathcal{G}_{0}(k):=\mathcal{G}(k)\backslash\mathcal{P}_{0}(k)$ denote the set of order-$k$ genograms with at least one positive vertex. Let $\mathcal{P}_{1}(k):=\{G\in \mathcal{G}(k): s_{j}\leq 0, 1\leq j\leq k-1, s_{k}\geq 1 \}$ be the set of order-$k$ genograms where $v[k]$ is the only positive vertex. Note that from the compatible conditions of identifiers, we know any genogram in $\mathcal{P}_{0}(k)$ or $\mathcal{P}_{1}(k)$ has only one branch.

Next we formally present a refined version of \cref{thm:wfwgraph12restate}.
\begin{corollary}\label{thm:wfwgraph12}
  Given $k\geq 1$, suppose that $\mathbb{E}[X_i]=0$ and $\mathbb{E}[\lvert X_i\rvert^{k+1}]<\infty$ for every $i\in T$. Then the equations below hold for any $f\in \mathcal{C}^{k}(\mathbb{R})$ for which the displayed expectations are finite:
  \begin{align}
    \mathbb{E} [Wf(W)]
    ={} & \sum_{j=1}^{k}\frac{\kappa_{j+1}(W)}{j !}\mathbb{E}[\partial^{j}f(W)] \notag\\
        &{}+\sum_{H\in \mathcal{G}(k+2)}b_{H}\mathcal{U}_{f}(H)
        \label{eq:wfwgraph1} \\
    ={} & \sum_{j=1}^{k-1}\frac{\kappa_{j+1}(W)}{j !}\mathbb{E}[\partial^{j}f(W)] \notag\\
        &{}+\frac{\widetilde{\kappa}_{k+1}}{k!}\mathbb{E} [\partial^{k}f(W)]
        +
    \sum_{\substack{
    H\in \mathcal{P}_{0}(k+2)\sqcup                                                                                                                   \\
        \mathcal{P}_{1}(k+2)\sqcup \mathcal{G}_{0}(k+1)
      }}
    b_{H}\mathcal{U}_{f}(H),\label{eq:wfwgraph2}
  \end{align}
  where, for the fixed neighborhood radius $m$, $\widetilde{\kappa}_{k+1}$ and $b_{H}$ are defined as
  \begin{align}
    \widetilde{\kappa}_{k+1}:= & \kappa_{k+1}(W)+\!\!\sum_{H\in \mathcal{G}_{0}(k+1)}\frac{k!\,b_{H}}{k+2-u(k+1)}\mathcal{S}(H),                                        \\
    b_{H}:=                    & \begin{cases}(-1)^{\gamma_{H}+\tau_{H}}                                                  & \text{if } \lvert H\rvert =2     \\
             (-1)^{\gamma_{H}+\tau_{H}}\prod_{j=2}^{\lvert H \rvert-1}\frac{1}{j+1-u(j)} & \text{if } \lvert H\rvert \geq 3\end{cases}.
  \end{align}
  Here $\gamma_{H}$ is the number of leaves on $H$, and $\tau_{H}$ is the number of negative vertices on $H$.
\end{corollary}

\subsection*{Properties of the coefficients $b_{H}$}

We remark that from the definition of $b_{H}$, it is straightforward that for any $H$, we have the bound $\lvert b_{H} \rvert\leq 1$. Moreover, if two genograms $H_{1},H_{2}$ share the same tree structure $(V,E)$ and the set of negative vertices (i.e., $\{ j:s_{j}=-1 \}$), then $b_{H_{1}}=b_{H_{2}}$.

\subsection{Controlling the remainders}\label{sec:remaindercontrol}

From now on, we consider the case when $(X_{i})_{i\in T}$ is an $\alpha$-mixing random field of mean-zero random variables (see \cref{thm:amix}), and proceed to control the terms $\mathcal{S}(H)$ and $\mathcal{U}_{f}(H)$ in \cref{thm:wfwgraph12} to obtain \cref{thm:wfwgraph3} and \cref{thm:wfwgraph4}. \cref{thm:wfwgraph3} provides the remainder control essential for the proof of \cref{thm:hugehugethm}, and \cref{thm:wfwgraph4} is a more refined version of remainder control that has been used in showing \cref{THM:BARBOURGRAPH2}.

\begin{proposition}\label{thm:wfwgraph3}
Let $\bigl(X_{i}\bigr)_{i\in T}$ be a mean-zero random field indexed by $T\subset \mathbb{Z}^{d}$ with $\alpha$-mixing coefficients $(\alpha_{\ell} )$. Given $p>0$ with $k=\lceil p\rceil$ and $\omega=p+1-k$, suppose that \eqref{rr} holds.  If \eqref{moment} holds for some $r\in(k+1+\omega,\infty]$, then for any $f\in \mathcal{C}^{k, \omega }(\mathbb{R})$ and $m\in\mathbb{N}_{+}$ we have
  \begin{equation}\label{eq:wfwgraph3}
    \begin{aligned}
    \mathbb{E} [Wf(W)]=  &\sum_{j=1}^{k}\frac{\kappa_{j+1}(W)}{j !}\mathbb{E}[\partial^{j}f(W)]+
    \mathcal{O}\biggl(\lvert f \rvert_{k, \omega }\lvert T\rvert^{-(k+\omega-1)/2}\Bigl(m^{d(k+\omega)}\\
    &\ +\sum_{\ell=m+1}^{m+1+\lfloor\frac{\lvert T\rvert^{1/d}}{2}\rfloor}\ell^{d(k+\omega) -\omega}\alpha_{\ell}^{(r-k- 1-\omega)/r}\Bigr)\biggr).
    \end{aligned}
  \end{equation}
  Without the preceding strict-moment assumption, there is also an endpoint
  version.  If \eqref{moment} holds with $r=k+1+\omega$, then for any
  $f\in\mathcal C^{k,\omega}(\mathbb R)$ and $m\in\mathbb N_+$ such that
  $\alpha_{m+1}=0$,
  \begin{equation}\label{eq:wfwgraph3endpoint}
    \mathbb E[Wf(W)]
    =\sum_{j=1}^{k}\frac{\kappa_{j+1}(W)}{j!}
      \mathbb E[\partial^jf(W)]
      +\mathcal O\bigl(\lvert f\rvert_{k,\omega}
      \lvert T\rvert^{-(k+\omega-1)/2}m^{d(k+\omega)}\bigr).
  \end{equation}
\end{proposition}

\clearpage
\begin{proposition}\label{thm:wfwgraph4}
  Let $(X_{i})_{i\in T}$ be a mean-zero random field indexed by $T\subset \mathbb{Z}^{d}$ with $\alpha$-mixing coefficients $(\alpha_{\ell})$. Given an integer $k\geq 1$, suppose that there exists $r\in[k+2,\infty]$ such that \eqref{rr} and \eqref{moment} hold. Then for any $f\in \mathcal{C}^{k-1, 1 }(\mathbb{R})\cap \mathcal{C}^{k, 1 }(\mathbb{R})$, $m\in\mathbb{N}_{+}$ and real number $\delta\in [0,1]\cap[0,r-1-k)$ we have
  \begin{align}
      \mathbb{E} [Wf(W)]
      ={} & \sum_{j=1}^{k-1}\frac{\kappa_{j+1}(W)}{j !}
             \mathbb{E}[\partial^{j}f(W)]
        +\frac{\widetilde{\kappa}_{k+1}}{k!}
             \mathbb{E} [\partial^{k}f(W)] \notag\\
      &{}+\mathcal{O}\biggl(
        \lvert f \rvert_{k,1}\lvert T \rvert^{-k/2}m^{d(k+1)}\notag\\
      &\qquad{}+\lvert f \rvert_{k-1,1}^{1-\delta}
        \lvert f \rvert_{k,1}^{\delta}
        \lvert T \rvert^{-(k-1+\delta)/2}m^{dk}\notag\\
      &\qquad\quad{}\times
        \sum_{\ell=m+1}^{m+1+\lfloor\frac{\lvert T\rvert^{1/d}}{2}\rfloor}
        \ell^{d\delta-\delta}\alpha_{\ell}^{(r-k-1-\delta)/r}
        \notag\\
      &\qquad{}+\lvert f \rvert_{k-1,1}\lvert T \rvert^{-(k-1)/2}\notag\\
      &\qquad\quad{}\times
        \sum_{\ell=m+1}^{m+1+\lfloor\frac{\lvert T\rvert^{1/d}}{2}\rfloor}
        \ell^{dk-1}\alpha_{\ell}^{(r-k-1)/r}\biggr).
        \label{eq:wfwgraph4}
  \end{align}
  where $\widetilde{\kappa}_{k+1}$ is the quantity in \cref{thm:wfwgraph12}. It depends on the indexed joint distribution of $(X_{i})_{i\in T}$ and on the chosen radius $m$, but not on $f$ or $\delta$, and it satisfies
  \begin{equation*}
    \bigl\lvert \widetilde{\kappa}_{k+1}-\kappa_{k+1}(W) \bigr\rvert \lesssim\lvert T \rvert^{-(k-1)/2}\sum_{\ell=m+1}^{m+1+\lfloor\frac{\lvert T\rvert^{1/d}}{2}\rfloor}\ell^{dk-1}\alpha_{\ell}^{(r-k- 1)/r}.
  \end{equation*}
\end{proposition}
\section{Proofs of the results in \cref{sec:mixingmainpart}}\label{sec:tech1}

First we introduce the following lemma and prove \cref{thm:growingvertex}.

\begin{lemma}\label{thm:dfslabellemma}
  Let $(V,E)$ be a rooted tree whose vertices are ordered from a depth-first traversal (i.e., the labels satisfy \ref{itm:dft} of \cref{thm:rqmofgenogram}). Suppose $i\notin A(j)$ and $j\notin A(i)$ for some $1\leq i <j\leq k$. Then for any $t$ such that $i\in A(t)$, we have $t<j$.
\end{lemma}

\begin{proof}[Proof of \cref{thm:dfslabellemma}.]

  First let $\prec$ be the strict total order on $V$ such that for any $u,w\in V$, $u\prec w$ if and only if the label of $u$ is smaller than the label of $w$. For any two vertex subsets $U,W\subseteq V$, we denote $U\prec W$ if and only if $u\prec w$ for any $u\in U$ and $w\in W$.
  \begin{claim}
    Let $\prec$ be defined as above. Suppose $u,w\in V$ ($u\neq w$) are siblings in $(V,E)$. Let 
    \begin{align*}
      U&:=\{v:v=u\text{ or }u\text{ is an ancestor of }v\},\\
      W&:=\{v:v=w\text{ or }w\text{ is an ancestor of }v\}.
    \end{align*}
    Then either $U\prec W$ or $W\prec U$.
  \end{claim}

 To prove this we can perform induction on $k=\lvert U \rvert+\lvert W \rvert$. If $k=2$, this is true because $U$ and $W$ each contain one element and $\prec$ is a strict total order on $V$. Now suppose the claim is true for $k$ and consider the case for $k+1$. Without loss of generality, assume $\lvert U \rvert\geq 2$, and thus, $u$ is not a leaf.

  If there is only one leaf in $U$, denoted by $u_{0}$ ($u_{0}\neq u$), then restricting on $V':=V\backslash \{ u_{0}\}$, we get a new rooted tree $(V',E')$ and the order on $V'$ is induced by that on $V$. Then $(V',E')$ also satisfies \ref{itm:dft} of \cref{thm:rqmofgenogram} and $u,w$ are still siblings in $(V',E')$. By inductive hypothesis, either $U\backslash \{ u_{0} \}\prec W$ or $W\prec U\backslash \{ u_{0} \}$. Since $u_{0}$ is the only leaf in $U$, the parent of $u_{0}$, denoted by $u_{1}$ has only one child, which implies that $u_{0}:=\min_{\prec} \{ v: u_{1}\prec v \}$. Thus, we have $U\prec W$ or $W\prec U$.

  If there are at least two leaves in $U$, two of which are denoted by $u_{1}$ and $u_{2}$ ($u,u_{1},u_{2}$ are mutually different). Let $V_{1}':=V\backslash \{ u_{1}\}$ and $V_{2}':=V\backslash \{ u_{2}\}$. Restricting on $V_{1}'$ or $V_{2}'$, we get a new rooted tree $(V_{1}',E_{1}')$ or $(V_{2}',E_{2}')$ with the vertex order on $V_{1}'$ or $V_{2}'$ induced by $\prec$, respectively. By inductive hypothesis on $(V_{1}',E_{1}')$, we get $U\backslash \{ u_{1} \}\prec W$ or $W\prec U\backslash \{ u_{1} \}$. Similarly, we have $U\backslash \{ u_{2} \}\prec W$ or $W\prec U \backslash \{ u_{2} \}$. Since $u\in (U\backslash \{ u_{1}\})\cap (U\backslash \{ u_{2} \})\neq \emptyset$ and $U=(U\backslash \{ u_{1}\})\cup (U\backslash \{ u_{2} \})$, we conclude that $U\prec W$ or $W\prec U$.

  By induction the claim is true. Let $h:=\max A(i)\cap A(j)$, which is the closest common ancestor of $v[i]$ and $v[j]$. Since $i\notin A(j)$, we have $h<i$ and similarly $h<j$. Let $v[i_{0}]$ be the vertex such that $i_{0}\in A(i)$ and $p(i_{0})=h$, and $v[j_{0}]$ be the vertex such that $j_{0}\in A(j)$ and $p(j_{0})=h$. The definition of $h$ implies that $i_{0}\neq j_{0}$. 
  
  Let $U:=\{ v[t]: t=i_{0}\text{ or }i_{0}\in A(t)\}$ and $W:=\{ v[t]:t=j_{0}\text{ or }j_{0}\in A(t) \}$. Note that $v[i_{0}]$ and $v[j_{0}]$ are siblings. Using the claim above we get $U\prec W$ or $W\prec U$. Noticing that $v[i]\in U$, $v[j]\in W$ and $i<j$, we conclude that $U\prec W$. If $i\in A(t)$, then we have $v[t]\in U$ and $t<j$.
  \hfill
\end{proof}

\begin{proof}[Proof of \cref{thm:growingvertex}.]
  Firstly, if $p(j+1)\neq j$, then $v[j]$ is a leaf. Otherwise, suppose $v[j]$ has a child $v[h]$ where $h>j+1$. Then $p(h)>p(j+1)$ contradicts \ref{itm:dft} of \cref{thm:rqmofgenogram}.

  Suppose $p(j+1)\neq j$ and $p(j+1)\notin A(j)$. \ref{itm:dft} of \cref{thm:rqmofgenogram} implies that $p(j+1)\leq j$. Thus, we know $p(j+1)<j$. Since $v[j]$ is a leaf, we have $j\notin A(p(j+1))$. Thus, by \cref{thm:dfslabellemma}, for any $t$ such that $p(j+1)\in A(t)$, we have $t<j$. In particular, let $t=j+1$ and we get $j+1<j$. This is a contradiction. Thus, either $p(j+1)= j$ or $v[j]$ is a leaf and $p(j+1)\in A(j)$.\hfill
\end{proof}

From the derivation of \eqref{eq:split2}, we note that the expansion of $\mathcal{T}_{f}(G)$ is typically achieved by {constructing a telescoping sum} followed by the Taylor expansions. We will see later that the first operation of growing genograms corresponds to the {telescoping sum argument}, and the second one corresponds to the Taylor expansion, which will be formalized in \cref{THM:STEP1GRAPH,THM:STEP2GRAPH}.

\begin{lemma}\label{THM:STEP1GRAPH}
  Given an integer $\ell\geq 1$, an order-$\ell$ genogram $G$, and a function $f\in\mathcal{C}^{\ell-1}(\mathbb{R})$, suppose that all the displayed expectations are finite. Let $\Omega[j,s](G)$ ($1\leq j\leq \ell$, $s\geq 0$) be the genogram obtained by growing a child from the vertex $v[j]$, as defined in \cref{itm:grow}. Then we have
  \begin{align*}
    &\mathcal{T}_{f}(G)-\mathcal{S}(G)\ \mathbb{E} \bigl[\partial^{\ell-1}f(W)\bigr]\\
    = & -\sum_{s\geq 0}\mathcal{U}_{f}\bigl(\Omega [\ell, s](G)\bigr)
    +\sum_{\substack{j\in A(\ell)\cup\{\ell\}:\\ s_{j}\geq 1}}\sum_{s=0}^{s_{j}-1}\mathcal{U}_{f}\bigl(\Omega[p(j),s](G)\bigr), \\
    = & -\sum_{\substack{
    H\in \mathcal{G}(\ell+1):                                                                          \\
    H\supseteq G,                                                                                      \\
    s_{\ell+1}\geq 0,                                                                                  \\
        p(\ell+1,H)=\ell
      }}
    \mathcal{U}_{f}(H)
    +\sum_{\substack{
    H\in \mathcal{G}(\ell+1):                                                                          \\
    H\supseteq G,                                                                                      \\
    s_{\ell+1}\geq 0,                                                                                  \\
        p(\ell+1,H)<\ell
      }}
    \mathcal{U}_{f}(H),
  \end{align*}
  where $p(j)$ is the label of the parent of $v[j]$, $s_{j}$ is the identifier, and $\mathcal{G}(k)$ is the set of all order-$k$ genograms.
\end{lemma}
Note that \cref{THM:STEP1GRAPH} generalizes the idea of expanding to a telescoping sum while deriving \eqref{eq:split2} and will be proven in \cref{sec:lemma4}.

\begin{lemma}\label{THM:STEP2GRAPH}
  Given two integers $\ell\geq 2$, $k\geq 0$, an order-$\ell$ genogram $G$, and a function $f\in\mathcal{C}^{\ell+k-1}(\mathbb{R})$, suppose that all the displayed expectations are finite. Set $\Lambda[0](G):=G$, and for $1\leq j\leq k+1$ let $\Lambda [j](G)$ be the genogram obtained by gluing a path of $j$ negative vertices to $v[\ell]$, as defined in \cref{itm:glue}. Then we have
  \begin{align}\label{eq:step2graph}
    \mathcal{U}_{f}(G)=&\sum_{j=0}^{k}(-1)^{j+1}\frac{(\ell-u(\ell))!}{(j+1+ \ell-u(\ell))!}\mathcal{T}_{f}\bigl(\Lambda[j](G)\bigr)\\*
    &\ +(-1)^{k+1}\frac{(\ell- u(\ell))!}{(k+1+\ell-u(\ell))!}\mathcal{U}_{f}\bigl(\Lambda[k+1](G)\bigr),\nonumber
  \end{align}
  where $u(j)$ is given by \eqref{eq:defuj}.
\end{lemma}
As we have pointed out, \cref{THM:STEP2GRAPH} is a direct consequence of the Taylor expansion with the integral form of remainders, and the proof is also provided in \cref{sec:lemma4}.

\begin{proof}[Proof of \cref{thm:expansion1graph}.]
  For convenience, denote $\ell:=\lvert G \rvert$. We prove by performing induction on $k$ ($k\geq \ell$).

  If $k=\ell$, then we note that the set $\{H\supseteq G: |H|\le k\}=\{G\}$ only contains the genogram $G$. Moreover, in \eqref{eq:ahg} we set $a_{G,G}=1$. Therefore, we obtain that
  \begin{align*}
     & \sum_{\substack{
    H\supseteq G:                                                               \\
    \lvert H \rvert\leq k,                                                      \\
        s_{\lvert G \rvert+1}\geq 0
      }}
    a_{H,G}\,\mathcal{S}(H)\ \mathbb{E}
      \bigl[\partial^{\lvert H \rvert-1}f(W)\bigr]\nonumber\\
      &\quad+\sum_{\substack{
    H\in \mathcal{G}(k+1):                                                      \\
    H\supseteq G,                                                               \\
        s_{\lvert G \rvert+1}\geq 0
      }}
    b_{H,G}\,\mathcal{U}_{f}(H)                                                 \\
      =&\mathcal{S}(G)\ \mathbb{E} \bigl[\partial^{\ell-1}f(W)\bigr]
    +\sum_{\substack{
    H\in \mathcal{G}(\ell+1):                                                   \\
    H\supseteq G,                                                               \\
        s_{\ell+1}\geq 0
      }}
    b_{H,G}\,\mathcal{U}_{f}(H).
  \end{align*}

  Next let $H\in \mathcal{G}(\ell+1)$ be a genogram of order $\ell+1$ such that $H\supseteq G$ and such that  $s_{\ell+1}\geq 0$. Note that $s_{\ell+1}\geq 0$ implies that $\tau_{H}=\tau_{G}$. To calculate $\gamma_{H}-\gamma_{G}$, we check that the parent of $v[\ell+1]$ is either $v[\ell]$ (i.e., $p(\ell+1,H)=\ell$) or an ancestor of $v[\ell]$. If its parent is $v[\ell]$ then the number of leaves in $H$, denoted by $\gamma_{H}$, is the same as that of $G$, and hence by \eqref{eq:bhg} $b_{H,G}=-1$. Otherwise, the number of leaves increases by exactly one: $\gamma_{H}=\gamma_{G}+1$, which implies $b_{H,G}=1$.
  Thus, \eqref{eq:expansion1graph} reduces to
  \begin{equation*}
    \sum_{\substack{
        H\in \mathcal{G}(\ell+1):\\
        H\supseteq G,\\
        s_{\ell+1}\geq 0
      }}
    b_{H,G}\,\mathcal{U}_{f}(H)
    =-    \sum_{\substack{
        H\in \mathcal{G}(\ell+1):\\
        H\supseteq G,\\
        s_{\ell+1}\geq 0,\\
        p(\ell+1,H)=\ell
      }}
    \mathcal{U}_{f}(H)
    + \sum_{\substack{
        H\in \mathcal{G}(\ell+1):\\
        H\supseteq G,\\
        s_{\ell+1}\geq 0,\\
        p(\ell+1,H)<\ell
      }}
    \mathcal{U}_{f}(H),
  \end{equation*}
  Using \cref{THM:STEP1GRAPH} we directly obtain that
  \begin{equation*}
    \mathcal{T}_{f}(G)
    =\mathcal{S}(G)\ \mathbb{E} \bigl[\partial^{\ell-1}f(W)\bigr]
    -\sum_{\substack{
        H\in \mathcal{G}(\ell+1):\\
        H\supseteq G,\\
        s_{\ell+1}\geq 0,\\
        p(\ell+1,H)=\ell
      }}
    \mathcal{U}_{f}(H)
    +\sum_{\substack{
        H\in \mathcal{G}(\ell+1):\\
        H\supseteq G,\\
        s_{\ell+1}\geq 0,\\
        p(\ell+1,H)<\ell
      }}
    \mathcal{U}_{f}(H),
  \end{equation*}
  and the desired result is proven.

  Now supposing the statement is true for $k$, we will establish that the desired result also holds for $k+1$.
  By inductive hypothesis we have
  \begin{equation}\label{eq:109put1}
    \mathcal{T}_{f}(G)
    =\sum_{\substack{
        H\supseteq G:\\
        \lvert H \rvert\leq k,\\
        s_{\ell+1}\geq 0
      }}
    a_{H,G}\,\mathcal{S}(H)\ \mathbb{E} \bigl[\partial^{\lvert H \rvert-1}f(W)\bigr]
    +\sum_{\substack{
        H\in \mathcal{G}(k+1):\\
        H\supseteq G,\\
        s_{\ell+1}\geq 0
      }}
    b_{H,G}\,\mathcal{U}_{f}(H).
  \end{equation}
  For any $H\in\mathcal{G}(k+1)$, by \cref{THM:STEP2GRAPH} we have
  \begin{equation}\label{eq:109put2}
    \mathcal{U}_{f}(H)=-\frac{1}{k+2-u(k+1,H)}\bigl(\mathcal{T}_{f}(H)+\mathcal{U}_{f}\bigl(\Lambda[1](H)\bigr)\bigr).
  \end{equation}

  Combining \eqref{eq:109put2} with \eqref{eq:109put1}, we get
  \begin{align}\label{eq:109put3}
    \mathcal{T}_{f}(G)
    = & \sum_{\substack{
    H\supseteq G:                                                                 \\
    \lvert H \rvert\leq k,                                                        \\
        s_{\ell+1}\geq 0
      }}
    a_{H,G}\,\mathcal{S}(H)\ \mathbb{E} \bigl[\partial^{\lvert H \rvert-1}f(W)\bigr]
    -\sum_{\substack{
    H\in \mathcal{G}(k+1):                                                        \\
    H\supseteq G,                                                                 \\
        s_{\ell+1}\geq 0
      }}
    \frac{b_{H,G}}{k+2-u(k+1,H)}\mathcal{T}_{f}(H)                       \\
      & \ -\sum_{\substack{
    H\in \mathcal{G}(k+1):                                                        \\
    H\supseteq G,                                                                 \\
        s_{\ell+1}\geq 0
      }}
    \frac{b_{H,G}}{k+2-u(k+1,H)}\mathcal{U}_{f}\bigl(\Lambda[1](H)\bigr)\nonumber \\
 \overset{(a)}{   =} & \sum_{\substack{
    H\supseteq G:                                                                 \\
    \lvert H \rvert\leq k,                                                        \\
        s_{\ell+1}\geq 0
      }}
    a_{H,G}\,\mathcal{S}(H)\ \mathbb{E}
      \bigl[\partial^{\lvert H \rvert-1}f(W)\bigr]\nonumber\\
      &\quad+\!\!\!\sum_{\substack{
    H\in \mathcal{G}(k+1):                                                        \\
    H\supseteq G,                                                                 \\
        s_{\ell+1}\geq 0
      }}\!\!
    a_{H,G}\mathcal{T}_{f}(H)\nonumber\\
      &\quad+\!\!\!\!\sum_{\substack{
    K\in \mathcal{G}(k+2):                                                        \\
    K\supseteq G,                                                                 \\
    s_{\ell+1}\geq 0,                                                             \\
        s_{k+2}=-1
      }}\!\!
    b_{K,G}\mathcal{U}_{f}\bigl(K\bigr).\nonumber
  \end{align}
  Note that $s_{k+2}=-1$ implies that $\tau_{K}=\tau_{H}+1$ and $\gamma_{K}=\gamma_{H}$ (since $p(k+2,K)=k+1$). Moreover, equality $(a)$ is due to the fact that for any $H\in \mathcal{G}(k+1)$ we have $a_{H,G}=-\frac{b_{H,G}}{k+2-u(k+1,H)}$, and that for any $K\in\mathcal{G}(k+2)$ with $s_{k+2}=-1$ we have
  $b_{K,G}=(-1)^{\gamma_{K}-\gamma_{H}+\tau_{K}-\tau_{H}}\frac{b_{H,G}}{k+2-u(k+1,H)}=-\frac{b_{H,G}}{k+2-u(k+1,H)}$.

Next by another application of \cref{THM:STEP1GRAPH}, we get that
  \begin{equation}
    \mathcal{T}_{f}(H)-\mathcal{S}(H)\ \mathbb{E} \bigl[\partial^{k}f(W)\bigr]
    =-\sum_{\substack{
        K\in \mathcal{G}(k+2):\\
        K\supseteq H,\\
        s_{k+2}\geq 0,\\
        p(k+2,K)=k+1
      }}
    \mathcal{U}_{f}(K)
    +\sum_{\substack{
        K\in \mathcal{G}(k+2):\\
        K\supseteq H,\\
        s_{k+2}\geq 0,\\
        p(k+2,K)<k+1}}
    \mathcal{U}_{f}(K).
  \end{equation}
  Combining this with \eqref{eq:109put3}, we get
  \begin{equation}\label{eq:greatputin}
    \begin{aligned}
      \mathcal{T}_{f}(G)
      = & \sum_{\substack{
      H\supseteq G:                                        \\
      \lvert H \rvert\leq k,                               \\
          s_{\ell+1}\geq 0
        }}
      a_{H,G}\,\mathcal{S}(H)\ \mathbb{E} \bigl[\partial^{\lvert H \rvert-1}f(W)\bigr]
      +\sum_{\substack{
      H\in \mathcal{G}(k+1):                               \\
      H\supseteq G,                                        \\
          s_{\ell+1}\geq 0
        }}
      a_{H,G}\mathcal{S}(H)\ \mathbb{E} [\partial^{k}f(W)] \\
        & \ -\!\!\sum_{\substack{
      K\in \mathcal{G}(k+2):                               \\
      K\supseteq H\supseteq G,                             \\
      s_{\ell+1}\geq 0,                                    \\
      s_{k+2}\geq 0,                                       \\
          p(k+2,K)=k+1
        }}\!\!
      a_{H,G}\mathcal{U}_{f}(K)
      +\!\!\sum_{\substack{
      K\in \mathcal{G}(k+2):                               \\
      K\supseteq H\supseteq G,                             \\
      s_{\ell+1}\geq 0,                                    \\
      s_{k+2}\geq 0,                                       \\
          p(k+2,K)<k+1}}\!\!
      a_{H,G}\mathcal{U}_{f}(K)+\!\!\sum_{\substack{
      K\in \mathcal{G}(k+2):                               \\
      K\supseteq G,                                        \\
      s_{\ell+1}\geq 0,                                    \\
          s_{k+2}=-1
        }}\!\!
      b_{K,G}\mathcal{U}_{f}\bigl(K\bigr),
    \end{aligned}
  \end{equation}
  where $H$ is the order-$(k+1)$ sub-genogram of $K$.

  Now we simplify \eqref{eq:greatputin}.
  For $K\in \mathcal{G}(k+2)$ such that $p(k+2,K)=k+1$, we have $\gamma_{K}=\gamma_{H}$. And $s_{k+2}\geq 0$ implies that $\tau_{K}=\tau_{H}$. Thus, $b_{K,G}=(-1)^{\gamma_{K}-\gamma_{H}+\tau_{K}-\tau_{H}}\frac{b_{H,G}}{k+2-u(k+1,H)}=\frac{b_{H,G}}{k+2-u(k+1,H)}=-a_{H,G}$. For $K\in \mathcal{G}(k+2)$ such that $p(k+2,K)<k+1$, we have $\gamma_{K}=\gamma_{H}+1$.
  Thus, $b_{K,G}=(-1)^{\gamma_{K}-\gamma_{H}+\tau_{K}-\tau_{H}}\frac{b_{H,G}}{k+2-u(k+1,H)}=-\frac{b_{H,G}}{k+2-u(k+1,H)}=a_{H,G}$.
  \eqref{eq:greatputin} reduces to
  \begin{align}\label{eq:finalkplus1}
    \mathcal{T}_{f}(G)
    = & \sum_{\substack{
    H\supseteq G:                                                 \\
    \lvert H \rvert\leq k,                                        \\
        s_{\ell+1}\geq 0
      }}
    a_{H,G}\,\mathcal{S}(H)\ \mathbb{E}
      \bigl[\partial^{\lvert H \rvert-1}f(W)\bigr]\nonumber\\
      &\quad+\sum_{\substack{
    H\in \mathcal{G}(k+1):                                        \\
    H\supseteq G,                                                 \\
        s_{\ell+1}\geq 0
      }}
    a_{H,G}\mathcal{S}(H)\ \mathbb{E} [\partial^{k}f(W)] \\
      & \ +\sum_{\substack{
    K\in \mathcal{G}(k+2):                                        \\
    K\supseteq G,                                                 \\
    s_{\ell+1}\geq 0,                                             \\
    s_{k+2}\geq 0,                                                \\
        p(k+2,K)=k+1
      }}
    b_{K,G}\mathcal{U}_{f}(K)\nonumber\\
      &\quad+\sum_{\substack{
    K\in \mathcal{G}(k+2):                                        \\
    K\supseteq G,                                                 \\
    s_{\ell+1}\geq 0,                                             \\
    s_{k+2}\geq 0,                                                \\
        p(k+2,K)<k+1}}
    b_{K,G}\mathcal{U}_{f}(K)\nonumber\\
      &\quad+\sum_{\substack{
    K\in \mathcal{G}(k+2):                                        \\
    K\supseteq G,                                                 \\
    s_{\ell+1}\geq 0,                                             \\
        s_{k+2}=-1
      }}
    b_{K,G}\mathcal{U}_{f}\bigl(K\bigr)\nonumber                  \\
    = & \sum_{\substack{
    H\supseteq G:                                                 \\
    \lvert H \rvert\leq k+1,                                      \\
        s_{\ell+1}\geq 0
      }}
    a_{H,G}\,\mathcal{S}(H)\ \mathbb{E}
      \bigl[\partial^{\lvert H \rvert-1}f(W)\bigr]\nonumber\\
      &\quad+\sum_{\substack{
    K\in \mathcal{G}(k+2):                                        \\
    K\supseteq G,                                                 \\
        s_{\ell+1}\geq 0
      }}
    b_{K,G}\,\mathcal{U}_{f}(K).
  \end{align}
  The last equality is from the observation that
  \begin{align*}
    &\{ H\supseteq G: \lvert H \rvert\leq k+1,s_{\ell+1}\geq 0 \}\\
    &\qquad= \{ H\supseteq G:\lvert H \rvert\leq k,s_{\ell+1}\geq 0 \}\sqcup \{ H\in \mathcal{G}(k+1):H\supseteq G,s_{\ell+1}\geq 0 \}, \\
    &\{ K\in \mathcal{G}(k+2):K\supseteq G,s_{\ell+1}\geq 0 \}\\
    &\qquad= \{ K\in \mathcal{G}(k+2):K\supseteq G,s_{\ell+1}\geq 0,s_{k+2}\geq 0,p(k+2,K)=k+1 \}\sqcup                                 \\
     & \qquad\quad\  \{ K\in \mathcal{G}(k+2):K\supseteq G,s_{\ell+1}\geq 0,s_{k+2}\geq 0,p(k+2,K)<k+1 \}\sqcup                              \\
    & \qquad\quad\ \{ K\in \mathcal{G}(k+2):K\supseteq G,s_{\ell+1}\geq 0,s_{k+2}=-1\},\nonumber
  \end{align*}
  where $\sqcup$ denotes the disjoint union of sets.
  Note that \eqref{eq:finalkplus1} is precisely \eqref{eq:expansion1graph} for the case $k+1$. By induction \cref{thm:expansion1graph} is proven.\hfill
\end{proof}

Before we show \cref{thm:wfwgraph12}, the following lemma is needed:
\begin{lemma}\label{thm:expanduniquegraph}
  Given a genogram $G$ and an integer $k\geq \lvert G \rvert$, suppose there exist two sets of constants that only depend on $G$ and the joint distribution of $(X_{i})_{i\in T}$, $(Q_{\lvert G \rvert},\cdots,Q_{k})$ and $(Q_{\lvert G \rvert}',\cdots,Q_{k}')$, which satisfy that for any polynomial $f$ of degree at most $k-1$,
  \begin{equation*}
    \mathcal{T}_{f}(G)=\sum_{j=\lvert G \rvert}^{k}Q_{j}\mathbb{E} [\partial^{j-1}f(W)]=\sum_{j=\lvert G \rvert}^{k}Q_{j}'\mathbb{E} [\partial^{j-1}f(W)].
  \end{equation*}
  Then $Q_{j}=Q_{j}'$ for any $\lvert G \rvert\leq j\leq k$.
\end{lemma}

\begin{proof}[Proof of \cref{thm:expanduniquegraph}.]
  We prove the lemma by contradiction.

  Let $j$ be the smallest number such that $Q_{j}\neq Q_{j}'$. Since the coefficients do not depend on $f$, choose $f(x)=c x^{j-1}$ with $c\neq 0$. The terms indexed by numbers smaller than $j$ cancel by the minimality of $j$, while all terms indexed by numbers larger than $j$ vanish. Hence
  \[
    c(j-1)!Q_j=c(j-1)!Q_j',
  \]
  which is a contradiction. Therefore, $Q_{j}= Q_{j}'$ for any $\lvert G \rvert\leq j\leq k$.
\end{proof}

Now we proceed with the proof of \cref{thm:wfwgraph12}.
\begin{proof}[Proof of \cref{thm:wfwgraph12}.]
  We write $G_0=(1,\emptyset, 0)$ to be the genogram consisting only of the root. By definition $$\mathcal{T}_{f}(G_0)=\mathbb{E} [Wf(W)].$$ As $G_0$ is a genogram of order $1$, applying \cref{thm:expansion1graph} we have that for any $f\in \mathcal{C}^{k}(\mathbb{R})$,
  \begin{equation}\label{1234}
    \mathcal{T}_{f}(G_0)=\mathbb{E} [Wf(W)]=\sum_{j=1}^{k+1}Q_{j}\mathbb{E} [\partial^{j-1} f(W)]+\sum_{H\in \mathcal{G}(k+2)}b_{H,G_0}\mathcal{U}_{f}(H),
  \end{equation}
  where $Q_1=\mathcal{S}(G_0)=0$ and for all $j\ge 2$ we take
  \begin{equation}\label{eq:exprqj}
    Q_{j}=\sum_{H\in \mathcal{G}(j):H\supseteq G_0, s_{2}\geq 0}a_{H,G_0}\mathcal{S}(H)\overset{(a)}{=}-\sum_{H\in \mathcal{G}(j)}\frac{b_{H,G_0}}{j+1-u(j,H)}\mathcal{S}(H),
  \end{equation}
  where $a_{H,G_0}$ and $b_{H,G_0}$ are defined in \cref{thm:expansion1graph}, and where to obtain $(a)$ we used the fact that for all genograms $H$ of order larger than $2$, we have $a_{H,G_0}=-\frac{b_{H,G_0}}{|H|+1-u(|H|,H)}$. The expression above shows that $Q_{j}$ only depends on the joint distribution of $(X_{i})_{i\in T}$. Furthermore, we see that by definition $b_{H,G_{0}}$ is identical to the $b_{H}$ defined in the corollary.

  For any $H\in \mathcal{G}(k+2)$ and polynomial $f$ of degree at most $k$, $\partial^{k}f$ is a constant. Therefore, $\Delta_{f}(H)=0$, which directly implies that
  \begin{equation*}
    \begin{aligned}
      \mathcal{U}_{f}(H)
      ={}&\sigma^{-(k+1)}
        \sum_{i_{1}\in B_{1}\backslash D_{1}}\cdots
        \sum_{i_{k+1}\in B_{k+1}\backslash D_{k+1}}\\
      &\quad\times
        \mathcal{E}_{H}\bigl(X_{i_{1}},\ldots,X_{i_{k+1}},
        \Delta_{f}(H)\bigr)=0.
    \end{aligned}
  \end{equation*}
  Therefore, by combining this with \eqref{1234} we obtain that
  \begin{equation}\label{eq:bracketf}
    \mathbb{E} [Wf(W)]=\sum_{j=1}^{k+1}Q_{j}\mathbb{E} [\partial^{j-1} f(W)].
  \end{equation}
  On the other hand, for any random variable, the moments $(\mu_{j})_{j\geq 0}$ and cumulants $(\kappa_{j})_{j\geq 0}$, provided that they exist, are connected through the following relations \citep{smith1995recursive}:
  \begin{equation}\label{eq:lemcumueq}
    \mu_{n}=\sum_{j=1}^{n}\binom{n-1}{j-1}\kappa_{j}\mu_{n-j}.
  \end{equation} 
  If we choose $f(x)= x^{j}$ where $j\in [k]$, then by \eqref{eq:lemcumueq} we have
  \begin{align*}
      & \mathbb{E} [Wf(W)]=\mu_{j+1}(W)
    =\sum_{h=1}^{j+1}\binom{j}{h-1}\kappa_{h}(W)\mu_{j+1-h}(W) \\
    = & \sum_{h=0}^{j}\binom{j}{h}\kappa_{h+1}(W)\mu_{j-h}(W)
    =\sum_{h=1}^{k+1}\frac{\kappa_{h}(W)}{(h-1)!}\mathbb{E} [\partial^{h-1} f(W)].
  \end{align*}
  Any polynomial $f$ of degree at most $k$ can be written as $f(x)=\sum_{j=0}^{k}a_{j}x^{j}$. By linearity of expectations, we know
  \begin{equation*}
    \mathbb{E} [Wf(W)]=\sum_{j=1}^{k+1}\frac{\kappa_{j}(W)}{(j-1) !}\mathbb{E} [\partial^{j-1} f(W)].
  \end{equation*}
  Compare this to \eqref{eq:bracketf} and apply \cref{thm:expanduniquegraph}. We conclude that $Q_{j}=\kappa_{j}(W)/(j-1)!$ for any $j\in [k+1]$. 
  Thus, for any $f\in \mathcal{C}^{k}(\mathbb{R})$, we have shown
  \begin{align*}
    \mathbb{E} [Wf(W)]= & \sum_{j=2}^{k+1}Q_{j}\mathbb{E} [\partial^{j-1} f(W)]+\sum_{H\in \mathcal{G}(k+2)}b_{H}\mathcal{U}_{f}(H)                  \\
    =                   & \sum_{j=1}^{k}Q_{j+1}\mathbb{E} [\partial^{j} f(W)]+\sum_{H\in \mathcal{G}(k+2)}b_{H}\mathcal{U}_{f}(H)                    \\
    =                   & \sum_{j=1}^{k}\frac{\kappa_{j+1}(W)}{j !}\mathbb{E}[\partial^{j}f(W)]+\sum_{H\in \mathcal{G}(k+2)}b_{H}\mathcal{U}_{f}(H).
  \end{align*}
  Since $f\in\mathcal{C}^{k}(\mathbb{R})\subseteq \mathcal{C}^{k-1}(\mathbb{R})$, we also have
  \begin{equation}\label{eq:1011putin1}
    \mathbb{E} [Wf(W)]
    =\sum_{j=1}^{k-1}\frac{\kappa_{j+1}(W)}{j !}\mathbb{E}[\partial^{j}f(W)]+\sum_{H\in \mathcal{G}(k+1)}b_{H}\mathcal{U}_{f}(H).
  \end{equation}
  As $\mathcal{P}_0(k+1)\subseteq \mathcal{G}(k+1)$ we can decompose $\sum_{H\in \mathcal{G}(k+1)}b_{H}\mathcal{U}_{f}(H)$ into two sums as
  \begin{equation}\label{eq:1011putin2}
    \sum_{H\in \mathcal{G}(k+1)}b_{H}\mathcal{U}_{f}(H)=\sum_{H\in \mathcal{P}_{0}(k+1)}b_{H}\mathcal{U}_{f}(H)+\sum_{H\in \mathcal{G}_{0}(k+1)\setminus \mathcal{P}_0(k+1)}b_{H}\mathcal{U}_{f}(H).
  \end{equation}

  For $H\in \mathcal{P}_{0}(k+1)$, applying \cref{THM:STEP2GRAPH} we obtain
  \begin{align}\label{eq:afterfix}
     \sum_{H\in \mathcal{P}_{0}(k+1)}b_{H}\mathcal{U}_{f}(H)
     =&-\sum_{H\in \mathcal{P}_{0}(k+1)}\frac{b_{H}}{k+2-u(k+1)}\mathcal{T}_{f}(H)\\
     &\ -\sum_{H\in \mathcal{P}_{0}(k+1)}\sum_{\substack{K\in \mathcal{G}(k+2):\\K\supsetneq H,\\ s_{k+2}=-1}}\frac{b_{H}}{k+2-u(k+1)}\mathcal{U}_{f}(K) \nonumber\\
     \overset{(*)}{=}&-\sum_{H\in \mathcal{P}_{0}(k+1)}\frac{b_{H}}{k+2-u(k+1)}\mathcal{T}_{f}(H)+\sum_{\substack{K\in \mathcal{P}_{0}(k+2):\\ s_{k+2}=-1}}b_{K}\mathcal{U}_{f}(K)\nonumber
  \end{align}
  where in $(*)$ we use the fact that $b_{K}=(-1)^{\gamma_{K}-\gamma_{H}+\tau_{K}-\tau_{H}}\frac{b_{H}}{k+2-u(k+1)}=-\frac{b_{H}}{k+2-u(k+1)}$ since $s_{k+2}=-1$ implies that $\gamma_{K}=\gamma_{H}$ and $\tau_{K}=\tau_{H}+1$.

  Noting that for any $H\in \mathcal{P}_{0}(k+1)$, if an order-$(k+2)$ genogram $K$ satisfies that $K\supsetneq H$, then we have $p(k+2,K)=k+1$ since we have that $s_{j}\leq 0$ for any $j\in [k+1]$. Thus, \cref{THM:STEP1GRAPH} implies that
  \begin{align*}
    &\sum_{H\in \mathcal{P}_{0}(k+1)}\bigl(\mathcal{T}_{f}(H)-\mathcal{S}(H)\ \mathbb{E} [\partial^{k}f(W)]\bigr)
    = -\sum_{H\in \mathcal{P}_{0}(k+1)}\sum_{\substack{K\in \mathcal{G}(k+2):\\ K\supsetneq H,\\ s_{k+2}\geq 0}}\mathcal{U}_{f}(K)\\
    =&-\sum_{\substack{K\in \mathcal{P}_{0}(k+2):\\ s_{k+2}=0}}\mathcal{U}_{f}(K)-\sum_{K\in \mathcal{P}_{1}(k+2)}\mathcal{U}_{f}(K).
  \end{align*}
  Combining this with \eqref{eq:afterfix} we get
  \begin{align}\label{eq:1011putin3}
    &\sum_{H\in \mathcal{P}_{0}(k+1)}b_{H}\mathcal{U}_{f}(H)\\
   =&-\sum_{H\in \mathcal{P}_{0}(k+1)}\frac{b_{H}}{k+2-u(k+1)}\mathcal{S}(H)\ \mathbb{E} \bigl[\partial^{k}f(W)\bigr]+\sum_{\substack{K\in \mathcal{P}_{0}(k+2):\\ s_{k+2}=-1}}b_{K}\mathcal{U}_{f}(K) \nonumber\\*
   &\ +\sum_{\substack{K\in \mathcal{P}_{0}(k+2):\\ s_{k+2}=0}}
      \frac{b_{K[k+1]}}{k+2-u(k+1,K[k+1])}\mathcal{U}_{f}(K)\nonumber\\*
   &\ +\sum_{K\in \mathcal{P}_{1}(k+2)}
      \frac{b_{K[k+1]}}{k+2-u(k+1,K[k+1])}\mathcal{U}_{f}(K)\nonumber\\
   \overset{(a)}{=}&-\sum_{H\in \mathcal{P}_{0}(k+1)}\frac{b_{H}}{k+2-u(k+1)}\mathcal{S}(H)\ \mathbb{E} \bigl[\partial^{k}f(W)\bigr]+\sum_{\substack{K\in \mathcal{P}_{0}(k+2):\\ s_{k+2}=-1}}b_{K}\mathcal{U}_{f}(K)\nonumber \\*
   &\ +\sum_{\substack{K\in \mathcal{P}_{0}(k+2):\\ s_{k+2}=0}}b_{K}\mathcal{U}_{f}(K)+\sum_{K\in \mathcal{P}_{1}(k+2)}b_{K}\mathcal{U}_{f}(K)\nonumber\\
   =& -\sum_{H\in \mathcal{P}_{0}(k+1)}\frac{b_{H}}{k+2-u(k+1)}\mathcal{S}(H)\ \mathbb{E} \bigl[\partial^{k}f(W)\bigr]+\sum_{K\in \mathcal{P}_{0}(k+2)\sqcup \mathcal{P}_{1}(k+2)}b_{K}\mathcal{U}_{f}(K).\nonumber
 \end{align}
 Again $(a)$ is obtained by taking $H=K[k+1]$, the unique order-$(k+1)$ sub-genogram of $K$, and checking that $\gamma_{K}=\gamma_{H}$ and $\tau_{K}=\tau_{H}$. Thus,
 $b_{K}=b_{K[k+1]}/\bigl(k+2-u(k+1,K[k+1])\bigr)$.

  Since we have already established that $Q_{j}=\kappa_{j}(W)/(j-1)!$, by \eqref{eq:exprqj}, the following equation holds that
  \begin{equation*}
    \frac{\kappa_{k+1}(W)}{k!}=Q_{k+1}=-\sum_{H\in \mathcal{G}(k+1)}\frac{b_{H}}{k+2-u(k+1)}\mathcal{S}(H).
  \end{equation*}
  Thus, we get 
  \begin{equation}\label{eq:1011putin4}
    \begin{aligned}
    \frac{\widetilde{\kappa}_{k+1}}{k!}=&\frac{\kappa_{k+1}(W)}{k!}+\sum_{H\in \mathcal{G}_{0}(k+1)}\frac{b_{H}}{k+2-u(k+1)}\mathcal{S}(H)\\
    =&-\sum_{H\in \mathcal{P}_{0}(k+1)}\frac{b_{H}}{k+2-u(k+1)}\mathcal{S}(H).
    \end{aligned}
  \end{equation}
  Combining \eqref{eq:1011putin2}, \eqref{eq:1011putin3}, and \eqref{eq:1011putin4} with \eqref{eq:1011putin1}, we obtain \eqref{eq:wfwgraph2}.
\hfill
\end{proof}

\begin{lemma}\label{THM:REMAINDERCTRL1234}
  Let $(X_{i})_{i\in T}$ be a mean-zero random field indexed by $T\subsetneq \mathbb{Z}^{d}$ ($\lvert T\rvert<\infty$). Given $k\geq 1$, assume that the non-degeneracy condition holds that $\liminf_{\lvert T \rvert\to \infty}\sigma^{2}/\lvert T \rvert>0$. For each assertion below, suppose that $\sup_{i\in T}\lVert X_i\rVert_r\leq M<\infty$ at the indicated moment order $r$, where $M$ is uniform over the collection of fields as $\lvert T\rvert\to\infty$.

  For any $r\in(k+1,\infty]$ and $m\in\mathbb{N}_{+}$, we have
    \begin{equation}
      \!\!\!\!\!\!\sum_{H \in \mathcal{G}_{0}(k+1)}\!\!\!\!\bigl\lvert\mathcal{S}(H) \bigr\rvert
      \lesssim \lvert T \rvert^{-(k-1)/2}\sum_{\ell=m+1}^{m+1+\lfloor\frac{\lvert T\rvert^{1/d}}{2}\rfloor}\ell^{dk-1}\alpha_{\ell}^{(r-k-1)/r}.\label{eq:remainderctrl1}
    \end{equation}

  Given $\omega\in (0,1]$, for any $r\in(k+\omega,\infty]$, $f\in \mathcal{C}^{k-1, \omega}(\mathbb{R})$, and $m\in\mathbb{N}_{+}$, we have
    \begin{equation}
      \biggl\lvert\sum_{H \in \mathcal{G}_{0}(k+1)}b_{H}\mathcal{U}_{f}(H) \biggr\rvert
      \lesssim \lvert f \rvert_{k-1,\omega} \lvert T \rvert^{-(k+\omega-2)/2}\sum_{\ell=m+1}^{m+1+\lfloor\frac{\lvert T\rvert^{1/d}}{2}\rfloor}\ell^{d(k+\omega-1)-\omega }\alpha_{\ell}^{(r-k- \omega  )/r}.\label{eq:remainderctrl2}
    \end{equation}
  At the endpoint $r=k+\omega$, if $\alpha_{m+1}=0$, the left-hand side of
  \eqref{eq:remainderctrl2} is zero.
  For any $r\in[k+\omega,\infty]$, $f\in \mathcal{C}^{k-1, \omega}(\mathbb{R})$, and $m\in\mathbb{N}_{+}$, we have
    \begin{equation}
      \biggl\lvert\sum_{H\in \mathcal{P}_{0}(k+1)}b_{H}\mathcal{U}_{f}(H) \biggr\rvert
      \lesssim \lvert f \rvert_{k-1, \omega } \lvert T \rvert^{-(k+\omega-2)/2}m^{d(k+ \omega -1 )}.\label{eq:remainderctrl3}
    \end{equation}
    Moreover, if $k\geq 2$, for any $\delta\in [0,1]$, $r\in(k+\delta,\infty]$, $f\in \mathcal{C}^{k-2,1}(\mathbb{R})\cap \mathcal{C}^{k-1,1}(\mathbb{R})$, and $m\in\mathbb{N}_{+}$, we have
    \begin{equation}
      \begin{aligned}
        \biggl\lvert\sum_{H\in \mathcal{P}_{1}(k+1)}
          b_{H}\mathcal{U}_{f}(H) \biggr\rvert
        \lesssim{}& \lvert f \rvert_{k-2,1}^{1-\delta}
          \lvert f \rvert_{k-1,1}^{\delta}
          \lvert T \rvert^{-(k-2+\delta)/2}m^{d(k-1)}\\
        &\times\sum_{\ell=m+1}^{m+1+\lfloor\frac{\lvert T\rvert^{1/d}}{2}\rfloor}
          \ell^{d\delta-\delta}\alpha_{\ell}^{(r-k-\delta)/r}.
      \end{aligned}\label{eq:remainderctrl4}
    \end{equation}
\end{lemma}
  
  The proof of this lemma is long and technical, so we move it to \cref{sec:lemma5}. Now combining \cref{thm:wfwgraph12} and \cref{THM:REMAINDERCTRL1234}, we have \cref{thm:wfwgraph3,thm:wfwgraph4}.

\begin{proof}[Proofs of \cref{thm:wfwgraph3,thm:wfwgraph4}.]
  Using the preceding corollary, we have the expansion
  \begin{equation*}
    \mathbb{E} [Wf(W)]
    = \sum_{j=1}^{k}\frac{\kappa_{j+1}(W)}{j !}\mathbb{E}[\partial^{j}f(W)]+\sum_{H\in \mathcal{G}(k+2)}b_{H}\mathcal{U}_{f}(H).
  \end{equation*}
  At the endpoint $r=k+1+\omega$, if $\alpha_{m+1}=0$, the
  $\mathcal G_0(k+2)$ sum is zero by the endpoint assertion in
  \cref{THM:REMAINDERCTRL1234}, while the $\mathcal P_0(k+2)$ sum is
  bounded by \eqref{eq:remainderctrl3}.  Hence
  \begin{equation*}
    \mathbb E[Wf(W)]
    =\sum_{j=1}^{k}\frac{\kappa_{j+1}(W)}{j!}
      \mathbb E[\partial^jf(W)]
      +\mathcal O\bigl(\lvert f\rvert_{k,\omega}
      \lvert T\rvert^{-(k+\omega-1)/2}m^{d(k+\omega)}\bigr),
  \end{equation*}
  which proves the endpoint version of \cref{thm:wfwgraph3}.
  In the strict range $r>k+1+\omega$, \cref{THM:REMAINDERCTRL1234} with
  $m\in\mathbb{N}_{+}$ gives the upper bound
  \begin{align*}
    &\biggl\lvert\sum_{H\in \mathcal{G}(k+2)}b_{H}\mathcal{U}_{f} (H)\biggr\rvert\\
    \leq &
    \biggl\lvert\sum_{H\in \mathcal{P}_{0}(k+2)}b_{H}\mathcal{U}_{f}(H) \biggr\rvert + \biggl\lvert\sum_{H \in \mathcal{G}_{0}(k+2)}b_H\mathcal{U}_{f}(H) \biggr\rvert                                                                                                                                              \\
    \lesssim                                                                         & \lvert f \rvert_{k, \omega } \lvert T \rvert^{-(k+\omega-1)/2}\Bigl(m^{d(k+\omega)}+\sum_{\ell=m+1}^{m+1+\lfloor\frac{\lvert T\rvert^{1/d}}{2}\rfloor}\ell^{d(k+\omega)-\omega }\alpha_{\ell}^{(r-k- 1-\omega )/r}\Bigr).
  \end{align*}
  Again by \cref{thm:wfwgraph12}, we have
  \begin{equation*}
    \begin{aligned}
    \mathbb{E} [Wf(W)]
    = &\sum_{j=1}^{k-1}\frac{\kappa_{j+1}(W)}{j !}\mathbb{E}[\partial^{j}f(W)] \\
    &\ +\frac{\widetilde{\kappa}_{k+1}}{k!}\mathbb{E} [\partial^{k}f(W)]+
    \sum_{H\in \mathcal{G}_{0}(k+1)}b_{H}\mathcal{U}_{f}(H)+\!\!\!\!\!\!\!\!\sum_{H\in \mathcal{P}_{0}(k+2)\cup \mathcal{P}_{1}(k+2)}\!\!\!\!\!\!\!\!b_{H}\mathcal{U}_{f}(H),
    \end{aligned}
  \end{equation*}
  where $\widetilde{\kappa}_{k+1}$ satisfies that
  \begin{align*}
    \widetilde{\kappa}_{k+1}:=\kappa_{k+1}(W)+\sum_{H\in \mathcal{G}_{0}(k+1)}\frac{k!\,b_{H}}{k+2-u(k+1)}\mathcal{S}(H).
  \end{align*}
  Note that $\lvert b_{H} \rvert\leq 1$. By \cref{THM:REMAINDERCTRL1234} we get
  \begin{align*}
    \bigl\lvert \widetilde{\kappa}_{k+1}-\kappa_{k+1}(W) \bigr\rvert \leq & k!
    \sum_{H \in \mathcal{G}_{0}(k+1)}\!\!\!\!\bigl\lvert\mathcal{S}(H) \bigr\rvert\lesssim \lvert T \rvert^{-(k-1)/2}\sum_{\ell=m+1}^{m+1+\lfloor\frac{\lvert T\rvert^{1/d}}{2}\rfloor}\ell^{dk-1}\alpha_{\ell}^{(r-k-1)/r}.
  \end{align*}
  By \cref{THM:REMAINDERCTRL1234} with $\omega=1$ we have
  \begin{align*}
             & \biggl\lvert \sum_{H\in \mathcal{P}_{0}(k+2)}\!\!\!\!\!\!b_{H}\mathcal{U}_{f}(H)\biggr\rvert+\biggl\lvert \sum_{H\in  \mathcal{P}_{1}(k+2)}\!\!\!\!\!\!b_{H}\mathcal{U}_{f}(H) \biggr\rvert+\biggl\lvert\sum_{H\in \mathcal{G}_{0}(k+1)}\!\!\!\!\!\!b_{H}\mathcal{U}_{f}(H)\biggr\rvert                                                                                                                \\
    \lesssim & \lvert f \rvert_{k, 1  } \lvert T \rvert^{-k/2}m^{d(k+ 1  )}+ \lvert f \rvert_{k-1, 1  }^{1-\delta}\lvert f \rvert_{k,1}^{\delta} \lvert T \rvert^{-(k-1+\delta)/2}m^{dk}\!\!\!\!\!\!\sum_{\ell=m+1}^{m+1+\lfloor\frac{\lvert T\rvert^{1/d}}{2}\rfloor}\!\!\!\!\!\!\ell^{d\delta-\delta }\alpha_{\ell}^{(r-k-1-\delta)/r}\\*
    &\ +\lvert f \rvert_{k-1, 1  }\lvert T \rvert^{-(k-1)/2}\!\!\!\!\!\!\sum_{\ell=m+1}^{m+1+\lfloor\frac{\lvert T\rvert^{1/d}}{2}\rfloor}\!\!\!\!\!\!\ell^{dk-1 }\alpha_{\ell}^{(r-k- 1  )/r}.
  \end{align*}
  Therefore, \eqref{eq:wfwgraph4} is proven.\hfill
\end{proof}
\section{Proofs of Lemmas~\ref{THM:STEP1GRAPH} and~\ref{THM:STEP2GRAPH}}\label{sec:lemma4}

\begin{proof}[Proof of \cref{THM:STEP1GRAPH}.]
  Consider the label set of the positive vertices among $v[\ell]$ and its ancestors, and let
  \[
    z:=\bigl\lvert\{j\in A(\ell)\cup\{\ell\}:s_j\geq1\}\bigr\rvert.
  \]
  If $z\geq1$, we write
  \begin{equation*}
    \{j\in A(\ell)\cup\{\ell\}:s_j\geq1\}=\{r_1,r_2,\ldots,r_z\},
  \end{equation*}
  where $2\leq r_1<r_2<\cdots<r_z\leq\ell$.

  Then consider all possible genograms constructed by adding a non-negative child to a vertex of $G$. On one hand, \ref{itm:fourthrq} of \cref{thm:rqmofgenogram} implies that it is only possible to add such a child to $v[\ell]$ or $v[p(r_{j})]$ for some $1\leq j\leq z$ (provided that $z\geq 1$), and if a genogram is obtained by adding $v[\ell+1]$ as a non-negative child of $v[p(r_{j})]$, then $s_{\ell+1}\leq s_{r_{j}}-1$ because $\ell+1>r_{j}$.

  On the other hand, we show that for $1\leq j\leq z$,
  \begin{center}
    $s_{r_{j}}\leq s_{t}$ for any $t\leq \ell$ such that $p(t)=p(r_{j})$.
  \end{center}
  Supposing this is true and $v[\ell+1]$ is added as a child of $v[p(r_j)]$ with $s_{\ell+1}\leq s_{r_j}-1$, then $s_{\ell+1}$ is smaller than the identifiers of all $v[\ell+1]$'s siblings. Hence \cref{thm:rqmofgenogram} holds and $\Omega[p(r_j),s_{\ell+1}](G)$ is indeed a genogram with the compatible labeling.

  In fact, $p(t)=p(r_{j})$ ($t\leq \ell$) implies that $t\notin A(r_{j})$ and $r_{j}\notin A(t)$. If $t>r_{j}$, by \cref{thm:dfslabellemma}, for any $t'$ such that $r_{j}\in A(t')$, we have $t'<t$. In particular, let $t'=\ell$. Then $\ell<t$, contradicting $t\leq \ell$. Thus, $t\leq r_{j}$. The \cref{itm:fourthrq} of \cref{thm:rqmofgenogram} implies that $s_{t}\geq s_{r_{j}}$.

  Thus, we have shown
  \begin{equation}\label{eq:defmcah}
    \begin{aligned}
      \mathcal{H}:= & \bigl\{ H\in \mathcal{G}(\ell+1): H\supseteq G, s_{\ell+1}\geq 0 \bigr\}                                                     \\
      =             & \bigl\{ \Omega[p(r_{j}),s](G):1\leq j\leq z,0\leq s\leq  s_{r_{j}}-1 \bigr\}\sqcup \bigl\{ \Omega[\ell,s](G):s\geq 0\bigr\},
    \end{aligned}
  \end{equation}
  where $\sqcup$ denotes the disjoint union of sets.
  This directly implies that
  \begin{equation}
  \begin{aligned} &-\sum_{s\geq 0}\mathcal{U}_{f}\bigl(\Omega [\ell, s](G)\bigr)
    +\sum_{\substack{j\in A(\ell)\cup\{\ell\}:\\s_{j}\geq 1}}\sum_{s=0}^{s_{j}-1}\mathcal{U}_{f}\bigl(\Omega[p(j),s](G)\bigr)\\
    = & -\sum_{\substack{
    H\in \mathcal{G}(\ell+1): \\
    H\supseteq G,             \\
    s_{\ell+1}\geq 0,         \\
        p(\ell+1)=\ell
      }}
    \mathcal{U}_{f}(H)
    +\sum_{\substack{
    H\in \mathcal{G}(\ell+1): \\
    H\supseteq G,             \\
    s_{\ell+1}\geq 0,         \\
        p(\ell+1)<\ell
      }}
    \mathcal{U}_{f}(H).
  \end{aligned}
\end{equation}

  Therefore, to conclude the proof we only need to show that this is also equal to $$\mathcal{T}_{f}(G)-\mathcal{S}(G)\ \mathbb{E} \bigl[\partial^{\ell-1}f(W)\bigr].$$ 
  For convenience, we list all elements of $\mathcal{H}$ (see \eqref{eq:defmcah}) in a sequence. If $z=0$, we write
  \begin{equation*}
    H_{t}:=\Omega[\ell,t-1](G).
  \end{equation*}
  Otherwise, let
  \[
    H_t:=
    \begin{cases}
      \Omega[p(r_j),s-1](G)
        & \text{if }t=\displaystyle\sum_{i=1}^{j-1}s_{r_i}+s,
          \ 1\leq j\leq z,\ 1\leq s\leq s_{r_j}, \\
      \Omega[\ell,s-1](G)
        & \text{if }t=\displaystyle\sum_{i=1}^{z}s_{r_i}+s,
          \ s\geq1.
    \end{cases}
  \]
  In other words, $H_t$ is a genogram of order $\ell+1$ and the sequence $(H_t)_{t\geq 1}$ can be enumerated as
  \begin{align*}
    (H_{t})_{t\geq 1}:\quad & \Omega [p(r_{1}),0](G)\ ,\quad\Omega [p(r_{1}),1](G)\ ,\quad\cdots,\quad\Omega [p(r_{1}),s_{r_{1}}-1](G)\ , \\
                            & \Omega [p(r_{2}),0](G)\ ,\quad\Omega [p(r_{2}),1](G)\ ,\quad\cdots,\quad\Omega [p(r_{2}),s_{r_{2}}-1](G)\ , \\
                            & \cdots\cdots,                                                                                               \\
                            & \Omega [p(r_{z}),0](G)\ ,\quad\Omega [p(r_{z}),1](G)\ ,\quad\cdots,\quad\Omega [p(r_{z}),s_{r_{z}}-1](G)\ , \\
                            & \Omega [\ell,0](G)\ ,\quad\Omega [\ell,1](G)\ ,\quad\cdots
  \end{align*}
  {We write $$c_{0}:=\begin{cases} 1 &\text{ if } z=0\\s_{r_{1}}+\cdots+s_{r_{z}}+1&\text{ if }z\geq 1\end{cases},$$ and remark that $H_{t}=\Omega [\ell,t-c_{0}](G)$ for $t\geq c_{0}$.}

  For any $t\ge 1$, we note that $H_t$ is an order-($\ell+1$) genogram. We write $B_j(H_t)$ and $D_j(H_t)$ respectively the outer and inner constraints of $i_{j}$ with respect to $H_t$. We remark that as $H_t[\ell]=G$, this directly implies that $B_j(H_t)=B_j(G)=B_{j}$ and $D_j(H_t)=D_j(G)=D_j$ for all $j\le \ell$.

  Let $i_1,\cdots, i_{\ell}$ be indices in $i_{j}\in B_{j}\backslash D_{j}$ for $1\leq j\leq \ell$. We note that the sets $B_{\ell+1}(H_t)$ and $D_{\ell+1}(H_t)$ will depend on the value of $t\ge 1$.

  Firstly, in $H_1$ we remark that by definition \eqref{eq:defgj}, if $z\geq 1$, the vertices $v[\ell+1]$ and $v[r_{1}]$ have the same parent and that the identifier of $\ell+1$ is $0$. {This implies that if $z\geq 1$,
  $$
    D_{\ell+1}(H_{1})=D_{g(\ell+1,H_{1})}=D_{g(r_{1},H_{1})}=D_{g(r_{1},G)}=D_{g(r_{1})}=D_{1}=\emptyset.
  $$
  If $z=0$, the progenitor of $v[\ell+1]$ is either $v[\ell]$ or $v[1]$. In both cases, the inner constraint is empty, i.e., $D_{\ell}=D_{1}=\emptyset$. Thus,
  \begin{equation*}
    D_{\ell+1}(H_{1})=D_{g(\ell+1,H_{1})}=\emptyset.
  \end{equation*}
  Note that here we have used $g(j,H)$ to denote the progenitor label of $v[j]$ with respect to the genogram $H$. Throughout the proof, $H$ is omitted only if $H=G$.}

  Next we establish that for all $t\ge 1$ the following holds: $B_{\ell+1}(H_{t})=D_{\ell+1}(H_{t+1})$.

    If $z=0$, the result is directly implied by the definitions \eqref{eq:defbj} and \eqref{eq:defdj} as we have
    \begin{align*}
      B_{\ell+1}(H_{t})
       & =N^{(t-1)}\bigl(i_{h}:h\in A(\ell+1,H_{t})\bigr)\cup D_{g(\ell+1,H_{t})}     \\
       & =N^{(t-1)}\bigl(i_{h}:h\in A(\ell+1,H_{t+1})\bigr)\cup D_{g(\ell+1,H_{t+1})} \\
       & =D_{\ell+1}(H_{t+1}).
    \end{align*}
    Note that $A(j,H)$ is used to denote the label set of $v[j]$'s ancestors with respect to the genogram $H$. Again $H$ is omitted only if $H=G$.

  If $z\geq 1$, we remark that according to the values of $t$ the relationship between the genograms $H_{t+1}$ and $H_t$ will be different. To make this clear, we distinguish 4 different cases according to the values of $t$:
  \begin{enumerate}[label=(\alph*)]
    \item \label{itm:lemmacase1} $t=\sum_{i=1}^{j-1}s_{r_i}+k+1$ for $1\leq j\leq z$ and $0\leq k\leq s_{r_{j}}-2$, in which case we observe that
    $$
      B_{\ell+1}(H_t)=B_{\ell+1}\bigl(\Omega[p(r_{j}),k](G)\bigr),\quad D_{\ell+1}(H_{t+1})=D_{\ell+1}\bigl(\Omega[p(r_{j}),k+1](G)\bigr);
    $$
    \item \label{itm:lemmacase2} $t=\sum_{i=1}^{z}s_{r_i}+k+1$ for $k\ge 0$, where we have
    $$
      B_{\ell+1}(H_t)=B_{\ell+1}\bigl(\Omega[\ell,k](G)\bigr),\quad D_{\ell+1}(H_{t+1})=D_{\ell+1}\bigl(\Omega[\ell,k+1](G)\bigr);
    $$
    \item \label{itm:lemmacase3} $t=\sum_{i=1}^{j}s_{r_i}$ for $1\leq j\leq z-1$, where we observe that
    $$
      B_{\ell+1}(H_t)=B_{\ell+1}\bigl(\Omega[p(r_{j}),s_{r_{j}}-1](G)\bigr),\quad D_{\ell+1}(H_{t+1})=D_{\ell+1}\bigl(\Omega[p(r_{j+1}),0](G)\bigr);
    $$
    \item \label{itm:lemmacase4} $t=\sum_{i=1}^{z}s_{r_i}$, in which case we have
    $$
      B_{\ell+1}(H_t)=B_{\ell+1}\bigl(\Omega[p(r_{z}),s_{r_{z}}-1](G)\bigr),\quad D_{\ell+1}(H_{t+1})=D_{\ell+1}\bigl(\Omega[\ell,0](G)\bigr).
    $$
  \end{enumerate}
 
    Again cases \ref{itm:lemmacase1} and \ref{itm:lemmacase2} are directly implied by the definitions \eqref{eq:defbj} and \eqref{eq:defdj}. Indeed, for all $1\leq j\leq z$, we have
    \begin{align*}
      B_{\ell+1}(H_{t})
       & =N^{(k)}\bigl(i_{h}:h\in A(\ell+1,H_{t})\bigr)\cup D_{g(\ell+1,H_{t})}     \\
       & =N^{(k)}\bigl(i_{h}:h\in A(\ell+1,H_{t+1})\bigr)\cup D_{g(\ell+1,H_{t+1})} \\
       & =D_{\ell+1}(H_{t+1}).
    \end{align*}
    We now prove cases \ref{itm:lemmacase3} and \ref{itm:lemmacase4}. For case \ref{itm:lemmacase3}, let $t=\sum_{i=1}^j s_{r_i}$ for a given $j\le z-1$. Since in $H_t$, the vertices $v[\ell+1]$ and $v[r_j]$ are siblings, we have $A(\ell+1,H_t)=A(r_j)$ and $g(\ell+1,H_t)=g(r_j)$. In $H_{t+1}$, the vertices $v[\ell+1]$ and $v[r_{j+1}]$ are siblings, and the closest positive ancestor of $v[r_{j+1}]$ is $v[r_j]$. Hence $g(\ell+1,H_{t+1})=r_j$.

    For case \ref{itm:lemmacase4}, take $j=z$ and $t=\sum_{i=1}^z s_{r_i}$. The vertex $v[r_z]$ is the closest positive vertex among $v[\ell]$ and its ancestors, so $g(\ell+1,H_{t+1})=r_z$. Equipped with these observations, we remark in either case that
    \begin{align*}
      B_{\ell+1}(H_t)=&N^{(s_{r_j}-1)}\bigl(i_{t}:t\in A(\ell+1,H_{t})\bigr)\cup D_{g(\ell+1,H_{t})}\\
      =&N^{(s_{r_{j}}-1)}\bigl(i_{t}:t\in A(r_{j})\bigr)\cup D_{g(r_{j})},
    \end{align*}
    On the other hand, we also remark that as the identifier of $\ell+1$ in $H_{t+1}$ is $0$, we have
    \begin{gather*}
      D_{\ell+1}(H_{t+1})=D_{g(\ell+1,H_{t+1})}=D_{r_{j}}=N^{(s_{r_{j}}-1)}\bigl(i_{t}:t\in A(r_{j})\bigr)\cup D_{g(r_{j})}.
    \end{gather*}
    Thus, $B_{\ell+1}(H_t)=D_{\ell+1}(H_{t+1})$.

    Therefore, we have established that $B_{\ell+1}(H_{t})=D_{\ell+1}(H_{t+1})$ for any $t\geq 1$.

    Since the index set $T$ of the random field is finite, there exists a finite number $c_1>c_0$ such that $N^{(c_1-c_0-1)}\bigl(i_t:t\in A(\ell)\text{ or }t=\ell\bigr)=T$. Then
    \begin{align*}
      B_{\ell+1}(H_{c_1-1})
      &\supseteq N^{(c_1-c_0-1)}
        \bigl(i_t:t\in A(\ell+1,H_{c_1-1})\bigr)\\
      &=N^{(c_1-c_0-1)}
        \bigl(i_{t}:t\in A(\ell)\text{ or }t=\ell\bigr)=T.
    \end{align*}
    Thus, $B_{\ell+1}(H_{c_1-1})=D_{\ell+1}(H_{c_1})=T$. The definitions also give $D_{\ell+1}(H_{c_0})=D_\ell$, and hence
    \begin{align*}
      \emptyset=D_{\ell+1}(H_{1})
      &\subseteq B_{\ell+1}(H_{1})=D_{\ell+1}(H_2)
        \subseteq\cdots\subseteq D_{\ell+1}(H_{c_0})=D_\ell\\
      &\subseteq\cdots\subseteq B_{\ell+1}(H_{c_1-1})
        =D_{\ell+1}(H_{c_1})=T.
    \end{align*}

  We prove that for $c_{0}\leq t\leq c_{1}-1$,
  \begin{equation}\label{eq:egphi1}
  \begin{aligned}
    &-\mathcal{E}_{H_{t}}\bigl(X_{i_{1}},\cdots,X_{i_{\ell}}\ ,\ \Delta_{f}(H_{t})\bigr)\\
    = & \mathcal{E}_{G}\Bigl(X_{i_{1}},\cdots\ ,\ X_{i_{\ell}}\bigl(\partial^{\ell-1}f\bigl(W(D_{\ell+1}(H_{t}))\bigr)-\partial^{\ell-1}f\bigl(W(B_{\ell+1}(H_{t}))\bigr)\bigr)\Bigr) \\
      & -\mathcal{E}_{G}\bigl(X_{i_{1}},\cdots,X_{i_{\ell}}\bigr)\Bigl(\mathbb{E} \bigl[\partial^{\ell-1}f\bigl(W(D_{\ell+1}(H_{t}))\bigr)\bigr]-\mathbb{E} \bigl[\partial^{\ell-1}f\bigl(W(B_{\ell+1}(H_{t}))\bigr)\bigr]\Bigr);
  \end{aligned}
\end{equation}
  and similarly that for $1\leq t\leq c_{0}-1$,
  \begin{equation}\label{eq:egphi2}
    \begin{aligned}
    &\mathcal{E}_{H_{t}}\bigl(X_{i_{1}},\cdots,X_{i_{\ell}}\ ,\ \Delta_{f}(H_{t})\bigr)\\
    =&\mathcal{E}_{G}\bigl(X_{i_{1}},\cdots,X_{i_{\ell}}\bigr)\Bigl(\mathbb{E} \bigl[\partial^{\ell-1}f\bigl(W(B_{\ell+1}(H_{t}))\bigr)\bigr]-\mathbb{E} \bigl[\partial^{\ell-1}f\bigl(W(D_{\ell+1}(H_{t}))\bigr)\bigr]\Bigr).
    \end{aligned}
  \end{equation}

  Note that by definition $\mathcal{E}_{G}$ is the product of $\mathcal{D}^{*}$ factors. Let
  \begin{equation*}
    q_{0}:=\sup \bigl\{ j: j=1\text{ or }p(j)\neq j-1\text{ for }2\leq j\leq \ell \bigr\}.
  \end{equation*}
  
    Intuitively, $v[q_{0}]$ is the starting vertex of the last branch of $(V,E)$. Now set $w=\bigl\lvert\{t:q_{0}+1\leq t\leq \ell\ \&\  s_t\ge 0\}\bigr\rvert$. If $w\geq 1$, write
    $\{q_{1},\cdots,q_{w}\}=\{t:q_{0}+1\leq t\leq \ell\ \&\  s_t\ge 0\}$ with $q_{0}+1\leq q_{1}<\cdots<q_{w}\leq \ell$. Define
    \begin{equation*}
      q_{*}:=\begin{cases}q_{0}&\text{if }w=0,\\q_{w}&\text{if }w\geq 1.\end{cases}
    \end{equation*}
    By definition, the last factor in $\mathcal{E}_{G}(X_{i_{1}},\cdots,X_{i_{\ell}})$ is given by
    \begin{equation*}
      \begin{cases}
        \mathcal{D}^{*}\bigl(X_{i_{q_{0}}}\cdots X_{i_{\ell}}\bigr)                                                    & \text{ if }w=0     \\
        \mathcal{D}^{*}\bigl(X_{i_{q_{0}}}\cdots X_{i_{q_{1}-1}}\ ,\ \cdots\ ,\ X_{i_{q_{w}}}\cdots X_{i_{\ell}}\bigr) & \text{ if }w\geq 1
      \end{cases}.
    \end{equation*}
    And the last factor in $\mathcal{E}_{G}(X_{i_{1}},\cdots,X_{i_{\ell-1}},X_{i_{\ell}}\partial^{\ell-1}f\bigl(W(D_{\ell})\bigr))$ is
    \begin{equation*}
      \begin{cases}
        \mathcal{D}^{*}\bigl(X_{i_{q_{0}}}\cdots X_{i_{\ell}}\partial^{\ell-1}f\bigl(W(D_{\ell})\bigr)\bigr)                                                                                                & \text{ if }w=0     \\
        \begin{aligned}
          \mathcal{D}^{*}\bigl(&X_{i_{q_{0}}}\cdots X_{i_{q_{1}-1}},\ldots,
          X_{i_{q_{(w-1)}}}\cdots X_{i_{q_{w}-1}},\\
          &X_{i_{q_{w}}}\cdots X_{i_{\ell}}
          \partial^{\ell-1}f\bigl(W(D_{\ell})\bigr)\bigr)
        \end{aligned} & \text{ if }w\geq 1
      \end{cases}.
    \end{equation*}
    For convenience, in this proof we temporarily denote
    $$
      \mathfrak{D}^{*}(\,\cdot\,):=
      \begin{cases}
        \mathcal{D}^{*}(\,\cdot\,)                                                                                                        & \text{ if }w=0     \\
        \mathcal{D}^{*}\bigl(X_{i_{q_{0}}}\cdots X_{i_{q_{1}-1}}\ ,\ \cdots\ ,\ X_{i_{q_{(w-1)}}}\cdots X_{i_{q_{w}-1}}\ ,\ \cdot\ \bigr) & \text{ if }w\geq 1
      \end{cases}.
    $$

  Since $s_{\ell+1}\geq 0$, we have
  \begin{equation*}
    \Delta_{f}(H_{t})=\partial^{\ell-1}f\bigl(W(B_{\ell+1}(H_{t}))\bigr)-\partial^{\ell-1}f\bigl(W(D_{\ell+1}(H_{t}))\bigr).
  \end{equation*}

  For $c_{0}\leq t\leq c_{1}-1$, we derive that
  \begin{align*}
    &-\mathfrak{D}^{*}\Bigl(X_{i_{q_{*}}}\cdots X_{i_{\ell}}\ ,\ \Delta_{f}(H_{t})\Bigr)\\
    =                & \mathfrak{D}^{*}\Bigl(X_{i_{q_{*}}}\cdots X_{i_{\ell}}\ ,\ \partial^{\ell-1}f\bigl(W(D_{\ell+1}(H_{t}))\bigr)-\partial^{\ell-1}f\bigl(W(B_{\ell+1}(H_{t}))\bigr)\Bigr)                                                          \\
    \overset{(*)}{=} & \mathfrak{D}^{*}\Bigl(X_{i_{q_{*}}}\cdots X_{i_{\ell}}\mathcal{D}\bigl(\partial^{\ell-1}f\bigl(W(D_{\ell+1}(H_{t}))\bigr)-\partial^{\ell-1}f\bigl(W(B_{\ell+1}(H_{t}))\bigr)\bigr)\Bigr)                                        \\
    =                & \mathfrak{D}^{*}\Bigl(X_{i_{q_{*}}}\cdots X_{i_{\ell}}\bigl(\partial^{\ell-1}f\bigl(W(D_{\ell+1}(H_{t}))\bigr)-\partial^{\ell-1}f\bigl(W(B_{\ell+1}(H_{t}))\bigr)\bigr)\Bigr)                                                   \\
                     & \ -\mathfrak{D}^{*}\bigl(X_{i_{q_{*}}}\cdots X_{i_{\ell}}\bigr)\Bigl(\mathbb{E} \bigl[\partial^{\ell-1}f\bigl(W(D_{\ell+1}(H_{t}))\bigr)\bigr]-\mathbb{E} \bigl[\partial^{\ell-1}f\bigl(W(B_{\ell+1}(H_{t}))\bigr)\bigr]\Bigr).
  \end{align*}
  Here the equality $(*)$ is implied by \cref{thm:dastaddterm}.
  Noticing that $v[\ell+1]$ is the child of $v[\ell]$, we combine this with the definition of the $\mathcal{E}_{G}$ operator and obtain \eqref{eq:egphi1}.

  For $1\leq t\leq c_{0}-1$, we note that
  \begin{align*}
    &\mathfrak{D}^{*}\bigl(X_{i_{q_{*}}}\cdots X_{i_{\ell}}\bigr)\ \mathcal{D}^{*}(\Delta_{f}(H_{t}))\\
    =&\mathfrak{D}^{*}\bigl(X_{i_{q_{*}}}\cdots X_{i_{\ell}}\bigr)\Bigl(\mathbb{E} \bigl[\partial^{\ell-1}f\bigl(W(B_{\ell+1}(H_{t}))\bigr)\bigr]-\mathbb{E} \bigl[\partial^{\ell-1}f\bigl(W(D_{\ell+1}(H_{t}))\bigr)\bigr]\Bigr).
  \end{align*}
  Since $v[\ell+1]$ is added as a child of $v[p(r_j)]$ for some $1\leq j\leq z$, the last vertex $v[\ell]$ remains a leaf in $H_t$. We obtain \eqref{eq:egphi2} by applying the definition of the $\mathcal{E}_{G}$ operator again.

  Finally, using the definitions of $\mathcal{T}_{f}(G)$ and $\mathcal{S}(G)$, we derive
  \begin{align*}
     & \mathcal{T}_{f}(G)-\mathcal{S}(G)\ \mathbb{E} \bigl[\partial^{\ell-1}f(W)\bigr]=                                                                     \\
     & \qquad\sigma^{-\ell}\sum_{i_{1}\in B_{1}\backslash D_{1}}\sum_{i_{2}\in B_{2}\backslash  D_{2}}\cdots\sum_{i_{\ell}\in B_{\ell}\backslash  D_{\ell}}
    \Bigl(\mathcal{E}_{G} \bigl(X_{i_{1}},\cdots ,X_{i_{\ell-1}},\ X_{i_{\ell}}\partial^{\ell-1}f\bigl(W(D_{\ell})\bigr)\bigr)\\*
    &\qquad\qquad-\mathcal{E}_{G}\bigl(X_{i_{1}},\cdots ,X_{i_{\ell-1}},\ X_{i_{\ell}}\bigr)\ \mathbb{E} \bigl[\partial^{\ell-1}f(W)\bigr]\Bigr).
  \end{align*}

  Using $D_{\ell+1}(H_1)=\emptyset$, $D_{\ell+1}(H_{c_0})=D_\ell$, $D_{\ell+1}(H_{c_1})=T$, and $B_{\ell+1}(H_t)=D_{\ell+1}(H_{t+1})$, equations \eqref{eq:egphi1} and \eqref{eq:egphi2} telescope to
  \begin{align*}
      & \mathcal{E}_{G}\bigl(X_{i_{1}},\cdots,X_{i_{\ell-1}},X_{i_{\ell}}\partial^{\ell-1}f\bigl(W(D_{\ell})\bigr)\bigr)-\mathcal{E}_{G}\bigl(X_{i_{1}},\cdots,X_{i_{\ell}}\bigr)\ \mathbb{E} \bigl[\partial^{\ell-1}f(W)\bigr] \\
    = & -\sum_{t=c_{0}}^{c_{1}-1}\mathcal{E}_{H_{t}}\bigl(X_{i_{1}},\cdots,X_{i_{\ell}}\ ,\ \Delta_{f}(H_{t})\bigr)+\sum_{t=1}^{c_{0}-1}\mathcal{E}_{H_{t}}\bigl(X_{i_{1}},\cdots,X_{i_{\ell}}\ ,\ \Delta_{f}(H_{t})\bigr).
  \end{align*}
  For $t\geq c_1$, both $D_{\ell+1}(H_t)$ and $B_{\ell+1}(H_t)$ equal $T$, so $\mathcal{U}_f(H_t)=0$.
  Taking the sums over $i_{j}\in B_{j}\backslash D_{j}$ where $1\leq j\leq \ell$, we obtain
  \begin{align*}
    &\mathcal{T}_{f}(G)-\mathcal{S}(G)\ \mathbb{E} \bigl[\partial^{\ell-1}f(W)\bigr]\\
    = & -\sum_{t=c_{0}}^{c_{1}-1}\mathcal{U}_{f}\bigl(H_{t}\bigr)
    +\sum_{t=1}^{c_{0}-1}\mathcal{U}_{f}\bigl(H_{t}\bigr)             \\
    = & -\sum_{s\geq 0}\mathcal{U}_{f}\bigl(\Omega [\ell, s](G)\bigr)
    +\sum_{\substack{j\in A(\ell)\cup\{\ell\}:\\s_{j}\geq 1}}\sum_{s=0}^{s_{j}-1}\mathcal{U}_{f}\bigl(\Omega[p(j),s](G)\bigr).
  \end{align*}
  \hfill
\end{proof}

\begin{proof}[Proof of \cref{THM:STEP2GRAPH}.]
  If $u(\ell)=\ell$ then by definition,
  \begin{equation*}
    \Delta_{f}(G)=\partial^{\ell-2}f\bigl(W(B_{\ell})\bigr)-\partial^{\ell-2}f\bigl(W(D_{\ell})\bigr).
  \end{equation*}
  Applying the Taylor expansion with integral-form remainders, we get
  \begin{align}\label{eq:12}
    &\Delta_{f}(G)\\
    ={}  & \sum_{j=1}^{k+1}\frac{1}{j!}
            \bigl(W(B_{\ell})-W(D_{\ell})\bigr)^{j}
            \partial^{\ell-2+j}f\bigl(W(D_{\ell})\bigr)\nonumber\\
         &+\frac{1}{(k+1)!}
            \bigl(W(B_{\ell})-W(D_{\ell})\bigr)^{k+1}\cdot\nonumber\\*
                      &\quad\int_{0}^{1}(k+1)v^{k}
                        \begin{aligned}[t]
                        \Bigl(&\partial^{k+\ell-1}f\bigl(vW(D_{\ell})+(1-v)W(B_{\ell})\bigr)\\[-2pt]
                              &-\partial^{k+\ell-1}f\bigl(W(D_{\ell})\bigr)\Bigr)
                        \end{aligned}\dif v\nonumber\\
    ={}               & \sum_{j=1}^{k+1}\frac{(-1)^{j}}{j!\sigma^{j}}
                         \Bigl(\sum_{i\in B_{\ell}\backslash D_{\ell}}X_{i}\Bigr)^{j}
                         \partial^{\ell-2+j}f\bigl(W(D_{\ell})\bigr)\nonumber\\
                      &+\frac{(-1)^{k+1}}{(k+1)!\sigma^{k+1}}
                         \Bigl(\sum_{i\in B_{\ell}\backslash D_{\ell}}X_{i}\Bigr)^{k+1}\cdot\\*\nonumber
                      &\quad\int_{0}^{1}(k+1)v^{k}
                        \begin{aligned}[t]
                        \Bigl(&\partial^{k+\ell-1}f\bigl(vW(D_{\ell})+(1-v)W(B_{\ell})\bigr)\\[-2pt]
                              &-\partial^{k+\ell-1}f\bigl(W(D_{\ell})\bigr)\Bigr)
                        \end{aligned}\dif v\nonumber\\
    \overset{(*)}{ =} & \sum_{j=1}^{k+1}\frac{(-1)^{j}\sigma^{-j}}{j!}
                        \sum_{i_{\ell}\in B_{\ell}\backslash D_{\ell}}\cdots
                        \sum_{i_{\ell+j-1}\in B_{\ell+j-1}\backslash D_{\ell+j-1}}\nonumber\\
                      &\quad\times X_{i_{\ell}}\cdots X_{i_{\ell+j-1}}
                        \partial^{\ell+j-2}f\bigl(W(D_{\ell+j-1})\bigr)\nonumber\\
                      &+\frac{(-1)^{k+1}\sigma^{-(k+1)}}{(k+1)!}
                        \sum_{i_{\ell}\in B_{\ell}\backslash D_{\ell}}\cdots
                        \sum_{i_{\ell+k}\in B_{\ell+k}\backslash D_{\ell+k}}\nonumber\\
                      &\quad\times X_{i_{\ell}}\cdots X_{i_{\ell+k}}
                        \int_{0}^{1}(k+1)v^{k}
                        \begin{aligned}[t]
                          \Bigl(&\partial^{\ell+k-1}f\bigl(vW(D_{\ell+k+1})
                            +(1-v)W(B_{\ell+k+1})\bigr)\\[-2pt]
                          &-\partial^{\ell+k-1}f\bigl(W(D_{\ell+k+1})\bigr)\Bigr)
                        \end{aligned}\dif v\nonumber\\
    =                 & -\sigma^{-1}\sum_{i_{\ell}\in B_{\ell}\backslash D_{\ell}}X_{i_{\ell}}\partial^{\ell-1}f\bigl(W(D_{\ell})\bigr)                                                                                                                                                            \\*\nonumber
                      &+\sum_{j=2}^{k+1}\frac{(-1)^{j}\sigma^{-j}}{j!}
                        \sum_{i_{\ell}\in B_{\ell}\backslash D_{\ell}}\cdots
                        \sum_{i_{\ell+j-1}\in B_{\ell+j-1}\backslash D_{\ell+j-1}}\nonumber\\
                      &\quad\times X_{i_{\ell}}\cdots X_{i_{\ell+j-1}}
                        \partial^{\ell+j-2}f\bigl(W(D_{\ell+j-1})\bigr)\nonumber\\
                      &+\frac{(-1)^{k+1}\sigma^{-(k+1)}}{(k+1)!}
                        \sum_{i_{\ell}\in B_{\ell}\backslash D_{\ell}}\cdots
                        \sum_{i_{\ell+k}\in B_{\ell+k}\backslash D_{\ell+k}}\nonumber\\
                      &\quad\times X_{i_{\ell}}\cdots X_{i_{\ell+k}}
                        \Delta_{f}\bigl(\Lambda[k+1](G)\bigr),
  \end{align}
  where to obtain $(*)$ we have used the fact that $v[\ell+1],\cdots,v[\ell+k+1]$ are negative vertices in $\Lambda[k+1](G)$ and $\Lambda[j](G)\subseteq \Lambda[k+1](G)$ for $0\leq j\leq k$, which implies that $B_{\ell}=B_{\ell+1}=\cdots=B_{\ell+k+1}$ and $D_{\ell}=D_{\ell+1}=\cdots=D_{\ell+k+1}$ from the constructions of $B_{j}$'s and $D_{j}$'s.

  In \eqref{eq:defvq0}, we defined $q_0$ to be
  \begin{equation*}
    q_{0}:=\sup \bigl\{ j: j=1\text{ or }p(j)\neq j-1\text{ for }2\leq j\leq \ell \bigr\}.
  \end{equation*}
  
    We write $w:=\bigl\lvert\{ t: q_{0}+1\leq t\leq \ell~\&~s_{t}\geq 0\}\bigr\rvert$. Since $\ell=u(\ell)$, we know $s_{\ell}\geq 0$. If $q_{0}<\ell$, write the elements of $\{ t: q_{0}+1\leq t\leq \ell~\&~s_{t}\geq 0\}=\{ q_{1},\cdots,q_{w} \}$ in increasing order, so that $q_{0}+1\leq q_{1}<\cdots<q_{w}=\ell$ and $w\geq1$. If $q_{0}=\ell$, then $w=0$. By definition, the term $\mathcal{E}_{G}(X_{i_{1}},\cdots,X_{i_{\ell-1}},\Delta_{f}(G))$ is the product of $\mathcal{D}^{*}$ factors, and its last factor is
    \begin{equation*}
      \begin{cases}
        \mathcal{D}^{*}\bigl(\Delta_{f}(G)\bigr) & \text{ if }w=0, \\
        \mathcal{D}^{*}\bigl(X_{i_{q_{0}}}\cdots  X_{i_{\ell-1}}\ ,\ \Delta_{f}(G)\bigr)                                                        & \text{ if }w=1     \\
        \mathcal{D}^{*}\bigl(X_{i_{q_{0}}}\cdots X_{i_{q_{1}-1}}\ ,\ \cdots\ ,\ X_{i_{q_{(w-1)}}}\cdots X_{i_{\ell-1}}\ ,\  \Delta_{f}(G)\bigr) & \text{ if }w\geq 2
      \end{cases}.
    \end{equation*}
    For convenience, in this proof we temporarily denote
    $$
      \mathfrak{D}^{*}(\,\cdot\,):=
      \begin{cases}
        \mathcal{D}^{*}(\,\cdot\,) & \text{ if }w=0, \\
        \mathcal{D}^{*}\bigl(X_{i_{q_{0}}}\cdots  X_{i_{\ell-1}}\ ,\ \cdot\ \bigr)                                                       & \text{ if } w=1     \\
        \mathcal{D}^{*}\bigl(X_{i_{q_{0}}}\cdots X_{i_{q_{1}-1}}\ ,\ \cdots\ ,\ X_{i_{q_{(w-1)}}}\cdots X_{i_{\ell-1}}\ ,\ \cdot\ \bigr) & \text{ if } w\geq 2
      \end{cases}.
    $$
  Let $\mathfrak{E}$ denote the product of the remaining $\mathcal{D}^{*}$ factors, with the empty product interpreted as $1$. In particular, when $w=0$ we have
  \begin{equation*}
    \mathfrak{E}=\mathcal{E}_{G[\ell-1]}(X_{i_{1}},\ldots,X_{i_{\ell-1}}).
  \end{equation*}
  Then
  \begin{equation*}
    \mathcal{E}_{G}(X_{i_{1}},\cdots,X_{i_{\ell-1}},\Delta_{f}(G))=\mathfrak{E}\cdot \mathfrak{D}^{*}(\Delta_{f}(G)).
  \end{equation*}
  Combining this with \eqref{eq:12} we obtain that
  \begin{align*}
    &\mathcal{E}_{G}(X_{i_{1}},\ldots,X_{i_{\ell-1}},\Delta_{f}(G))\\
    &\qquad=\mathfrak{E}\cdot\mathfrak{D}^{*}(\Delta_{f}(G))\\
    ={}&-\sigma^{-1}\sum_{i_{\ell}\in B_{\ell}\backslash D_{\ell}}
       \mathfrak{E}\cdot\mathfrak{D}^{*}\bigl(
       X_{i_{\ell}}\partial^{\ell-1}f\bigl(W(D_{\ell})\bigr)\bigr)\\
      &+\sum_{j=2}^{k+1}\frac{(-1)^{j}}{j!\sigma^{j}}
       \sum_{i_{\ell}\in B_{\ell}\backslash D_{\ell}}\cdots
       \sum_{\substack{i_{\ell+j-1}\in\\B_{\ell+j-1}\backslash D_{\ell+j-1}}}\\
      &\qquad\times\mathfrak{E}\cdot\mathfrak{D}^{*}\bigl(
       X_{i_{\ell}}\cdots X_{i_{\ell+j-1}}
       \partial^{\ell+j-2}f\bigl(W(D_{\ell+j-1})\bigr)\bigr)\\
      &+\frac{(-1)^{k+1}}{(k+1)!\sigma^{k+1}}
       \sum_{i_{\ell}\in B_{\ell}\backslash D_{\ell}}\cdots
       \sum_{i_{\ell+k}\in B_{\ell+k}\backslash D_{\ell+k}}\\
      &\qquad\times\mathfrak{E}\cdot\mathfrak{D}^{*}\bigl(
       X_{i_{\ell}}\cdots X_{i_{\ell+k}}
       \Delta_{f}\bigl(\Lambda[k+1](G)\bigr)\bigr)\\
    \overset{(*)}{={}}&-\sigma^{-1}
       \sum_{i_{\ell}\in B_{\ell}\backslash D_{\ell}}
       \mathcal{E}_{G}\bigl(X_{i_{1}},\ldots,X_{i_{\ell-1}},
       X_{i_{\ell}}\partial^{\ell-1}f\bigl(W(D_{\ell})\bigr)\bigr)\\
      &+\sum_{j=2}^{k+1}\frac{(-1)^{j}}{j!\sigma^{j}}
       \sum_{i_{\ell}\in B_{\ell}\backslash D_{\ell}}\cdots
       \sum_{\substack{i_{\ell+j-1}\in\\B_{\ell+j-1}\backslash D_{\ell+j-1}}}\\
      &\qquad\times\mathcal{E}_{\Lambda[j-1](G)}\bigl(
       X_{i_{1}},\ldots,X_{i_{\ell+j-2}},
       X_{i_{\ell+j-1}}\partial^{\ell+j-2}f\bigl(W(D_{\ell+j-1})\bigr)\bigr)\\
      &+\frac{(-1)^{k+1}}{(k+1)!\sigma^{k+1}}
       \sum_{i_{\ell}\in B_{\ell}\backslash D_{\ell}}\cdots
       \sum_{i_{\ell+k}\in B_{\ell+k}\backslash D_{\ell+k}}\\
      &\qquad\times\mathcal{E}_{\Lambda[k+1](G)}\bigl(
       X_{i_{1}},\ldots,X_{i_{\ell+k}},
       \Delta_{f}\bigl(\Lambda[k+1](G)\bigr)\bigr),
  \end{align*}
  where to get $(*)$ we have used the condition that $v[\ell+1],\cdots,v[\ell+k+1]$ are negative vertices in $\Lambda[k+1](G)$. Indeed, the factorization stays the same because all these vertices are added to the same branch, and the possibly empty list $q_{1},\cdots,q_{w}$ remains unchanged because all the new vertices are negative. Taking the sum over $i_{j}\in B_{j}\backslash D_{j}$ for $1\leq j\leq \ell-1$, we have
  \begin{align}\label{eq:ufgexpansion}
    \mathcal{U}_{f}(G)
    = & \sum_{j=1}^{k+1}(-1)^{j}\frac{1}{j!}\mathcal{T}_{f}\bigl(\Lambda[j-1](G)\bigr)+(-1)^{k+1}\frac{1}{(k+1)!}\mathcal{U}_{f}\bigl(\Lambda [k+1](G)\bigr)                 \\
    = & \sum_{j=0}^{k}(-1)^{j+1}\frac{1}{(j+1)!}\mathcal{T}_{f}\bigl(\Lambda[j](G)\bigr)+(-1)^{k+1}\frac{1}{(k+1)!}\mathcal{U}_{f}\bigl(\Lambda [k+1](G)\bigr).\nonumber   
  \end{align}
  Now consider the case $u(\ell)<\ell$.
  Let $G[u(\ell)]:=(V',E',s_{1:u(\ell)})\subseteq G$ be the order-$u(\ell)$ sub-genogram of $G$ as defined in the last paragraph of \cref{sec:genogram}. Now by \eqref{eq:ufgexpansion}, we have
  \begin{align*}
    \mathcal{U}_{f}(G[u(\ell)])
    =&\sum_{j=0}^{k}(-1)^{j+1}\frac{1}{(j+1)!}\mathcal{T}_{f}\bigl(\Lambda[j](G[u(\ell)])\bigr)\\
    &\ +(-1)^{k+1}\frac{1}{(k+1)!}\mathcal{U}_{f}\bigl(\Lambda [k+1](G[u(\ell)])\bigr).
  \end{align*}
  Replacing $k$ by $\ell-u(\ell)-1$ and $\ell-u(\ell)+k$ respectively, we get that
  \begin{align}
    \label{eq:23} &
    \begin{aligned}
      \mathcal{U}_{f}(G[u(\ell)])
      ={}& \sum_{j=0}^{\ell-u(\ell)-1}\frac{(-1)^{j+1}}{(j+1)!}
        \mathcal{T}_{f}\bigl(\Lambda[j](G[u(\ell)])\bigr)\\
        &{}+\frac{(-1)^{\ell-u(\ell)}}{(\ell-u(\ell))!}\\[-2pt]
        &\qquad\times\mathcal{U}_{f}\bigl(
        \Lambda[\ell-u(\ell)](G[u(\ell)])\bigr),
    \end{aligned} \\
    \label{eq:34} &
    \begin{aligned}
      \mathcal{U}_{f}(G[u(\ell)])
      ={}& \sum_{j=0}^{\ell-u(\ell)+k}\frac{(-1)^{j+1}}{(j+1)!}
        \mathcal{T}_{f}\bigl(\Lambda[j](G[u(\ell)])\bigr)\\
        &{}+\frac{(-1)^{\ell-u(\ell)+k+1}}{(\ell-u(\ell)+k+1)!}\\[-2pt]
        &\qquad\times\mathcal{U}_{f}\bigl(
        \Lambda[\ell-u(\ell)+k+1](G[u(\ell)])\bigr).
    \end{aligned}
  \end{align}
  By taking the difference of \eqref{eq:23} and \eqref{eq:34} we obtain that
  \begin{align*}
      & (-1)^{\ell-u(\ell)}\frac{1}{(\ell-u(\ell))!}\mathcal{U}_{f}\bigl(\Lambda [\ell-u(\ell)](G[u(\ell)])\bigr)                 \\
    = & \sum_{j=\ell-u(\ell)}^{\ell-u(\ell)+k}(-1)^{j+1}\frac{1}{(j+1)!}\mathcal{T}_{f}\bigl(\Lambda[j](G[u(\ell)])\bigr)         \\*
      & \ +(-1)^{\ell-u(\ell)+k+1}\frac{1}{(\ell-u(\ell)+k+1)!}\mathcal{U}_{f}\bigl(\Lambda [\ell-u(\ell)+k+1](G[u(\ell)])\bigr).
  \end{align*}
  Thus, we have
  \begin{align*}
      & \mathcal{U}_{f}(G)=\mathcal{U}_{f}\bigl(\Lambda[\ell-u(\ell)](G[u(\ell)])\bigr)                                                                                                                               \\
    = & \sum_{j=0}^{k}(-1)^{j+1}\frac{(\ell-u(\ell))!}{(j+1+\ell-u(\ell))!}\mathcal{T}_{f}\bigl(\Lambda[j+\ell-u(\ell)](G[u(\ell)])\bigr)                                                                             \\
      & \ +(-1)^{k+1}\frac{(\ell-u(\ell))!}{(k+1+\ell-u(\ell))!}\mathcal{U}_{f}\bigl(\Lambda[k+1+\ell-u(\ell)](G[u(\ell)])\bigr)                                                                                      \\
    = & \sum_{j=0}^{k}(-1)^{j+1}\frac{(\ell-u(\ell))!}{(j+1+\ell-u(\ell))!} \mathcal{T}_{f}\bigl(\Lambda[j](G)\bigr)\\*
    &\ +(-1)^{k+1}\frac{(\ell-u(\ell))!}{(k+1+\ell-u(\ell))!}\mathcal{U}_{f}\bigl(\Lambda[k+1](G)\bigr).
  \end{align*}
  \hfill
\end{proof}
\section{Proof of Lemma~\ref{THM:REMAINDERCTRL1234}}\label{sec:lemma5}

Before we proceed to Lemma~\ref{THM:REMAINDERCTRL1234}, we need a generalization of Young's inequality.

For $t\in\mathbb N_+$, let
\begin{equation*}
  C(t):=\bigl\{(\ell,\eta_{1:\ell}):
  \ell,\eta_{1:\ell}\in\mathbb N_+,\ \textstyle\sum_{j=1}^{\ell}\eta_j=t\bigr\}.
\end{equation*}
For a composition $(\ell,\eta_{1:\ell})\in C(t)$, write
$n_{0}:=0$ and $n_{a}:=\sum_{b=1}^{a}\eta_{b}$ for $1\leq a\leq \ell$.
We use the ordinary bracket notation
\begin{equation*}
  [\eta_{1},\ldots,\eta_{\ell}]\triangleright(Y_{1},\ldots,Y_{t})
  :=\prod_{a=1}^{\ell}\mathbb{E}\biggl[\prod_{i=n_{a-1}+1}^{n_{a}}Y_{i}\biggr].
\end{equation*}

\begin{lemma}\label{thm:lemmayoung}
  Given $t\in\mathbb{N}_{+}$, let $(Y_i)_{i=1}^{t}$ be a sequence of random variables, and real numbers $p_{1},\cdots, p_{t}>1$ satisfy that $1/p_{1}+\cdots+1/p_{t}=1$. Then for any $(\ell, \eta_{1:\ell})\in C(t)$, we have that
  \begin{equation}\label{eq:lemmayoung3}
    [\eta_{1},\cdots, \eta_{\ell}]\triangleright (\lvert Y_{1} \rvert,\cdots,\lvert Y_{t} \rvert)\leq
    \frac{1}{p_{1}}\mathbb{E} [\lvert Y_{1} \rvert^{p_{1}}]+\cdots+ \frac{1}{p_{t}}\mathbb{E} [\lvert Y_{t} \rvert^{p_{t}}].
  \end{equation}
\end{lemma}

\begin{proof}[Proof of \cref{thm:lemmayoung}.]
  First, we prove
  \begin{align}
     & \mathbb{E} [\lvert Y_{1}\cdots Y_{t} \rvert]\leq \frac{1}{p_{1}}\mathbb{E} [\lvert Y_{1} \rvert^{p_{1}}]+\cdots+ \frac{1}{p_{t}}\mathbb{E} [\lvert Y_{t} \rvert^{p_{t}}],\label{eq:lemmayoung1}                            \\
     & \mathbb{E} [\lvert Y_{1} \rvert]\cdots\mathbb{E} [\lvert Y_{t} \rvert] \leq \frac{1}{p_{1}}\mathbb{E} [\lvert Y_{1} \rvert^{p_{1}}]+\cdots+ \frac{1}{p_{t}}\mathbb{E} [\lvert Y_{t} \rvert^{p_{t}}].\label{eq:lemmayoung2}
  \end{align}
  To achieve this goal, note that Young's inequality is stated as follows: For any $a_{1},\cdots,a_{t}\geq 0$, and $p_{1},\cdots,p_{t}>1$ such that $1/p_{1}+\cdots+1/p_{t}=1$, we have
  \begin{equation*}
    a_{1}\cdots a_{t}\leq \frac{1}{p_{1}}a_{1}^{p_{1}}+\cdots+\frac{1}{p_{t}}a_{t}^{p_{t}}.
  \end{equation*}
  Thus, by Young's inequality we know that
  \begin{equation*}
    \lvert Y_{1}\cdots Y_{t}\rvert\leq \frac{1}{p_{1}}\lvert Y_{1} \rvert^{p_{1}}+\cdots+\frac{1}{p_{t}}\lvert Y_{t} \rvert^{p_{t}}.
  \end{equation*}
  Taking the expectation, we have
  \begin{equation*}
    \mathbb{E} [\lvert Y_{1}\cdots Y_{t} \rvert]\leq \frac{1}{p_{1}}\mathbb{E} [\lvert Y_{1} \rvert^{p_{1}}]+\cdots+\frac{1}{p_{t}}\mathbb{E} [\lvert Y_{t} \rvert^{p_{t}}].
  \end{equation*}
  Again by Young's inequality, we obtain that
  \begin{equation*}
    \mathbb{E} [\lvert Y_{1} \rvert]\cdots \mathbb{E} [\lvert Y_{t} \rvert]\leq \frac{1}{p_{1}}\mathbb{E} [\lvert Y_{1} \rvert]^{p_{1}}+\cdots+\frac{1}{p_{t}}\mathbb{E} [\lvert Y_{t} \rvert]^{p_{t}}.
  \end{equation*}
  By Jensen's inequality, $\mathbb{E} [\lvert Y_{i} \rvert]^{p_{i}}\leq \mathbb{E} [\lvert Y_{i} \rvert^{p_{i}}]$ for $i\in [t]$.
  This implies that
  \begin{equation*}
    \mathbb{E} [\lvert Y_{1} \rvert]\cdots \mathbb{E} [\lvert Y_{t} \rvert]\leq \frac{1}{p_{1}}\mathbb{E} [\lvert Y_{1} \rvert^{p_{1}}]+\cdots+\frac{1}{p_{t}}\mathbb{E} [\lvert Y_{t} \rvert^{p_{t}}].
  \end{equation*}

  Finally, we prove \eqref{eq:lemmayoung3}. If $\ell=1$, the result is
  exactly \eqref{eq:lemmayoung1}. Suppose that $\ell\geq 2$, and define
  \begin{equation*}
    \frac{1}{q_{a}}:=\sum_{i=n_{a-1}+1}^{n_{a}}\frac{1}{p_{i}},
    \qquad 1\leq a\leq \ell.
  \end{equation*}
  Then $q_{a}>1$ for every $a$ and $\sum_{a=1}^{\ell}1/q_{a}=1$. Hence
  \begin{align*}
    &[\eta_{1},\ldots,\eta_{\ell}]\triangleright
      (\lvert Y_{1}\rvert,\ldots,\lvert Y_{t}\rvert)\\
    &\quad\overset{\eqref{eq:lemmayoung2}}{\leq}
      \sum_{a=1}^{\ell}\frac{1}{q_{a}}
      \mathbb{E}\biggl[\biggl(\prod_{i=n_{a-1}+1}^{n_{a}}
      \lvert Y_{i}\rvert\biggr)^{q_{a}}\biggr]\\
    &\quad\overset{\eqref{eq:lemmayoung1}}{\leq}
      \sum_{a=1}^{\ell}\sum_{i=n_{a-1}+1}^{n_{a}}
      \frac{1}{p_{i}}\mathbb{E}[\lvert Y_{i}\rvert^{p_{i}}]
      =\sum_{i=1}^{t}\frac{1}{p_{i}}\mathbb{E}[\lvert Y_{i}\rvert^{p_{i}}].
  \end{align*}
  In the second inequality, a singleton block gives equality directly; on
  every block of size at least two, \eqref{eq:lemmayoung1} is applied with
  exponents $p_i/q_a>1$.
\end{proof}

Secondly, let's show some properties of the $\mathcal{D}$ and $\mathcal{D}^{*}$ operators that will be useful later.

\begin{lemma}\label{thm:dastaddterm}Let $(Y_i)_{i=1}^{t}$ be a sequence of random variables. Suppose for any $i,j\in\mathbb{N}_{+}$ such that $i\leq j\leq t$, we have $\mathbb{E}\bigl[ |Y_{i}\cdots Y_{j}|\bigr]<\infty$.
  Then the following holds for all $t\in\mathbb{N}_{+}$ such that $t\geq 2$ and for any $j=1,\cdots, t-1$:
  \begin{align}
     & \mathcal{D}^{*}\bigl(Y_{1},\cdots,Y_{j}\mathcal{D}(Y_{j+1},\cdots,Y_{t})\bigr)=\mathcal{D}^{*}(Y_{1},\cdots,Y_{t}),\label{eq:dastaddterm1} \\
     & \mathcal{D}\bigl(Y_{1},\cdots,Y_{j}\mathcal{D}(Y_{j+1},\cdots,Y_{t})\bigr)=\mathcal{D}(Y_{1},\cdots,Y_{t}).\label{eq:dastaddterm3}
  \end{align}
  In particular,
  \begin{equation}\label{eq:dastaddterm2}
    \begin{aligned}
        & \mathcal{D}^{*}(Y_{1},\cdots,Y_{t})=\mathcal{D}^{*}(Y_{1},\cdots, Y_{t-2}\ ,\ Y_{t-1} \mathcal{D}(Y_{t}))           \\
      = & \mathcal{D}^{*}(Y_{1},\cdots, Y_{t-2}\ ,\ Y_{t-1} Y_{t})-\mathcal{D}^{*}(Y_{1},\cdots,Y_{t-1})\ \mathbb{E} [Y_{t}].
    \end{aligned}
  \end{equation}
  Moreover, we know that $$\mathbb{E}\bigl[\mathcal{D}(Y_1,\cdots, Y_t)\bigr]=0.$$
\end{lemma}

\begin{proof}[Proof of \cref{thm:dastaddterm}.]
  We perform induction on $j$ to prove that
  \begin{equation*}
    \mathcal{D}^*\bigl(Y_{1},\cdots,Y_{j}\mathcal{D}(Y_{j+1},\cdots,Y_{t})\bigr)=\mathcal{D}^*(Y_{1},\cdots,Y_{t}).
  \end{equation*}

  If $j=1$, this is precisely the definition. Supposing the lemma holds for $j$ ($j\leq t-1$), consider the case for $j+1$.
  By definition,
  \begin{align*}
      & \mathcal{D}\bigl(Y_{1},\cdots,Y_{j+1}\mathcal{D}(Y_{j+2},\cdots,Y_{t})\bigr)                             \\
    = & \mathcal{D}\Bigl(Y_{1}\mathcal{D}\bigl(Y_{2},\cdots,Y_{j+1}\mathcal{D}(Y_{j+2},\cdots,Y_{t})\bigr)\Bigr) \\
    = & \mathcal{D} \bigl(Y_{1}\mathcal{D}(Y_{2},\cdots,Y_{t})\bigr)=\mathcal{D}(Y_{1},\cdots,Y_{t}).
  \end{align*}
  Note that we have used the inductive hypothesis in the second equation.
  By induction, \eqref{eq:dastaddterm3} is proven.

  Now for any $j=1\cdots t-1$,
  \begin{align*}
      & \mathcal{D}^{*}\bigl(Y_{1},\cdots,Y_{j}\mathcal{D}(Y_{j+1},\cdots,Y_{t})\bigr)                         \\
    = & \mathbb{E} \Bigl[Y_{1}\mathcal{D}\bigl(Y_{2},\cdots,Y_{j}\mathcal{D}(Y_{j+1},\cdots,Y_{t})\bigr)\Bigr] \\
    = & \mathbb{E} \bigl[Y_{1}\mathcal{D}(Y_{2},\cdots,Y_{t})\bigr]=\mathcal{D}^{*}(Y_{1},\cdots,Y_{t}).
  \end{align*}
  Finally, we remark that
  \begin{align*}
    \mathbb{E}\bigl[\mathcal{D}(Y_1,\dots,Y_t)\bigr] & =\mathbb{E}\bigl[\mathcal{D}(Y_1\mathcal{D}(Y_2,\dots,Y_t))\bigr]
    \\&=\mathbb{E}\bigl[Y_1\mathcal{D}(Y_2,\dots,Y_t)\bigr]-\mathbb{E}\bigl[Y_1\mathcal{D}(Y_2,\dots,Y_t)\bigr]=0.
  \end{align*}\hfill
\end{proof}

\begin{lemma}\label{thm:dastexpr}
  Let $(Y_i)_{i=1}^{t}$ be a sequence of random variables. Suppose that
  $\mathbb{E}\bigl[|Y_i\cdots Y_j|\bigr]<\infty$ whenever
  $1\leq i\leq j\leq t$. Then $\mathcal{D}^{*}$ and $\mathcal{D}$ admit the
  following expressions:
  \begin{align}
    \mathcal{D}^{*}(Y_{1},\ldots,Y_{t})
      &=\sum_{(\ell,\eta_{1:\ell})\in C(t)}(-1)^{\ell-1}
        [\eta_{1},\ldots,\eta_{\ell}]\triangleright
        (Y_{1},\ldots,Y_{t}),\label{eq:dastexpr} \\
    \mathcal{D}(Y_{1},\ldots,Y_{t})
      &=Y_{1}\cdots Y_{t}-\mathcal{D}^{*}(Y_{1},\ldots,Y_{t})\nonumber\\
      &\quad-\sum_{j=1}^{t-1}Y_{1}\cdots Y_{j}\,
        \mathcal{D}^{*}(Y_{j+1},\ldots,Y_{t}).\label{eq:dexpr}
  \end{align}
\end{lemma}

\begin{proof}[Proof of \cref{thm:dastexpr}.]
  We perform induction on $t$.

  If $t=1$, then by definition $\mathcal{D}^{*}(Y_{1})=\mathbb{E} [Y_{1}]=[1]\triangleright (Y_{1})$ and $\mathcal{D}(Y_{1})=Y_{1}-\mathbb{E}[Y_{1}]=Y_{1}-\mathcal{D}^{*}(Y_{1})$.

  Supposing the results hold for $1,2,\cdots,t-1$, we consider the case $t$. Suppose that we have $\mathbb{E}[|Y_1\cdots Y_t|]<\infty$. By the inductive hypothesis, we have
  \begin{align*}
                      & \mathcal{D}^{*}(Y_{1},Y_{2},\cdots,Y_{t})=\mathbb{E} [Y_{1}\mathcal{D}(Y_{2},\cdots,Y_{t})]                                                                                       \\
    \overset{(a)}{=}  & \mathbb{E} [Y_{1}Y_{2}\cdots Y_{t}]- \mathbb{E}(Y_1)\mathcal{D}^*(Y_2,\dots,Y_t)-\sum_{j=3}^{t}\mathbb{E} [Y_{1}\cdots Y_{j-1}\ \mathcal{D}^{*}(Y_{j},\cdots,Y_{t})]              \\  \overset{}{=} & \mathbb{E} [Y_{1}Y_{2}\cdots Y_{t}]-\sum_{j=2}^{t}\mathbb{E} [Y_{1}\cdots Y_{j-1}\ \mathcal{D}^{*}(Y_{j},\cdots,Y_{t})]                                    \\
    \overset{(b)} {=} & \mathbb{E} [Y_{1}Y_{2}\cdots Y_{t}]-\sum_{j=2}^{t}\sum_{ C (t-j+1)}\mathbb{E}[Y_{1}\cdots Y_{j-1}\ (-1)^{\ell-1}[\eta_{1},\cdots,\eta_{\ell}]\triangleright (Y_{j},\cdots,Y_{t})] \\
    =                 & \mathbb{E} [Y_{1}Y_{2}\cdots Y_{t}]+\sum_{j=2}^{t}\sum_{  C  (t-j+1)}(-1)^{\ell}[j-1,\eta_{1},\cdots,\eta_{\ell}]\triangleright (Y_{1},\cdots,Y_{t})                              \\
    \overset{(c)}{=}  & \mathbb{E} [Y_{1}Y_{2}\cdots Y_{t}]+\sum_{ C  (t)\setminus \{\ell=1,~\eta_1=t\} }(-1)^{\ell-1}[\eta_{1},\cdots,\eta_{\ell}]\triangleright (Y_{1},\cdots,Y_{t})                    \\
    =                 & \sum_{ C  (t)}(-1)^{\ell-1}[\eta_{1},\cdots,\eta_{\ell}]\triangleright (Y_{1},\cdots,Y_{t}),
  \end{align*}
  where to get $(a)$ we have used the fact that by inductive hypothesis \eqref{eq:dexpr} holds for $t-1$, and to get $(b)$ we have used the fact that we assumed that \eqref{eq:dastexpr} hold for $t-1$. Finally, to get $(c)$ we have used the fact that
  \begin{align*}
      & C (t)={\textstyle \bigl\{\ell,\eta_{1:\ell}\in \mathbb{N}_{+}: \sum_{j=1}^{\ell}\eta_j=t\bigr\}}                                                         \\
    = & {\textstyle \{\ell=1,~\eta_1=t\}\cup\bigcup_{a=1}^{t-1}\bigl\{\ell,\eta_{1:\ell}\in\mathbb{N}_{+}:\ell\ge 2,~\eta_1=a, ~\sum_{j=1}^{\ell}\eta_j=t\bigr\}.}
  \end{align*}

  Moreover, we also have
  \begin{align*}
      & \mathcal{D}(Y_{1},Y_{2},\cdots,Y_{t})=Y_{1}\mathcal{D}(Y_{2},\cdots,Y_{t})-\mathcal{D}^{*}(Y_{1},\cdots,Y_{t})                                                             \\
    = & Y_{1}Y_{2}\cdots Y_{t}-Y_{1}\mathcal{D}^{*}(Y_{2},\cdots,Y_{t}) \\
      &\quad-\sum_{j=2}^{t-1}Y_{1}\cdots Y_{j}\mathcal{D}^{*}(Y_{j+1},\cdots,Y_{t})
       -\mathcal{D}^{*}(Y_{1},\cdots,Y_{t}) \\
    = & Y_{1}Y_{2}\cdots Y_{t}-\mathcal{D}^*(Y_1,\cdots,Y_t) \\
      &\quad-\sum_{j=1}^{t-1}Y_{1}\cdots Y_{j}\mathcal{D}^{*}(Y_{j+1},\cdots,Y_{t}).
  \end{align*}
  Thus, the results also hold for $t$. And the proof is complete by induction.\hfill
\end{proof}

Next we need the following notions of \emph{compositional} $\mathcal{D}^{*}$ and $\mathcal{D}$ \emph{operators}.
\begin{align}
  [\eta_{1},\ldots,\eta_{\ell}]\triangleright \mathcal{D}^{*}(Y_{1},\ldots,Y_{t})
    &:={\mathcal{D}}^{*}\bigl(Y_{1}\cdots Y_{\eta_{1}},
      Y_{\eta_{1}+1}\cdots Y_{\eta_{1}+\eta_{2}},\nonumber\\
    &\hspace{7em}\ldots,
      Y_{\eta_{1}+\cdots+\eta_{\ell-1}+1}\cdots Y_{t}\bigr),
      \label{eq:compositionaldast} \\
  [\eta_{1},\ldots,\eta_{\ell}]\triangleright \mathcal{D}(Y_{1},\ldots,Y_{t})
    &:={\mathcal{D}}\bigl(Y_{1}\cdots Y_{\eta_{1}},
      Y_{\eta_{1}+1}\cdots Y_{\eta_{1}+\eta_{2}},\nonumber\\
    &\hspace{7em}\ldots,
      Y_{\eta_{1}+\cdots+\eta_{\ell-1}+1}\cdots Y_{t}\bigr).
      \label{eq:compositionald}
\end{align}
Note that a compositional $\mathcal{D}$ term is a random variable while a compositional $\mathcal{D}^{*}$ operator gives a deterministic value. We remark that
\begin{align}\label{mean0cd}
  &\mathbb{E}\bigl[[\eta_1,\dots, \eta_l]\triangleright \mathcal{D}(Y_1,\dots,Y_t)\bigr]\nonumber\\
  &\quad=\mathbb{E}\Bigl[\mathcal{D}\bigl(Y_1\cdots Y_{\eta_1},\ldots,\nonumber\\
  &\hspace{8em}Y_{\eta_1+1}\cdots
  Y_{\eta_1+\cdots+\eta_{\ell-1}+1}\cdots Y_t\bigr)\Bigr]=0.
\end{align}
Moreover, by definition and \eqref{eq:dastaddterm1}, we can directly check that
\begin{equation}\label{eq:strangeeq}
  \begin{aligned}
      & [\eta_{1},\cdots,\eta_{\ell}]\triangleright \mathcal{D}^{*}\bigl(Y_{1},\cdots,Y_{t}\bigr)                                                                                                                                                       \\
    = & [\eta_{1},\cdots,\eta_{s-1},\eta_{s}+1]
        \triangleright \mathcal{D}^{*}\Bigl(
        Y_{1},\cdots,Y_{\eta_{1}+\cdots+\eta_{s}},\\
      &\qquad [\eta_{s+1},\cdots,\eta_{\ell}]\triangleright
        \mathcal{D}\bigl(Y_{\eta_{1}+\cdots+\eta_{s}+1},\cdots,Y_{t}\bigr)
        \Bigr).
  \end{aligned}
\end{equation}

The following lemma shows some upper bounds on their norms.

\begin{lemma}\label{thm:sec11lemma1}Let $(Y_i)_{i=1}^{t}$ be random variable such that for all $i,j\in \mathbb{N}_{+}$ such that $i\leq j\leq t$ we have $\mathbb{E}\bigl[ \lvert Y_i Y_{i+1}\cdots Y_j\rvert \bigr]<\infty$. Then for any $q\geq 1$ the following holds
  \begin{align}
     & \bigl\lvert [\eta_{1},\cdots,\eta_{\ell}]\triangleright
       \mathcal{D}^{*}(Y_{1},\cdots,Y_{t}) \bigr\rvert\nonumber\\
     &\qquad\leq \sum_{(s,\zeta_{1:s})\in C(t)}
       [\zeta_{1},\cdots,\zeta_{s}]\triangleright
       \bigl(\lvert Y_{1} \rvert,\cdots,\lvert Y_{t} \rvert\bigr),
       \label{eq:dastbound1} \\
     & \begin{aligned}
         \bigl\lvert [\eta_{1},\cdots,\eta_{\ell}] & \triangleright \mathcal{D}^{*}(Y_{1},\cdots,Y_{t}) \bigr\rvert\leq \bigl\lvert\,\mathbb{E} [ Y_{1}\cdots Y_{t} ]\,\bigr\rvert
         \\
                                                   & +\sum_{j=1}^{t-1}\sum_{(s,\zeta_{1:s})\in C (j)}\ [\zeta_{1},\cdots,\zeta_{s}]\triangleright \bigl(\lvert Y_{1} \rvert,\cdots,\lvert Y_{j} \rvert\bigr)\cdot \bigl\lvert\,\mathbb{E} [ Y_{j+1}\cdots Y_{t} ]\,\bigr\rvert,
       \end{aligned}\label{eq:dastbound2} \\
     & \bigl\lVert [\eta_{1},\cdots,\eta_{\ell}]\triangleright
       \mathcal{D}(Y_{1},\cdots,Y_{t}) \bigr\rVert _{q}\nonumber\\
     &\qquad\leq 2\sum_{(s,\zeta_{1:s})\in C(t)}
       \Bigl([\zeta_{1},\cdots,\zeta_{s}]\triangleright
       \bigl(\lvert Y_{1} \rvert^{q},\cdots,\lvert Y_{t}\rvert^{q}\bigr)\Bigr)^{1/q}.
       \label{eq:dbound}
  \end{align}
\end{lemma}

\begin{proof}[Proof of \cref{thm:sec11lemma1}.]
  Applying \cref{thm:dastexpr}, we get
  \begingroup\footnotesize
  \begin{align*}
         & \bigl\lvert [\eta_{1},\cdots,\eta_{\ell}]\triangleright \mathcal{D}^{*}(Y_{1},\cdots,Y_{t}) \bigr\rvert                                                                                                                                                         \\
    =    & \Bigl\lvert \mathcal{D}^{*}\bigl(Y_{1}\cdots Y_{\eta_{1}},Y_{\eta_{1}+1}\cdots Y_{\eta_{1}+\eta_{2}},\ldots,Y_{\eta_{1}+\cdots+\eta_{\ell-1}+1}\cdots Y_{t}\bigr) \Bigr\rvert \\
    =    & \Bigl\lvert \sum_{(s,\zeta_{1:s})\in C(\ell)}
             (-1)^{s-1}[\zeta_{1},\ldots,\zeta_{s}]\triangleright\\
         &\hspace{4em}\bigl(Y_{1}\cdots Y_{\eta_{1}},
             Y_{\eta_{1}+1}\cdots Y_{\eta_{1}+\eta_{2}},\ldots,
             Y_{\eta_{1}+\cdots+\eta_{\ell-1}+1}\cdots Y_{t}\bigr)
           \Bigr\rvert \\
    \leq & \sum_{(s,\zeta_{1:s})\in C(\ell)}
             [\zeta_{1},\ldots,\zeta_{s}]\triangleright\\
         &\hspace{4em}\bigl(\lvert Y_{1}\cdots Y_{\eta_{1}}\rvert,
             \lvert Y_{\eta_{1}+1}\cdots Y_{\eta_{1}+\eta_{2}}\rvert,\ldots,
             \lvert Y_{\eta_{1}+\cdots+\eta_{\ell-1}+1}\cdots Y_{t}\rvert\bigr) \\
    \leq & \sum_{(s,\lambda_{1:s})\in C(t)}
             [\lambda_{1},\ldots,\lambda_{s}]\triangleright
             \bigl(\lvert Y_{1}\rvert,\ldots,\lvert Y_{t}\rvert\bigr),
  \end{align*}
  \endgroup
  where in the last inequality we have used the following regrouping. For
  every $(s,\zeta_{1:s})\in C(\ell)$, let $\rho_{0}:=0$,
  $\rho_{a}:=\sum_{b=1}^{a}\zeta_{b}$, and
  \begin{equation*}
    \lambda_{a}:=\sum_{h=\rho_{a-1}+1}^{\rho_{a}}\eta_{h},
    \qquad 1\leq a\leq s.
  \end{equation*}
  Then $(s,\lambda_{1:s})\in C(t)$ and
  \begin{align*}
      & [\lambda_1,\cdots,\lambda_s]\triangleright (|Y_1|,|Y_2|,\cdots,|Y_t|)                                                                                                                            \\
    = & [\zeta_{1},\cdots,\zeta_{s}]\triangleright \bigl(|Y_{1}\cdots Y_{\eta_{1}}|\ ,\ |Y_{\eta_{1}+1}\cdots Y_{\eta_{1}+\eta_{2}}|\ ,\cdots,\ |Y_{\eta_{1}+\cdots+\eta_{\ell-1}+1}\cdots Y_{t}|\bigr).
  \end{align*}
  Using similar ideas, we observe that
  \begin{align*}
                                     & \bigl\lvert [\eta_{1},\cdots,\eta_{\ell}]\triangleright \mathcal{D}^{*}(Y_{1},\cdots,Y_{t}) \bigr\rvert                                                                                                                                                                                                 \\
    =\ \                             & \bigl\lvert \mathcal{D}^{*}\bigl(Y_{1}\cdots Y_{\eta_{1}}\ ,\ Y_{\eta_{1}+1}\cdots Y_{\eta_{1}+\eta_{2}}\ ,\cdots,\ Y_{\eta_{1}+\cdots+\eta_{\ell-1}+1}\cdots Y_{t}\bigr) \bigr\rvert                                                                                                                   \\
    \overset{\eqref{eq:dastexpr}}{=} & \Bigl\lvert \sum_{(s,\zeta_{1:s})\in C (\ell)}(-1)^{s-1}[\zeta_{1},\cdots,\zeta_{s}]\triangleright \bigl(Y_{1}\cdots Y_{\eta_{1}}\ ,\cdots,\ Y_{\eta_{1}+\cdots+\eta_{\ell-1}+1}\cdots Y_{t}\bigr) \Bigr\rvert                                                                                             \\
    =\ \                             & \Bigl\lvert \mathbb{E} [Y_{1}\cdots Y_{t}] + \sum_{j=1}^{\ell-1}\sum_{(s,\zeta_{1:s})\in C (j)}(-1)^{s}[\zeta_{1},\cdots,\zeta_{s}]\triangleright                                                                                                                                                       \\*
                                     & \quad \bigl(Y_{1}\cdots Y_{\eta_{1}},\cdots,\ Y_{\eta_{1}+\cdots+\eta_{j-1}+1}\cdots Y_{\eta_{1}+\cdots+\eta_{j}}\bigr)\cdot
    \mathbb{E} [Y_{\eta_{1}+\cdots+\eta_{j}+1}\cdots Y_{t}] \Bigr\rvert                                                                                                                                                                                                                                                                        \\
    \leq\ \                          & \bigl\lvert \mathbb{E} [Y_{1}\cdots Y_{t}]\bigr\rvert+\sum_{j=1}^{\ell-1}\sum_{ (s,\zeta_{1:s})\in C (j)}[\zeta_{1},\cdots,\zeta_{s}]\triangleright                                                                                                                                                     \\*
                                     & \quad\bigl(\bigl\lvert Y_{1}\cdots Y_{\eta_{1}}\bigr\rvert\ ,\cdots,\bigl\lvert Y_{\eta_{1}+\cdots+\eta_{j-1}+1}\cdots Y_{\eta_{1}+\cdots+\eta_{j}}\bigr\rvert\bigr)\cdot\bigl\lvert\mathbb{E} [Y_{\eta_{1}+\cdots+\eta_{j}+1}\cdots Y_{t}] \bigr\rvert                                                 \\
    \overset{(*)}{\leq} \ \          & \bigl\lvert \mathbb{E} [Y_{1}\cdots Y_{t}]\bigr\rvert +\sum_{j=1}^{\ell-1}\sum_{(s,\lambda_{1:s})\in C (\eta_{1}+\cdots+\eta_{j})}[\lambda_{1},\cdots,\lambda_{s}]\triangleright                                                                                                                            \\*
                                     & \quad \bigl(\lvert Y_{1}\rvert\ ,\ \lvert Y_{2}\rvert\ ,\cdots,\ \bigl\lvert Y_{\eta_{1}+\cdots+\eta_{j}}\bigr\rvert\bigr)\cdot \bigl\lvert\mathbb{E} [Y_{\eta_{1}+\cdots+\eta_{j}+1}\cdots Y_{t}] \bigr\rvert                                                                                          \\
    \leq \ \                         & \bigl\lvert \mathbb{E} [Y_{1}\cdots Y_{t}]\bigr\rvert
      +\sum_{h=1}^{t-1}\sum_{(s,\lambda_{1:s})\in C(h)}
       [\lambda_{1},\cdots,\lambda_{s}]\triangleright\\
                                      &\quad \bigl(\lvert Y_{1}\rvert,\ldots,
       \lvert Y_{h}\rvert\bigr)\cdot
       \bigl\lvert\mathbb{E} [Y_{h+1}\cdots Y_{t}] \bigr\rvert,
  \end{align*}
  where to obtain $(*)$ we used the same cumulative regrouping above, now
  applied to $(s,\zeta_{1:s})\in C(j)$ and $\eta_{1:j}$, so that
  \begin{align*}
      & [\lambda_{1},\cdots,\lambda_{s}]\triangleright \bigl(\lvert Y_{1}\rvert,\cdots,\bigl\lvert Y_{\eta_{1}+\cdots+\eta_{j}}\bigr\rvert\bigr)                                                                   \\
    = & [\zeta_{1},\cdots,\zeta_{s}]\triangleright \bigl(\bigl\lvert Y_{1}\cdots Y_{\eta_{1}} \bigr\rvert,\cdots,\bigl\lvert Y_{\eta_{1}+\cdots+\eta_{j-1}+1}\cdots Y_{\eta_{1}+\cdots+\eta_{j}}\bigr\rvert\bigr).
  \end{align*}
  Now let's prove \eqref{eq:dbound}. By \cref{thm:dastexpr}, we observe that
  \begingroup\footnotesize
  \begin{align*}
                                  & \bigl\lVert [\eta_{1},\cdots,\eta_{\ell}]\triangleright \mathcal{D}(Y_{1},\cdots,Y_{t}) \bigr\rVert _{q}                                                                                                                                                    \\
    =\ \                          & \bigl\lVert \mathcal{D}\bigl(Y_{1}\cdots Y_{\eta_{1}},Y_{\eta_{1}+1}\cdots Y_{\eta_{1}+\eta_{2}},\ldots,Y_{\eta_{1}+\cdots+\eta_{\ell-1}+1}\cdots Y_{t}\bigr) \bigr\rVert _{q} \\
    \overset{\eqref{eq:dexpr}}{=} & \biggl\lVert Y_{1}\cdots Y_{t}
      -\mathcal{D}^{*}\bigl(Y_{1}\cdots Y_{\eta_{1}},\ldots,
       Y_{\eta_{1}+\cdots+\eta_{\ell-1}+1}\cdots Y_{t}\bigr) \\
                                  &\quad -\sum_{j=1}^{\ell-1}Y_{1}\cdots Y_{\eta_{1}+\cdots+\eta_{j}}\cdot \mathcal{D}^{*}\Bigl(Y_{\eta_{1}+\cdots +\eta_{j}+1}\cdots Y_{\eta_{1}+\cdots+\eta_{j+1}},\ldots,Y_{\eta_{1}+\cdots+\eta_{\ell-1}+1}\cdots Y_{t}\Bigr)\biggr\rVert _{q} \\
    = \ \                         & \biggl\lVert Y_{1}Y_{2}\cdots Y_{t}
      -[\eta_{1},\cdots,\eta_{\ell}]\triangleright \mathcal{D}^{*}(Y_{1},\cdots,Y_{t})\\
                                  &\quad
      -\sum_{j=1}^{\ell-1}Y_{1}\cdots Y_{\eta_{1}+\cdots+\eta_{j}}\cdot
      \sum_{(s,\zeta_{1:s})\in C(\ell-j)}(-1)^{s-1}
      [\zeta_{1},\cdots,\zeta_{s}]\triangleright\\
                                  &\hspace{5em}\bigl(Y_{\eta_{1}+\cdots +\eta_{j}+1}\cdots Y_{\eta_{1}+\cdots+\eta_{j+1}},\ldots,Y_{\eta_{1}+\cdots+\eta_{\ell-1}+1}\cdots Y_{t}\bigr)\biggr\rVert _{q}
  \end{align*}
  \endgroup
  We upper-bound this using the triangle inequality. Indeed, we obtain that
  \begingroup\small
  \begin{align*}
                                          & \bigl\lVert [\eta_{1},\cdots,\eta_{\ell}]\triangleright \mathcal{D}(Y_{1},\cdots,Y_{t}) \bigr\rVert _{q}
    \\
    \leq \ \                              & \bigl\lvert [\eta_{1},\cdots,\eta_{\ell}]\triangleright \mathcal{D}^{*}(Y_{1},\cdots,Y_{t}) \bigr\rvert+\lVert Y_{1}Y_{2}\cdots Y_{t} \rVert_{q}\\*
    &\ +\sum_{j=1}^{\ell-1}\sum_{(s,\zeta_{1:s}) \in C (\ell-j)}\lVert Y_{1}\cdots Y_{\eta_{1}+\cdots+\eta_{j}} \rVert_{q}\cdot                             \\*
                                          & \ \ \Bigl\lvert[\zeta_{1},\cdots,\zeta_{s}]\triangleright
    \bigl(Y_{\eta_{1}+\cdots +\eta_{j}+1}\cdots Y_{\eta_{1}+\cdots+\eta_{j+1}},\cdots,Y_{\eta_{1}+\cdots+\eta_{\ell-1}+1}\cdots Y_{t}\bigr)\Bigr\rvert                                                                                                                                                                                           \\
    \overset{(*)}{\leq}\ \                & \bigl\lvert [\eta_{1},\cdots,\eta_{\ell}]\triangleright \mathcal{D}^{*}(Y_{1},\cdots,Y_{t}) \bigr\rvert +\lVert Y_{1}Y_{2}\cdots Y_{t} \rVert_{q}\\*
    &\ +\sum_{j=1}^{\ell-1}\sum_{(s,\zeta_{1:s})\in C (\ell-j)}\lVert Y_{1}\cdots Y_{\eta_{1}+\cdots+\eta_{j}} \rVert_{q}\cdot                             \\*
                                          & \ \ \Bigl([\zeta_{1},\cdots,\zeta_{s}]\triangleright
    \bigl(\lvert Y_{\eta_{1}+\cdots +\eta_{j}+1}\cdots Y_{\eta_{1}+\cdots+\eta_{j+1}}\rvert^{q},\cdots,\lvert Y_{\eta_{1}+\cdots+\eta_{\ell-1}+1}\cdots Y_{t}\rvert^{q}\bigr)\Bigr)^{1/q}                                                                                                                                                        \\
    \overset{(**)}{\leq}\                 & \bigl\lvert [\eta_{1},\cdots,\eta_{\ell}]\triangleright \mathcal{D}^{*}(Y_{1},\cdots,Y_{t}) \bigr\rvert+\lVert Y_{1}Y_{2}\cdots Y_{t} \rVert_{q}\\*
    &\ +\sum_{j=1}^{\ell-1}\sum_{(s,\lambda_{1:s})\in C (\eta_{j+1}+\cdots+\eta_{\ell})}\lVert Y_{1}\cdots Y_{\eta_{1}+\cdots +\eta_{j}} \rVert_{q}\cdot    \\*
                                          & \ \  \bigl([\lambda_{1},\cdots,\lambda_{s}]\triangleright \bigl(\lvert Y_{\eta_{1}+\cdots+\eta_{j}+1}\rvert^{q},\cdots,\lvert Y_{t}\rvert^{q}\bigr)\bigr)^{1/q}                                                                                                                                     \\
    \leq \ \                              & \bigl\lvert [\eta_{1},\cdots,\eta_{\ell}]\triangleright \mathcal{D}^{*}(Y_{1},\cdots,Y_{t}) \bigr\rvert+\lVert Y_{1}Y_{2}\cdots Y_{t} \rVert_{q}+\sum_{h=1}^{t-1}\sum_{(s,\lambda_{1:s})\in C (t-h)}\lVert Y_{1}\cdots Y_{h} \rVert_{q}\cdot                                                         \\*
                                          & \quad \bigl([\lambda_{1},\cdots,\lambda_{s}]\triangleright \bigl(\lvert Y_{h+1}\rvert^{q},\cdots,\lvert Y_{t}\rvert^{q}\bigr)\bigr)^{1/q}                                                                                                                                                            \\
    =\ \                                  & \bigl\lvert [\eta_{1},\cdots,\eta_{\ell}]\triangleright \mathcal{D}^{*}(Y_{1},\cdots,Y_{t}) \bigr\rvert+\sum_{(s,\lambda_{1:s})\in C (t)}\Bigl([\lambda_{1},\cdots,\lambda_{s}]\triangleright \bigl(\lvert Y_{1} \rvert^{q},\cdots,\lvert Y_{t}\rvert^{q}\bigr)\Bigr)^{1/q}                          \\
    \overset{\eqref{eq:dastbound1}}{\leq} & \sum_{(s,\zeta_{1:s})\in C (t)}\ [\zeta_{1},\cdots,\zeta_{s}]\triangleright \bigl(\lvert Y_{1} \rvert,\cdots,\lvert Y_{t} \rvert\bigr)\\*
    &\ +\sum_{(s,\zeta_{1:s})\in C (t)}\Bigl([\zeta_{1},\cdots,\zeta_{s}]\triangleright \bigl(\lvert Y_{1} \rvert^{q},\cdots,\lvert Y_{t}\rvert^{q}\bigr)\Bigr)^{1/q} \\
    \leq \ \                              & 2\sum_{(s,\zeta_{1:s})\in C(t)}\Bigl([\zeta_{1},\cdots,\zeta_{s}]\triangleright\bigl(\lvert Y_{1} \rvert^{q},\cdots,\lvert Y_{t}\rvert^{q}\bigr)\Bigr)^{1/q},
  \end{align*}
  \endgroup
  where to obtain $(*)$ we have used the fact that by Jensen inequality for any random variable $X\in \mathcal{L}^q(\mathbb{R})$ we have $\bigl\lvert \mathbb{E} [X] \bigr\rvert\leq \bigl(\mathbb{E} [\lvert X \rvert^{q}]\bigr)^{1/q}$; and where to obtain $(**)$ we have used the fact that
  \begin{align*}
      & [\lambda_{1},\cdots,\lambda_{s}]\triangleright \bigl(\lvert Y_{\eta_{1}+\cdots+\eta_{j}+1}\rvert^{q},\cdots,\lvert Y_{t}\rvert^{q}\bigr)                                                                              \\
    = & [\zeta_{1},\cdots,\zeta_{s}]\triangleright \bigl(\lvert Y_{\eta_{1}+\cdots +\eta_{j}+1}\cdots Y_{\eta_{1}+\cdots+\eta_{j+1}}\rvert^{q},\cdots,\lvert Y_{\eta_{1}+\cdots+\eta_{\ell-1}+1}\cdots Y_{t}\rvert^{q}\bigr).
  \end{align*}\hfill
\end{proof}

\begin{lemma}\label{thm:bounddstarmix}
  Let $(X_{i})_{i\in T}$ be mean-zero random fields of random variables, and let $S_{T_{0},\omega}$ be a random variable that satisfies $\lvert S_{T_{0},\omega} \rvert\leq \bigl\lvert \sum_{i\in T_{0}}X_{i} \bigr\rvert^{\omega}$, where $T_{0}\subseteq T$ and $0\leq \omega\leq 1$. When $\omega=0$, this condition is interpreted as $|S_{T_0,0}|\leq1$, and factors $|T_0|^\omega$ below are interpreted as $1$. Fix $t\in \mathbb{N}_{+}$ and $i_{1:t}\in T$. In each assertion below, assume $\sup_i\lVert X_i\rVert_r\leq M<\infty$ at the indicated moment order.

  For any $\ell,\eta_{1:\ell}\in \mathbb{N}_{+}$ such that $\eta_{1}+\cdots+\eta_{\ell}=t+1$, the first bound below holds for $r\in[t+1,\infty]$, and the second holds for $r\in[t+\omega,\infty]$:
  \begin{align}
     & \bigl\lvert [\eta_{1},\cdots,\eta_{\ell}]\triangleright \mathcal{D}^{*}(X_{i_{1}},\cdots,X_{i_{t}},X_{i_{t+1}}) \bigr\rvert\leq 2^{t}M^{t+1},\label{eq:dstarnomix1}                                       \\
     & \bigl\lvert [\eta_{1},\cdots,\eta_{\ell}]\triangleright \mathcal{D}^{*}(X_{i_{1}},\cdots,X_{i_{t}},S_{T_{0},\omega}) \bigr\rvert\leq 2^{t}\lvert T_{0} \rvert^{\omega}M^{t+\omega}.\label{eq:dstarnomix2}
  \end{align}

  We further let $j$ be an integer that satisfies $2\leq j\leq t+1$ and define the $\sigma$-algebras $\mathcal{F}_{1}$ and $\mathcal{F}_{2}$ by
  \begin{gather*}
    \mathcal{F}_{1}:=\sigma \bigl(X_{i_{1}},\cdots,X_{i_{j-1}}\bigr),\quad \mathcal{F}_{2}:=\sigma \bigl(X_{i_{j}},\cdots,X_{i_{t+1}}\bigr),\\
    \mathcal{F}_{3}:=
    \begin{cases}
      \sigma \bigl(X_{i_{j}},\cdots,X_{i_{t}},S_{T_{0},\omega}\bigr) & \text{ if }j\leq t \\
      \sigma (S_{T_{0},\omega})                                      & \text{ if }j=t+1
    \end{cases}.
  \end{gather*}
  For any $s,\ell,\eta_{1:\ell}\in \mathbb{N}_{+}$ such that $\eta_{1}+\cdots+\eta_{s}=j-1$ and $\eta_{1}+\cdots+\eta_{\ell}=t+1$, the first bound below holds for $r\in(t+1,\infty]$, and the second holds for $r\in(t+\omega,\infty]$:
  \begin{align}
     & \bigl\lvert [\eta_{1},\cdots,\eta_{\ell}]\triangleright \mathcal{D}^{*}(X_{i_{1}},\cdots,X_{i_{t}},X_{i_{t+1}}) \bigr\rvert\leq 2^{t+3}\bigl(\alpha (\mathcal{F}_{1},\mathcal{F}_{2})\bigr)^{(r-t-1)/r}M^{t+1},\label{eq:dstarmix1}                                           \\
     & \bigl\lvert [\eta_{1},\cdots,\eta_{\ell}]\triangleright \mathcal{D}^{*}(X_{i_{1}},\cdots,X_{i_{t}},S_{T_{0},\omega}) \bigr\rvert\leq 2^{t+3}\lvert T_{0} \rvert^{\omega}\bigl(\alpha(\mathcal{F}_{1},\mathcal{F}_{3})\bigr)^{(r-t-\omega)/r}M^{t+\omega},\label{eq:dstarmix2}
  \end{align}
  where $\alpha (\mathcal{F}_{1},\mathcal{F}_{2})$ is the $\alpha$-mixing coefficients between $\mathcal{F}_{1}$ and $\mathcal{F}_{2}$, and $\alpha (\mathcal{F}_{1},\mathcal{F}_{3})$ is the $\alpha$-mixing coefficients between $\mathcal{F}_{1}$ and $\mathcal{F}_{3}$.
  At the endpoint moment orders, there is also an exact-independence
  conclusion: the left-hand side of \eqref{eq:dstarmix1} is zero when
  $r=t+1$ and $\mathcal F_1$ and $\mathcal F_2$ are independent, and the
  left-hand side of \eqref{eq:dstarmix2} is zero when $r=t+\omega$ and
  $\mathcal F_1$ and $\mathcal F_3$ are independent.
\end{lemma}

\begin{proof}[Proof of \cref{thm:bounddstarmix}.]
  To prove \eqref{eq:dstarnomix1}, we remark that
  \begin{align*}
                             & \Bigl\lvert [\eta_{1},\cdots,\eta_{\ell}]\triangleright \mathcal{D}^{*}\bigl(X_{i_{1}},\cdots,X_{i_{t+1}}\bigr)\Bigr\rvert                                   \\
    \overset{\eqref{eq:dastbound1}}{\leq}
                             & \sum_{(\ell',\zeta_{1:\ell'})\in C (t+1)}\ [\zeta_{1},\cdots,\zeta_{\ell'}]\triangleright \bigl(\lvert X_{i_{1}}\rvert,\cdots,\lvert X_{i_{t+1}}\rvert\bigr) \\
    \overset{(*)}{\leq} \ \  & 2^{t}\frac{1}{t+1}\bigl(\bigl\lVert X_{i_{1}} \bigr\rVert _{t+1}^{t+1}+\cdots+\bigl\lVert X_{i_{t+1}} \bigr\rVert _{t+1}^{t+1}\bigr)\leq 2^{t}M^{t+1},
  \end{align*}
  where $(*)$ is implied by \cref{thm:lemmayoung} and the fact that $\lvert C(t+1) \rvert=2^{t}$ \citep{heubach2009combinatorics}.

  For \eqref{eq:dstarnomix2}, we have
  \begin{align*}
                                          & \Bigl\lvert [\eta_{1},\cdots,\eta_{\ell}]\triangleright \mathcal{D}^{*}\bigl(X_{i_{1}},\cdots,X_{i_{t}}, S_{T_{0},\omega}\bigr)\Bigr\rvert                                                                   \\
    \overset{\eqref{eq:dastbound1}}{\leq} & \sum_{(\ell',\zeta_{1:\ell'})\in C  (t+1)}[\zeta_{1},\cdots,\zeta_{\ell'}]\triangleright \bigl(\lvert X_{i_{1}} \rvert,\cdots,\lvert X_{i_{t}} \rvert,\lvert S_{T_{0},\omega}\rvert\bigr)                    \\
  \end{align*}
  If $\omega=0$, omit the bounded factor $|S_{T_0,0}|\leq1$ from each
  bracket. Applying \cref{thm:lemmayoung} to the remaining $t$ variables with
  all exponents equal to $t$ (and arguing directly when $t=1$), each summand
  is at most $M^t$. Since $|C(t+1)|=2^t$, this proves
  \eqref{eq:dstarnomix2} in this case.

  Suppose now that $\omega>0$. If $T_0=\varnothing$, then
  $S_{T_0,\omega}=0$ and all the asserted bounds are immediate. Hence, in all
  positive-$\omega$ arguments below, we may assume that $T_0\neq\varnothing$.
  Then
  \begin{align*}
    &\Bigl\lvert [\eta_{1},\cdots,\eta_{\ell}]\triangleright
      \mathcal{D}^{*}\bigl(X_{i_{1}},\cdots,X_{i_{t}},S_{T_{0},\omega}\bigr)
      \Bigr\rvert\\
    &\quad\leq \lvert  T_{0}\rvert^{\omega}
    \sum_{(\ell',\zeta_{1:\ell'})\in C(t+1)}
    [\zeta_{1},\cdots,\zeta_{\ell'}]\triangleright\\*
    &\quad\biggl(\lvert X_{i_{1}} \rvert,\cdots,\lvert X_{i_{t}} \rvert,
    \biggl\lvert\frac{1}{\lvert T_{0}\rvert}
    \sum_{i\in T_{0}}X_{i} \biggr\rvert^{\omega}\biggr) \\
    \overset{(*)}{\leq}\ \                & \lvert T_{0}\rvert^{\omega}\!\!\!
    \sum_{(\ell',\zeta_{1:\ell'})\in C(t+1)}
    \Biggl(\frac{1}{t+\omega}\Bigl(
    \bigl\lVert X_{i_{1}} \bigr\rVert _{t+\omega}^{t+\omega}+\cdots
    +\bigl\lVert X_{i_{t}} \bigr\rVert _{t+\omega}^{t+\omega}\Bigr)\\*
    &\quad+\frac{\omega}{t+\omega}
    \biggl\lVert \frac{1}{\lvert T_{0} \rvert}
    \sum_{i\in T_{0}}X_{i} \biggr\rVert _{t+\omega}^{t+\omega}\Biggr) \\
    \overset{(**)}{\leq}\                 & 2^{t}\lvert T_{0} \rvert^{\omega}
    \Biggl(\frac{1}{t+\omega}\bigl(
    \bigl\lVert X_{i_{1}} \bigr\rVert _{t+\omega}^{t+\omega}+\cdots
    +\bigl\lVert X_{i_{t}} \bigr\rVert _{t+\omega}^{t+\omega}\bigr)\\*
    &\quad+\frac{\omega}{t+\omega}\frac{1}{\lvert T_{0} \rvert}
    \sum_{i\in T_{0}}\bigl\lVert X_{i} \bigr\rVert _{t+\omega}^{t+\omega}
    \Biggr) \\
    \leq \ \                              & 2^{t}\lvert  T_{0} \rvert^{\omega}\cdot M^{t+\omega}.
  \end{align*}
  Here, for $\omega>0$, we have used \cref{thm:lemmayoung} in $(*)$, and $(**)$ is implied by $\lvert C(t+1) \rvert\leq 2^{t}$ and Jensen's inequality as
  \begin{equation*}
    \begin{aligned}
      &\textstyle \bigl\lVert \frac{1}{\lvert T_{0} \rvert} \sum_{i\in T_{0}}X_{i}\bigr\rVert _{t+\omega}^{t+\omega}=\mathbb{E} \bigl[\bigl\lvert\frac{1}{\lvert T_{0} \rvert}\sum_{i\in T_{0}}X_{i}\bigr\rvert^{t+\omega}\bigr]\\
    \leq &\textstyle \mathbb{E} \bigl[\frac{1}{\lvert T_{0} \rvert}\sum_{i\in T_{0}}\lvert X_{i} \rvert^{t+\omega}\bigr]=\frac{1}{\lvert T_{0} \rvert}\sum_{i\in T_{0}}\mathbb{E} [\lvert X_{i} \rvert^{t+\omega}]=\frac{1}{\lvert T_{0} \rvert}\sum_{i\in T_{0}}\lVert X_{i} \rVert_{t+\omega}^{t+\omega}.
    \end{aligned}
  \end{equation*}

  To show \eqref{eq:dstarmix1}, we remark that assumption we have that $s$ is such that  $\eta_{1}+\cdots+\eta_{s}=j-1$. Therefore, according to \eqref{eq:strangeeq}, we get that
  \begin{align}\label{eq:thinkofaname1}
      & [\eta_{1},\cdots,\eta_{\ell}]\triangleright \mathcal{D}^{*}\bigl(X_{i_{1}},\cdots, X_{i_{t+1}}\bigr)                                                                                                       \\
    = & [\eta_{1},\cdots,\eta_{s},\cdots,\eta_{\ell}]\triangleright \mathcal{D}^{*}\bigl(X_{i_{1}},\cdots,X_{i_{j}},\cdots,X_{i_{t+1}}\bigr)\nonumber                                                                       \\
    = & [\eta_{1},\cdots,\eta_{s-1},\eta_{s}+1]\triangleright \mathcal{D}^{*}\Bigl(X_{i_{1}},\cdots,X_{i_{j-1}}\ ,\ [\eta_{s+1},\cdots,\eta_{\ell}]\triangleright\mathcal{D}\bigl(X_{i_{j}},\cdots,X_{i_{t+1}}\bigr)\Bigr).\nonumber
  \end{align}

  Moreover, by exploiting \eqref{eq:dastbound2}, we obtain that 
  \begingroup\small
  \begin{equation}\label{eq:thinkofaname2}
    \begin{aligned}
           & \biggl\lvert [\eta_{1},\cdots,\eta_{s-1},\eta_{s}+1]\triangleright \mathcal{D}^{*}\Bigl(X_{i_{1}},\cdots,X_{i_{j-1}}\ ,\ [\eta_{s+1},\cdots,\eta_{\ell}]\triangleright\mathcal{D}\bigl(X_{i_{j}},\cdots,X_{i_{t+1}}\bigr)\Bigr)\biggr\rvert                                                                                                                       \\
      \leq & \Bigl\lvert\mathbb{E} \bigl[ X_{i_{1}}\cdots X_{i_{j-1}} \!\!\cdot [\eta_{s+1},\cdots,\eta_{\ell}]\triangleright\mathcal{D}\bigl(X_{i_{j}},\cdots,X_{i_{t+1}}\bigr)\bigr] \Bigr\rvert                                                                                                                                                                             \\
           & \ +\sum_{w=1}^{j-2}\sum_{(s',\zeta_{1:s'})\in  C (w)}\ [\zeta_{1},\cdots,\zeta_{s'}]\triangleright \bigl(\lvert X_{i_{1}} \rvert,\cdots,\lvert X_{i_{w}} \rvert\bigr)\cdot \\
           &\qquad\Bigl\lvert\mathbb{E} \bigl[ X_{i_{w+1}}\cdots X_{i_{j-1}} \!\!\cdot [\eta_{s+1},\cdots,\eta_{\ell}]\triangleright\mathcal{D}\bigl(X_{i_{j}},\cdots,X_{i_{t+1}}\bigr)\bigr]\Bigr\rvert
    \end{aligned}
  \end{equation}
  \endgroup
  For $w=1$, the following bound is immediate from
  $\mathbb E|X_{i_1}|\leq M$; for $w\geq2$, it follows from
  \cref{thm:lemmayoung} with all exponents equal to $w$:
  \begin{equation}\label{eq:dontlikelabel7}
    [\zeta_{1},\cdots,\zeta_{s'}]\triangleright \bigl(\lvert X_{i_{1}} \rvert,\cdots,\lvert X_{i_{w}} \rvert\bigr)\leq \frac{1}{w}\bigl(\lVert X_{i_{1}} \rVert_{w}^{w}+\cdots+\lVert X_{i_{w}} \rVert_{w}^{w}\bigr)\leq M^{w}.
  \end{equation}
  Combining \eqref{eq:thinkofaname1}, \eqref{eq:thinkofaname2}, and \eqref{eq:dontlikelabel7}, we get
  \begin{align}\label{eq:someneedslabel3}
                     & \Bigl\lvert [\eta_{1},\cdots,\eta_{\ell}]\triangleright \mathcal{D}^{*}\bigl(X_{i_{1}},\cdots, X_{i_{t+1}}\bigr)\Bigr\rvert                                                                                                                        \\
    \leq             & \Bigl\lvert \mathbb{E} \bigl[ X_{i_{1}}\cdots X_{i_{j-1}} \!\!\cdot [\eta_{s+1},\cdots,\eta_{\ell}]\triangleright\mathcal{D}\bigl(X_{i_{j}},\cdots,X_{i_{t+1}}\bigr)\bigr]\Bigr\rvert\nonumber                                                              \\*
                     & \ +\sum_{w=1}^{j-2}\sum_{(s',\zeta_{1:s'})\in C (w)}\ M^{w} \Bigl\lvert\mathbb{E} \bigl[ X_{i_{w+1}}\cdots X_{i_{j-1}} \!\!\cdot [\eta_{s+1},\cdots,\eta_{\ell}]\triangleright\mathcal{D}\bigl(X_{i_{j}},\cdots,X_{i_{t+1}}\bigr)\bigr]\Bigr\rvert\nonumber \\
    \overset{(*)}{=} & \sum_{w=0}^{j-2}2^{(w-1)\vee 0}M^{w}\Bigl\lvert\mathbb{E} \bigl[ X_{i_{w+1}}\cdots X_{i_{j-1}} \!\!\cdot [\eta_{s+1},\cdots,\eta_{\ell}]\triangleright\mathcal{D}\bigl(X_{i_{j}},\cdots,X_{i_{t+1}}\bigr)\bigr]\Bigr\rvert\nonumber                         \\
    =                & \sum_{w=1}^{j-1}2^{(w-2)\vee 0}M^{w-1}\Bigl\lvert\mathbb{E} \bigl[ X_{i_{w}}\cdots X_{i_{j-1}} \!\!\cdot [\eta_{s+1},\cdots,\eta_{\ell}]\triangleright\mathcal{D}\bigl(X_{i_{j}},\cdots,X_{i_{t+1}}\bigr)\bigr]\Bigr\rvert.\nonumber
  \end{align}
  where $(*)$ is due to the fact that $\lvert C  (w)\rvert= 2^{w-1}$.

  As mentioned in \eqref{mean0cd}, by definition of the compositional $\mathcal{D}$ operator we know that 
  $$\mathbb{E} \bigl[[\eta_{s+1},\cdots,\eta_{\ell}]\triangleright\mathcal{D}\bigl(X_{i_{j}},\cdots,X_{i_{t+1}}\bigr)\bigr]=0.$$ Thus, we can apply \cref{thm:covineq}:
  \begin{equation}\label{eq:lemma125pf}
    \begin{aligned}
           & \Bigl\lvert \mathbb{E} \bigl[ X_{i_{w}}\cdots X_{i_{j-1}} \!\!\cdot [\eta_{s+1},\cdots,\eta_{\ell}]\triangleright\mathcal{D}\bigl(X_{i_{j}},\cdots,X_{i_{t+1}}\bigr)\bigr] \Bigr\rvert                                                                                        \\
      \leq & 8\bigl(\alpha (\mathcal{F}_{1},\mathcal{F}_{2})\bigr)^{(r-t+w-2)/r}\bigl\lVert X_{i_{w}}\cdots X_{i_{j-1}} \bigr\rVert _{r/(j-w)}\cdot \\*
      &\ \bigl\lVert [\eta_{s+1},\cdots,\eta_{\ell}]\triangleright\mathcal{D}\bigl(X_{i_{j}},\cdots,X_{i_{t+1}}\bigr) \bigr\rVert _{r/(t-j+2)}.
    \end{aligned}
  \end{equation}

  By \cref{thm:lemmayoung} when $j-w\geq2$ (and directly when $j-w=1$),
  we get
  \begin{equation}\label{eq:dontlikelabel2}
    \begin{aligned}
      \bigl\lVert X_{i_{w}}\cdots X_{i_{j-1}} \bigr\rVert _{r/(j-w)}= & \Bigl(\mathbb{E} \Bigl[\bigl\lvert X_{i_{w}}\cdots X_{i_{j-1}} \bigr\rvert^{r/(j-w)}\Bigr]\Bigr)^{(j-w)/r}                        \\
      \leq                                                            & \Bigl(\frac{1}{j-w}\bigl(\lVert X_{i_{w}} \rVert_{r}^{r}+\cdots+\lVert X_{i_{j-1}} \rVert_{r}^{r}\bigr)\Bigr)^{(j-w)/r}\leq M^{j-w}.
    \end{aligned}
  \end{equation}

  Moreover, remark that
  \begingroup\small
  \begin{align}\label{eq:someneedslabel2}
                                      & \bigl\lVert [\eta_{s+1},\cdots,\eta_{\ell}]\triangleright\mathcal{D}\bigl(X_{i_{j}},\cdots,X_{i_{t+1}}\bigr)\bigr\rVert _{r/(t-j+2)}\nonumber                                                                          \\
    \overset{\eqref{eq:dbound}}{\leq} & 2\sum_{(\ell',\zeta_{1:\ell'})\in C(t-j+2)}\Bigl([\zeta_{1},\cdots,\zeta_{\ell'}]\triangleright \bigl(\lvert X_{i_{j}} \rvert^{r/(t-j+2)},\cdots,\lvert X_{i_{t+1}} \rvert^{r/(t-j+2)}\bigr)\Bigr)^{(t-j+2)/r}\nonumber \\
    \overset{(*)}{\leq}\ \            & 2^{t-j+2}\Bigl(\frac{1}{t-j+2}\bigl(\lVert X_{i_{j}} \rVert_{r}^{r}+\cdots+\lVert X_{i_{t+1}} \rVert_{r}^{r}\bigr)\Bigr)^{(t-j+2)/r}\leq 2^{t-j+2} M^{t-j+2}.
  \end{align}
  \endgroup
  Note that $(*)$ is implied by the fact that
  $\lvert C(t-j+2)\rvert=2^{t-j+1}$ and \cref{thm:lemmayoung} when
  $t-j+2\geq2$; the one-variable case is immediate.

  Substituting \eqref{eq:dontlikelabel2} and \eqref{eq:someneedslabel2} into \eqref{eq:lemma125pf}, we get
  \begin{align*}
    &\Bigl\lvert \mathbb{E} \bigl[ X_{i_{w}}\cdots X_{i_{j-1}} \!\!\cdot [\eta_{s+1},\cdots,\eta_{\ell}]\triangleright\mathcal{D}\bigl(X_{i_{j}},\cdots,X_{i_{t+1}}\bigr)\bigr] \Bigr\rvert\\
    \leq & 2^{t-j+5}\bigl(\alpha (\mathcal{F}_{1},\mathcal{F}_{2})\bigr)^{(r-t+w-2)/r} M^{t-w+2}.
  \end{align*}
  Combining this and \eqref{eq:someneedslabel3}, we get
  \begin{align*}
         & \Bigl\lvert [\eta_{1},\cdots,\eta_{\ell}]\triangleright \mathcal{D}^{*}\bigl(X_{i_{1}},\cdots,X_{i_{t+1}}\bigr)\Bigr\rvert \\
    \leq & \sum_{w=1}^{j-1}2^{(w-2)\vee 0}\cdot 2^{t-j+5}\bigl(\alpha (\mathcal{F}_{1},\mathcal{F}_{2})\bigr)^{(r-t+w-2)/r}M^{t+1}    \\
    \leq & 2^{t+3}\bigl(\alpha (\mathcal{F}_{1},\mathcal{F}_{2})\bigr)^{(r-t-1)/r}M^{t+1}.
  \end{align*}

  Lastly, we prove \eqref{eq:dstarmix2}. We consider two cases, $2\leq j\leq t$ and $j=t+1$.

  If $2\leq j\leq t$, by \eqref{eq:strangeeq}, we have
  \begingroup\small
  \begin{align}\label{eq:dontlikelabel5}
      & [\eta_{1},\cdots,\eta_{\ell}]\triangleright \mathcal{D}^{*}\bigl(X_{i_{1}},\cdots,X_{i_{t}},S_{T_{0},\omega}\bigr)                                                                                                                             \\
    = & [\eta_{1},\cdots,\eta_{s},\cdots,\eta_{\ell}]\triangleright \mathcal{D}^{*}\bigl(X_{i_{1}},\cdots,X_{i_{j}},\cdots,X_{i_{t}},S_{T_{0},\omega}\bigr)\nonumber                                                                                                \\
    = & [\eta_{1},\cdots,\eta_{s-1},\eta_{s}+1]\triangleright \mathcal{D}^{*}\Bigl(X_{i_{1}},\cdots,X_{i_{j-1}}\ ,\ [\eta_{s+1},\cdots,\eta_{\ell}]\triangleright\mathcal{D}\bigl(X_{i_{j}},\cdots,X_{i_{t}},S_{T_{0},\omega}\bigr)\Bigr).\nonumber
  \end{align}
  \endgroup
  By \eqref{eq:dastbound2}, we get
  \begingroup\small
  \begin{equation}\label{eq:dontlikelabel6}
    \begin{aligned}
           & \biggl\lvert [\eta_{1},\cdots,\eta_{s-1},\eta_{s}+1]\triangleright \mathcal{D}^{*}\Bigl(X_{i_{1}},\cdots,X_{i_{j-1}}\ ,\ [\eta_{s+1},\cdots,\eta_{\ell}]\triangleright\mathcal{D}\bigl(X_{i_{j}},\cdots,X_{i_{t}},S_{T_{0},\omega}\bigr)\Bigr)\biggr\rvert                                                                                                                          \\
      \leq & \Bigl\lvert\mathbb{E} \bigl[ X_{i_{1}}\cdots X_{i_{j-1}} \!\!\cdot [\eta_{s+1},\cdots,\eta_{\ell}]\triangleright\mathcal{D}\bigl(X_{i_{j}},\cdots,X_{i_{t}},S_{T_{0},\omega}\bigr)\bigr]\Bigr\rvert                                                                                                                                                                              \\*
           & \ +\sum_{w=1}^{j-2}\sum_{(s',\zeta_{1:s'})\in C (w)}\ [\zeta_{1},\cdots,\zeta_{s'}]\triangleright \bigl(\lvert X_{i_{1}} \rvert,\cdots,\lvert X_{i_{w}} \rvert\bigr)\cdot \\
           &\qquad \Bigl\lvert\mathbb{E} \bigl[ X_{i_{w+1}}\cdots X_{i_{j-1}} \!\!\cdot [\eta_{s+1},\cdots,\eta_{\ell}]\triangleright\mathcal{D}\bigl(X_{i_{j}},\cdots,X_{i_{t}},S_{T_{0},\omega}\bigr)\bigr]\Bigr\rvert.
    \end{aligned}
  \end{equation}
  \endgroup

  Combining \eqref{eq:dontlikelabel5}, \eqref{eq:dontlikelabel6}, and \eqref{eq:dontlikelabel7}, we have
  \begingroup\small
  \begin{align}\label{eq:dontlikelabel4}
                     & \Bigl\lvert [\eta_{1},\cdots,\eta_{\ell}]\triangleright \mathcal{D}^{*}\bigl(X_{i_{1}},\cdots,X_{i_{t}},S_{T_{0},\omega}\bigr)\Bigr\rvert                                                                                                                      \\
    \leq             & \Bigl\lvert\mathbb{E} \bigl[ X_{i_{1}}\cdots X_{i_{j-1}} \!\!\cdot [\eta_{s+1},\cdots,\eta_{\ell}]\triangleright\mathcal{D}\bigl(X_{i_{j}},\cdots,X_{i_{t}},S_{T_{0},\omega}\bigr)\bigr]\Bigr\rvert\nonumber                                                              \\*
                     & \ +\sum_{w=1}^{j-2}\sum_{(s',\zeta_{1:s'})\in C (w)}\ M^{w}\Bigl\lvert\mathbb{E} \bigl[ X_{i_{w+1}}\cdots X_{i_{j-1}} \!\!\cdot [\eta_{s+1},\cdots,\eta_{\ell}]\triangleright\mathcal{D}\bigl(X_{i_{j}},\cdots,X_{i_{t}},S_{T_{0},\omega}\bigr)\bigr]\Bigr\rvert\nonumber \\
    \overset{(*)}{=} & \sum_{w=1}^{j-1}2^{(w-2)\vee 0}M^{w-1}\Bigl\lvert\mathbb{E} \bigl[ X_{i_{w}}\cdots X_{i_{j-1}} \!\!\cdot [\eta_{s+1},\cdots,\eta_{\ell}]\triangleright\mathcal{D}\bigl(X_{i_{j}},\cdots,X_{i_{t}},S_{T_{0},\omega}\bigr)\bigr]\Bigr\rvert,\nonumber  
  \end{align}
  \endgroup
  where $(*)$ is due to the fact that $\lvert C(w) \rvert=2^{w-1}$.

  We apply \cref{thm:covineq} and obtain
  \begin{equation}\label{eq:dontlikelabel1}
    \begin{aligned}
           & \Bigl\lvert \mathbb{E} \bigl[ X_{i_{w}}\cdots X_{i_{j-1}} \!\!\cdot [\eta_{s+1},\cdots,\eta_{\ell}]\triangleright\mathcal{D}\bigl(X_{i_{j}},\cdots,X_{i_{t}},S_{T_{0},\omega}\bigr)\bigr]\Bigr\rvert                                                                                                       \\
      \leq & 8\bigl(\alpha (\mathcal{F}_{1},\mathcal{F}_{3})\bigr)^{(r-t+w-1-\omega)/r}\bigl\lVert X_{i_{w}}\cdots X_{i_{j-1}} \bigr\rVert _{r/(j-w)}\cdot\\*
      &\  \bigl\lVert [\eta_{s+1},\cdots,\eta_{\ell}]\triangleright\mathcal{D}\bigl(X_{i_{j}},\cdots,X_{i_{t}},S_{T_{0},\omega}\bigr) \bigr\rVert _{r/(t+1+\omega-j)}.
    \end{aligned}
  \end{equation}

  If $\omega=0$, then $|S_{T_0,0}|\leq1$. Applying \eqref{eq:dbound},
  omitting this bounded factor in every bracket, and applying
  \cref{thm:lemmayoung} to the $t+1-j$ remaining variables with equal
  exponents (directly when $t+1-j=1$), gives
  \begin{equation*}
    \bigl\lVert [\eta_{s+1},\cdots,\eta_{\ell}]\triangleright
    \mathcal{D}\bigl(X_{i_j},\ldots,X_{i_t},S_{T_0,0}\bigr)\bigr\rVert_{r/(t+1-j)}
    \leq 2^{t+2-j}M^{t+1-j}.
  \end{equation*}

  For $\omega>0$, we observe that
  \begin{align}\label{eq:dontlikelabel3}
                                      & \bigl\lVert [\eta_{s+1},\cdots,\eta_{\ell}]\triangleright\mathcal{D}\bigl(X_{i_{j}},\cdots,X_{i_{t}},S_{T_{0},\omega}\bigr) \bigr\rVert _{r/(t+1+\omega-j)}                                                                                                                                                                                                                              \\
    \overset{\eqref{eq:dbound}}{\leq} & 2\lvert T_{0} \rvert^{\omega}\!\!\sum_{(\ell',\lambda_{1:\ell'})\in C  (t+2-j)}\Biggl([\lambda_{1},\cdots,\lambda_{\ell'}]\triangleright \nonumber\\*
    &\quad\biggl(\lvert X_{i_{j}} \rvert^{\frac{r}{t+1+\omega-j}},\cdots,\lvert X_{i_{t}} \rvert^{\frac{r}{t+1+\omega-j}},\biggl\lvert \frac{1}{\lvert T_{0}\rvert}\sum_{i\in T_{0}}X_{i} \biggr\rvert^{\frac{r\omega}{t+1+\omega-j}}\biggr)\Biggr)^{\frac{t+1+\omega-j}{r}}\nonumber \\
    \overset{(*)}{\leq} \ \           & 2\lvert T_{0} \rvert^{\omega}\!\!\sum_{(\ell',\lambda_{1:\ell'}) \in C  (t+2-j)}\biggl(\frac{1}{t+1+\omega-j}\bigl(\lVert X_{i_{j}} \rVert_{r}^{r}+\cdots+\lVert X_{i_{t}} \rVert_{r}^{r}\bigr)\nonumber\\*
    &\quad+\frac{\omega}{t+1+\omega-j}\biggl\lVert \frac{1}{\lvert T_{0} \rvert}\sum_{i\in T_{0}}X_{i} \biggr\rVert _{r}^{r}\biggr)^{\frac{t+1+\omega-j}{r}}\nonumber                                                     \\
    \overset{(**)}{\leq} \            & 2^{t+2-j}\lvert T_{0}\rvert^{\omega}\biggl(\frac{1}{t+1+\omega-j}\bigl(\lVert X_{i_{j}} \rVert_{r}^{r}+\cdots+\lVert X_{i_{t}} \rVert_{r}^{r}\bigr)\nonumber\\*
    &\quad+\frac{\omega}{t+1+\omega-j}\frac{1}{\lvert T_{0} \rvert}\sum_{i\in T_{0}}\lVert X_{i} \rVert_{r}^{r}\biggr)^{\frac{t+1+\omega-j}{r}}\nonumber                                                                                                           \\
    \leq\ \                           & 2^{t+2-j}\lvert T_{0} \rvert^{\omega}M^{t+1+\omega-j},\nonumber
  \end{align}
  where we have used \cref{thm:lemmayoung} in $(*)$, and $(**)$ is implied by $\lvert C(t+2-j) \rvert= 2^{t+1-j}$ and Jensen's inequality as
  \begin{align*}
    & \biggl\lVert \frac{1}{\lvert T_{0} \rvert} \sum_{i\in T_{0}}X_{i}\biggr\rVert _{r}^{r}=\mathbb{E} \biggl[\bigl\lvert\frac{1}{\lvert T_{0} \rvert}\sum_{i\in T_{0}}X_{i}\bigr\rvert^{r}\biggr]\\
    \leq &\mathbb{E} \biggl[\frac{1}{\lvert T_{0} \rvert}\sum_{i\in T_{0}}\lvert X_{i} \rvert^{r}\biggr]=\frac{1}{\lvert T_{0} \rvert}\sum_{i\in T_{0}}\mathbb{E} [\lvert X_{i} \rvert^{r}]=\frac{1}{\lvert T_{0} \rvert}\sum_{i\in T_{0}}\lVert X_{i} \rVert_{r}^{r}.
  \end{align*}

  Thus, the final bound in \eqref{eq:dontlikelabel3} holds for every
  $\omega\in[0,1]$. Substituting \eqref{eq:dontlikelabel2} and
  \eqref{eq:dontlikelabel3} into \eqref{eq:dontlikelabel1}, we have
  \begin{align*}
         & \Bigl\lvert \mathbb{E} \bigl[ X_{i_{w}}\cdots X_{i_{j-1}} \!\!\cdot [\eta_{s+1},\cdots,\eta_{\ell}]\triangleright\mathcal{D}\bigl(X_{i_{j}},\cdots,X_{i_{t}},S_{T_{0},\omega}\bigr)\bigr]\Bigr\rvert \\
    \leq & 2^{t+5-j}\lvert T_{0} \rvert^{\omega}\bigl(\alpha (\mathcal{F}_{1},\mathcal{F}_{3})\bigr)^{(r-t+w-1-\omega)/r} M^{t+1+\omega-w}.
  \end{align*}

  Combining this with \eqref{eq:dontlikelabel4}, we obtain
  \begin{align*}
         & \Bigl\lvert [\eta_{1},\cdots,\eta_{\ell}]\triangleright \mathcal{D}^{*}\bigl(X_{i_{1}},\cdots,X_{i_{t}},S_{T_{0},\omega}\bigr)\Bigr\rvert                                   \\
    \leq & \sum_{w=1}^{j-1}2^{(w-2)\vee 0}M^{w-1}\cdot 2^{t+5-j}\lvert T_{0} \rvert^{\omega}\bigl(\alpha (\mathcal{F}_{1},\mathcal{F}_{3})\bigr)^{(r-t+w-1-\omega)/r} M^{t+1+\omega-w} \\
    \leq & 2^{t+3}\lvert T_{0} \rvert^{\omega}\bigl(\alpha (\mathcal{F}_{1},\mathcal{F}_{3})\bigr)^{(r-t-\omega)/r}M^{t+\omega}.
  \end{align*}

  If $j=t+1$, then $\eta_{1}+\cdots+\eta_{s}=j-1=t$ implies that $s=\ell-1$ and $\eta_{\ell}=1$. By \eqref{eq:strangeeq} we have
  \begin{align}\label{eq:hatelabel1}
      & [\eta_{1},\cdots,\eta_{\ell}]\triangleright \mathcal{D}^{*}\bigl(X_{i_{1}},\cdots,X_{i_{t}},S_{T_{0},\omega}\bigr)              \\
    = & [\eta_{1},\cdots,\eta_{s},1]\triangleright \mathcal{D}^{*}\bigl(X_{i_{1}},\cdots,X_{i_{t}},S_{T_{0},\omega}\bigr)\nonumber                 \\
    = & [\eta_{1},\cdots,\eta_{s-1},\eta_{s}+1]\triangleright \mathcal{D}^{*}\bigl(X_{i_{1}},\cdots,X_{i_{t}},\mathcal{D}(S_{T_{0},\omega})\bigr).\nonumber  
  \end{align}

  By \eqref{eq:dastbound2}, we get
  \begingroup\small
  \begin{align}\label{eq:hatelabel2}
                        & \Bigl\lvert [\eta_{1},\cdots,\eta_{s-1},\eta_{s}+1]\triangleright \mathcal{D}^{*}\bigl(X_{i_{1}},\cdots,X_{i_{t}},\mathcal{D}(S_{T_{0},\omega})\bigr)\Bigr\rvert                                                                                                                         \\
    \leq                & \Bigl\lvert\mathbb{E} \bigl[ X_{i_{1}}\cdots X_{i_{t}} \!\!\cdot \mathcal{D}(S_{T_{0},\omega})\bigr]\Bigr\rvert\nonumber                                                                                                                                                                           \\*
                        & \ +\sum_{w=1}^{t}\sum_{(s',\zeta_{1:s'})\in C (w)}\ [\zeta_{1},\cdots,\zeta_{s'}]\triangleright \bigl(\lvert X_{i_{1}} \rvert,\cdots,\lvert X_{i_{w}} \rvert\bigr)\cdot \Bigl\lvert\mathbb{E} \bigl[ X_{i_{w+1}}\cdots X_{i_{t}} \!\!\cdot \mathcal{D}(S_{T_{0},\omega})\bigr]\Bigr\rvert\nonumber \\
    \overset{(*)}{\leq} & \sum_{w=1}^{t}2^{(w-2)\vee 0}M^{w-1}\Bigl\lvert\mathbb{E} \bigl[ X_{i_{w}}\cdots X_{i_{t}} \!\!\cdot \mathcal{D}(S_{T_{0},\omega})\bigr]\Bigr\rvert,\nonumber 
  \end{align}
  \endgroup
  where $(*)$ is due to \eqref{eq:dontlikelabel7} and the fact that $\lvert C(w) \rvert=2^{w-1}$.

  Combining \eqref{eq:hatelabel1} and \eqref{eq:hatelabel2}, we have
  \begingroup\small
  \begin{equation}\label{eq:hatelabel3}
    \Bigl\lvert [\eta_{1},\cdots,\eta_{\ell}]\triangleright \mathcal{D}^{*}\bigl(X_{i_{1}},\cdots,X_{i_{t}},S_{T_{0},\omega}\bigr) \Bigr\rvert \leq \sum_{w=1}^{t}2^{(w-2)\vee 0}M^{w-1}\Bigl\lvert\mathbb{E} \bigl[ X_{i_{w}}\cdots X_{i_{t}} \!\!\cdot \mathcal{D}(S_{T_{0},\omega})\bigr]\Bigr\rvert.
  \end{equation}
  \endgroup

  If $\omega=0$, then $\lVert\mathcal{D}(S_{T_{0},0})\rVert_{\infty}\leq 2$.
  Applying \cref{thm:covineq} with exponents
  $r/(t+1-w)$ and $\infty$, and then using \eqref{eq:dontlikelabel2}, gives
  \begin{equation*}
    \Bigl\lvert \mathbb{E} \bigl[ X_{i_{w}}\cdots X_{i_{t}} \!\!\cdot
    \mathcal{D}(S_{T_{0},0})\bigr]\Bigr\rvert
    \leq 16\bigl(\alpha (\mathcal{F}_{1},\mathcal{F}_{3})\bigr)^{(r-t+w-1)/r}
    M^{t-w+1}.
  \end{equation*}

  We now suppose $\omega>0$. Applying \cref{thm:covineq}, we obtain
  \begin{equation}\label{eq:hatelabel4}
  \begin{aligned}
    &\Bigl\lvert \mathbb{E} \bigl[ X_{i_{w}}\cdots X_{i_{t}} \!\!\cdot \mathcal{D}(S_{T_{0},\omega})\bigr]\Bigr\rvert\\
    \leq & 8\bigl(\alpha (\mathcal{F}_{1},\mathcal{F}_{3})\bigr)^{(r-t+w-1-\omega)/r}\bigl\lVert X_{i_{w}}\cdots X_{i_{t}} \bigr\rVert _{r/(t+1-w)}\cdot \bigl\lVert \mathcal{D}(S_{T_{0},\omega}) \bigr\rVert _{r/\omega}.
  \end{aligned}
\end{equation}

  Note that
  \begin{align}\label{eq:hatelabel5}
         & \bigl\lVert \mathcal{D}(S_{T_{0},\omega}) \bigr\rVert _{r/\omega}\leq  2\bigl\lVert S_{T_{0},\omega} \bigr\rVert _{r/\omega}
    \leq 2\biggl\lVert \sum_{i\in T_{0}}X_{i} \biggr\rVert _{r}^{\omega}                                                      \\
    \leq & 2\lvert T_{0} \rvert ^{\omega}\biggl\lVert \frac{1}{\lvert T_{0} \rvert}\sum_{i\in T_{0}}X_{i} \biggr\rVert _{r}^{\omega}
    \overset{(*)}{\leq }2\lvert T_{0} \rvert  ^{\omega} \biggl(\frac{1}{\lvert T_{0} \rvert}\sum_{i\in T_{0}}\lVert X_{i}\rVert _{r}^{r}\biggr)^{\omega/r}
    \leq 2\lvert T_{0} \rvert^{\omega}M^{\omega}.\nonumber 
  \end{align}
  Once again we have used Jensen's inequality to get $(*)$.

  Substituting \eqref{eq:dontlikelabel2} and \eqref{eq:hatelabel5} into \eqref{eq:hatelabel4}, and combining the result with the preceding $\omega=0$ bound, we get for every $\omega\in[0,1]$
  \begin{equation*}
    \Bigl\lvert \mathbb{E} \bigl[ X_{i_{w}}\cdots X_{i_{t}} \!\!\cdot \mathcal{D}(S_{T_{0},\omega})\bigr]\Bigr\rvert\leq  16\bigl(\alpha (\mathcal{F}_{1},\mathcal{F}_{3})\bigr)^{(r-t+w-1-\omega)/r}\lvert T_{0} \rvert^{\omega}M^{t-w+\omega+1}.
  \end{equation*}

  Combining this and \eqref{eq:hatelabel3}, we obtain
  \begin{align*}
         & \Bigl\lvert [\eta_{1},\cdots,\eta_{\ell}]\triangleright \mathcal{D}^{*}\bigl(X_{i_{1}},\cdots,X_{i_{t}},S_{T_{0},\omega}\bigr) \Bigr\rvert                        \\
    \leq & \sum_{w=1}^{t}2^{(w-2)\vee 0}M^{w-1}\cdot 16\bigl(\alpha (\mathcal{F}_{1},\mathcal{F}_{3})\bigr)^{(r-t+w-1-\omega)/r}\lvert T_{0} \rvert^{\omega}M^{t-w+\omega+1} \\
    \leq & 2^{t+3}\lvert T_{0} \rvert^{\omega}\bigl(\alpha (\mathcal{F}_{1},\mathcal{F}_{3})\bigr)^{(r-t-\omega)/r}M^{t+\omega}.
  \end{align*}\hfill

  It remains to verify the exact-independence conclusions.  In
  \eqref{eq:someneedslabel3}, every term on the right is the absolute value
  of $\mathbb E[UV]$, where $U$ is $\mathcal F_1$-measurable, $V$ is
  $\mathcal F_2$-measurable, and $\mathbb EV=0$ by \eqref{mean0cd}.  Hence
  every such term is zero when $\mathcal F_1$ and $\mathcal F_2$ are
  independent.  The same argument applies to \eqref{eq:dontlikelabel4}
  when $j\leq t$, with $V$ measurable with respect to $\mathcal F_3$, and
  to \eqref{eq:hatelabel3} when $j=t+1$, where
  $V=\mathcal D(S_{T_0,\omega})$ has mean zero.  Thus the complete
  compositional $\mathcal D^*$ block vanishes before any absolute-value or
  mixing bound is used.  The endpoint moment assumptions make all these
  products integrable by the already proved bounds
  \eqref{eq:dstarnomix1}--\eqref{eq:dstarnomix2}.
\end{proof}

Before providing the upper bound for the sums of $\mathcal{S}(H)$ and $\mathcal{U}_{f}(H)$, we first control the $\mathcal{E}_{G}$ operator defined in \cref{sec:summationterms} using $\alpha$-mixing coefficients.

\begin{lemma}\label{thm:mixingsplit12}
  Let $(X_{i})_{i\in T}$ be a mean-zero random field, fix the neighborhood
  radius $m\in\mathbb N_+$, and let $G$ be an order-$(k+1)$ genogram. In each assertion below, assume $\sup_i\lVert X_i\rVert_r\leq M<\infty$ at the indicated moment order $r$. Then the following bounds hold:
  In the interpolated bounds, a factor with exponent zero is interpreted as
  one, including when its base is zero.
  \begin{enumerate}
    \item \label{itm:nmbound1} For any $r\in[k+1,\infty]$, we have
    \begin{equation}\label{eq:nomixingsplit1}
      \bigl\lvert\mathcal{E}_{G}\bigl(X_{i_{1}},\cdots,X_{i_{k}},X_{i_{k+1}}\bigr)\bigr\rvert \leq 2^{k} M^{k+1}.
    \end{equation}
    \item \label{itm:nmbound2} For any $\omega\in(0,1]$, $r\in[k+\omega,\infty]$, and $f\in \mathcal{C}^{k-1,\omega }(\mathbb{R})$, we have
    \begin{equation}\label{eq:nomixingsplit2}
      \bigl\lvert\mathcal{E}_{G}\bigl(X_{i_{1}},\cdots,X_{i_{k}},\Delta_{f}(G)\bigr)\bigr\rvert \leq 2^{k} \sigma^{-\omega}\bigl\lvert B_{k+1}\backslash D_{k+1} \bigr\rvert^{\omega}\cdot \lvert f \rvert_{k-1,\omega}M^{k+\omega},
    \end{equation}
    where $\sigma^{2}:=\operatorname{Var} \left(\sum_{i\in T}X_{i}\right)$.
 
    \par If $k\geq 2$, then for any $\omega\in[0,1]$, $r\in[k+\omega,\infty]$, and $f\in \mathcal{C}^{k-2,1}(\mathbb{R})\cap \mathcal{C}^{k-1,1}(\mathbb{R})$, we have
    \begin{equation}\label{eq:nomixingsplit2new}
      \begin{aligned}
      &\bigl\lvert\mathcal{E}_{G}\bigl(X_{i_{1}},\cdots,X_{i_{k}},\Delta_{f}(G)\bigr)\bigr\rvert\\
       \leq &2^{k+1} \sigma^{-\omega}\bigl\lvert B_{k+1}\backslash D_{k+1} \bigr\rvert^{\omega}\cdot \lvert f \rvert_{k-2,1}^{1-\omega}\lvert f \rvert_{k-1,1}^{\omega}M^{k+\omega}.
      \end{aligned}
    \end{equation}
 
    \item \label{itm:nmbound3} For any $\omega\in(0,1]$, $r\in[k+\omega,\infty]$, $f\in\mathcal{C}^{k-1,\omega}(\mathbb{R})$, and $c_{1},c_{2}\in\mathbb{N}_{+}$ such that $c_{1}<c_{2}$, we have
    \begin{equation*}
      \biggl\lvert\sum_{c_{1}\leq s_{k+1}<c_{2}}\mathcal{E}_{G}\bigl(X_{i_{1}},\cdots,X_{i_{k}},\Delta_{f}(G)\bigr)\biggr\rvert \leq 2^{k} \sigma^{-\omega } (c_{2}-c_{1})^{\omega} \lvert f \rvert_{k-1,\omega }M^{k+\omega },
    \end{equation*}
    where the sum is taken over genograms whose $c_{1}\leq s_{k+1}<c_{2}$ with the vertex set, the edge set, and $s_{1:k}$ fixed.
   
    \par If $k\geq 2$, then for any $\omega\in[0,1]$, $r\in[k+\omega,\infty]$, $f\in \mathcal{C}^{k-2,1}(\mathbb{R})\cap \mathcal{C}^{k-1,1}(\mathbb{R})$, and $c_{1},c_{2}\in\mathbb{N}_{+}$ such that $c_{1}<c_{2}$, we have
    \begin{equation*}
      \begin{aligned}
      &\biggl\lvert\sum_{c_{1}\leq s_{k+1}<c_{2}}\!\!\!\!\!\!\mathcal{E}_{G}\bigl(X_{i_{1}},\cdots,X_{i_{k}},\Delta_{f}(G)\bigr)\biggr\rvert \\
      \leq &2^{k+1} \sigma^{-\omega } (c_{2}-c_{1})^{\omega} \lvert f \rvert_{k-2, 1}^{1-\omega}\lvert f \rvert_{k-1,1}^{\omega}M^{k+\omega }.
      \end{aligned}
    \end{equation*}

  \end{enumerate}
  Now suppose there exists $1<j\leq k+1$ such that $s_{j}\geq 1$. Then the following holds:
  \begin{enumerate}
    \setcounter{enumi}{3}
    \item \label{itm:egbound1} For any $r\in(k+1,\infty]$, we have
    \begin{equation}\label{eq:mixingsplit1}
      \bigl\lvert\mathcal{E}_{G}\bigl(X_{i_{1}},\cdots,X_{i_{k}},X_{i_{k+1}}\bigr)\bigr\rvert \leq 2^{k+3} \alpha_{\ell_{0}}^{(r-k-1)/r} M^{k+1},
    \end{equation}
    where $\ell_{0}$ is the smallest integer $\ell$ such that
    \begin{equation*}
      k(2\ell+1)^{d}\geq  \max_{1\leq j\leq k+1} s_{j}+k(2m+1)^{d}.
    \end{equation*}
    \item \label{itm:egbound2} For any $\omega\in(0,1]$, $r\in(k+\omega,\infty]$, and $f\in \mathcal{C}^{k-1,\omega }(\mathbb{R})$, we have
    \begin{equation}\label{eq:mixingsplit2}
      \begin{aligned}
      &\bigl\lvert\mathcal{E}_{G}\bigl(X_{i_{1}},\cdots,X_{i_{k}},\Delta_{f}(G)\bigr)\bigr\rvert\\ 
      \leq & 2^{k+3} \sigma^{-\omega}\bigl\lvert B_{k+1}\backslash D_{k+1} \bigr\rvert^{\omega}\cdot \lvert f \rvert_{k-1,\omega}\alpha_{\ell_{0}}^{(r-k-\omega)/r}M^{k+\omega},
      \end{aligned}
    \end{equation}
    where $\ell_{0}$ is defined as above.
      \par If $k\geq 2$, then for any $\omega\in[0,1]$, $r\in(k+\omega,\infty]$, and $f\in \mathcal{C}^{k-2,1}(\mathbb{R})\cap \mathcal{C}^{k-1,1}(\mathbb{R})$, we have
      \begin{equation}\label{eq:mixingsplit2new}
        \begin{aligned}
         &\bigl\lvert\mathcal{E}_{G}\bigl(X_{i_{1}},\cdots,X_{i_{k}},\Delta_{f}(G)\bigr)\bigr\rvert\\
          \leq &2^{k+4} \sigma^{-\omega}\bigl\lvert B_{k+1}\backslash D_{k+1} \bigr\rvert^{\omega}\cdot \lvert f \rvert_{k-2,1}^{1-\omega}\lvert f \rvert_{k-1,1}^{\omega}\alpha_{\ell_{0}}^{(r-k-\omega)/r}M^{k+\omega},
        \end{aligned}
      \end{equation}
    \item \label{itm:egbound3} For any $\omega\in(0,1]$, $r\in(k+\omega,\infty]$, $f\in\mathcal{C}^{k-1,\omega}(\mathbb{R})$, and $c_{1},c_{2}\in\mathbb{N}_{+}$ such that $c_{1}<c_{2}$, we have
    \begin{equation}\label{eq:thirdstatement}
      \begin{aligned}
      &\biggl\lvert\sum_{c_{1}\leq s_{k+1}<c_{2}}\!\!\!\!\!\!\mathcal{E}_{G}\bigl(X_{i_{1}},\cdots,X_{i_{k}},\Delta_{f}(G)\bigr)\biggr\rvert \\
      \leq & 2^{k+3} \sigma^{-\omega } (c_{2}-c_{1})^{\omega} \lvert f \rvert_{k-1,\omega }\alpha_{\ell_{0}}^{(r-k-\omega)/r}M^{k+\omega },
      \end{aligned}
    \end{equation}
    where $\ell_{0}$ is the smallest integer $\ell$ such that
    \begin{equation*}
      k(2\ell+1)^{d}\geq c_{1}\!\vee\!\max_{1\leq j\leq k}s_{j}+k(2m+1)^{d}.
    \end{equation*}
      \par If $k\geq 2$, then for any $\omega\in[0,1]$, $r\in(k+\omega,\infty]$, $f\in \mathcal{C}^{k-2,1}(\mathbb{R})\cap \mathcal{C}^{k-1,1}(\mathbb{R})$, and $c_{1},c_{2}\in\mathbb{N}_{+}$ such that $c_{1}<c_{2}$, we have
      \begin{equation}\label{eq:thirdstatementnew}
        \begin{aligned}
        &\biggl\lvert\sum_{c_1\leq s_{k+1}<c_2}\mathcal{E}_{G}\bigl(X_{i_{1}},\cdots,X_{i_{k}},\Delta_{f}(G)\bigr)\biggr\rvert\\
         \leq &2^{k+4} \sigma^{-\omega}(c_{2}-c_{1})^{\omega}\cdot \lvert f \rvert_{k-2,1}^{1-\omega}\lvert f \rvert_{k-1,1}^{\omega}\alpha_{\ell_{0}}^{(r-k-\omega)/r}M^{k+\omega}.
        \end{aligned}
      \end{equation}
  \end{enumerate}
  If $\alpha_{m+1}=0$, every left-hand side in
  \cref{itm:egbound1,itm:egbound2,itm:egbound3} is in fact zero.  This
  exact-zero conclusion remains valid at $r=k+1$ in
  \cref{itm:egbound1} and at $r=k+\omega$ in
  \cref{itm:egbound2,itm:egbound3}, under the corresponding function
  assumptions stated above.
\end{lemma}

\begin{proof}[Proof of \cref{thm:mixingsplit12}.]
  We will perform induction on $k$ to prove this lemma. But before that, we will present some preliminary results.

  Firstly, we observe that if $s_{k+1}\geq 1$, then $B_{k+1}\setminus D_{k+1}$ is a singleton and therefore, $\bigl\lvert B_{k+1}\setminus D_{k+1}\bigr\rvert=1$ which implies that \cref{itm:nmbound2} is a special case of \cref{itm:nmbound3} and that \cref{itm:egbound2} is a special case of \cref{itm:egbound3} by setting $c_{2}=c_{1}+1$. For notational convenience, we combine the two cases by denoting
  \begin{equation*}
    \widetilde{\Delta}_{f}:=\begin{cases}
      \Delta_{f}(G)                               & \text{ if }s_{k+1}\leq 0 \\
      \sum_{c_{1}\leq s_{k+1}<c_{2}}\Delta_{f}(G) & \text{ if }s_{k+1}\geq 1
    \end{cases}.
  \end{equation*}

  Then by definition of $\mathcal{E}_{G}$, we have
  \begin{equation*}
    \mathcal{E}_{G}\bigl(X_{i_{1}},\cdots,X_{i_{k}},\widetilde{\Delta}_{f}\bigr)=
    \begin{cases}
      \mathcal{E}_{G}\bigl(X_{i_{1}},\cdots,X_{i_{k}},\Delta_{f}(G)\bigr)                               & \text{ if }s_{k+1}\leq 0 \\
      \sum_{c_{1}\leq s_{k+1}<c_{2}}\mathcal{E}_{G}\bigl(X_{i_{1}},\cdots,X_{i_{k}},\Delta_{f}(G)\bigr) & \text{ if }s_{k+1}\geq 1
    \end{cases}.
  \end{equation*}
  Further let
  \begin{align*}
     & \widetilde{B}_{k+1}:=\begin{cases}
                              B_{k+1}                                                           & \text{ if }s_{k+1}\leq 0 \\
                              N^{(c_{2}-1)}\bigl( i_{\ell}: \ell\in A(k+1)\bigr)\cup D_{g(k+1)} & \text{ if }s_{k+1}\geq 1
                            \end{cases}, \\
     & \widetilde{D}_{k+1}:=\begin{cases}
                              D_{k+1}                                                         & \text{ if }s_{k+1}\leq 0 \\
                              N^{(c_{1}-1)}\bigl(i_{\ell}:\ell\in A(k+1)\bigr)\cup D_{g(k+1)} & \text{ if }s_{k+1}\geq 1
                            \end{cases}
  \end{align*}
  If $s_{k+1}\geq c_{1}\geq 1$, we have $u(k+1)=k+1$. The following holds due to a telescoping sum argument:
  \begin{equation*}
    \widetilde{\Delta}_{f}= \partial^{k-1}f\bigl(W(\widetilde{B}_{k+1})\bigr)-\partial^{k-1}f\bigl(W(\widetilde{D}_{k+1})\bigr).
  \end{equation*}
  Thus, we get that
  \begin{align*}
    \widetilde{\Delta}_{f} =
    \left\{\begin{aligned}
      &\partial^{k-1}f(W(\widetilde{B}_{k+1}))-\partial^{k-1}f(W(\widetilde{D}_{k+1})) && \text{if }u(k+1)=k+1    \\
      &\begin{aligned}
         & \int_{0}^{1}(k+1-u(k+1))v^{k-u(k+1)}\cdot                                                                               \\
         & \quad \bigl(\partial^{k-1}f(vW(\widetilde{D}_{k+1})+(1-v)W(\widetilde{B}_{k+1}))\\
         &\ \quad -\partial^{k-1}f(W(\widetilde{D}_{k+1}))\bigr)\dif v
      \end{aligned}
   && \text{if } u(k+1)\leq k
    \end{aligned}\right..
  \end{align*}
  If $u(k+1)=k+1$, we have
  \begin{equation*}
    \lvert \widetilde{\Delta}_{f} \rvert \leq \bigl\lvert \partial^{k-1}f\bigl(W(\widetilde{B}_{k+1})\bigr)-\partial^{k-1}f\bigl(W(\widetilde{D}_{k+1})\bigr) \bigr\rvert\leq \sigma^{-\omega}\lvert f \rvert_{k-1,\omega}\biggl\lvert \sum_{i\in \widetilde{B}_{k+1}\backslash \widetilde{D}_{k+1}}X_{i} \biggr\rvert^{\omega}.
  \end{equation*}
  If $u(k+1)\leq k$, we have
  \begin{align*}
    \lvert \widetilde{\Delta}_{f} \rvert\leq & \int_{0}^{1}(k+1-u(k+1))v^{k-u(k+1)}\cdot                                                                                                                                               \\
                                             & \ \Bigl\lvert \partial^{k-1}f\bigl(vW(\widetilde{D}_{k+1})+(1-v)W(\widetilde{B}_{k+1})\bigr)-\partial^{k-1}f\bigl(W(\widetilde{D}_{k+1})\bigr) \Bigr\rvert\dif v                        \\
    \leq                                     & \sigma^{-\omega}\lvert f \rvert_{k-1,\omega}
      \biggl\lvert \sum_{i\in \widetilde{B}_{k+1}\backslash \widetilde{D}_{k+1}}X_{i} \biggr\rvert^{\omega}
      \int_{0}^{1}(k+1-u(k+1))v^{k-u(k+1)}(1-v)^{\omega}\dif v \\
    \leq                                     & \sigma^{-\omega}\lvert f \rvert_{k-1,\omega}\biggl\lvert \sum_{i\in \widetilde{B}_{k+1}\backslash \widetilde{D}_{k+1}}X_{i} \biggr\rvert^{\omega}.
  \end{align*}
  Therefore, in both cases, we can write
  \begin{equation*}
    \widetilde{\Delta}_{f}=\sigma^{-\omega}\lvert f \rvert_{k-1,\omega}S_{T_{0},\omega},
  \end{equation*}
  where $T_{0}=\widetilde{B}_{k+1}\backslash \widetilde{D}_{k+1}$ and $S_{T_{0},\omega}$ {(which depends on $f$ by definition)} is measurable with respect to
  $\sigma(X_i:i\in T\backslash\widetilde D_{k+1})$ and satisfies
  \begin{equation*}
    \textstyle\lvert S_{T_{0},\omega} \rvert\leq \bigl\lvert \sum_{i\in T_{0}}X_{i} \bigr\rvert^{\omega}.
  \end{equation*}

    If $k\geq 2$ and $f\in \mathcal{C}^{k-2,1}(\mathbb{R})\cap \mathcal{C}^{k-1,1}(\mathbb{R})$, $f\in\mathcal{C}^{k-1,1}(\mathbb{R})$ implies that
    \begin{equation*}
      \textstyle\lvert \widetilde{\Delta}_{f} \rvert\leq \sigma^{-1}\lvert f \rvert_{k-1,1}\bigl\lvert \sum_{i\in T_{0}}X_{i} \bigr\rvert.
    \end{equation*}
    On the other hand, $f\in\mathcal{C}^{k-2,1}(\mathbb{R})$ implies that
    \begin{equation*}
      \lvert \widetilde{\Delta}_{f} \rvert\leq 2 \lvert f \rvert_{k-2,1}.
    \end{equation*}
    Thus, for any $\omega \in [0,1]$, we have
    \begin{equation*}
      \textstyle\lvert \widetilde{\Delta}_{f} \rvert\leq 2 \sigma^{-\omega}\lvert f \rvert_{k-2,1}^{1-\omega}\lvert f \rvert_{k-1,1}^{\omega}\bigl\lvert \sum_{i\in T_{0}}X_{i} \bigr\rvert^{\omega}.
    \end{equation*}
    In this setting, we can write 
    \begin{equation*}
      \textstyle\widetilde{\Delta}_{f}=2\sigma^{-\omega}\lvert f \rvert_{k-2,1}^{1-\omega}\lvert f \rvert_{k-1,1}^{\omega}S_{T_{0},\omega},
    \end{equation*}
    where $T_{0}=\widetilde{B}_{k+1}\backslash \widetilde{D}_{k+1}$ and $S_{T_{0},\omega}$ {(which depends on $f$ by definition)} is measurable with respect to
    $\sigma(X_i:i\in T\backslash\widetilde D_{k+1})$ and satisfies
    \begin{equation*}
      \textstyle\lvert S_{T_{0},\omega} \rvert\leq \bigl\lvert \sum_{i\in T_{0}}X_{i} \bigr\rvert^{\omega}.
    \end{equation*}
  Consequently, both the H\"older and the interpolated versions in
  \cref{itm:nmbound2,itm:nmbound3} reduce to the following estimate (with the
  corresponding prefactor displayed above):
  \begin{enumerate}
    \setcounter{enumi}{6}
    \item \label{itm:nmbound4} 
    \begin{equation}\label{eq:nomixingsplit3}
      \bigl\lvert\mathcal{E}_{G}\bigl(X_{i_{1}},\cdots,X_{i_{k}},S_{T_{0},\omega}\bigr)\bigr\rvert \leq 2^{k} \bigl\lvert T_{0}\bigr\rvert^{\omega}M^{k+\omega}.
    \end{equation}
  \end{enumerate}
  Likewise, both versions in \cref{itm:egbound2,itm:egbound3} reduce to
  \begin{enumerate}
    \setcounter{enumi}{7}
    \item \label{itm:egbound4} 
    \begin{equation}\label{eq:mixingsplit3}
      \bigl\lvert\mathcal{E}_{G}\bigl(X_{i_{1}},\cdots,X_{i_{k}},S_{T_{0},\omega}\bigr)\bigr\rvert \leq 2^{k+3} \bigl\lvert T_{0} \bigr\rvert^{\omega}\alpha_{\ell_{0}}^{(r-k-\omega)/r}M^{k+\omega},
    \end{equation}
    where $\ell_{0}$ is the smallest integer $\ell$ such that
    \begin{equation*}
      k(2\ell+1)^{d}\geq c_{1}\!\vee\!\max_{1\leq j\leq k}s_{j}+k(2m+1)^{d}.
    \end{equation*}
  \end{enumerate}

  Secondly, if $s_{j}\geq 1$ for some $1<j\leq k+1$, we denote the $\sigma$-algebras $\mathcal{F}_{j-}:=\sigma(X_{i_{t}}:t\in A(j))$ and $$\mathcal{F}_{j+}:=\begin{cases}\sigma (X_{i}:i\in T\backslash D_{j})&\text{ if }2\leq j\leq k\\\sigma (X_{i}: i\in T\backslash \widetilde{D}_{k+1})&\text{ if }j=k+1\end{cases}.$$ We will establish that $\alpha (\mathcal{F}_{j-},\mathcal{F}_{j+})\leq \alpha_{\ell_{0}}$ where $\alpha(\,\cdot\,,\,\cdot\,)$ is defined in \cref{thm:amix}.

  To achieve this goal, we will write $\mathbb{B}_{\|\cdot\|}(i,b):=\{z\in \mathbb{Z}^d: \|i-z\|\le b\}$ for any $i\in \mathbb{Z}^{d}$, where $\|\cdot\|$ is the maximum norm on $\mathbb{Z}^d$. This is the set of elements at a distance at most $b$ from $i$. Similarly, if $I\subset\mathbb{Z}^d$ we write $\mathbb{B}_{\|\cdot\|}(I,b):=\{z\in \mathbb{Z}^d: \min_{i\in I}\|i-z\|\le b\}$. Set
  \begin{equation*}
    I_j:=\{i_t:t\in A(j)\}.
  \end{equation*}
  We denote by $\ell$ ($\ell\geq m+1$) the distance between $i_{j}$ and $I_j$ in $\mathbb{Z}^{d}$, and by $q$ the number of indices whose distance from $I_j$ is at least $m+1$ and at most $\ell$; that is, we set
  $$
  \begin{aligned}
  q:=&\bigl|\{s\in T: d(I_j,s)\in [m+1,\ell]\}\bigr| \\
  =&\bigl|T\cap \mathbb{B}_{\|\cdot\|}(I_j, \ell)\backslash
       \mathbb{B}_{\|\cdot\|}(I_j, m)\bigr|.
  \end{aligned}
  $$ 
  To bound $q$, we note that for any $i\in \mathbb{Z}^d$ and $b\in \mathbb{N}$ we have exactly $(2b+1)^d$ elements in $\mathbb{B}_{\|\cdot \|}(i,b)$.
  Thus, we have
  \begin{equation*}
    q\leq \lvert A(j)\rvert\bigl((2\ell+1)^{d}-(2m+1)^{d}\bigr)\leq k\bigl((2\ell+1)^{d}-(2m+1)^{d}\bigr).
  \end{equation*}
  Moreover, by definition, $N^{(s_{j})}(i_{t}:t\in A (j))\setminus N(i_{t}:t\in A(j))$ contains the smallest $s_{j}$ indices (with respect to the strict order on $\mathbb{Z}^{d}$) in $T\setminus N(i_{t}:t\in A(j))$. We remark that all the elements in $N^{(s_{j})}(i_{t}:t\in A (j))\setminus N(i_{t}:t\in A (j))$ have distance at least $m+1$ and at most $\ell$ from $\{ i_{t}:t\in A(j) \}$, meaning that
  \begin{align*}
    &N^{(s_{j})}(i_{t}:t\in A(j))\backslash N(i_{t}:t\in A(j))\\
    &\qquad\subseteq \mathbb{B}_{\|\cdot\|}(I_j,\ell)
      \backslash \mathbb{B}_{\|\cdot\|}(I_j,m).
  \end{align*}
  Thus, we have $q\geq s_{j}$.
  As a result,
  \begin{align*}
    k(2\ell+1)^{d}\geq s_{j}+k(2m+1)^{d}.
  \end{align*}
  As $\ell_{0}:=\min_{\ell}\{\ell: k(2\ell+1)^{d}\geq s_{j}+k(2m+1)^{d}\}$, we have $\ell\geq \ell_{0}$. Thus, we obtain that  $\alpha(\mathcal{F}_{j-},\mathcal{F}_{j+})\leq\alpha_{\ell}\leq \alpha_{\ell_{0}}$.  Since also $\ell\geq m+1$, if $\alpha_{m+1}=0$ then
  $\alpha(\mathcal F_{j-},\mathcal F_{j+})=0$, so these two sigma-fields
  are independent.

  Now we finish the proof of \cref{itm:nmbound1,itm:nmbound4,itm:egbound1,itm:egbound4} of \cref{thm:mixingsplit12} by performing induction on $k$.

  If $k=1$, by definition we have
  \begin{equation*}
    \mathcal{E}_{G}(X_{i_{1}},X_{i_{2}})=\mathcal{D}^{*}(X_{i_{1}},X_{i_{2}}),\quad \mathcal{E}_{G}(X_{i_{1}},S_{T_{0},\omega})=\mathcal{D}^{*}(X_{i_{1}},S_{T_{0},\omega}).
  \end{equation*}
  By \eqref{eq:dstarnomix1} and \eqref{eq:dstarnomix2}, we have
  \begin{align*}
    \bigl\lvert \mathcal{D}^{*}(X_{i_{1}},X_{i_{2}}) \bigr\rvert\leq 2 M^{2},\quad \bigl\lvert \mathcal{D}^{*}(X_{i_{1}},S_{T_{0},\omega}) \bigr\rvert\leq 2 \lvert T_{0} \rvert^{\omega}M^{1+\omega}.
  \end{align*}
  Thus, \cref{itm:nmbound1,itm:nmbound4} of \cref{thm:mixingsplit12} hold for $k=1$.
  Now supposing $s_{2}\geq 1$, by \eqref{eq:dstarmix1} and \eqref{eq:dstarmix2}, we get
  \begin{align*}
     & \bigl\lvert \mathcal{D}^{*}(X_{i_{1}},X_{i_{2}}) \bigr\rvert\leq 2^{4}\bigl(\alpha (\mathcal{F}_{1},\mathcal{F}_{2})\bigr)^{(r-2)/r}M^{2},                                                   \\
     & \bigl\lvert \mathcal{D}^{*}(X_{i_{1}},S_{T_{0},\omega}) \bigr\rvert \leq 2^{4}\lvert T_{0} \rvert^{\omega}\bigl(\alpha (\mathcal{F}_{1},\mathcal{F}_{3})\bigr)^{(r-1-\omega)/r}M^{1+\omega},
  \end{align*}
  where $\mathcal{F}_{1}:=\sigma (X_{i_{1}})=\mathcal{F}_{2-}$, $\mathcal{F}_{2}:=\sigma (X_{i_{2}})\subseteq \mathcal{F}_{2+}$, and
  $$\textstyle\mathcal{F}_{3}:=\sigma(S_{T_{0},\omega})\subseteq
  \sigma (X_{i}:i\in T\backslash \widetilde{D}_{2})=\mathcal{F}_{2+}.$$
  As we have shown $\alpha (\mathcal{F}_{2-},\mathcal{F}_{2+})\leq \alpha_{\ell_{0}}$, we obtain
  \begin{equation*}
    \textstyle \alpha(\mathcal{F}_{1},\mathcal{F}_{2})\leq \alpha_{\ell_{0}},\quad \alpha(\mathcal{F}_{1},\mathcal{F}_{3})\leq \alpha_{\ell_{0}}.
  \end{equation*}
  Thus, \cref{itm:egbound1,itm:egbound4} of \cref{thm:mixingsplit12} also hold for $k=1$.

  Suppose \cref{itm:nmbound1,itm:egbound1,itm:nmbound4,itm:egbound4} of \cref{thm:mixingsplit12} are true for $\lvert G \rvert\leq k$. Consider the case where $\lvert G \rvert=k+1$. Let
  \begin{equation}
    q_{0}:=\sup \{ j: j=1\text{ or }p(j)\neq j-1\text{ for }2\leq j\leq k+1 \},
  \end{equation}
  We remark that $q_0$ is the first vertex in the branch of $G$ with the highest indices. We set $w:=\bigl\lvert \{t:q_{0}+1\leq t\leq k+1\ \&\ s_t\ge 0\}\bigr\rvert$ to be the number of all indices $q_{0}+1\leq t\leq k+1$ such that the identifier $s_{t}\geq 0$. If $\max_{1\leq j\leq k+1}s_{j}\geq 1$, we let $j_{0}$ be an integer that satisfies $s_{j_{0}}=\max_{1\leq j\leq k+1}s_{j}\geq 1$. We remark that such an index always exists.

  We will first propose a simplified formulation for $\mathcal{E}_G$ that will hold irrespective of the value of $w$. Then we will distinguish two main cases in our analysis namely (i) when $q_0=1$ and (ii) when $q_0\ge 2$.

  To achieve this goal, we first remark that if  $w=0$, by definition we know that for any random variables $Y_{1},\cdots,Y_{k+1}$ the following holds
  \begin{equation*}
    \mathcal{E}_{G}(Y_{1},\cdots,Y_{k+1})=
    \begin{cases}
      \mathcal{D}^{*}\bigl(Y_{1}Y_{2}\cdots Y_{k+1}\bigr)                                                                             & \text{ if }q_{0}=1     \\
      \mathcal{E}_{G[q_{0}-1]}\bigl(Y_{1},\cdots,Y_{q_{0}-1}\bigr)\cdot \mathcal{D}^{*}\bigl(Y_{q_{0}}Y_{q_{0}+1}\cdots Y_{k+1}\bigr) & \text{ if }q_{0}\geq 2
    \end{cases},
  \end{equation*}
  where $G[q_{0}-1]\subseteq G$ is the unique order-($q_{0}-1$) sub-genogram of $G$ as defined in \cref{sec:genogram}.

  For $w\geq 1$, we write $\{t:q_{0}+1\leq t\leq k+1\ \&\ s_t\ge 0\}=\{ q_{1},\cdots,q_{w} \}$.  Without loss of generality, we suppose that the sequence $q_{0}+1\leq q_{1}<\cdots<q_{w}\leq k+1$ is increasing. By definition
  \begingroup\small
  \begin{equation*}
    \mathcal{E}_{G}(Y_{1},\cdots,Y_{k+1})=
    \begin{cases}
      \mathcal{D}^{*}\bigl(Y_{1}\cdots Y_{q_{1}-1}\ ,\ Y_{q_{1}}\cdots Y_{q_{2}-1}\ ,\ \cdots\ ,\ Y_{q_{w}}\cdots Y_{k+1}\bigr)                                                                       & \!\!\!\!\text{ if }q_{0}=1     \\
      \begin{aligned}
        &\mathcal{E}_{G[q_{0}-1]}\bigl(Y_{1},\cdots,Y_{q_{0}-1}\bigr)\cdot\\
        &\  \mathcal{D}^{*}\bigl(Y_{q_{0}}\cdots Y_{q_{1}-1}\ ,\ Y_{q_{1}}\cdots Y_{q_{2}-1}\ ,\ \cdots\ ,\ Y_{q_{w}}\cdots Y_{k+1}\bigr)
      \end{aligned} & \!\!\!\!\text{ if }q_{0}\geq 2
    \end{cases}\!\!.
  \end{equation*}
  \endgroup

  Set $q_{w+1}:=k+2$, then by exploiting the definition of compositional $\mathcal{D}^{*}$ operators, we remark that $\mathcal{E}_G$ will take the following form irrespectively of the fact that $w\ge 1$ or not:
  \begingroup\small
  \begin{equation}
    \mathcal{E}_{G}(Y_{1},\cdots,Y_{k+1}):=
    \begin{cases}
      [q_{1}-q_{0},\cdots,q_{w+1}-q_{w}]\triangleright\mathcal{D}^{*}\bigl(Y_{1},\cdots , Y_{k+1}\bigr)                                                                       & \text{ if }q_{0}=1     \\
      \begin{aligned}
       &\mathcal{E}_{G[q_{0}-1]}\bigl(Y_{1},\cdots,Y_{q_{0}-1}\bigr)\cdot\\
       &\quad [q_{1}-q_{0},\cdots,q_{w+1}-q_{w}]\triangleright\mathcal{D}^{*}\bigl(Y_{q_{0}},\cdots , Y_{k+1}\bigr)
      \end{aligned} & \text{ if }q_{0}\geq 2
    \end{cases}.
  \end{equation}
  \endgroup
  In particular, we know that
  \begin{align}
     & \mathcal{E}_{G}(X_{i_{1}},\cdots,X_{i_{k+1}})=                                                                                                                                                                                                                                                                                                                                                                                                   \\
     & \quad\left\{\begin{aligned}
              &[q_{1}-q_{0},\cdots,q_{w+1}-q_{w}]\triangleright\mathcal{D}^{*}\bigl(X_{i_{1}},\cdots, X_{i_{k+1}}\bigr)                                                                               && \text{ if }q_{0}=1     \\
              &\begin{aligned}
                &\mathcal{E}_{G[q_{0}-1]}\bigl(X_{i_{1}},\cdots,X_{i_{q_{0}-1}}\bigr)\cdot \\
                &\quad [q_{1}-q_{0},\cdots,q_{w+1}-q_{w}]\triangleright\mathcal{D}^{*}\bigl(X_{i_{q_{0}}},\cdots, X_{i_{k+1}}\bigr)
              \end{aligned} && \text{ if }q_{0}\geq 2
            \end{aligned}\right., \nonumber                                                                                                 \\
     & \mathcal{E}_{G}(X_{i_{1}},\cdots,X_{i_{k}},S_{T_{0},\omega})=                                                                                                                                                                                                                                                                                                                                                                                       \\
     & \quad \left\{\begin{aligned}
               &[q_{1}-q_{0},\cdots,q_{w+1}-q_{w}]\triangleright\mathcal{D}^{*}\bigl(X_{i_{1}},\cdots, X_{i_{k}},S_{T_{0},\omega}\bigr)                                                                                               && \text{ if }q_{0}=1     \\
               &\begin{aligned}
               &\mathcal{E}_{G[q_{0}-1]}\bigl(X_{i_{1}},\cdots,X_{i_{q_{0}-1}}\bigr)\cdot \\
               &\quad [q_{1}-q_{0},\cdots,q_{w+1}-q_{w}]\triangleright\mathcal{D}^{*}\bigl(X_{i_{q_{0}}},X_{i_{q_{0}+1}},\cdots ,X_{i_{k}},S_{T_{0},\omega}\bigr)
               \end{aligned} && \text{ if }q_{0}\geq 2
             \end{aligned}\right..\nonumber  
  \end{align}

  We will use this simplified representation to prove the desired result. If $q_{0}=1$, by \eqref{eq:dstarnomix1} and \eqref{eq:dstarnomix2} we remark that
  \begin{align*}
     & \bigl\lvert [q_{1}-q_{0},\cdots,q_{w+1}-q_{w}]\triangleright\mathcal{D}^{*}(X_{i_{1}},\cdots,X_{i_{k+1}}) \bigr\rvert\leq 2^{k+3}M^{k+1},                                                 \\
     & \bigl\lvert [q_{1}-q_{0},\cdots,q_{w+1}-q_{w}]\triangleright\mathcal{D}^{*}(X_{i_{1}},\cdots,X_{i_{k}},S_{T_{0},\omega}) \bigr\rvert\leq 2^{k+3}\lvert T_{0} \rvert^{\omega}M^{k+\omega}.
  \end{align*}
  Therefore, \cref{itm:nmbound1,itm:nmbound4} are true when $q_{0}=1$.

  Supposing $s_{j_{0}}=\max_{1\leq j\leq k+1}s_{j}\geq 1$ ($j_{0}\geq 2$ since $s_{1}=0$), by definition of $q_{1},\cdots,q_{w}$ we know there is some $1\leq w'\leq w$ such that $q_{w'}=j_{0}$. Hence
  \begin{equation*}
    (q_{1}-q_{0})+\cdots+(q_{w'}-q_{w'-1})=j_{0}-1.
  \end{equation*}
  By \eqref{eq:dstarmix1} and \eqref{eq:dstarmix2} we have
  \begin{equation*}
    \begin{aligned}
     & \bigl\lvert [q_{1}-q_{0},\cdots,q_{w+1}-q_{w}]\triangleright\mathcal{D}^{*}(X_{i_{1}},\cdots,X_{i_{k+1}}) \bigr\rvert\\
     \leq &2^{k+3}\bigl(\alpha (\mathcal{F}_{1},\mathcal{F}_{2})\bigr)^{(r-k-1)/r}M^{k+1},                                                      
    \end{aligned}
  \end{equation*}
  and
  \begin{equation*}
    \begin{aligned}
     & \bigl\lvert [q_{1}-q_{0},\cdots,q_{w+1}-q_{w}]\triangleright\mathcal{D}^{*}(X_{i_{1}},\cdots,X_{i_{k}},S_{T_{0},\omega}) \bigr\rvert\\
     \leq &2^{k+3}\lvert T_{0} \rvert^{\omega}\bigl(\alpha (\mathcal{F}_{1},\mathcal{F}_{3})\bigr)^{(r-k-\omega)/r}M^{k+\omega},
    \end{aligned}
  \end{equation*}
  where
  \begingroup\small
  \begin{align*}
     & \mathcal{F}_{1}:=\sigma (X_{i_{1}},\cdots,X_{i_{j_{0}-1}})=\mathcal{F}_{j_{0}-},                                               \\
     & \mathcal{F}_{2}:=\sigma (X_{i_{j_{0}}},\cdots,X_{i_{k+1}})\overset{(*)}{\subseteq}
    \begin{cases}
      \sigma (X_{i}:i\in T\backslash D_{j_{0}})= \mathcal{F}_{j_{0}+}          & \text{ if }j_{0}\leq k \\
      \sigma (X_{i}:i\in T\backslash \widetilde{D}_{k+1})=\mathcal{F}_{j_{0}+} & \text{ if }j_{0}=k+1
    \end{cases}, \\
     & \mathcal{F}_{3}:=
    \begin{cases}
      \begin{aligned}
        \sigma (X_{i_{j_{0}}},\cdots,X_{i_{k}},S_{T_{0},\omega})
        \subseteq \               & \sigma (X_{i}:i=j_{0},\cdots k,\text{ or }i\in T\backslash \widetilde{D}_{k+1}) \\
        \overset{(**)}{\subseteq} & \sigma (X_{i}:i\in T\backslash D_{j_{0}})=\mathcal{F}_{j_{0}+}
      \end{aligned}
                                                                                                                  & \text{ if }j_{0}\leq k \\
      \sigma (S_{T_{0},\omega})\subseteq \sigma (X_{i}:i\in T\backslash \widetilde{D}_{k+1})=\mathcal{F}_{j_{0}+} & \text{ if }j_{0}=k+1
    \end{cases}.
  \end{align*}
  \endgroup
  Here $(*)$ and $(**)$ are due to the fact that $p(j)=j-1$ for any $2\leq j\leq k+1$ since it implies that $1,\cdots,j_{0}-1\in A(j)$ for any $j_{0}\leq j\leq k+1$.
  Thus, we have
  \begin{align*}
     & \bigl\lvert [q_{1}-q_{0},\cdots,q_{w+1}-q_{w}]\triangleright\mathcal{D}^{*}(X_{i_{1}},\cdots,X_{i_{k+1}}) \bigr\rvert\leq 2^{k+3}\alpha_{\ell_{0}}^{(r-k-1)/r}M^{k+1},                                                       \\
     & \bigl\lvert [q_{1}-q_{0},\cdots,q_{w+1}-q_{w}]\triangleright\mathcal{D}^{*}(X_{i_{1}},\cdots,X_{i_{k}},S_{T_{0},\omega}) \bigr\rvert\leq 2^{k+3}\lvert T_{0} \rvert^{\omega}\alpha_{\ell_{0}} ^{(r-k-\omega)/r}M^{k+\omega}.
  \end{align*}
  Therefore, \cref{itm:egbound1,itm:egbound4} of \cref{thm:mixingsplit12} are true when $q_{0}=1$.

  If $q_{0}\geq 2$, note that
  \begin{equation}\label{eq:inductivenomix}
    \bigl\lvert \mathcal{E}_{G[q_{0}-1]}\bigl(X_{i_{1}},\cdots,X_{i_{q_{0}-1}}\bigr) \bigr\rvert \leq 2^{q_{0}-2}M^{q_{0}-1},
  \end{equation}
  which is true for $q_{0}=2$ since $\mathcal{E}_{G[1]}(X_{i_{1}})=\mathbb{E} [X_{i_{1}}]$ and is {precisely the inductive hypothesis for $q_{0}\geq 3$}.

  {Thus, we have}
  \begin{align*}
                        & \bigl\lvert \mathcal{E}_{G}(X_{i_{1}},\cdots,X_{i_{k+1}}) \bigr\rvert                                                                                                                                                                                \\
    \leq                & \bigl\lvert \mathcal{E}_{G[q_{0}-1]}\bigl(X_{i_{1}},\cdots,X_{i_{q_{0}-1}}\bigr)\bigr\rvert\cdot\\
    &\quad\bigl\lvert [q_{1}-q_{0},\cdots,q_{w+1}-q_{w}]\triangleright\mathcal{D}^{*}\bigl(X_{i_{q_{0}}},X_{i_{q_{0}+1}},\cdots ,X_{i_{k+1}}\bigr) \bigr\rvert \\
    \overset{(*)}{\leq} & 2^{q_{0}-2}M^{q_{0}-1}\cdot 2^{k+1-q_{0}}M^{k+2-q_{0}}\leq 2^{k}M^{k+1},
  \end{align*}
  where $(*)$ is due to \eqref{eq:inductivenomix} and \eqref{eq:dstarnomix1}.
  Similarly, we have
  \begin{align*}
         & \bigl\lvert \mathcal{E}_{G}(X_{i_{1}},\cdots,X_{i_{k}},S_{T_{0},\omega}) \bigr\rvert                                                                                                                                                                                \\
    \leq & \bigl\lvert \mathcal{E}_{G[q_{0}-1]}\bigl(X_{i_{1}},\cdots,X_{i_{q_{0}-1}}\!\bigr)\bigr\rvert\cdot\\
    &\quad \bigl\lvert [q_{1}\!-\!q_{0},\cdots,q_{w+1}\!-\!q_{w}]\triangleright\mathcal{D}^{*}\bigl(X_{i_{q_{0}}},X_{i_{q_{0}+1}},\cdots, X_{i_{k}},S_{T_{0},\omega}\bigr) \bigr\rvert \\
    \leq & 2^{q_{0}-2}M^{q_{0}-1}\cdot 2^{k+1-q_{0}}\lvert T_{0} \rvert^{\omega}M^{k+1+\omega-q_{0}}\leq 2^{k}\lvert T_{0} \rvert^{\omega}M^{k+\omega}.
  \end{align*}
  Therefore, \cref{itm:nmbound1,itm:nmbound4} are true when $q_{0}\geq 2$.

  Supposing $s_{j_{0}}=\max_{1\leq j\leq k+1}s_{j}\geq 1$, we claim that $p(j_{0})=j_{0}-1$. In fact, if $p(j_{0})<j_{0}-1$, set $j_{0}'=p(j_{0})+1$. Since $j_{0}'<j_{0}$ and $v[j_{0}']$ and $v[j_{0}]$ are siblings, by \cref{itm:fourthrq} of \cref{thm:rqmofgenogram} we have $s_{j_{0}'}>s_{j_{0}}$, which contradicts the definition of $j_{0}$. Therefore, we have shown $p(j_{0})=j_{0}-1$, and thus, $j_{0}\neq q_{0}$ by definition of $q_{0}$.

  If $j_{0}\geq q_{0}+1$, by definition of $q_{1},\cdots,q_{w}$ we know there is some $1\leq w'\leq w$ such that $q_{w'}=j_{0}$. Hence
  \begin{equation*}
    (q_{1}-q_{0})+\cdots+(q_{w'}-q_{w'-1})=j_{0}-q_{0}.
  \end{equation*}
  Thus, by \eqref{eq:dstarmix1} and \eqref{eq:dstarmix2} we have
  \begin{align*}
     & \bigl\lvert [q_{1}-q_{0},\cdots,q_{w+1}-q_{w}]\triangleright\mathcal{D}^{*}(X_{i_{q_{0}}},\cdots,X_{i_{k+1}}) \bigr\rvert\\
     &\qquad\quad\leq 2^{k-q_{0}+4}\bigl(\alpha (\mathcal{F}_{1},\mathcal{F}_{2})\bigr)^{(r-k+q_{0}-2)/r}M^{k-q_{0}+2},                                                          \\
     & \bigl\lvert [q_{1}-q_{0},\cdots,q_{w+1}-q_{w}]\triangleright\mathcal{D}^{*}(X_{i_{q_{0}}},\cdots,X_{i_{k}},S_{T_{0},\omega}) \bigr\rvert\\
     &\qquad\quad\leq 2^{k-q_{0}+4}\lvert T_{0} \rvert^{\omega}\bigl(\alpha (\mathcal{F}_{1},\mathcal{F}_{3})\bigr)^{(r-k+q_{0}-1-\omega)/r}M^{k-q_{0}+1+\omega},
  \end{align*}
  where
  \begingroup\small
  \begin{align*}
     & \mathcal{F}_{1}:=\sigma (X_{i_{q_{0}}},\cdots,X_{i_{j_{0}-1}})=\mathcal{F}_{j_{0}-},                                           \\
     & \mathcal{F}_{2}:=\sigma (X_{i_{j_{0}}},\cdots,X_{i_{k+1}})\overset{(*)}{\subseteq}
    \begin{cases}
      \sigma (X_{i}:i\in T\backslash D_{j_{0}})= \mathcal{F}_{j_{0}+}          & \text{ if }j_{0}\leq k \\
      \sigma (X_{i}:i\in T\backslash \widetilde{D}_{k+1})=\mathcal{F}_{j_{0}+} & \text{ if }j_{0}=k+1
    \end{cases}, \\
     & \mathcal{F}_{3}:=
    \begin{cases}
      \begin{aligned}
        \sigma (X_{i_{j_{0}}},\cdots,X_{i_{k}},S_{T_{0},\omega})
        \subseteq \               & \sigma (X_{i}:i=j_{0},\cdots k,\text{ or }i\in T\backslash \widetilde{D}_{k+1}) \\
        \overset{(**)}{\subseteq} & \sigma (X_{i}:i\in T\backslash D_{j_{0}})=\mathcal{F}_{j_{0}+}
      \end{aligned}
                                                                                                                  & \text{ if }j_{0}\leq k \\
      \sigma (S_{T_{0},\omega})\subseteq \sigma (X_{i}:i\in T\backslash \widetilde{D}_{k+1})=\mathcal{F}_{j_{0}+} & \text{ if }j_{0}=k+1
    \end{cases}.
  \end{align*}
  \endgroup
  Here $(*)$ and $(**)$ are due to the fact that $p(j)=j-1$ for any $q_{0}+1\leq j\leq k+1$ since it implies that $q_{0},\cdots,j_{0}-1\in A(j)$ for any $j_{0}\leq j\leq k+1$.
  Thus, we have
  \begin{align*}
         & \bigl\lvert \mathcal{E}_{G}(X_{i_{1}},\cdots,X_{i_{k+1}}) \bigr\rvert                                                                                                                                                                                \\
    \leq & \bigl\lvert \mathcal{E}_{G[q_{0}-1]}\bigl(X_{i_{1}},\cdots,X_{i_{q_{0}-1}}\bigr)\bigr\rvert\cdot\\
    &\quad\bigl\lvert [q_{1}-q_{0},\cdots,q_{w+1}-q_{w}]\triangleright\mathcal{D}^{*}\bigl(X_{i_{q_{0}}},X_{i_{q_{0}+1}},\cdots, X_{i_{k+1}}\bigr) \bigr\rvert \\
    \leq & 2^{q_{0}-2}M^{q_{0}-1}\cdot 2^{k-q_{0}+4}\alpha_{\ell_{0}}^{(r-k+q_{0}-2)/r}M^{k-q_{0}+2}\\
    \leq & 2^{k+3}\alpha_{\ell_{0}}^{(r-k-1)/r}M^{k+1},
  \end{align*}
  and
  \begin{align*}
         & \bigl\lvert \mathcal{E}_{G}(X_{i_{1}},\cdots,X_{i_{k}},S_{T_{0},\omega}) \bigr\rvert                                                                                                                                                                                \\
    \leq & \bigl\lvert \mathcal{E}_{G[q_{0}-1]}\bigl(X_{i_{1}},\cdots,X_{i_{q_{0}-1}}\!\bigr)\bigr\rvert\cdot\\
    &\quad\bigl\lvert [q_{1}\!-\!q_{0},\cdots,q_{w+1}\!-\!q_{w}]\triangleright\mathcal{D}^{*}\bigl(X_{i_{q_{0}}},X_{i_{q_{0}+1}},\cdots, X_{i_{k}},S_{T_{0},\omega}\bigr) \bigr\rvert \\
    \leq & 2^{q_{0}-2}M^{q_{0}-1}\cdot 2^{k-q_{0}+4}\lvert T_{0} \rvert^{\omega}\alpha_{\ell_{0}}^{(r-k+q_{0} -1-\omega)/r}M^{k-q_{0}+1+\omega}\\
    \leq & 2^{k+3}\lvert T_{0} \rvert^{\omega}\alpha_{\ell_{0}}^{(r-k-\omega)/r}M^{k+\omega}.
  \end{align*}
  If $2\leq j_{0}\leq q_{0}-1$, by inductive hypothesis we have
  \begin{equation}\label{eq:inductivemix}
    \bigl\lvert \mathcal{E}_{G[q_{0}-1]}(X_{i_{1}},\cdots,X_{i_{q_{0}-1}}) \bigr\rvert\leq 2^{q_{0}+1}\alpha_{\ell_{0}}^{(r-q_{0}+1)/r}M^{q_{0}-1}.
  \end{equation}
  Thus, we have
  \begin{align*}
                        & \bigl\lvert \mathcal{E}_{G}(X_{i_{1}},\cdots,X_{i_{k+1}}) \bigr\rvert                                                                                                                                                                                \\
    \leq                & \bigl\lvert \mathcal{E}_{G[q_{0}-1]}\bigl(X_{i_{1}},\cdots,X_{i_{q_{0}-1}}\bigr)\bigr\rvert\cdot\\
    &\quad \bigl\lvert [q_{1}-q_{0},\cdots,q_{w+1}-q_{w}]\triangleright\mathcal{D}^{*}\bigl(X_{i_{q_{0}}},X_{i_{q_{0}+1}},\cdots ,X_{i_{k+1}}\bigr) \bigr\rvert \\
    \overset{(*)}{\leq} & 2^{q_{0}+1}\alpha_{\ell_{0}}^{(r-q_{0}+1)/r}M^{q_{0}-1}\cdot 2^{k-q_{0}+1}M^{k-q_{0}+2}\\
    \leq & 2^{k+3}\alpha_{\ell_{0}}^{(r-k-1)/r}M^{k+1},
  \end{align*}
  where $(*)$ is implied by \eqref{eq:inductivemix} and \eqref{eq:nomixingsplit1}.
  Similarly,
  \begin{align*}
         & \bigl\lvert \mathcal{E}_{G}(X_{i_{1}},\cdots,X_{i_{k}},S_{T_{0},\omega}) \bigr\rvert                                                                                                                                                                                \\
    \leq & \bigl\lvert \mathcal{E}_{G[q_{0}-1]}\bigl(X_{i_{1}},\cdots,X_{i_{q_{0}-1}}\!\bigr)\bigr\rvert\cdot\\
    &\quad \bigl\lvert [q_{1}\!-\!q_{0},\cdots,q_{w+1}\!-\!q_{w}]\triangleright\mathcal{D}^{*}\bigl(X_{i_{q_{0}}},X_{i_{q_{0}+1}},\cdots ,X_{i_{k}},S_{T_{0},\omega}\bigr) \bigr\rvert \\
    \leq & 2^{q_{0}+1}\alpha_{\ell_{0}}^{(r-q_{0}+1)/r}M^{q_{0}-1}\cdot 2^{k-q_{0}+1}\lvert T_{0} \rvert^{\omega}M^{k-q_{0}+1+\omega}\\
    \leq & 2^{k+3}\lvert T_{0} \rvert^{\omega}\alpha_{\ell_{0}}^{(r-k-\omega)/r}M^{k+\omega}.
  \end{align*}
  Therefore, \cref{itm:egbound1,itm:egbound4} of \cref{thm:mixingsplit12} are true when $q_{0}\geq 2$.

  We finally prove the exact-zero assertion, using the same factorization
  but keeping each complete compositional block intact.  For $k=1$ it is
  the exact-independence conclusion of \cref{thm:bounddstarmix}.  In the
  induction step, if $q_0=1$, the positive maximal identifier $j_0$ gives
  precisely the independent split in the displayed $\mathcal D^*$ block,
  so that block is zero.  If $q_0\geq2$, then $j_0\neq q_0$ by the argument
  above.  When $j_0\geq q_0+1$, the terminal $\mathcal D^*$ factor is zero
  by exact independence; when $j_0\leq q_0-1$, the prefix
  $\mathcal E_{G[q_0-1]}$ is zero by the inductive hypothesis.  This proves
  the assertion simultaneously for an $X$ in the final slot and for the
  complete variable $S_{T_0,\omega}$.  The representation of
  $\widetilde\Delta_f$ by $S_{T_0,\omega}$ above---obtained only after the
  entire positive terminal block has been telescoped---then gives the
  claimed conclusions for $\Delta_f(G)$ and for that terminal block.
  By induction the proof is complete.\hfill
\end{proof}

Equipped with the tools in \cref{thm:mixingsplit12}, we are able to show the proof of \cref{THM:REMAINDERCTRL1234}.

\begin{proof}[Proof of \cref{THM:REMAINDERCTRL1234}.]
  For each of \eqref{eq:remainderctrl1}--\eqref{eq:remainderctrl4}, we conduct the sum in two steps:
  \begin{enumerate}
    \item \label{itm:step1lastlemma} Fixing an ordered tree $(V,E,\prec)$ with the compatible labeling. Here $V=\bigl\{ v[1],\cdots,v[k+1] \bigr\}$ denotes the vertex set and $E$ denotes the edge set. We take the sum of $\mathcal{S}(H)$ or $\mathcal{U}_{f}(H)$ over all possible values of $s_{1},\cdots, s_{k+1}$ such that $H$ is a genogram that induces $(V,E,\prec)$;
    \item \label{itm:step2lastlemma} Sum over all possible ordered trees $(V,E,\prec)$ of order ($k+1$).
  \end{enumerate}
  In all identifier sums below, an identifier tuple that does not give a
  compatible genogram for the fixed ordered tree is understood to contribute
  zero. This convention lets us write rectangular identifier ranges without
  adding the same compatibility condition to every display.
  Note that an ordered tree corresponds to infinitely many genograms (as there is an infinite number of possible identifiers). However, when the index set $T$ of the random field is finite, only finitely many genograms give non-zero values of $\mathcal{S}(H)$ and $\mathcal{U}_{f}(H)$.

  For the second step, we observe that the total number of ordered trees of order ($k+1$) solely depends on $k$. (In fact, this is exactly the $k$-th Catalan number \citep{roman2015introduction}.) Hence summing over all such trees only contributes to the constant in the bounds. As for the first step, the following statement will be crucial to our proof.

  \begin{claim}
    Fix a positive integer $s\geq 1$. For any $2\leq t \leq k+1$, given a sequence $i_{1},\cdots,i_{t-1}$, the sum of $\bigl\lvert B_{t}\backslash D_{t} \bigr\rvert$ over $-1\leq s_{t}\leq s$ is smaller or equal to $2(k(2m+1)^{d}+ s)$.
  \end{claim}

  To see this we will consider the following three cases: 
  \begin{enumerate}
    \item When $s_t=-1$ and $s_{u(t)}=0$;
    \item When $s_t=-1$ and $s_{u(t)}\ge 1$;
    \item When $0\le s_t\le s$.
  \end{enumerate} Firstly, if $s_{t}=-1$ and $s_{u(t)}=0$, then we note that
  \begin{align*}
    B_{t}\backslash D_{t}
      &=B_{u(t)}\backslash D_{u(t)}
        \subseteq N(i_{h}:h\in A(u(t)))\\
      &\subseteq N(i_{h}:h\in A(t))
        \subseteq N^{(s)}(i_{h}:h\in A(t)).
  \end{align*}
  If $s_{t}=-1$ and $s_{u(t)}\geq 1$, then by definition, $B_t\backslash D_t=B_{u(t)}\backslash D_{u(t)}$ has at most one element namely $i_{u(t)}$. Thus, $B_{t}\backslash D_{t}\subseteq N(i_{h}:h\in A(t))\subseteq N^{(s)}(i_{h}:h\in A(t))$.

  Finally, if $0\leq s_{t}\leq s$, by definition, $B_{t}\backslash D_{t}\subseteq N^{(s)}(i_{h}:h\in A(t))$.

  To bound $\sum_{s_t\le s}\bigl\lvert B_t\backslash D_t\bigr\rvert$ we remark that the sets $ B_{t}\backslash D_{t}$ are disjoints for different values of $0\le s_t\le s$. Thus, this implies that $$\sum_{s_t\le s}\ \bigl\lvert B_t\backslash D_t\bigr\rvert\le 2\bigl\lvert N^{(s)}(i_h: h\in A(t))\bigr\rvert.$$

  To further bound this, note that for any index $i$, the indices with distance from $i$ at most $m$ lie in the $d$-dimensional hypercube centered at $i$ with the sides of length $2m+1$. Thus, $\bigl\lvert N(i_{h}:h\in A(t)) \bigr\rvert\leq k(2m+1)^{d}$. By noticing that for any subset $J$ (by definition of $N^{(s)}(\cdot)$) there is at most $s$ elements in $N^{(s)}(J)\setminus N(J)$ we obtain that
  \begin{equation*}
    \bigl\lvert N^{(s)}(i_{h}:h\in A(t)) \bigr\rvert\leq\bigl\lvert N(i_{h}:h\in A(t)) \bigr\rvert+s\leq k(2m+1)^{d}+s.
  \end{equation*}

  Next we establish \eqref{eq:remainderctrl1}.

  Suppose $v[j_{0}]$ is a vertex with the largest identifier among all vertices and $s_{j_{0}}=s$. By \eqref{eq:mixingsplit1} of \cref{thm:mixingsplit12}, we obtain that
  \begin{equation*}
    \bigl\lvert\mathcal{E}_{H}\bigl(X_{i_{1}},\cdots,X_{i_{k+1}}\bigr)\bigr\rvert \lesssim \alpha_{\ell_{s}}^{(r-k-1)/r} M^{k+1},
  \end{equation*}
  where $\ell_{s}$ is the smallest integer $\ell$ that satisfies
  \begin{equation}\label{eq:defofls}
    k(2\ell+1)^{d}\geq s+k(2m+1)^{d}.
  \end{equation}
  Thus, we obtain that
  \begin{align*}
    &\sum_{\substack{
    s_{1:(k+1)}:                                                                                                                                                                                              \\
    s_{j_{0}}=s,                                                                                                                                                                                              \\
        s_{h}\leq s,~\forall h\neq j_{0}
      }}
    \bigl\lvert\mathcal{S}(H) \bigr\rvert\\
    \leq & \sigma^{-(k+1)}\sum_{\substack{
    s_{1:(k+1)}:                                                                                                                                                                                              \\
    s_{j_{0}}=s,                                                                                                                                                                                              \\
        s_{h}\leq s,\forall h\neq j_{0}
      }}
    \sum_{i_{1}\in B_{1}\backslash D_{1}}\sum_{i_{2}\in B_{2}\backslash D_{2}}\cdots \sum_{i_{k+1}\in B_{k+1}\backslash D_{k+1}}\bigl\lvert\mathcal{E}_{G}\bigl(X_{i_{1}},\cdots,X_{i_{k+1}}\bigr)\bigr\rvert \\
    \leq & 2^{2k+3}\sigma^{-(k+1)}\lvert T \rvert\bigl(k(2m+1)^{d}+s\bigr)^{k-1}\alpha_{\ell_{s}}^{(r-k-1)/r}M^{k+1}.
  \end{align*}

  Since the set $\{ s_{1},\cdots,s_{k+1}:\max_{1\leq h \leq k+1} s_{h} =s\}$
  is the union (not necessarily disjoint) of $\{ s_{1},\cdots,s_{k+1}: s_{j}=s ,s_{h}\leq s\ \forall 1\leq h\leq k+1\}$ over $2\leq j\leq k+1$, we have
  \begin{align*}
    \sum_{\substack{
    s_{1:(k+1)}:                                                                                                       \\
    \max_{1\leq h\leq k+1}s_{h}=s
    }}
    \bigl\lvert \mathcal{S}(H) \bigr\rvert
    \leq & \sum_{j_{0}=2}^{k+1}\sum_{\substack{
    s_{1:(k+1)}:                                                                                                       \\
    s_{j_{0}}=s,                                                                                                       \\
        s_{h}\leq s,\forall h\neq j_{0}
      }}
    \bigl\lvert\mathcal{S}(H) \bigr\rvert                                                                              \\
    \leq & 2^{2k+3}k\sigma^{-(k+1)}\lvert T \rvert\bigl(k(2m+1)^{d}+s\bigr)^{k-1}\alpha_{\ell_{s}}^{(r-k-1)/r}M^{k+1}.
  \end{align*}

  Next we take the sum over all possible $|T|\ge s\geq 1$ and obtain that
  \begin{align}\label{n1}
    \sum_{\substack{
    s_{1:(k+1)}:                                                                                                                                    \\
        \max_{1\leq h\leq k+1}s_h\geq 1
      }}
    \bigl\lvert \mathcal{S}(H) \bigr\rvert
    \leq & \sum_{s=1}^{\lvert T \rvert} \sum_{\substack{
    s_{1:(k+1)}:                                                                                                                                    \\
    \max_{1\leq h\leq k+1}s_{h}=s
    }}
    \bigl\lvert \mathcal{S}(H) \bigr\rvert                                                                                                          \\
    \leq & 2^{2k+3}k\sigma^{-(k+1)}\lvert T \rvert\sum_{s= 1}^{\lvert T \rvert}\bigl(k(2m+1)^{d}+s\bigr)^{k-1}\alpha_{\ell_{s}}^{(r-k-1)/r}M^{k+1}.\nonumber
  \end{align}

  To further bound this it will be important to know what is the number of different possible values of $s$ will have the same $\ell_{s}=\ell$.
  To do so, we note that by definition \eqref{eq:defofls} of $\ell_{s}$ any such $s$ much satisfy 
  \begin{equation*}
    k(2\ell_{s}+1)^{d}\geq s+k(2m+1)^{d}\geq k(2\ell_{s}-1)^{d}+1,
  \end{equation*}
  which implies that for any $\ell\geq m+1$
  \begin{equation*}
    \bigl|\{s:\ell_s=\ell\}\bigr|\leq k(2\ell+1)^{d}-k(2\ell-1)^{d}\leq 2kd(2\ell+1)^{d-1}.
  \end{equation*}
  On the other hand, the minimality of $\ell_s$ and $s\leq\lvert T\rvert$
  imply that
  \begin{equation*}
    (2\ell_s-1)^d<\frac{s}{k}+(2m+1)^d
    \leq\frac{\lvert T\rvert}{k}+(2m+1)^d.
  \end{equation*}
  Using $(a^d+b^d)^{1/d}\leq a+b$, we obtain
  \begin{equation*}
    \ell_s<m+1+\frac{1}{2}\biggl(\frac{\lvert T\rvert}{k}\biggr)^{1/d}
    \leq m+1+\frac{\lvert T\rvert^{1/d}}{2}.
  \end{equation*}
  Thus, $\ell_{s}\leq m+1+\lfloor \frac{\lvert T \rvert^{1/d}}{2}\rfloor$ for any $s$. And we conclude that
  \begin{equation*}
    \bigl|\{s:\ell_s=\ell\}\bigr|\leq \begin{cases}2kd(2\ell+1)^{d-1} & \text{if } m+1+\lfloor\frac{\lvert T\rvert^{1/d}}{2}\rfloor\ge \ell\ge m+1\\0 & \text{otherwise }\end{cases}.
  \end{equation*}
  Therefore, by combining this with \eqref{n1} we obtain that
  \begin{align*}
    &\sum_{\substack{
    s_{1:(k+1)}:                                                                                                                                                                                             \\
        \max_{1\leq h\leq k+1}s_h\geq 1
      }}
    \bigl\lvert\mathcal{S}(H) \bigr\rvert\\
    \leq     & 2^{2k+3}k\sigma^{-(k+1)}\lvert T \rvert\sum_{\ell=m+ 1}^{m+1+\lfloor\frac{\lvert T\rvert^{1/d}}{2}\rfloor} 2kd(2\ell+1)^{d-1}\bigl(k(2\ell+1)^{d}\bigr)^{k-1}\alpha_{\ell}^{(r-k-1)/r}M^{k+1} \\
    \overset{(a)}{\lesssim} & \lvert T \rvert^{-(k-1)/2}\sum_{\ell=m+ 1}^{m+1+\lfloor\frac{\lvert T\rvert^{1/d}}{2}\rfloor} \ell^{dk-1}\alpha_{\ell}^{(r-k-1)/r},
  \end{align*}
  where to obtain $(a)$ we used the assumption that the asymptotic variance does not degenerate: 
  $$\liminf_{\lvert T \rvert\to \infty} \sigma^{2}/\lvert T \rvert>0.$$

  For any $G\in \mathcal{G}_{0}(k+1)$, there exists at least one positive vertex. Thus, $\max_{1\leq h\leq k+1}s_{h}\geq 1$. Now that the number of labeled rooted trees on $k+1$ vertices only depends on $k$, we conclude
  \begin{equation*}
    \begin{aligned}
      \sum_{H \in \mathcal{G}_{0}(k+1)}\bigl\lvert\mathcal{S}(H) \bigr\rvert
      =        & \sum_{\substack{
      (V,E,\prec):                                                                                                                                 \\
          \lvert V \rvert=k+1
        }}\
      \sum_{\substack{
      s_{1:(k+1)}:                                                                                                                                 \\
      H=(V,E,s_{1:(k+1)})                                                                                                                          \\
      \in \mathcal{G}_{0}(k+1)
      }}
      \bigl\lvert\mathcal{S}\bigl(H\bigr) \bigr\rvert                                                                                              \\
      \leq     & \sum_{\substack{
      (V,E,\prec):                                                                                                                                 \\
          \lvert V \rvert=k+1
        }}
      \sum_{\substack{
      s_{1:(k+1)}:                                                                                                                                 \\
      \max_{1\leq h\leq k+1}s_{h}\geq 1
      }}
      \bigl\lvert\mathcal{S}\bigl(H\bigr) \bigr\rvert                                                                                              \\
      \lesssim & \sum_{\substack{
      (V,E,\prec):                                                                                                                                 \\
          \lvert V \rvert=k+1
        }}
      \lvert T \rvert^{-(k-1)/2}\sum_{\ell=m+1}^{m+1+\lfloor\frac{\lvert T\rvert^{1/d}}{2}\rfloor}\ell^{dk-1}\alpha_{\ell}^{(r-k-1)/r}             \\
      \lesssim & \lvert T \rvert^{-(k-1)/2}\sum_{\ell=m+1}^{m+1+\lfloor\frac{\lvert T\rvert^{1/d}}{2}\rfloor}\ell^{dk-1}\alpha_{\ell}^{(r-k-1)/r}.
    \end{aligned}
  \end{equation*}

  Next we prove \eqref{eq:remainderctrl2}.

  At the endpoint $r=k+\omega$, suppose that $\alpha_{m+1}=0$.  Every
  $H\in\mathcal G_0(k+1)$ has a positive identifier, so the exact-zero
  conclusion of \cref{thm:mixingsplit12} gives
  $\mathcal E_H(X_{i_1},\ldots,X_{i_k},\Delta_f(H))=0$ for every admissible
  index tuple.  Consequently every $\mathcal U_f(H)$ in
  \eqref{eq:remainderctrl2} is zero.  This proves the endpoint assertion
  without replacing the complete generalized-covariance block by bounds on
  its individual pieces.  It remains to establish the displayed estimate
  in the strict range $r>k+\omega$.

  We partition the sum into two disjoint cases: (i) the largest identifier
  among the first $k$ vertices is strictly larger than $s_{k+1}$, and (ii)
  $s_{k+1}$ is at least as large as every preceding identifier. In the first
  case, set $s:=\max_{1\leq j\leq k}s_j$ and
  $j_0:=\min\{j\leq k:s_j=s\}$; then $s_{j_0}=s>s_{k+1}$.

  From the claim above, we know that for any $2\leq t \leq k$,
  $t\neq j_{0}$, given a sequence $i_{1},\cdots,i_{t-1}$, we have that
  $\sum_{s_t\le s}\bigl|B_t\setminus D_t\bigr|
  \leq2(k(2m+1)^{d}+s)$.

  Let $\ell_s$ be the smallest integer $\ell$ such that
  \begin{equation*}
    k(2\ell+1)^d\geq s+k(2m+1)^d.
  \end{equation*}
  For fixed $(V,E,\prec)$ and $s_{1:k}$, let $H_a$ denote the resulting
  genogram with $s_{k+1}=a$, whenever it is admissible; an inadmissible
  contribution is understood to be zero. The coefficient $b_H$ can change
  when $s_{k+1}=-1$, but it is constant as $s_{k+1}$ ranges over the positive
  integers. Since $|b_H|\leq 1$, \eqref{eq:mixingsplit2} and
  \eqref{eq:thirdstatement} give
  \begin{align*}
    &\biggl\lvert\sum_{s_{k+1}<s}b_H\mathcal{U}_{f}(H)\biggr\rvert\\
    &\quad\leq \sigma^{-k}
      \sum_{i_{1}\in B_{1}\backslash D_{1}}\cdots
      \sum_{i_{k}\in B_{k}\backslash D_{k}}
      \Biggl\{
      \sum_{a\in\{-1,0\}}
      \bigl\lvert\mathcal{E}_{H_a}
      \bigl(X_{i_{1}},\ldots,X_{i_{k}},\Delta_f(H_a)\bigr)\bigr\rvert\\
    &\hspace{15em}+
      \biggl\lvert\sum_{1\leq s_{k+1}<s}
      \mathcal{E}_{H}\bigl(X_{i_{1}},\ldots,X_{i_{k}},\Delta_f(H)\bigr)
      \biggr\rvert\Biggr\}\\
    &\quad\lesssim \lvert f\rvert_{k-1,\omega}\sigma^{-(k+\omega)}
      \sum_{i_{1}\in B_{1}\backslash D_{1}}\cdots
      \sum_{i_{k}\in B_{k}\backslash D_{k}}
      \bigl(k(2m+1)^d+s\bigr)^\omega
      \alpha_{\ell_s}^{(r-k-\omega)/r}M^{k+\omega}.
  \end{align*}
  Here we used $|B_{k+1}\backslash D_{k+1}|\leq k(2m+1)^d$ for the two
  nonpositive terminal identifiers. Thus the $-1$ and $0$ terms are bounded
  separately, while the cancellation in the positive block is retained.

  By this choice of $j_0$, the following conditions give a disjoint partition
  of the prefix-max case. Applying the preceding bound and the claim above, we have
  \begin{align*}
    &\biggl\lvert\sum_{\substack{
    s_{1:(k+1)}:                                                                                                                                                                                                                           \\
    \max_{1\leq h\leq k}s_h=s>s_{k+1}
    }}
    b_{H}\mathcal{U}_{f}(H) \biggr\rvert\\
    \leq    &  \sum_{j_{0}=1}^{k}\sum_{\substack{
    s_{1:k}:                                                                                                                                                                                                                               \\
    s_{j_{0}}=s,                                                                                                                                                                                                                           \\
    s_h<s,~1\leq h<j_0,\\
    s_h\leq s,~j_0<h\leq k
      }}
    \biggl\lvert \sum_{s_{k+1}<s}b_{H}\mathcal{U}_{f}(H) \biggr\rvert\\
    \lesssim & k\lvert f \rvert_{k-1,\omega}\sigma^{-(k+\omega)}\lvert T \rvert\bigl(k(2m+1)^{d}+s\bigr)^{k-2+\omega}\alpha_{\ell_{s}}^{(r-k-\omega)/r}M^{k+\omega}.
  \end{align*}

  Thus, we get that
  \begin{align*}
    &\biggl\lvert \sum_{\substack{
    s_{1:(k+1)}:                                                                                                                                                                                      \\
    \max_{1\leq h\leq k}s_h\geq1,\quad
    \max_{1\leq h\leq k}s_h>s_{k+1}
      }}
    b_{H}\mathcal{U}_{f}(H)  \biggr\rvert
    \leq     \sum_{s=1}^{\lvert T \rvert}\biggl\lvert\sum_{\substack{
    s_{1:(k+1)}:                                                                                                                                                                                      \\
    \max_{1\leq h\leq k}s_h=s>s_{k+1}
    }}
    b_{H}\mathcal{U}_{f}(H) \biggr\rvert                                                                                                                                                              \\
    \lesssim & k\lvert f \rvert_{k-1,\omega}\sigma^{-(k+\omega)}\lvert T \rvert\sum_{s= 1}^{\lvert T \rvert}\bigl(k(2m+1)^{d}+s\bigr)^{k-2+\omega}\alpha_{\ell_{s}}^{(r-k-\omega)/r}M^{k+\omega}      \\
   \overset{(a)}{ \lesssim }& \lvert f \rvert_{k-1,\omega}\lvert T \rvert^{-(k+\omega-2)/2}\sum_{\ell=m+ 1}^{m+1+\lfloor\frac{\lvert T\rvert^{1/d}}{2}\rfloor} \ell^{d(k+\omega-1)-1}\alpha_{\ell}^{(r-k-\omega)/r},
  \end{align*}
The step $(a)$ again uses the standing assumption $\limsup |T|/\sigma^2<\infty$.

  We now consider the terminal-max case. For fixed $s_{1:k}$, set
  \begin{equation*}
    \begin{aligned}
      a&:=\max\{1,s_1,\ldots,s_k\},\\
      L^-_\ell&:=k(2\ell-1)^d-k(2m+1)^d+1,\\
      L^+_\ell&:=k(2\ell+1)^d-k(2m+1)^d.
    \end{aligned}
  \end{equation*}
  and let $c_{1,\ell}:=a\vee L^-_\ell$ and
  $c_{2,\ell}:=L^+_\ell+1$. We retain only the nonempty blocks
  $c_{1,\ell}<c_{2,\ell}$. On each such block the terminal identifier is
  positive, the negative-vertex pattern and hence $b_H$ are fixed, and the
  separation radius in \eqref{eq:thirdstatement} is exactly $\ell$. Therefore,
  for fixed $s_{1:k}$ and $i_{1:k}$, write
  \begin{equation*}
    \mathfrak B_\ell:=\sum_{c_{1,\ell}\leq s_{k+1}<c_{2,\ell}}
      \mathcal{E}_{H}\bigl(X_{i_{1}},\ldots,X_{i_{k}},\Delta_{f}(H)\bigr).
  \end{equation*}
  By \eqref{eq:thirdstatement}, for every $\ell\geq m+1$,
  \begin{equation}\label{eq:epshsec}
    \begin{aligned}
      &\lvert\mathfrak B_\ell\rvert\\
      &\quad\lesssim \lvert f \rvert_{k-1,\omega}\sigma^{-\omega}
      (c_{2,\ell}-c_{1,\ell})^\omega
      \alpha_{\ell}^{(r-k-\omega)/r}M^{k+\omega}\\
      &\quad\lesssim \lvert f \rvert_{k-1,\omega}\sigma^{-\omega}
      \ell^{d\omega-\omega}\alpha_{\ell}^{(r-k-\omega)/r}M^{k+\omega}.
    \end{aligned}
  \end{equation}
  Indeed,
  $c_{2,\ell}-c_{1,\ell}\leq
  k\bigl((2\ell+1)^d-(2\ell-1)^d\bigr)\lesssim\ell^{d-1}$.

  Every positive terminal identifier satisfying
  $s_{k+1}\geq\max_{1\leq j\leq k}s_j$ belongs to exactly one of these
  truncated blocks. Taking the sum over $\ell$ and $s_{1:k}$, retaining the
  block cancellation before applying the claim, gives
  \begingroup\small
  \begin{align}\label{eq:needslabel101}
    &\biggl\lvert
    \sum_{\substack{
      s_{1:(k+1)}:\\
      s_{k+1}\geq\max\{1,s_1,\ldots,s_k\}
    }}b_H\mathcal U_f(H)\biggr\rvert\nonumber\\
    &\quad\leq \sigma^{-k}
    \sum_{\ell=m+1}^{m+1+\lfloor\lvert T\rvert^{1/d}/2\rfloor}
    \sum_{\substack{s_{1:k}:\ a\leq L^+_\ell}}
    \sum_{i_1\in B_1\backslash D_1}\cdots
    \sum_{i_k\in B_k\backslash D_k}
    \lvert\mathfrak B_\ell\rvert\nonumber\\
    &\quad\lesssim \lvert f\rvert_{k-1,\omega}\sigma^{-(k+\omega)}
    \sum_{\ell=m+1}^{m+1+\lfloor\lvert T\rvert^{1/d}/2\rfloor}
    \lvert T\rvert\bigl(k(2m+1)^d+L^+_\ell\bigr)^{k-1}
    \ell^{d\omega-\omega}\alpha_\ell^{(r-k-\omega)/r}M^{k+\omega}\nonumber\\
    &\quad\lesssim \lvert f\rvert_{k-1,\omega}\sigma^{-(k+\omega)}
    \lvert T\rvert
    \sum_{\ell=m+1}^{m+1+\lfloor\lvert T\rvert^{1/d}/2\rfloor}
    \ell^{d(k+\omega-1)-\omega}\alpha_\ell^{(r-k-\omega)/r}M^{k+\omega}\nonumber\\
    &\quad\lesssim \lvert f\rvert_{k-1,\omega}
    \lvert T\rvert^{-(k+\omega-2)/2}
    \sum_{\ell=m+1}^{m+1+\lfloor\lvert T\rvert^{1/d}/2\rfloor}
    \ell^{d(k+\omega-1)-\omega}\alpha_\ell^{(r-k-\omega)/r}.
  \end{align}
  \endgroup
  Therefore, we conclude that
  \begin{align*}
    &\biggl\lvert\sum_{H \in \mathcal{G}_{0}(k+1)} b_{H}\mathcal{U}_{f}(H) \biggr\rvert\\
    \leq     & \sum_{\substack{
    (V,E,\prec):                                                                                                                                                                                                                                                                       \\
        \lvert V \rvert=k+1
      }}
    \biggl\lvert
    \sum_{\substack{
    s_{1:(k+1)}:                                                                                                                                                                                                                                                                       \\
    H=(V,E,s_{1:(k+1)})                                                                                                                                                                                                                                                                \\
    \in \mathcal{G}_{0}(k+1)
    }}
    b_{H}\mathcal{U}_{f}(H) \biggr\rvert
    \leq   \sum_{\substack{
    (V,E,\prec):                                                                                                                                                                                                                                                                       \\
        \lvert V \rvert=k+1
      }}
    \biggl\lvert\sum_{\substack{
    s_{1:(k+1)}:                                                                                                                                                                                                                                                                       \\
    \max_{1\leq h\leq k+1}s_{h}\geq 1
    }}
    b_{H}\mathcal{U}_{f}(H) \biggr\rvert                                                                                                                                                                                                                                               \\
    \leq     & \sum_{\substack{
    (V,E,\prec):                                                                                                                                                                                                                                                                       \\
        \lvert V \rvert=k+1
      }}
    \biggl\lvert\sum_{\substack{
    s_{1:(k+1)}:                                                                                                                                                                                                                                                                       \\
    s_{j}\leq s_{k+1},\forall 1\leq j\leq k                                                                                                                                                                                                                                            \\
        s_{k+1}\geq 1
      }}
    b_{H}\mathcal{U}_{f}(H) \biggr\rvert
    +  \sum_{\substack{
    (V,E,\prec):                                                                                                                                                                                                                                                                       \\
        \lvert V \rvert=k+1
      }}
    \biggl\lvert \sum_{\substack{
    s_{1:(k+1)}:                                                                                                                                                                                                                                                                       \\
    \max_{1\leq h\leq k}s_h\geq1,\quad
    \max_{1\leq h\leq k}s_h>s_{k+1}
      }}
    b_{H}\mathcal{U}_{f}(H)  \biggr\rvert                                                                                                                                                                                                                                              \\
    \lesssim & \sum_{\substack{
    (V,E,\prec):                                                                                                                                                                                                                                                                       \\
        \lvert V \rvert=k+1
    }}                                                                                   \lvert f \rvert_{k-1,\omega} \lvert T \rvert^{-(k+\omega-2)/2}\sum_{\ell=m+1}^{m+1+\lfloor\frac{\lvert T\rvert^{1/d}}{2}\rfloor}\ell^{d(k+\omega-1)-\omega }\alpha_{\ell}^{(r-k- \omega  )/r} \\
    \lesssim & \lvert f \rvert_{k-1,\omega} \lvert T \rvert^{-(k+\omega-2)/2}\sum_{\ell=m+1}^{m+1+\lfloor\frac{\lvert T\rvert^{1/d}}{2}\rfloor}\ell^{d(k+\omega-1)-\omega }\alpha_{\ell}^{(r-k- \omega  )/r}.
  \end{align*}

  Next we prove \eqref{eq:remainderctrl3}. If $H\in \mathcal{P}_{0}(k+1)$, then for any $t\leq k$ we know that  $i_{t+1}\in N(i_{1:t})$. In other words, new indices lie in the $m$-neighborhood of previous ones. By \eqref{eq:nomixingsplit2}, we have
  \begin{equation*}
    \bigl\lvert\mathcal{E}_{H}\bigl(X_{i_{1}},\cdots,X_{i_{k}},\Delta_{f}(H)\bigr)\bigr\rvert \leq 2^{k} \sigma^{-\omega}\bigl\lvert B_{k+1}\backslash D_{k+1} \bigr\rvert^{\omega}\cdot \lvert f \rvert_{k-1,\omega}M^{k+\omega},
  \end{equation*}

  Taking the sums over $i_{j}\in B_{j}\backslash D_{j}$ for all $1\leq j\leq k$, we get
  \begin{align*}
    \bigl\lvert \mathcal{U}_{f}(H) \bigr\rvert
    \leq     & 2^{k} \sigma^{-\omega}\bigl\lvert B_{k+1}\backslash D_{k+1} \bigr\rvert^{\omega}\cdot \lvert f \rvert_{k-1,\omega}M^{k+\omega}\prod_{j=1}^{k}\bigl\lvert B_{j}\backslash D_{j} \bigr\rvert \\
    \lesssim & \lvert f \rvert_{k-1,\omega}\lvert T \rvert^{-(k+\omega-2)/2}m^{d(k+\omega-1)}.
  \end{align*}
  Therefore, we have
  \begingroup\small
  \begin{align*}
    \biggl\lvert\sum_{H \in \mathcal{P}_{0}(k+1)}\!\!\!\!b_{H}\mathcal{U}_{f}(H) \biggr\rvert
    \leq     & \sum_{\substack{
    (V,E,\prec):                                                                                  \\
        \lvert V \rvert=k+1
      }}
    \biggl\lvert
    \sum_{\substack{
    s_{1:(k+1)}:                                                                                  \\
    H=(V,E,s_{1:(k+1)})                                                                           \\
    \in \mathcal{P}_{0}(k+1)
    }}
    b_{H}\mathcal{U}_{f}(H) \biggr\rvert
    \leq  \sum_{\substack{
    (V,E,\prec):                                                                                  \\
        \lvert V \rvert=k+1
      }}
    \biggl\lvert\sum_{\substack{
    s_{1:(k+1)}:                                                                                  \\
    s_{j}=0\text{ or }-1,                                                                         \\
        \forall 1\leq j\leq k+1
      }}
    b_{H}\mathcal{U}_{f}(H) \biggr\rvert                                                          \\
    \leq     & \sum_{\substack{
    (V,E,\prec):                                                                                  \\
        \lvert V \rvert=k+1
      }}
    \sum_{\substack{
    s_{1:(k+1)}:                                                                                  \\
    s_{j}=0\text{ or }-1,                                                                         \\
        \forall 1\leq j\leq k+1
      }}
    \lvert b_{H}\rvert \bigl\lvert \mathcal{U}_{f}(H) \bigr\rvert
    \leq   \sum_{\substack{
    (V,E,\prec):                                                                                  \\
        \lvert V \rvert=k+1
      }}
    \sum_{\substack{
    s_{1:(k+1)}:                                                                                  \\
    s_{j}=0\text{ or }-1,                                                                         \\
        \forall 1\leq j\leq k+1
      }}
    \bigl\lvert \mathcal{U}_{f}(H) \bigr\rvert                                                    \\
    \lesssim & \sum_{\substack{
    (V,E,\prec):                                                                                  \\
        \lvert V \rvert=k+1
    }}  \lvert f \rvert_{k-1, \omega  }\lvert T \rvert^{-(k+\omega-2)/2}m^{d(k+\omega-1)}         \\
    \lesssim & \lvert f \rvert_{k-1, \omega  }\lvert T \rvert^{-(k+\omega-2)/2}m^{d(k+\omega-1)}.
  \end{align*}
  \endgroup

  Finally, to prove \eqref{eq:remainderctrl4}, we follow a derivation similar to \eqref{eq:needslabel101}. For a fixed ordered tree and a fixed compatible vector $s_{1:k}\in\{0,-1\}^k$, the coefficient $b_H$ is constant as the positive identifier $s_{k+1}$ varies; denote this common value by $b_{V,E,s_{1:k}}$. We group first by this fixed negative-vertex pattern. All sums over $s_{k+1}$ below are restricted to values for which the resulting $H$ is a genogram in $\mathcal{P}_1(k+1)$.
  \begingroup\small
  \begin{align*}
    &\biggl\lvert
    \sum_{H \in \mathcal{P}_{1}(k+1)}\!\!\!\!
    b_{H}\mathcal{U}_{f}(H)
    \biggr\rvert\\
    \leq  &                    
    \sum_{\substack{
    (V,E,\prec):                                                                                                                                                                                                                   \\
        \lvert V \rvert=k+1
      }}
    \biggl\lvert
    \sum_{\substack{
    s_{1:(k+1)}:                                                                                                                                                                                                                   \\
    H=(V,E,s_{1:(k+1)})                                                                                                                                                                                                            \\
    \in \mathcal{P}_{1}(k+1)
    }}
    b_{H}\mathcal{U}_{f}(H)
    \biggr\rvert
    \leq \sum_{\substack{
    (V,E,\prec):                                                                                                                                                                                                                   \\
        \lvert V \rvert=k+1
      }}
    \biggl\lvert
    \sum_{\substack{
    s_{1:(k+1)}:                                                                                                                                                                                                                   \\
    s_{j}\leq 0,\forall 1\leq j\leq k,                                                                                                                                                                                             \\
        s_{k+1}\geq 1
      }}
    b_{H}\mathcal{U}_{f}(H)
    \biggr\rvert                                                                                                                                                                                                                   \\
    \leq &
    \sum_{\substack{
    (V,E,\prec):\\
        \lvert V \rvert=k+1
      }}
    \sum_{\substack{
    s_{1:k}\in\{0,-1\}^{k}:\\
    s_{1:k}\text{ admits a compatible }\mathcal{P}_1(k+1)\text{ extension}
      }}
    \lvert b_{V,E,s_{1:k}}\rvert
    \biggl\lvert
    \sum_{s_{k+1}\geq1}\mathcal{U}_{f}(H)
    \biggr\rvert\\
    \leq &
    \sum_{\substack{
    (V,E,\prec):\\
        \lvert V \rvert=k+1
      }}
    \sum_{\substack{
    s_{1:k}\in\{0,-1\}^{k}:\\
    s_{1:k}\text{ admits a compatible }\mathcal{P}_1(k+1)\text{ extension}
      }}
    \biggl\lvert
    \sum_{s_{k+1}\geq1}\mathcal{U}_{f}(H)
    \biggr\rvert\\
    \leq                      & \sigma^{-k}\\
    &\quad\times
    \sum_{\substack{
    (V,E,\prec):                                                                                                                                                                                                                   \\
        \lvert V \rvert=k+1
      }}
    \sum_{\substack{
    s_{1:k}\in\{0,-1\}^{k}:\\
    s_{1:k}\text{ admits a compatible }\mathcal{P}_1(k+1)\text{ extension}
      }}\\
    &\quad\times
    \sum_{\ell=m+1}^{m+1+\lfloor\frac{\lvert T\rvert^{1/d}}{2}\rfloor}
    \left\lvert
      \sum_{\substack{
        s_{k+1}=k(2\ell-1)^{d}\\
        {}-k(2m+1)^{d}+1
      }}^{\substack{
        k(2\ell+1)^{d}\\
        {}-k(2m+1)^{d}
      }}
    \mathcal{E}_{H}\bigl(X_{i_{1}},\cdots,X_{i_{k}},\Delta_{f}(H)\bigr)
    \right\rvert                                                                                                                                                                                                                   \\
    \leq                      & \sigma^{-k}\\
    &\quad\times
    \sum_{\substack{
    (V,E,\prec):                                                                                                                                                                                                                   \\
        \lvert V \rvert=k+1
      }}
    \sum_{\substack{
    s_{1:k}\in\{0,-1\}^{k}:\\
    s_{1:k}\text{ admits a compatible }\mathcal{P}_1(k+1)\text{ extension}
      }}
    \sum_{\ell=m+1}^{m+1+\lfloor\frac{\lvert T\rvert^{1/d}}{2}\rfloor}\\
    &\quad\times
    \sum_{i_{1}\in B_{1}\backslash D_{1}}\!\!\!\!\cdots\!\!\!\!\sum_{i_{k}\in B_{k}\backslash D_{k}}
    \left\lvert
    \sum_{\substack{
      s_{k+1}=k(2\ell-1)^{d}\\
      {}-k(2m+1)^{d}+1
    }}^{\substack{
      k(2\ell+1)^{d}\\
      {}-k(2m+1)^{d}
    }}\!\!\!\!
    \mathcal{E}_{H}\bigl(X_{i_{1}},\cdots,X_{i_{k}},\Delta_{f}(H)\bigr)
    \right\rvert                                                                                                                                                                                                                   \\
    \overset{(*)}{\lesssim }  & \lvert f \rvert_{k-2,1}^{1-\delta}
      \lvert f \rvert_{k-1,1}^{\delta}\sigma^{-(k+\delta)}\\
    &\quad\times
    \sum_{\substack{
    (V,E,\prec):                                                                                                                                                                                                                   \\
        \lvert V \rvert=k+1
      }}
    \sum_{\substack{
    s_{1:k}\in\{0,-1\}^{k}:\\
    s_{1:k}\text{ admits a compatible }\mathcal{P}_1(k+1)\text{ extension}
      }}
    \sum_{\ell=m+1}^{m+1+\lfloor\frac{\lvert T\rvert^{1/d}}{2}\rfloor}\\
    &\quad\times
    \sum_{i_{1}\in B_{1}\backslash D_{1}}\!\!\!\!\cdots\!\!\!\!\sum_{i_{k}\in B_{k}\backslash D_{k}}\!\!
    \ell^{d\delta-\delta}\alpha_{\ell}^{(r-k-\delta)/r}                                                                                                                                                                            \\
    \overset{(**)}{\lesssim } &\lvert f \rvert_{k-2,1}^{1-\delta} \lvert f \rvert_{k-1,1}^{\delta}\sigma^{-(k+\delta)}2^{k-1}\lvert T \rvert\bigl(k(2m+1)^{d}\bigr)^{k-1}\\
                                &\quad\times\sum_{\ell=m+1}^{m+1+\lfloor\frac{\lvert T\rvert^{1/d}}{2}\rfloor}\ell^{d\delta-\delta}\alpha_{\ell}^{(r-k-\delta)/r}                                          \\
    \lesssim                  & \lvert f \rvert_{k-2,1}^{1-\delta}
      \lvert f \rvert_{k-1,1}^{\delta}
      \lvert T \rvert^{-(k+\delta-2)/2}m^{d(k-1)}\\
    &\quad\times\sum_{\ell=m+1}^{m+1+\lfloor\frac{\lvert T\rvert^{1/d}}{2}\rfloor}
      \ell^{d\delta-\delta}\alpha_{\ell}^{(r-k-\delta)/r}.
  \end{align*}
  \endgroup
  Note that the inequality $(*)$ is due to \eqref{eq:thirdstatementnew}. For $(**)$, we note that $s_{j}=0$ or $-1$ for $2\leq j\leq k$, and that the number of choices for $i_{j}$ ($2\leq j\leq k$) is upper-bounded by $k(2m+1)^{d}$ since $i_{j}$ lies in the $m$-neighborhood of $i_{1},\cdots,i_{j-1}$.\hfill
\end{proof}

  \section{Extended overview of the constructive graph approach}\label{sec:technical-overview}
  \subsection[Terms E4 to E9 rewritten using D-star notation]{Terms $E_4$ to $E_9$ rewritten using the $\mathcal{D}^*$ notation}\label{re-express}
  \begingroup\small
  \begin{align*}
  E_{4}
  =             & \frac{1}{\sigma_{n}^{2}}\sum_{i=1}^{n}\sum_{j=m+1}^{n-1}\mathcal{D}^{*} \biggl(X_{i}\,,\,X_{i-j}\int_{0}^{1}\Bigl(f'\bigl(\nu W_{i,j-1}+(1-\nu)W_{i,j}^{*}\bigr)-f'(W_{i,j-1})\Bigr)\dif\nu\biggr),                                                                         \\
                & \qquad + \frac{1}{\sigma_{n}^{2}}\sum_{i=1}^{n}\sum_{j=m+1}^{n-1}\mathcal{D}^{*} \biggl(X_{i}\,,\,X_{i+j}\int_{0}^{1}\Bigl(f'\bigl(\nu W_{i,j}^{*}+(1-\nu)W_{i,j}\bigr)-f'(W_{i,j}^{*})\Bigr)\dif\nu\biggr),                                                                     \\
  E_{5}=        & \frac{1}{\sigma_{n}^{2}}\sum_{i=1}^{n}\sum_{j=m+1}^{n-1}\Bigl(\mathcal{D}^{*} \bigl(X_{i}\,,\,X_{i-j}\,,\,f'(W_{i,j-1})\bigr)+\mathcal{D}^{*} \bigl(X_{i}\,,\,X_{i+j}\,,\,f'(W_{i,j}^{*})\bigr)\Bigr),                                                               \\
  E_{6}=        & \frac{1}{\sigma_{n}^{2}}\sum_{i=1}^{n}\sum_{j=m+1}^{n-1}\Bigl(\mathcal{D}^{*} \bigl(X_{i}\,,\,X_{i-j}\bigr)\ \mathcal{D}^{*} \bigl(f'(W_{i,j-1})-f'(W_{n})\bigr)\\*
  &\qquad\qquad\qquad\qquad\qquad\quad+\mathcal{D}^{*} \bigl(X_{i}\,,\,X_{i+j}\bigr)\ \mathcal{D}^{*} \bigl(f'(W_{i,j}^{*})-f'(W_{n})\bigr)\Bigr),
  \\
  E_{7}=&   \frac{1}{\sigma_{n}^{2}}\sum_{i=1}^{n}\sum_{j:\lvert j-i \rvert\leq m} \mathcal{D}^{\ast}(X_{i},X_{j})\ \mathcal{D}^{\ast}\bigl(f'(W_{i,j,m})-f'(W_{n})\bigr),      \\
  E_{8}=&   \frac{1}{\sigma_{n}^{2}}\sum_{i=1}^{n}\sum_{j:\lvert j-i \rvert\leq m}\mathcal{D}^{\ast}\bigl(X_{i},X_{j}, f'(W_{n})-f'(W_{i,j,m})\bigr),\\
  E_{9}=&   \frac{1}{\sigma_{n}^{2}}\sum_{i=1}^{n}\sum_{j:\lvert j-i \rvert\leq m} \mathcal{D}^{\ast}\bigl(X_{i},X_{j}, f'(W_{i,j,m})\bigr),\\
  F_{1}=&   \frac{1}{\sigma_{n}^{2}}\sum_{i=1}^{n}\sum_{j=i-m}^{i+m}\ \mathcal{D}^{*} (X_{i},X_{j})\ \mathbb{E} [f'(W_{n})],\\
  F_{2}=&   \frac{1}{\sigma_{n}^{2}}\sum_{i=1}^{n}\sum_{j=m+1}^{n-1}\mathcal{D}^{*} (X_{i}\,,\,X_{i-j}+X_{i+j})\ \mathbb{E} [f'(W_{n})].
  \end{align*}
  \endgroup
\subsection{Complete low-order derivation and conceptual guide}\label{sec:technical-overview-full}

This subsection supplies the complete low-order derivation and conceptual table
supporting the concise account in Section~1.3 of the main article.

The main methodological problem is to establish \eqref{eq:mysterious} at any
prescribed finite order. Although its principal part resembles the
moment--cumulant identity below, a direct expansion under long-range
dependence produces three coupled difficulties. The summation ranges are
nested and depend on previously chosen indices; Taylor and telescoping
operations generate products of generalized covariance terms; and these
terms must be recombined exactly to recover genuine cumulants. At the same
time, near and separated spatial configurations require different estimates.
In particular, taking absolute values term by term would discard
cancellations needed for the stated remainder bounds.

We call the resulting integer-decorated rooted trees \emph{spatial-expansion
genograms}, or simply \emph{genograms}. A vertex represents a summation index, the
rooted-tree structure records traversal and factorization, and the integer
identifiers record the associated spatial constraints. Adding a path
represents a Taylor expansion, whereas attaching a child with varying
identifiers represents a telescoping operation. The genealogy, identifiers,
and signed coefficients are therefore handled jointly: related genograms are
first grouped according to the exact algebraic identities, and only then are
moment, neighborhood-volume, and mixing estimates applied. This order of
operations makes it possible to use the same dependence-free expansion at
each prescribed order while confining the probabilistic input to the
remainder estimates.

The polynomial case illustrates why the desired principal terms should be
cumulants.
\begin{lemma}\label{thm:identitycumulant:technical}
  Let $X$ be a mean-zero random variable with finite moments
  $\mu_{1},\ldots,\mu_{k+1}$ and cumulants
  $\kappa_{1},\ldots,\kappa_{k+1}$. Then
  \begin{equation}
    \mu_{k+1}=\sum_{j=1}^{k-1}\frac{k!}{j!(k-j)!}\kappa_{j+1}\mu_{k-j}+\kappa_{k+1}.
  \end{equation}
\end{lemma}

Thus, \eqref{eq:mysterious} is exact with vanishing remainder when $f$ is a
polynomial. For a general H\"older-smooth function it need not hold for an
arbitrary $W_n$; our task is to prove it for the normalized field sum and to
control the remainder through the $\alpha$-mixing coefficients.

To illustrate the mechanism, we begin with $k=\omega=d=1$ for a stationary,
mean-zero random sequence with finite third moments. Our objective is to
derive the first-order form of \eqref{eq:mysterious} by expanding
$\mathbb{E}[W_n f(W_n)]-\mathbb{E}[f'(W_n)]$.
For any positive $i,m\geq 1$, $i+m\leq n$, we denote
\begin{equation}\label{eq:definewim:technical}
  W_{i,m}:=\frac{1}{\sigma_{n}}\sum_{\ell: \lvert\ell-i\rvert\geq m+1}X_{\ell},\quad W_{i,m}^{*}:=W_{i,m}+\frac{1}{\sigma_{n}}X_{i+m}.
\end{equation}
Intuitively for each $i$, 
\begin{itemize}
  \item we split $W_{n}$ into two parts, $W_{i,m}$ and $W_{n}-W_{i,m}$;
  \item $W_{n}-W_{i,m}$  is the sum of at most $2m+1$ $X_{j}$'s, and can converge to $0$ when $m\ll n\to\infty$;
  \item $W_{i,m}$ could include a large number of $X_{j}$'s in the sum but the dependence between $W_{i,m}$ and $X_{i}$ is low because they are at distance $m+1$ from each other; and
  \item $W_{i,m}^{\ast}$ contains one more term than $W_{i,m}$ and one fewer than $W_{i,m-1}$, serving as a bridge between $W_{i,m}$ and $W_{i,m-1}$ for technical reasons.
\end{itemize}

Now fixing $1\leq m<n$ we can check that the following expansion holds:
\begin{equation}\label{eq:split:technical}
\begin{aligned}
  \mathbb{E}[ W_{n}f (W_{n})]-\mathbb{E} [ f'(W_{n})]
  &=\frac{1}{\sigma_n}\sum_{i=1}^{n}\mathbb{E}[X_{i}f(W_n)] \\
  &\quad-\frac{1}{\sigma_{n}^{2}}\sum_{i=1}^{n}\sum_{j=1}^{n}
  \mathbb{E} [X_{i}X_{j}]\ \mathbb{E} [f'(W_{n})] \\
  &=E_{1}+E_{2}+E_{3}.
\end{aligned}
\end{equation}
where
\begin{align}
  &E_{1}=\frac{1}{\sigma_{n}}
  \sum_{i=1}^{n}\mathbb{E} \bigl[X_{i}\bigl(f (W_{n})-f (W_{i,m}) - f'(W_{n})(W_{n}-W_{i,m})\bigr)\bigr],\\
  &E_{2}=\frac{1}{\sigma_{n}}\sum_{i=1}^{n}\mathbb{E} [X_{i}f (W_{i,m})]-\frac{1}{\sigma_{n}^{2}}\sum_{i=1}^{n}\sum_{j: \lvert j-i \rvert\geq m+1}\mathbb{E} [X_{i}X_{j}]\ \mathbb{E} [f'(W_{n})],\\
  &E_{3}=\frac{1}{\sigma_{n}}  \sum_{i=1}^{n}\mathbb{E} \bigl[X_{i}(W_{n}-W_{i,m})f'(W_{n})\bigr]-\frac{1}{\sigma_{n}^{2}}\sum_{i=1}^{n}\sum_{j:\lvert j-i \rvert\leq m}\ \mathbb{E} [X_{i}X_{j}]\ \mathbb{E} [f'(W_{n})].
\end{align}
The separation here is designed so that $E_{1}$, $E_{2}$, and $E_{3}$ can be controlled using distinct properties of the random variables involved in the summations. Specifically, $E_{1}$ is bounded by applying the Taylor expansion in conjunction with the Hölder condition, contributing to the remainder term in \eqref{eq:mysterious}. In contrast, $E_{2}$ requires analysis leveraging the $\alpha$-mixing properties.

\textbf{Bounding $E_1$ using the Taylor expansion and Hölder conditions:}
\begin{equation}\label{eq:control1111:technical}
\begin{aligned}
  &\lvert E_{1} \rvert =  \frac{1}{\sigma_{n}}\biggl\lvert
  \sum_{i=1}^{n}\mathbb{E} \bigl[X_{i}\bigl(f (W_{n})-f (W_{i,m})
  - f'(W_{n})(W_{n}-W_{i,m})\bigr)\bigr]\biggr\rvert \\     
  \leq &                 \frac{\lVert f''\rVert}{2\sigma_{n}}\sum_{i=1}^{n}
  \mathbb{E}\bigl[\bigl\lvert X_{i}(W_{n}-W_{i,m})^{2}\bigr\rvert\bigr] 
= \frac{\lVert f''\rVert}{2\sigma_{n}^{3}}\sum_{i=1}^{n}
  \sum_{j:\lvert j-i \rvert\leq m}\sum_{\ell:\lvert \ell-i\rvert\leq m}
  \mathbb{E}[\lvert X_{i}X_{j}X_{\ell}\rvert]\\                            
    \leq &\frac{2(m+1)^{2}\lVert f''\rVert}{\sigma_{n}^{3}}\sum_{i=1}^{n}
  \mathbb{E}[\lvert X_{i}\rvert^{3}]\lesssim \lVert f'' \rVert m^{2} n^{-1/2}.
\end{aligned}
\end{equation}

\textbf{Bounding $E_2$ using a telescoping sum argument and $\boldsymbol{\alpha}$-mixing properties:}
\begin{align}
  E_{2}
  &= \frac{1}{\sigma_{n}}\sum_{i=1}^{n}\mathbb{E} [X_{i}f(W_{i,m})] \notag\\
  &\quad-\frac{1}{\sigma_{n}^{2}}\sum_{i=1}^{n}
  \sum_{j: \lvert j-i \rvert\geq m+1}\mathbb{E} [X_{i}X_{j}]\ \mathbb{E} [f'(W_{n})] \notag\\
  &=\frac{1}{\sigma_{n}}\sum_{i=1}^{n}\sum_{j\geq m+1}
  \mathbb{E} \bigl[X_{i}\bigl(f(W_{i,j-1})-f(W_{i,j})\bigr)\bigr] \notag\\
  &\quad-\frac{1}{\sigma_{n}^{2}}\sum_{i=1}^{n}
  \sum_{j: \lvert j-i \rvert\geq m+1}\mathbb{E} [X_{i}X_{j}]\ \mathbb{E} [f'(W_{n})] \notag\\
  &=E_{4}+E_{5}+E_{6}.\label{eq:split2:technical}
\end{align}
where
\begin{align}
E_{4}={}&\frac{1}{\sigma_{n}}\sum_{i=1}^{n}\sum_{j\geq m+1}
  \mathbb{E} \Bigl[X_{i}\Bigl(f(W_{i,j-1})-f(W_{i,j}) \notag\\
&\qquad-f'(W_{i,j-1})(W_{i,j-1}-W_{i,j}^{*})
  -f'(W_{i,j}^{*})(W_{i,j}^{*}-W_{i,j})\Bigr)\Bigr],\\
E_{5}={}&\frac{1}{\sigma_{n}^{2}}\sum_{i=1}^{n}\sum_{j\geq m+1}
  \mathbb{E} \Bigl[X_{i}\Bigl(
  X_{i-j}\bigl(f'(W_{i,j})-\mathbb{E} [f'(W_{i,j})]\bigr) \notag\\
&\qquad+X_{i+j}\bigl(f'(W_{i,j}^{*})-\mathbb{E} [f'(W_{i,j}^{*})]\bigr)
  \Bigr)\Bigr],\\
E_{6}={}&\frac{1}{\sigma_{n}^{2}}\sum_{i=1}^{n}\sum_{j\geq m+1}\Bigl(
  \mathbb{E} [X_{i}X_{i-j}]
  \bigl(\mathbb{E} [f'(W_{i,j})]-\mathbb{E}[f'(W_{n})]\bigr) \notag\\
&\qquad+\mathbb{E} [X_{i}X_{i+j}]
  \bigl(\mathbb{E} [f'(W_{i,j}^{*})]-\mathbb{E} [f'(W_{n})] \bigr)\Bigr).
\end{align}
The terms $E_{4}$ to $E_{6}$ are analyzed by exploiting the $\alpha$-mixing properties of $(X_{i})_{i=1}^{n}$. For instance, in $E_{5}$, the $\alpha$-mixing coefficient between the $\sigma$-algebra generated by $X_{i}$ and the $\sigma$-algebra generated by $X_{i-j}, X_{i+j}, W_{i,j}, W_{i,j}^{\ast}$ is at most $\alpha_{j}$. Consequently, we can proceed to bound $E_{5}$ by applying the covariance inequality as provided in \cref{thm:covineq}. 

Additionally, we note that $W_{i,j}^{\ast}$ is introduced to facilitate a more precise analysis, which is especially important when the dimension $d \geq 2$. In this case, the number of terms in $W_{i,j-1} - W_{i,j}$ grows with $j$ at the rate of $\mathcal{O}(j^{d-1})$, leading to significantly larger remainder terms if only coarse bounds are applied.

\textbf{Bounding $E_3$ using a mixture of the Hölder condition and $\boldsymbol{\alpha}$-mixing properties:}
\begin{align}\label{eq:split3:technical}
  E_{3}= & 
  \frac{1}{\sigma_{n}^{2}}\sum_{i=1}^{n}\sum_{j:\lvert j-i\rvert\leq m}\Bigl(\mathbb{E} [X_{i}X_{j}f'(W_{n})] -\mathbb{E} [X_{i}X_{j}]\ \mathbb{E} [f'(W_{n})]\Bigr) = E_{7}+E_{8}+E_{9},
\end{align}
where
\vspace*{-15pt}
\begin{align}
  & E_{7}=\frac{1}{\sigma_{n}^{2}}\sum_{i=1}^{n}\sum_{j:\lvert j-i\rvert\leq m}\mathbb{E} [X_{i}X_{j}]\bigl(\mathbb{E} [f'(W_{i,j,m})]-\mathbb{E} [f'(W_{n})]\bigr),\\
  & E_{8}=\frac{1}{\sigma_{n}^{2}}\sum_{i=1}^{n}\sum_{j:\lvert j-i\rvert\leq m}\mathbb{E} \bigl[X_{i}X_{j}\bigl(f'(W_{n})-f'(W_{i,j,m})\bigr)\bigr],\\
  & E_{9}=\frac{1}{\sigma_{n}^{2}}\sum_{i=1}^{n}\sum_{j:\lvert j-i\rvert\leq m}\mathbb{E} \bigl[X_{i}X_{j}\bigl(\mathbb{E} [f'(W_{i,j,m})]-f'(W_{i,j,m})\bigr)\bigr],
\end{align}
where we set
\vspace*{-5pt}
\begin{equation*}
  W_{i,j,m}:=\frac{1}{\sigma_{n}}\sum_{\ell: \lvert \ell-i \rvert\wedge \lvert \ell-j \rvert\geq m+1}X_{\ell}.
\end{equation*}
Here $E_{7}$ and $E_{8}$ can be bounded by exploiting the Hölder condition of $f'$ in a way similar to \eqref{eq:control1111:technical}, and $E_{9}$ can be controlled with the $\alpha$-mixing properties of $(X_{i})_{i=1}^{n}$. In summary, we obtain the expansion
\begin{equation}
  \mathbb{E} [W_{n}f(W_{n})]=\mathbb{E} [f'(W_{n})]+E_{1}+E_{4}+E_{5}+E_{6}+E_{7}+E_{8}+E_{9},
\end{equation}
where $E_{1}$ and $E_{4}$ to $E_{9}$ together give the remainder terms.

The preceding calculation is a low-order instance of a systematic finite
construction. A genogram records each summation index as a vertex, its nested
domain through traversal and ancestor-based identifiers, and each generalized
covariance factor as a branch. Path insertion and varying-identifier child
extensions implement Taylor expansion and telescoping. A term's order is its
number of vertices: in \eqref{eq:mysterious}, principal terms arise through
order $k+1$ and remainders have order $k+2$. Related genograms remain in signed
blocks until the exact cumulant and telescoping identities are formed.
\Cref{tab:concept:technical} summarizes these roles; the bijections are developed in
\cref{sec:genoroadmap} and defined formally in \cref{sec:mixingmainpart}.

The finite-sum identities do not use mixing. Their remainder estimates use
moments, variance, lattice volume growth, and a separated-block covariance
inequality, and already cover possibly nonstationary fields on general
increasing finite lattice subsets. Given analogous bounds, the construction
may adapt to graph- or metric-indexed fields, other weak-dependence notions,
and local statistics. These extensions are not proved here; multivariate and
functional versions also need new Stein and Edgeworth input.

\begin{table}[t]
  \centering
  \caption{Roles of the genogram construction in the high-order expansion.}
  \label{tab:concept:technical}
  \resizebox{\linewidth}{!}{%
    \begin{minipage}{1.30\linewidth}
      \small
      \setlength{\tabcolsep}{4pt}
      \renewcommand{\arraystretch}{1.12}
      \begin{tabularx}{\linewidth}{
        @{}
        >{\raggedright\arraybackslash}p{0.21\linewidth}
        >{\raggedright\arraybackslash}X
        >{\raggedright\arraybackslash}p{0.28\linewidth}
        @{}}
        \toprule
        Difficulty & Genogram mechanism & Consequence \\
        \midrule
        Nested relative sums
          & Traversal and ancestor-based identifiers encode the admissible
            domain of each successive index.
          & A finite systematic representation replaces term-by-term
            enumeration. \\
        Products of generalized covariances
          & Rooted-tree branches and identifier signs encode the successive
            factors.
          & Factorizations remain visible, and their signed sums recover
            genuine cumulants. \\
        Taylor expansions and finite telescopes
          & Path insertion represents Taylor expansion; varying-identifier
            child extensions represent telescoping.
          & Exact recursive identities reach any prescribed finite order. \\
        Near and separated configurations
          & Identifier magnitudes and neighborhood shells locate the first
            block separated from its ancestors.
          & Local terms use volume counts; separated terms use mixing
            estimates. \\
        Cancellation among related terms
          & Coefficients and blockwise identifier sums are kept intact before
            absolute values are taken.
          & Telescoping cancellation gives the powers in the stated remainder
            bounds. \\
        \bottomrule
      \end{tabularx}
    \end{minipage}%
  }
\end{table}
  
\ManuscriptAfterAppendices
 
\end{document}